\documentclass[12pt]{amsart}
\usepackage{amssymb,amsmath,amsfonts,amsthm,tikz,multirow,wasysym,color}
\usetikzlibrary{calc,arrows}

\usepackage[hidelinks, hyperindex]{hyperref}
\usepackage{subfig,overpic}
\usepackage[numbers]{natbib}
\usepackage{url}
\usepackage{doi}
\usepackage[normalem]{ulem}

\hyperref[sec:function]{}

\newcommand{\sub}{\subset}

\newcommand{\bb}{\mathbb}

\newcommand\aaa{\alpha}
\newcommand\bbb{\beta}
\newcommand\ccc{\gamma}
\newcommand\ddd{\delta}
\newcommand\eee{\epsilon}

\newcommand{\thin}{\hspace{0.1em}\rule{0.7pt}{0.8em}\hspace{0.1em}}
\newcommand{\thick}{\hspace{0.1em}\rule{1.5pt}{0.8em}\hspace{0.1em}}

\newtheorem{theorem}{Theorem}
\newtheorem{lemma}[theorem]{Lemma}

\newtheorem*{theorem*}{Theorem}

\theoremstyle{definition}
\newtheorem*{definition*}{Definition}

\newtheorem*{case*}{Case}

\newtheorem*{subcase*}{Subcase}

\newtheorem*{subsubcase*}{Subsubcase}

\def\Z{{\mathbb{Z}}}

\theoremstyle{remark}

\newtheorem{example}[equation]{Example}

\numberwithin{equation}{section}

\title[Tilings of the sphere by congruent pentagons V]{Tilings of the sphere by congruent pentagons V: Edge combination $a^{4}b$ with rational angles}

\subjclass[2020]{Primary 52C20, 05B45; Secondary 11R18, 11Y50, 14Q25.}
\keywords{spherical tiling, almost equilateral pentagon, classification, trigonometric Diophantine equation, root of unity.}

\author{Jinjin Liang}
\address{
School of Mathematical Sciences\\
Zhejiang Normal University\\
Jinhua 321004\\
P.R. China}
\email{liangjinjin@zjnu.edu.cn}
\author{Yixi Liao}
\address{
School of Mathematical Sciences\\
Zhejiang Normal University\\
Jinhua 321004\\
P.R. China}
\email{yixiliao@zjnu.edu.cn}
\author{Wenchuan Hu}
\address{
School of Mathematics\\
Sichuan University\\
Chengdu 610064\\
P.R. China}
\email{huwenchuan@gmail.com,wenchuan@scu.edu.cn}
\author{Erxiao Wang}
\address{
School of Mathematical Sciences\\
Zhejiang Normal University\\
Jinhua 321004\\
P.R. China}
\email{wang.eric@zjnu.edu.cn}

\begin{document}

\maketitle
\begin{abstract}
We classify edge-to-edge tilings of the sphere by congruent pentagons with the edge combination $a^4b$ and with rational angles in degree: they are a one-parameter family of symmetric $a^4b$-pentagonal subdivisions of the tetrahedron with $12$ tiles; a sequence of unique symmetric $a^4b$-pentagons admitting a symmetric $3$-layer earth map tiling by $4m$ tiles for any $m\ge4$, among which each odd $m$ case admits two standard flip modifications; and a unique non-symmetric and degenerate $a^4b$-pentagon admitting a non-symmetric $3$-layer earth map tiling and its standard flip modification with $20$ tiles. The full classification from this series and all induced non-edge-to-edge quadrilateral tilings from degenerate pentagons are summarized with their 3D pictures. 
\end{abstract}

\section{Introduction}

In an edge-to-edge tiling of the sphere by congruent pentagons, the pentagon admits five possible edge combinations \cite{wy1}: $a^{2}b^{2}c$, $a^{3}bc$, $a^{3}b^{2}$, $a^4b$, $a^{5}$. Tilings for $a^2b^2c$, $a^3bc$, $a^3b^2$, $a^5$, and for $a^4b$ with any irrational angle, were classified in \cite{wy1, wy2, wy3, slw}. The remarkable preprint \cite{cpy} of $174$  pages classified tilings for $a^4b$ to complete the full classification, while a related doctoral thesis of nearly a thousand pages provided more details. It is important to verify such complicated classifications by different methods. We will classify $a^4b$-tilings in a very different way: our recent preprint \cite{slw} classified such tilings with any irrational angle in degree (or simply called ``with general angles''), and this paper will classify such tilings with all angles being rational in degree (such pentagons will be simply called \textit{rational} hereafter), thereby completing the full classification.  Two papers together will not only be much shorter than \cite{cpy}, but also provide three new results: full geometric data for all prototiles, 3D pictures for all types of tilings, and counting the total number of different tilings for each prototile. As a byproduct, we also obtain many new non-edge-to-edge quadrilateral tilings from degenerate pentagons.
	\begin{figure}[htp]
		\centering
		\begin{tikzpicture}[>=latex,scale=1]
			
			%% a^4b
			
			\begin{scope}[xshift=3cm]
				
				\draw
				(-54:1) -- (18:1) -- (90:1) -- (162:1) -- (234:1);
				
				\draw[line width=1.5]
				(234:1) -- (-54:1);

				\node at (90:0.75) {\small $\alpha$};
				\node at (162:0.75) {\small $\beta$};
				\node at (15:0.75) {\small $\gamma$};
				\node at (234:0.75) {\small $\delta$};
				\node at (-54:0.75) {\small $\epsilon$};
				
			\end{scope}		
			
			\begin{scope}[xshift=6cm]
				
				\fill[gray!50]  (-54:1) -- (18:1) -- (90:1) -- (162:1) -- (234:1);
				
				\draw
				(-54:1) -- (18:1) -- (90:1) -- (162:1) -- (234:1);
				
				\draw[line width=1.5]
				(234:1) -- (-54:1);

				\node at (90:0.75) {\small $\alpha$};
				\node at (162:0.75) {\small $\ccc$};
				\node at (15:0.75) {\small $\bbb$};
				\node at (234:0.75) {\small $\eee$};
				\node at (-54:0.75) {\small $\ddd$};
				
			\end{scope}	
			
			\begin{scope}[xshift=8cm]
				
				\draw
				(0,0.4) -- node[above=-2] {\small $a$} ++(1,0);
				
				\draw[line width=1.5]
				(0,-0.1) -- node[above=-2] {\small $b$} ++(1,0);
				
			\end{scope}
		\end{tikzpicture}
		\caption{Pentagons with the edge combinations $a^4b$.}
		\label{pentagon1}
	\end{figure}
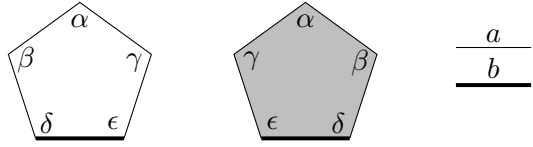
	
\begin{figure}[htp]
		\centering
		\includegraphics[scale=0.045]{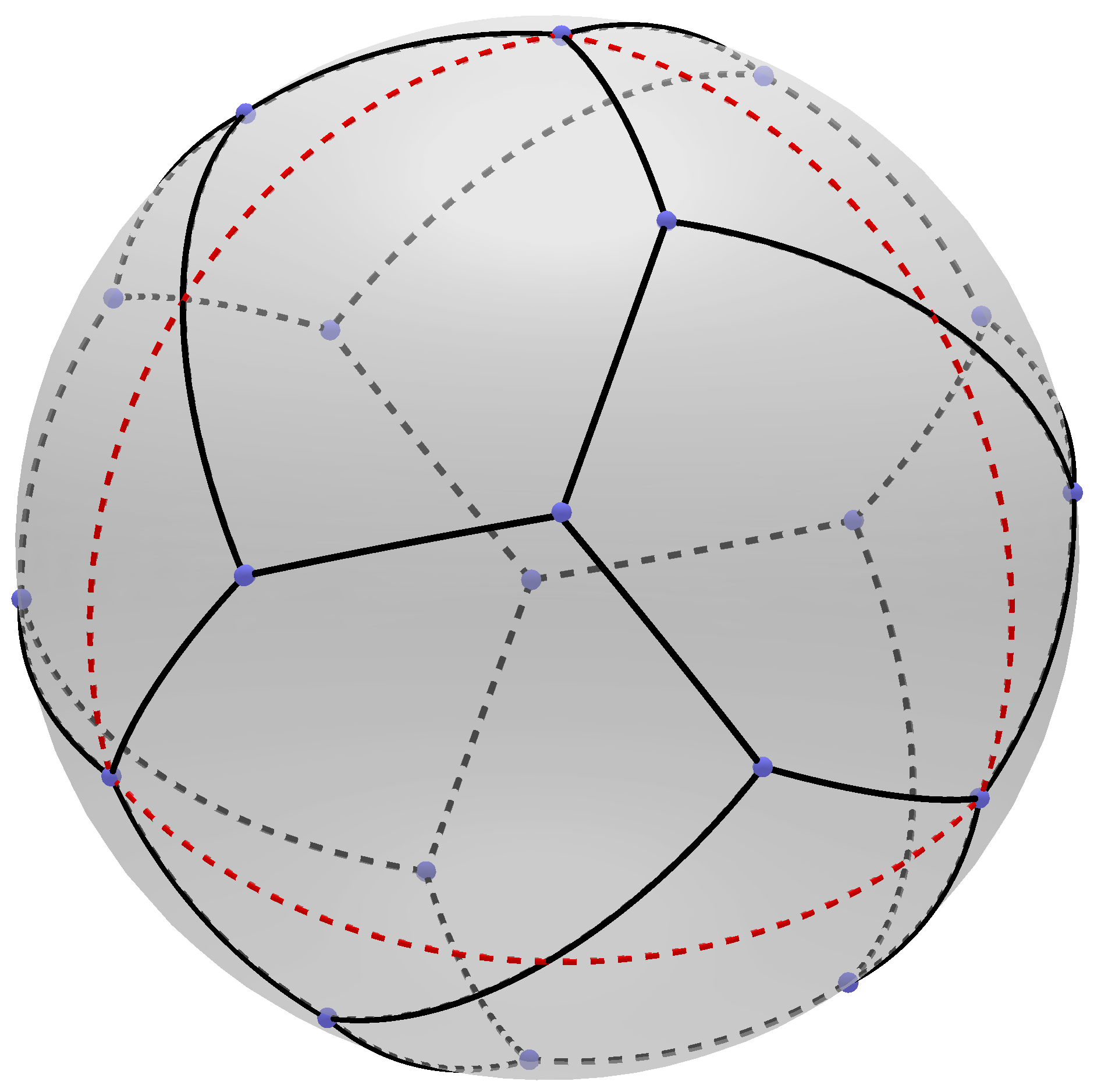}\hspace{20pt}
		 \includegraphics[scale=0.107]{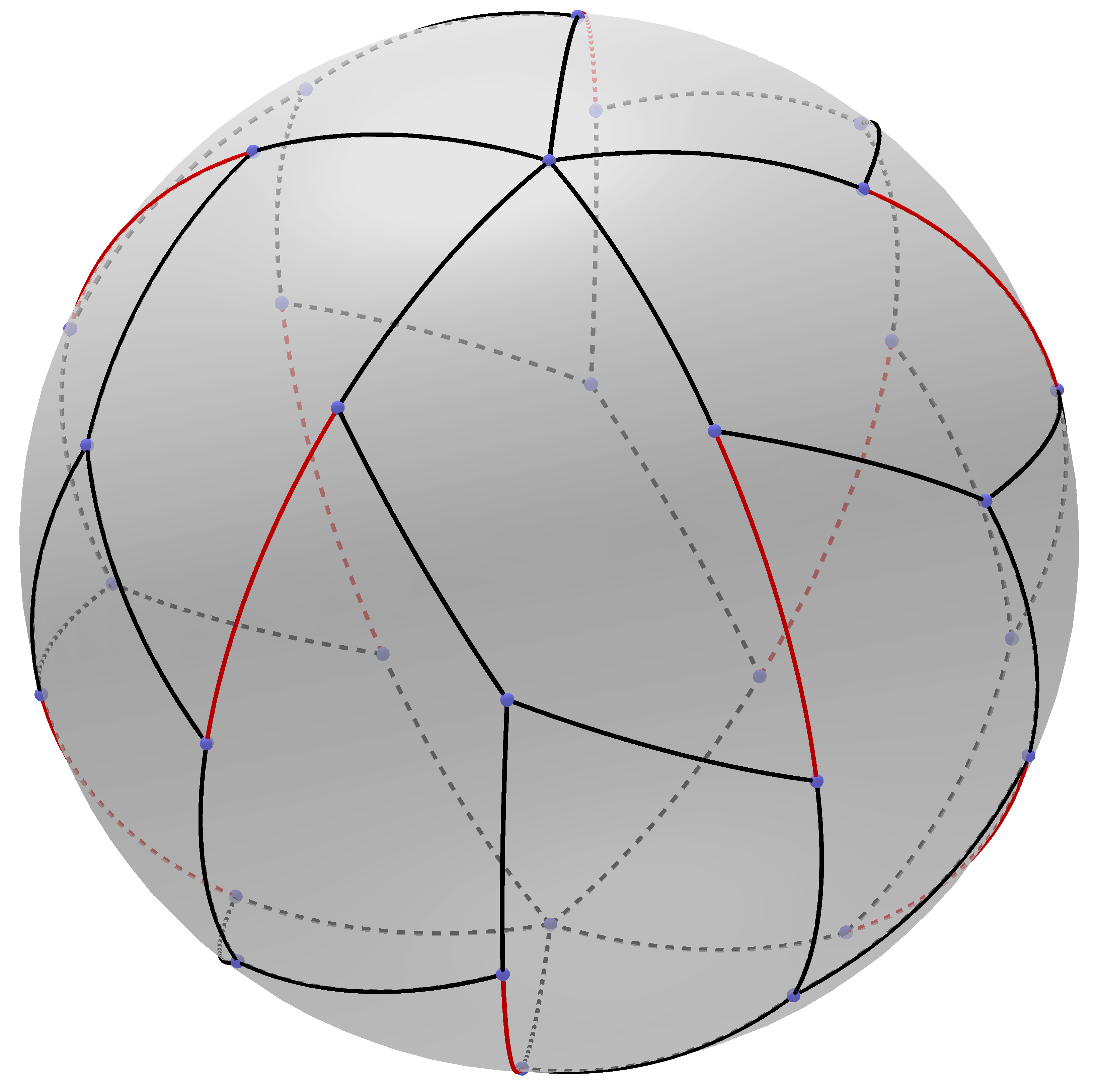}\hspace{20pt}
                 \includegraphics[scale=0.127]{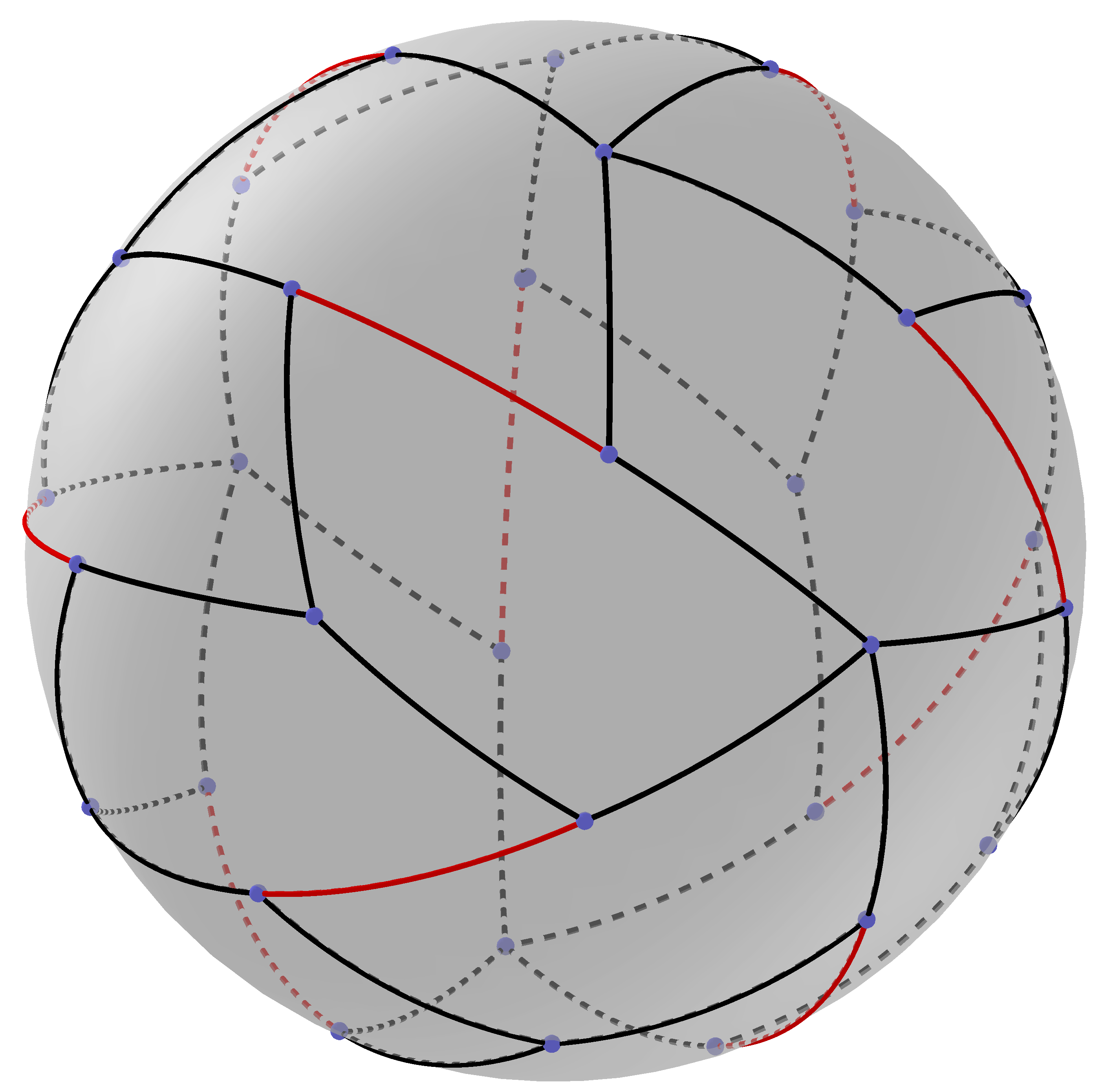}
  
		\caption{Pentagonal subdivision of the tetrahedron and two tilings of a non-symmetric $a^4b$-pentagon.} 
		\label{poufen}	
	\end{figure}

\begin{figure}[htp]
  \centering
  \begin{overpic}[scale=0.19]{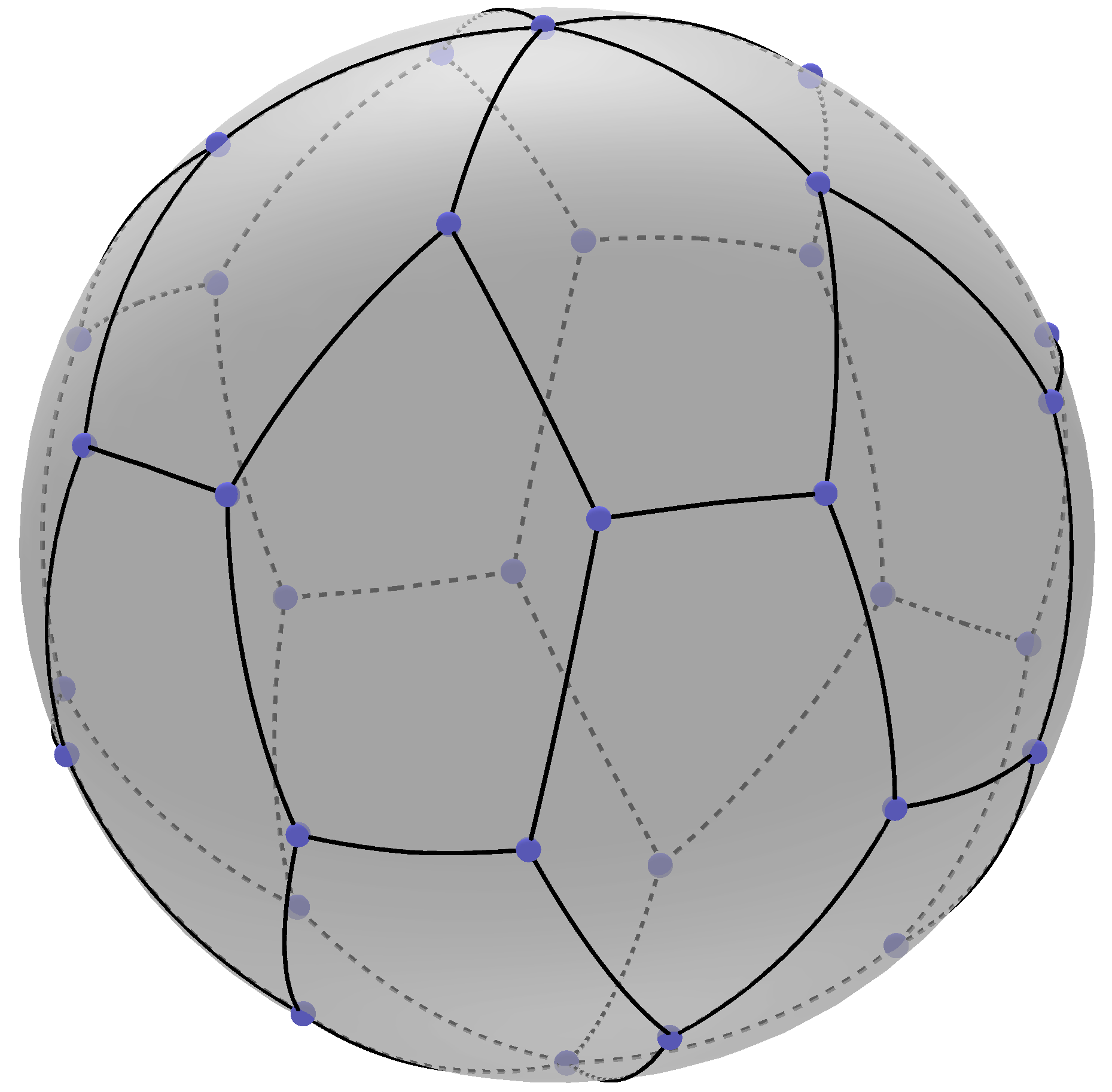}
   \put(90,5){\small $f=20$}
  \end{overpic}   \hspace{20pt}
  \includegraphics[scale=0.192]{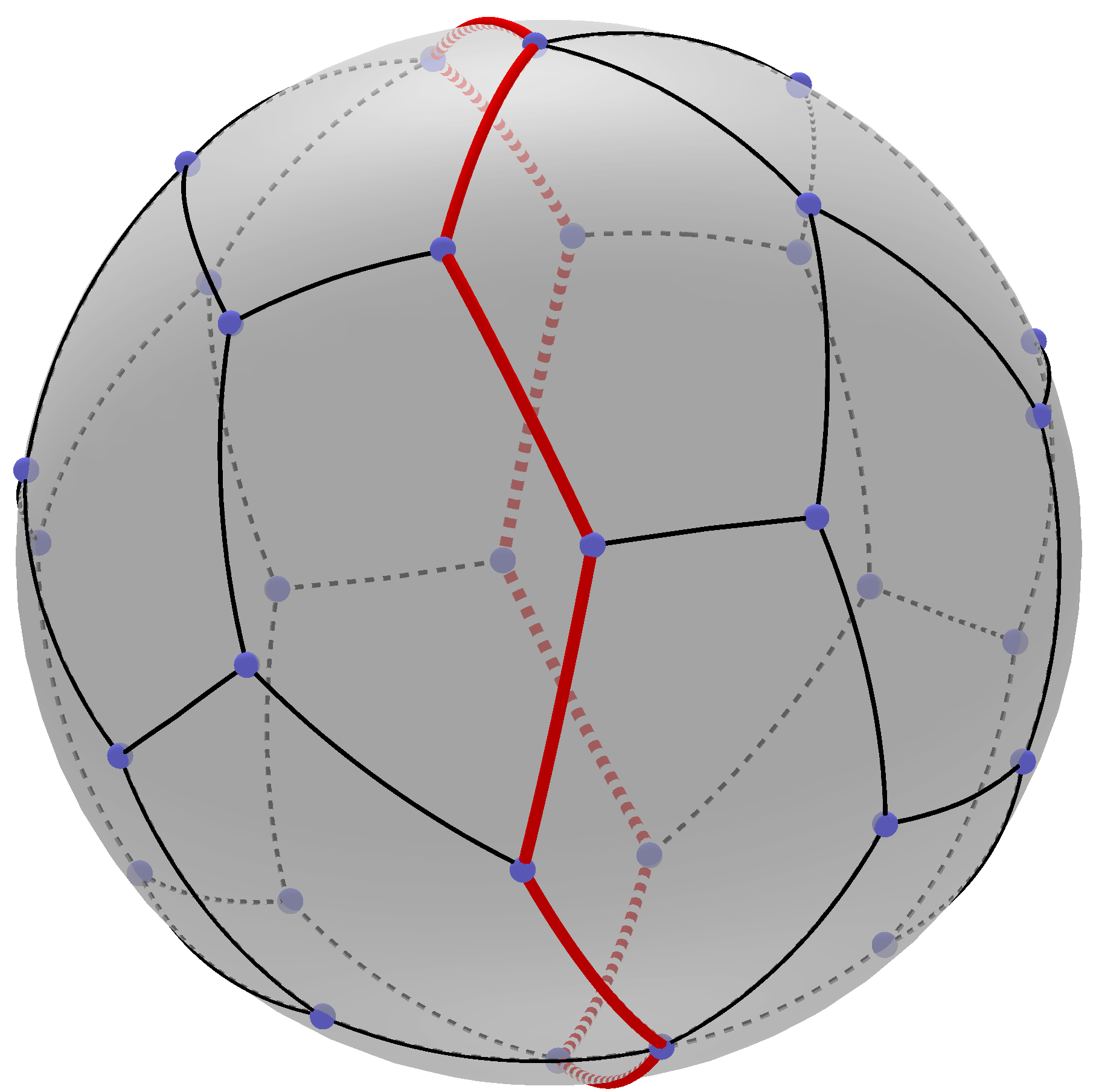}\hspace{20pt}
  \includegraphics[scale=0.21]{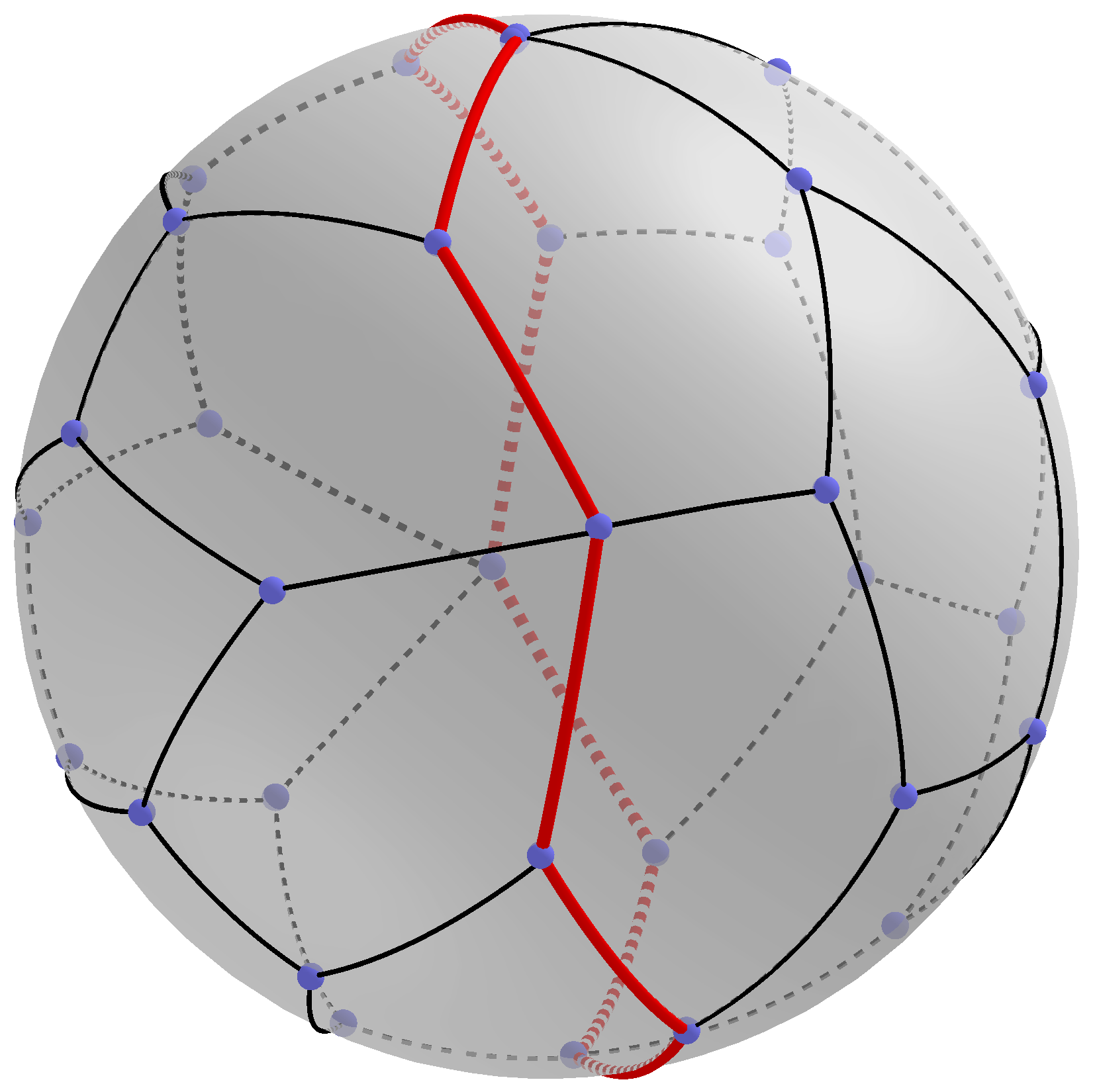}
  \caption{A symmetric $3$-layer earth map tiling and its two standard flip modifications.} 
  \label{sym3} 
 \end{figure}

Figure \ref{pentagon1} shows an $a^{4}b$-pentagon with edges $a$ (normal), $b$ (thick), and angles $\alpha$--$\epsilon$. The second image is its mirror (flip). The angles determine the orientation; conversely, the edge lengths and orientation determine angles uniquely. This allows tiling representation via shading instead of angle labels.
	\begin{theorem*}
		All  $a^4b$-tilings  with  rational angles are the following:
		\begin{enumerate}
			\item A one-parameter family of symmetric $a^4b$-pentagonal subdivisions of the tetrahedron with $12$ tiles (see the first picture of Figure \ref{poufen}).		
			\item A sequence of unique symmetric $a^4b$-pentagons with five angles $(\tfrac8f,\,1-\tfrac4f,\,1-\tfrac4f,\,\tfrac12+\tfrac2f,\,\tfrac12+\tfrac2f)\pi$ and		
			\begin{equation*}
				\cos a=1-2\left(\frac{\sqrt{5}-1}{4\cos{\frac{4\pi}f}}\right)^2,\,\,\,\, 
				\cos b=2 \left(\frac{(3-\sqrt5)\cos^2{\frac{4\pi}f}+\sqrt5-2}{\cos{\frac{4\pi}f}}\right)^2 -1,
			\end{equation*}
each admitting a symmetric $3$-layer earth map tiling $T(2m\,\aaa\bbb\ccc,\\
2m\,\bbb\ddd\eee,
			2m\,\ccc\ddd\eee,2\aaa^{m})$ by $f=4m$ tiles for any $m\ge4$, among which each odd $m=2k+1$ case admits two standard flip modifications (see Figure \ref{sym3}):
			\begin{itemize}
    \item $T(4k\,\aaa\bbb\ccc,\,(4k+2)\bbb\ddd\eee,\,(4k+2)\ccc\ddd\eee,\,2\aaa^{k+1}\bbb,\,2\aaa^{k+1}\ccc)$.
    \item $T((4k+2)\aaa\bbb\ccc,\,(4k+1)\bbb\ddd\eee,\,(4k+1)\ccc\ddd\eee,\,\aaa^{k+1}\bbb,\,\aaa^{k+1}\ccc,\,2\aaa^k\ddd\eee)$.
\end{itemize}

\item A unique non-symmetric $a^4b$-pentagon with angles $\left(10,12,6,5,15\right)\\
\pi/15$, where
\begin{equation*}
    \cos a = \frac{\sqrt{6}\,(5+\sqrt{5})^{3/2}}{60},\quad
    \cos b = \frac{\sqrt{5}}{3},
\end{equation*}
admitting two tilings:
\begin{itemize}
    \item A non-symmetric 3-layer earth map tiling $T(20\aaa\ddd\eee,\,10\bbb^2\ccc,\\
    \,2\ccc^5)$ with $f=20$ (see the second picture of Figure \ref{poufen}).
    \item A unique flip modification $T(20\aaa\ddd\eee,\,8\bbb^2\ccc,\,4\bbb\ccc^3)$ (see the third picture of Figure \ref{poufen}).
\end{itemize}					
		\end{enumerate}

	\end{theorem*}

The first two classes (Theorem 1 of \cite{slw}) consist of symmetric pentagons divisible into congruent quadrilaterals, which were classified in \cite{lw1, lw2, lw3}. The third class contains the sole rational pentagon from the 1-parameter families admitting non-symmetric 3-layer earth map tilings (Theorem 2 of \cite{slw}), which however does not admit new tilings. This is in sharp contrast to the quadrilateral case \cite{lw2}, where there are an abundance of rational prototiles  admitting new tilings.

\subsection*{Outline of the paper}
Section \ref{basic_facts} presents tiling fundamentals and classification techniques. Section \ref{region} introduces the method of solving algebraic equations in roots of unity with three examples. Section \ref{cases} applies this method to all cases of degree 3 $a^2b$-vertex types, completing the proof. The final section summarizes: (1) full classification from this five-paper series, and (2) all induced non-edge-to-edge quadrilateral tilings from degenerate pentagons, with 3D pictures.

\subsection*{Acknowledgements:}
Corresponding author: Erxiao Wang.  His research is partially supported by National Natural Science Foundation of China NSFC-RGC 12361161603 and Key Projects of Zhejiang Natural Science Foundation LZ22A010003.

\section{Basic Facts}\label{basic_facts}

For the basics of tilings we refer to the classic books \cite{gs, Adams}.

Let $v,e,f$ be the numbers of vertices, edges, and tiles respectively in a pentagonal tiling. Let $v_k$ represent the number of vertices of degree $k$. The following formulas are presented in \cite{wy1}.
\begin{align}
	f
	&=12+2\sum_{k\ge 4}(k-3)v_k=12+2v_4+4v_5+6v_6+\cdots, \label{vcountf} \\
	v_3
	&=20+\sum_{k\ge 4}(3k-10)v_k=20+2v_4+5v_5+8v_6+\cdots. \label{vcountv}
\end{align}
Thus, we have  $f\ge 12$, $v_3\ge 20$, and the majority of vertices have degree $3$.

\begin{lemma}[{\cite[Lemma 4]{wy1}}] \label{anglesum} 
	If all tiles in a tiling of the sphere by $f$ pentagons have the same five angles $\aaa,\bbb,\ccc,\ddd,\eee$, then 
	\[
	\aaa+\bbb+\ccc+\ddd+\eee = 3+\tfrac{4}{f} , 
	\]  
	ranging in $(3,\frac{10}3]$. In particular no vertex contains all five angles.
\end{lemma}

\begin{lemma}[{\cite[Lemma 1]{wy1}}] \label{base_tile}
	Any pentagonal tiling of the sphere has a ``special" tile with four vertices of degree $3$ and with the fifth vertex of degree $3$, $4$ or $5$.
\end{lemma}

\begin{lemma}[Parity Lemma, {\cite[Lemma 10]{wy2}}]\label{parity}
	In an $a^4b$-tiling, the total number of $ab$-angles (i.e. $\delta$ and $\epsilon$) at any vertex is even.    
\end{lemma}

In a tiling of the sphere by $f$ congruent tiles, each angle of the tile appears $f$ times in total. If one vertex has more $\aaa$ than $\bbb$, there must exist another vertex with more $\bbb$ than $\aaa$.
Such global counting induces many interesting and useful results. 

\begin{lemma}[Balance Lemma, {\cite[Lemma 1]{wy2}}]\label{balance}
	In an $a^4b$-tiling, if either $\ddd^2\cdots$ or $\eee^2\cdots$ is not a vertex, then any vertex either has no $\ddd,\eee$, or is of the form $\ddd\eee\cdots$ with no more $\ddd,\eee$ in the remainder.
\end{lemma}

We call a vertex with three $a$-edges an $a^{3}$-vertex, and a vertex with two $a$-edges and one $b$-edge an $a^{2}b$-vertex. 

\begin{lemma}\label{l2}
	In an $a^4b$-tiling, an $a^3$-vertex and an $a^2b$-vertex must both appear.
\end{lemma}

\begin{proof}
	By Lemma \ref{base_tile}, we know that either $\ddd\cdots$ or $\eee\cdots$ is a $a^2b$-vertex. If there is no $a^3$-vertex, then every degree 3 vertex is a $a^2b$-vertex. By the equations (\ref{vcountf}) and (\ref{vcountv}), we have $\#\ddd+\#\eee\ge2v_3>2f$ in the tiling, contradicting $\#\ddd+\#\eee=2f$.
\end{proof}

\begin{figure}[htp]
	\centering
	\begin{tikzpicture}[>=latex,scale=1] 
		
		\begin{scope}[xshift=-4cm]
			
			\draw
			(0,0) -- (0,0.8) -- (-0.6,1.2) -- (-1.2,0.6) -- (-0.8,0) -- (0,0);
			\draw
			(0,0) -- (0,0.8) -- (0.6,1.2) -- (1.2,0.6) -- (0.8,0) -- (0,0);
			\draw
			(0,0) -- (0,-0.8) -- (0.6,-1.2) -- (1.2,-0.6) -- (0.8,0) -- (0,0);
			\draw
			(0,0) -- (0,-0.8) -- (-0.6,-1.2) -- (-1.2,-0.6) -- (-0.8,0) -- (0,0);

			\draw[line width=1.5]
			(-1.2,0.6) -- (-0.8,0);
			\draw[line width=1.5]
			(0,0) -- (0.8,0);
			\draw[line width=1.5]
			(-1.2,-0.6) -- (-0.8,0);
			
			\node at (-0.7,0.2) {\small $\delta$}; 
			\node at (-1,0.6) {\small $\epsilon$};
			\node at (-0.6,0.95) {\small $\gamma$};
			\node at (-0.2,0.7) {\small $\alpha$};
			\node at (-0.2,0.2) {\small $\beta$};
			
			\node at (0.7,0.2) {\small $\eee$}; 
			\node at (1,0.6) {\small $\ccc$};
			\node at (0.6,0.95) {\small $\aaa$};
			\node at (0.2,0.7) {\small $\bbb$};
			\node at (0.2,0.2) {\small $\ddd$};
			
			\node at (0.7,-0.2) {\small $\eee$}; 
			\node at (1,-0.6) {\small $\ccc$};
			\node at (0.6,-0.95) {\small $\aaa$};
			\node at (0.2,-0.7) {\small $\bbb$};
			\node at (0.2,-0.2) {\small $\ddd$};
			
			\node at (-0.7,-0.2) {\small $\delta$}; 
			\node at (-1,-0.6) {\small $\epsilon$};
			\node at (-0.6,-0.95) {\small $\gamma$};
			\node at (-0.2,-0.7) {\small $\alpha$};
			\node at (-0.2,-0.2) {\small $\beta$};
			
		\end{scope}
		
		\begin{scope}[xshift=0cm]
			
			\draw
			(0,0) -- (0,0.8) -- (-0.6,1.2) -- (-1.2,0.6) -- (-0.8,0) -- (0,0);
			\draw
			(0,0) -- (0,0.8) -- (0.6,1.2) -- (1.2,0.6) -- (0.8,0) -- (0,0);
			\draw
			(0,0) -- (0,-0.8) -- (0.6,-1.2) -- (1.2,-0.6) -- (0.8,0) -- (0,0);
			\draw
			(0,0) -- (0,-0.8) -- (-0.6,-1.2) -- (-1.2,-0.6) -- (-0.8,0) -- (0,0);

			\draw[line width=1.5]
			(-1.2,0.6) -- (-0.8,0);
			\draw[line width=1.5]
			(0,0) -- (0.8,0);
			\draw[line width=1.5]
			(0,-0.8) -- (-0.6,-1.2);
			
			\node at (-0.7,0.2) {\small $\ddd$}; 
			\node at (-1,0.6) {\small $\eee$};
			\node at (-0.6,0.95) {\small $\ccc$};
			\node at (-0.2,0.7) {\small $\aaa$};
			\node at (-0.2,0.2) {\small $\beta$};
			
			\node at (0.7,0.2) {\small $\eee$}; 
			\node at (1,0.6) {\small $\ccc$};
			\node at (0.6,0.95) {\small $\aaa$};
			\node at (0.2,0.7) {\small $\bbb$};
			\node at (0.2,0.2) {\small $\ddd$};
			
			\node at (0.7,-0.2) {\small $\eee$}; 
			\node at (1,-0.6) {\small $\ccc$};
			\node at (0.6,-0.95) {\small $\aaa$};
			\node at (0.2,-0.7) {\small $\bbb$};
			\node at (0.2,-0.2) {\small $\ddd$};
			
			\node at (-0.7,-0.2) {\small $\aaa$}; 
			\node at (-1,-0.6) {\small $\ccc$};
			\node at (-0.6,-0.95) {\small $\eee$};
			\node at (-0.2,-0.7) {\small $\ddd$};
			\node at (-0.2,-0.2) {\small $\beta$};
			
		\end{scope}
		
		\begin{scope}[xshift=4cm]
			
			\draw
			(0,0) -- (0,0.8) -- (-0.6,1.2) -- (-1.2,0.6) -- (-0.8,0) -- (0,0);
			\draw
			(0,0) -- (0,0.8) -- (0.6,1.2) -- (1.2,0.6) -- (0.8,0) -- (0,0);
			\draw
			(0,0) -- (0,-0.8) -- (0.6,-1.2) -- (1.2,-0.6) -- (0.8,0) -- (0,0);
			\draw
			(0,0) -- (0,-0.8) -- (-0.6,-1.2) -- (-1.2,-0.6) -- (-0.8,0) -- (0,0);

			\draw[line width=1.5]
			(0,0.8) -- (-0.6,1.2);
			\draw[line width=1.5]
			(0,0) -- (0.8,0);
			\draw[line width=1.5]
			(0,-0.8) -- (-0.6,-1.2);
			
			\node at (-0.7,0.2) {\small $\aaa$}; 
			\node at (-1,0.6) {\small $\ccc$};
			\node at (-0.6,0.95) {\small $\eee$};
			\node at (-0.2,0.7) {\small $\ddd$};
			\node at (-0.2,0.2) {\small $\beta$};
			
			\node at (0.7,0.2) {\small $\eee$}; 
			\node at (1,0.6) {\small $\ccc$};
			\node at (0.6,0.95) {\small $\aaa$};
			\node at (0.2,0.7) {\small $\bbb$};
			\node at (0.2,0.2) {\small $\ddd$};
			
			\node at (0.7,-0.2) {\small $\eee$}; 
			\node at (1,-0.6) {\small $\ccc$};
			\node at (0.6,-0.95) {\small $\aaa$};
			\node at (0.2,-0.7) {\small $\bbb$};
			\node at (0.2,-0.2) {\small $\ddd$};
			
			\node at (-0.7,-0.2) {\small $\aaa$}; 
			\node at (-1,-0.6) {\small $\ccc$};
			\node at (-0.6,-0.95) {\small $\eee$};
			\node at (-0.2,-0.7) {\small $\ddd$};
			\node at (-0.2,-0.2) {\small $\beta$};
			
		\end{scope}
		
	\end{tikzpicture}
	\caption{Different Adjacent Angle Deductions of $\bbb^2\ddd^2$.} \label{anglecombo}
\end{figure}
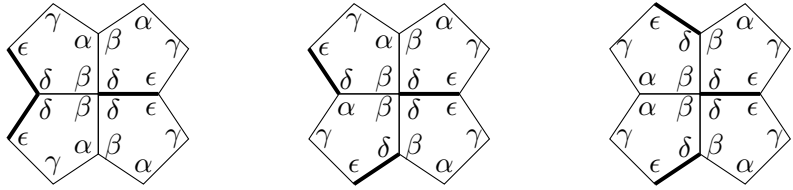

The very useful tool \textit{adjacent angle deduction} (abbreviated as \textbf{AAD}) has been introduced in \cite[Section 2.5]{wy1}. We give a quick review here using  Figure \ref{anglecombo}. Let ``$\thin$'' denote an a-edge and ``$\thick$'' denote a b-edge. Then we indicate the arrangements of angles and edges by denoting the vertices as $\thin\bbb\thin\bbb\thin\ddd\thick\ddd\thin$. The notation can be reversed, such as $\thin\bbb\thin\bbb\thin\ddd\thick\ddd\thin=\thin\ddd\thick\ddd\thin\bbb\thin\bbb\thin$; and it can be rotated, such as $\thin\bbb\thin\bbb\thin\ddd\thick\ddd\thin=\thin\bbb\thin\ddd\thick\ddd\thin\bbb\thin=\thick\ddd\thin\bbb\thin\bbb\thin\ddd\thick$. We also denote the first  vertex in Figure \ref{anglecombo} as $\ddd\thick\ddd\cdots$, $\thin\ddd\thick\ddd\thin\cdots,\bbb\thin\ddd\cdots$, $\thin\ddd\thick\ddd\thin\bbb\thin\cdots$, and denote the consecutive angle segments as $\ddd\thick\ddd$, $\thin\ddd\thick\ddd\thin,\bbb\thin\ddd$, $\thin\ddd\thick\ddd\thin\bbb\thin$.

The pictures of Figure \ref{anglecombo} have the same vertex $\thin\bbb\thin\bbb\thin\ddd\thick\ddd\thin$, but different arrangements of the four tiles. To indicate the difference, we write ${}^{\lambda}\theta^{\mu}$ to mean $\lambda,\mu$ are the two angles adjacent to $\theta$ in the pentagon. The first picture has the AAD $\thin^{\aaa}\bbb^{\ddd}\thin^{\ddd}\bbb^{\aaa}\thin^{\bbb}\ddd^{\eee}\thick^{\eee}\ddd^{\bbb}\thin$. The second and third have the AAD $\thin^{\aaa}\bbb^{\ddd}\thin^{\aaa}\bbb^{\ddd}\thin^{\bbb}\ddd^{\eee}\thick^{\eee}\ddd^{\bbb}\thin$ and $\thin^{\ddd}\bbb^{\aaa}\thin^{\aaa}\bbb^{\ddd}\thin^{\bbb}\ddd^{\eee}\thick^{\eee}\ddd^{\bbb}\thin$ respectively.  

\begin{lemma}[{\cite[Lemma 2]{wy2}}] \label{geometry1}
	In an $a^4b$-tiling, $\beta>\gamma$ is equivalent to $\delta<\epsilon$, and $\beta<\gamma$ is equivalent to $\delta>\epsilon$. 
\end{lemma}
	
\begin{lemma}[{\cite[Lemma 18]{cpy}}]\label{formula}
	If the angles $\alpha,\beta,\gamma,\delta,\epsilon$ and the edge length $a$ of an almost equilateral pentagon satisfy $a\not\in {\bb Z}\pi$, then
	\begin{align}
		&
		[
		((1-\cos\beta)\sin(\delta-\tfrac{1}{2}\alpha)
		-(1-\cos\gamma)\sin(\epsilon-\tfrac{1}{2}\alpha))\sin\tfrac{1}{2}(\delta-\epsilon) \nonumber \\
		&-(1-\cos(\beta-\gamma))\sin\tfrac{1}{2}\alpha\sin\tfrac{1}{2}(\delta+\epsilon)]\cos\tfrac{1}{2}(\delta+\epsilon-\alpha)=0, \label{coolsaet_eq1} \\
		& \sin\tfrac{1}{2}\alpha\sin\tfrac{1}{2}(\beta-\gamma)
		(\sin\tfrac{1}{2}\beta\sin\delta\cos a
		-\cos\tfrac{1}{2}\beta\cos\delta) \nonumber \\
		& +\sin\tfrac{1}{2}\gamma
		\sin\tfrac{1}{2}(\delta-\epsilon)
		\cos\tfrac{1}{2}(\delta+\epsilon-\alpha)=0, \label{coolsaet_eq2} \\
		& \sin\tfrac{1}{2}\alpha\sin\tfrac{1}{2}(\beta-\gamma)
		(\sin\tfrac{1}{2}\gamma\sin\epsilon\cos a
		-\cos\tfrac{1}{2}\gamma\cos\epsilon) \nonumber \\
		& +\sin\tfrac{1}{2}\beta
		\sin\tfrac{1}{2}(\delta-\epsilon)
		\cos\tfrac{1}{2}(\delta+\epsilon-\alpha)=0. \label{coolsaet_eq3}
	\end{align}
	Conversely, if $\alpha,\beta,\gamma,\delta,\epsilon,a$ satisfy the three equalities, and $\alpha,\beta-\gamma\not\in 2{\bb Z}\pi$, and one of $\delta,\epsilon\not\in {\bb Z}\pi$, then there is an almost equilateral pentagon with the given $\alpha,\beta,\gamma,\delta,\epsilon,a$.
\end{lemma}

The pentagon constructed in Lemma \ref{formula} may be self-intersecting. Specifically, the first two solutions in Table \ref{ade 2bc} yield pentagons via Lemma \ref{formula} that violate Lemma \ref{geometry1}, and thus are not simple.

The equation \eqref{coolsaet_eq1} only involves five angles and can be simplified further for any pentagon in $a^4b$-tilings as follows. 
\begin{lemma}
	In an $a^4b$-tiling, the following identity holds:
	\begin{align}
		&
		[(1-\cos\beta)\sin(\delta-\tfrac{1}{2}\alpha)
		-(1-\cos\gamma)\sin(\epsilon-\tfrac{1}{2}\alpha)] \, \sin\tfrac{1}{2}(\delta-\epsilon) \nonumber \\
		&-[1-\cos(\beta-\gamma)] \, \sin\tfrac{1}{2}\alpha \, \sin\tfrac{1}{2}(\delta+\epsilon)=0.\label{coolsaet_eq6} 
	\end{align}

\end{lemma}

\begin{proof}
     	In an $a^4b$-tiling, we have $0<a < \pi$ and then \eqref{coolsaet_eq1} holds. This implies \eqref{coolsaet_eq6} if $\cos\tfrac{1}{2}(\delta+\epsilon-\alpha)\neq 0$. 
     	
     	If $\cos\tfrac{1}{2}(\delta+\epsilon-\alpha)=0$, then $\ddd+\eee=\aaa \pm \pi$ by $\ddd+\eee<2\pi$. (If $\ddd+\eee\ge2\pi$, then w.l.o.g. let $\ddd\ge\pi$. There is no vertex $\ddd^2\cdots$ nor $\ddd\eee\cdots$, contradicting the balance lemma \ref{balance}.) Substituting $\ddd+\eee=\aaa \pm \pi$ into \eqref{coolsaet_eq2} and \eqref{coolsaet_eq3}, we derive the following equation:
\begin{equation}\label{ccoolsaet_eq1}
\cos\tfrac{1}{2}\beta\cos\delta\sin\tfrac{1}{2}\gamma\sin\epsilon
-\sin\tfrac{1}{2}\beta\sin\delta\cos\tfrac{1}{2}\gamma\cos\epsilon=0.
\end{equation}    
Then $\mathrm{LHS}\ \text{of}\ \eqref{coolsaet_eq6} = \pm 2\sin\left(\tfrac{1}{2}(\bbb - \ccc)\right) \cdot \mathrm{LHS}\ \text{of}\ \eqref{ccoolsaet_eq1} = 0$.
\end{proof}

\begin{figure}[htp]
		\centering
		\includegraphics[scale=0.158]{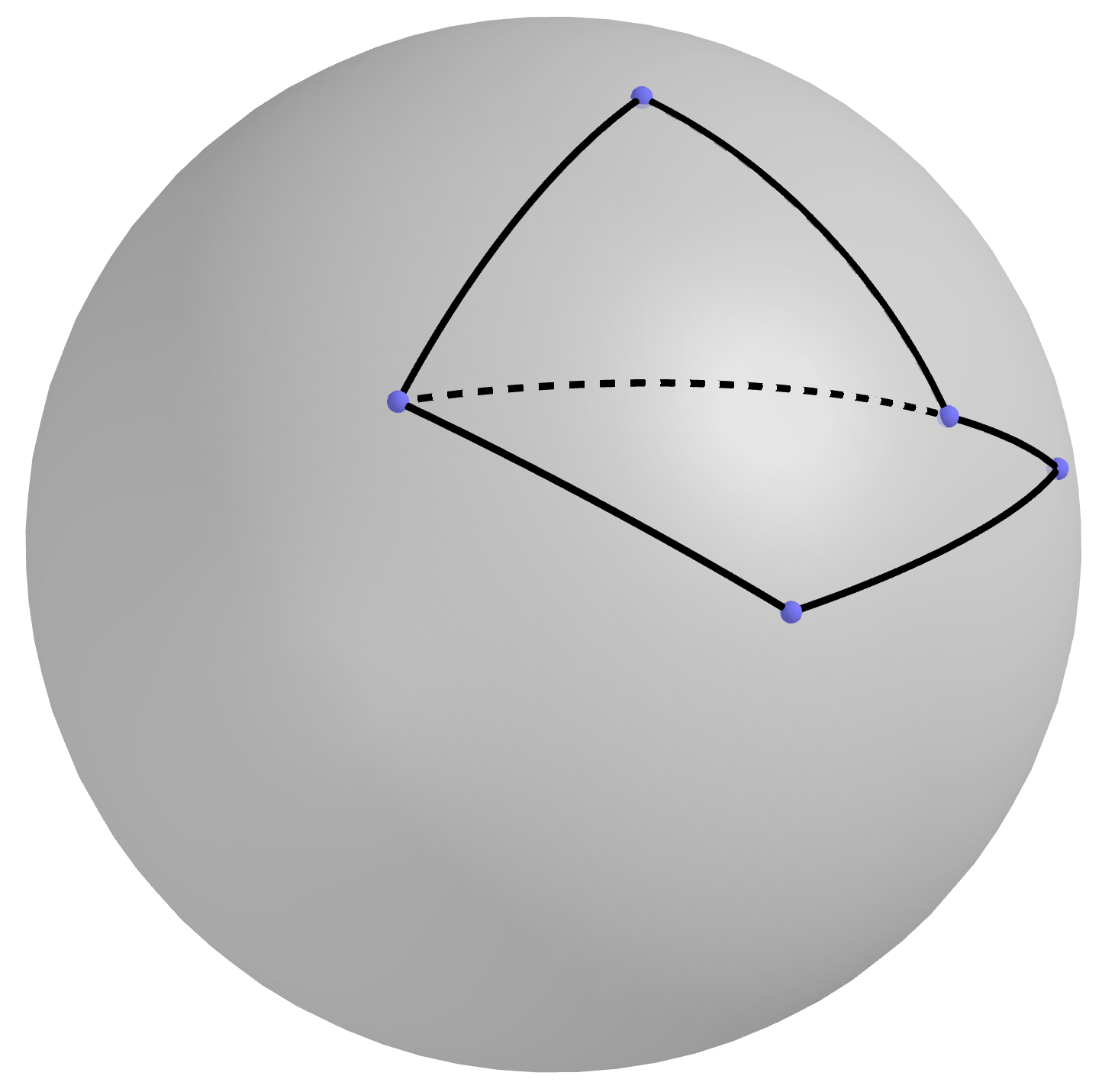}\hspace{60pt}
		 \includegraphics[scale=0.200]{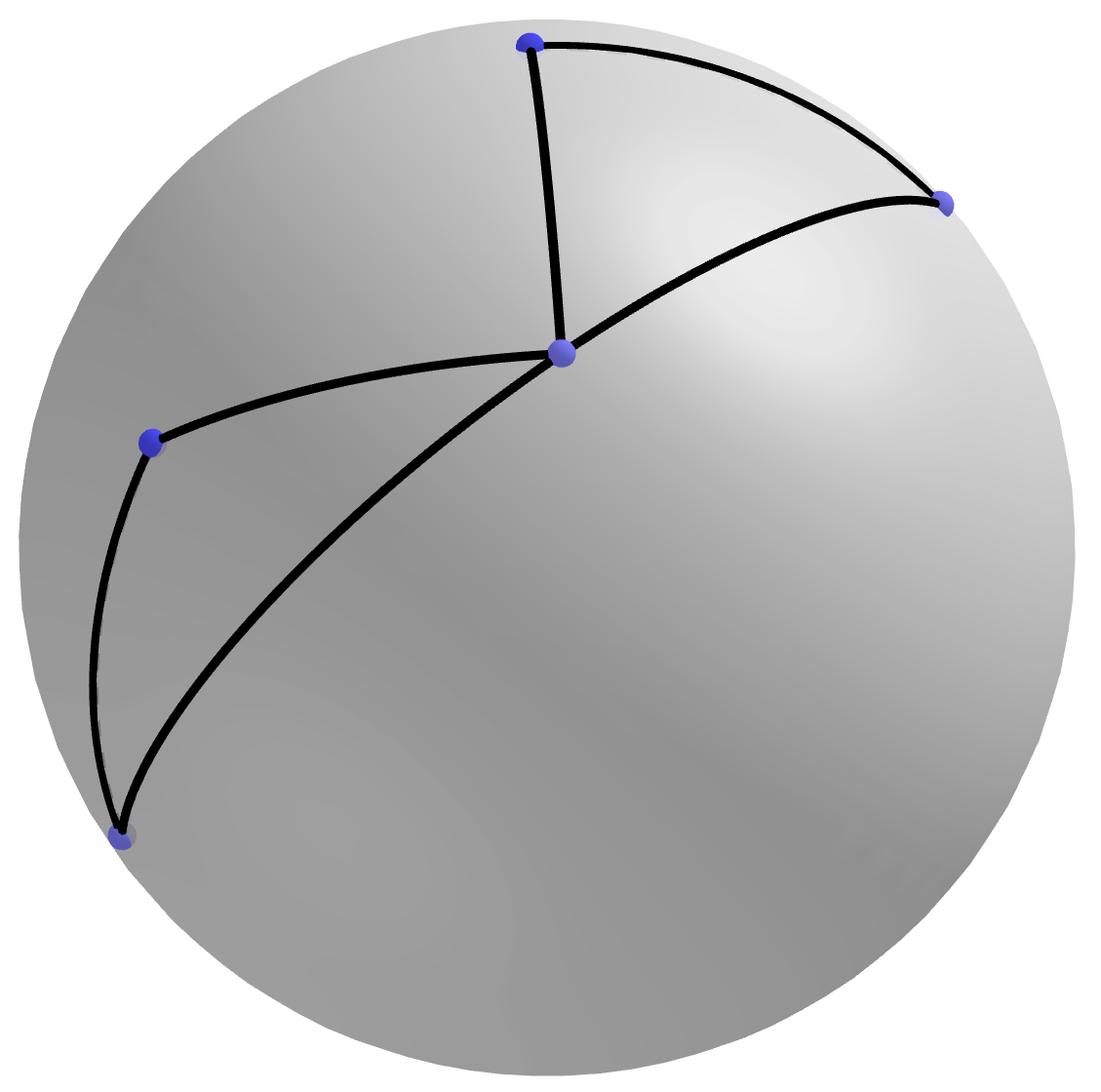}\hspace{20pt}

		\caption{ The geometric interpretation of $\ddd+\eee=\aaa \pm \pi$.} 
		\label{equation}	
	\end{figure}

The proof's last line is verifiable symbolically (e.g., Maple). GeoGebra 3D pictures reveal the geometric meaning of $\ddd+\eee=\aaa \pm \pi$: The $b$-edge and $\alpha$-vertex lie on a great circle (Figure \ref{equation}).

\medskip

Equation \eqref{coolsaet_eq6} is pivotal, and remarkably, we can find all  rational angle solutions to this trigonometric Diophantine equation via several algorithms, as detailed in the next section.

\section{Trigonometric Diophantine equations}
\label{region}
Using $x = e^{i\theta}$ for each angle variable, we convert trigonometric equations like \eqref{coolsaet_eq6} to algebraic form via $\cos \theta = \frac{1}{2}(x + x^{-1})$. Solving for rational angles is then equivalent to finding roots of unity solutions. Although \cite[\S 1.3]{kk} provides two algorithms for this, both have super-exponential complexity—limiting practicality to either $\leq 12$ monomials or $\leq 3
$ variables. Fortunately, linear angle constraints from vertices around special tiles reduce \eqref{coolsaet_eq6} to $2-3$ variables.

Let $K = \mathbb{Q}^{\text{ab}}$ be the maximal abelian extension of $\mathbb{Q}$ (containing all roots of unity). For cyclotomic field fundamentals, see \cite[Chapter VI]{Lang}. Our objective requires finding cyclotomic points (all coordinates roots of unity) on either:
\begin{itemize}
    \item An affine plane curve $V(f) \subset \mathbb{C}^2$ where $f \in K[x,y]$, or
    \item An affine surface $V(g) \subset \mathbb{C}^3$ where $g \in K[x,y,z]$.
\end{itemize}
We illustrate this by searching for rational $a^4b$-pentagons, first reviewing the univariate case.

\subsection{The cyclotomic roots of a one-variable polynomial} \label{subsection1.1}

Finding all roots of unity $\omega$ satisfying $f(\omega) = 0$ for $f \in \mathbb{Q}[x]$ is equivalent to finding all irreducible cyclotomic polynomial factors of $f$. Bradford and Davenport \cite{Bradford-Davenport} developed an efficient algorithm improving trial division, later extended to: 
\begin{itemize}
    \item Bivariate complex coefficients (Beukers--Smyth \cite{BS})
    \item Multivariate cases (Aliev--Smyth \cite{AS})
\end{itemize}
The method exploits a key property: For $\omega$ a root of unity of order $k \equiv r \pmod{4}$ ($r=0,1,2,3$), $\omega$ is conjugate to $-\omega$, $\omega^2$, $-\omega^2$, $\omega^2$ respectively. This yields a recursive algorithm computing $\gcd(f(x), f(-x))$, $\gcd(f(x), f(x^2))$, and $\gcd(f(x), f(-x^2))$.

In practice, we find it more convenient to first factor $f$ into irreducible parts  using symbolic software (e.g., SageMath, Maple). For fast factorization algorithms, see flintlib.org and libntl.org.

Let $g \in \mathbb{Z}[x]$ be irreducible. If $g(x) = g(-x)$ (equivalently, $g$ contains only even-degree monomials), set $f(x) = g(\sqrt{x})$. Iterate until obtaining $f(x) \neq f(-x)$. For this irreducible $f$:
\begin{itemize}
    \item If $\gcd\!\left(f(x),\, f(x^2)f(-x^2)\right) = 1$, then $f$ has no roots of unity zeros.
    \item If $\gcd\!\left(f(x),\, f(x^2)f(-x^2)\right) = f(x)$, then $f$ is cyclotomic and all its zeros are roots of unity.
\end{itemize}

For $f \in K[x]$, the coefficient field of $f$ is a subfield $\mathbb{Q}[\zeta_n]$ of $K$. Let $F(x) = \prod_{k \in (\mathbb{Z}/n\mathbb{Z})^\times} f_k(x)$, where $f_k$ is obtained by replacing $\zeta_n \mapsto \zeta_n^k$ in $f$'s coefficients. Then $F \in \mathbb{Q}[x]$, and the cyclotomic roots of $f$ are exactly the cyclotomic roots of $F$ that satisfy $f(\omega) = 0$. Thus, apply the algorithm to $F$ and filter solutions by evaluation in $f$.

For example, let $f(x) = x + \sqrt{2}$, then $f(x) = x + \zeta_8 + \zeta_8^{-1} \in \mathbb{Q}[\zeta_8][x]$. We have $\Z_8^\ast = \{1,3,5,7\} $ and $\zeta_8, \zeta_8^3, \zeta_8^5, \zeta_8^7$ are four primitive eighth roots of unity for the cyclotomic polynomial $\Phi_8(x) = x^4 + 1$. Then we construct $F \in \mathbb{Q}[x]$:
\[
(x + \zeta_8 + \zeta_8^{-1})(x + \zeta_8^3 + \zeta_8^{-3})(x + \zeta_8^5 + \zeta_8^{-5})(x + \zeta_8^7 + \zeta_8^{-7})=(x^2 - 2)^2.
\]

\subsection{The cyclotomic roots of a $2$-variable polynomial}\label{subsection1.2}

Let $f(x,y)$ $=$ $\sum_{i,j} a_{ij} x^i y^j \in \mathbb{C}[x,y]$ be a polynomial such that $f(0,0) \neq 0$. Cyclotomic points on $f(x,y) = 0$ (also called torsion points) can be found using the algorithm in  \cite{BS}.  We can also consider $f$ to be in $\mathbb{C}[x,x^{-1},y,y^{-1}]$, i.e., a Laurent polynomial.

Let $L(f)$ be the sublattice of $\mathbb{Z}^2$ generated by the set of differences $\{(i,j) - (i',j') \mid (i,j), (i',j') \text{ such that } a_{ij} \neq 0, a_{i'j'} \neq 0\}$.  If $f$ is not constant, the rank of this lattice is either 1 or 2. If $L(f) = \mathbb{Z}^2$, then $f$ is said to be \emph{full}.

\begin{enumerate}
    \item If $\operatorname{rank}(L(f)) = 1$, then $f$ is equivalent to a Laurent polynomial of the form $g(x^k y^l)$, where $k, l \in \mathbb{Z}$ and $g$ is a polynomial in one variable. Then any cyclotomic root of $g$ induces one rational parameter family of cyclotomic points  on $f(x,y) = 0$. 
    
    \item If $\operatorname{rank}(L(f)) = 2$, we can always make a full lattice $L(g)$ by the transformation
    \[
    \begin{cases}
        u = x^a y^c \\
        v = x^b y^d
    \end{cases}
    \]
    to the exponents such that $f(x,y) = g(u,v)$, where
    \[
    \begin{pmatrix}
        a & b \\
        c & d
    \end{pmatrix} \in \operatorname{GL}_2(\mathbb{Z}).
    \]
    %In practice, we do not necessarily need to proceed in this way.
\end{enumerate}

We now use two cases in the next section to illustrate the algorithm. 
%%%% ad2_bde_abc%%
\begin{example}\label{example1.2}
Case $\{\aaa\ddd^2,\bbb\ddd\eee,\aaa\bbb\ccc\}$ in the fifth row of Table \ref{two a2b}. 

By solving the linear system of angle sums, we get
 \[
 \aaa=2\pi-2\ddd, \, \bbb=\pi-\frac{4\pi}{f}, \, \ccc=2\ddd+\frac{4\pi}{f}-\pi,\, \eee=\pi-\ddd+\frac{4\pi}{f}, 
 \]
for $\ddd\in\left(\frac\pi2-\frac{2\pi}{f},\pi \right)$. In an $a^4b$-tiling by non-symmetric pentagons, two basic constraints are that $f\ge12$ must be an even integer (see the equation \eqref{vcountf}) and $ \bbb\neq \ccc$, $ \ddd\neq \eee$ (see Lemma \ref{geometry1}). The trigonometric equation \eqref{coolsaet_eq6}  reduces to 
\begin{align*}
	Q(\delta, f) &:= \sin \left(\delta +\frac{2 \pi}{f}\right)+\sin  \left(3 \delta -\frac{2 \pi}{f}\right)+2 \sin  \left(3 \delta +\frac{6 \pi}{f}\right)\\
	&-2 \sin  \left(\delta -\frac{6 \pi}{f}\right) + \sin \left(3 \delta +\frac{2 \pi}{f}\right) + \sin  \left(3 \delta -\frac{6 \pi}{f}\right)\\
	&-\sin  \left(\delta -\frac{2 \pi}{f}\right)+\sin  \left(3 \delta+\frac{10 \pi}{f} \right)=0	
\end{align*}
By setting $x={\mathrm e}^{2\mathrm{i} \delta},y={\mathrm e}^{\frac{4\mathrm{i} \pi}f}$,
the above trigonometric function $Q(\delta, f)$  is transformed to a polynomial
\begin{equation*}
P(x,y) := (y+1)(x^3y^4 + x^3y^3 + x^3y + x^2y^2 + 2xy^3 - 2x^2y - xy^2 - y^3 - y - 1).	
\end{equation*}

We consider each factor separately. The first factor gives $y=-1$ and $f=4$. This contradicts $f\ge 12$ in our pentagonal tilings.

Denote the second factor by $L(x,y)$, which is full. Then we compute the resultants of $L(x,y)$ with $L(-x,y)$, $L(x,-y)$, $L(-x,-y)$, $L(x^2,y^2)$, $L(-x^2,y^2)$, $L(x^2,-y^2)$, and $L(-x^2,-y^2)$ respectively. 
\medskip
      \begin{enumerate}
        \item $\mathrm{Resultant} (\mathit{L(x,y)} ,\mathit{L(-x,y)} ,x)$
 $ = 8 y^{3} (y^{3}+y^{2}+1) (y^{3}+y +1) (y^{6}+y^{5}+3 y^{4}-2 y^{3}+3 y^{2}+y +1)^{2}.$

It has no cyclotomic root by Bradford-Davenport's algorithm. 

\medskip
        \item $\mathrm{Resultant} (\mathit{L(x,y)} ,\mathit{L(x,-y)} ,x) $
 $ = -8 y^{3} (y^{2}+y +1) (y^{2}-y +1) (y^{8}+3 y^{6}+3 y^{2}+1) (y -1)^{3} (y +1)^{3}.$

The cyclotomic roots are $\zeta_6^k$ for $k = 0, 1, \dots, 5$.

        \item $ \mathrm{Resultant} (\mathit{L(x,y)} ,\mathit{L(-x,-y)} ,x)$
        $  = -8 y^{4} (8 y^{16}+11 y^{14}+11 y^{12}+y^{10}+2 y^{8}+y^{6}+11 y^{4}+11 y^{2}+8).$

It has no cyclotomic root. 

        \item  $\mathrm{Resultant} (\mathit{L(x,y)} ,\mathit{L(x^2,y^2)} ,x)$
         $ = -8 y^{7} (y^{2}-y +1) (8 y^{22}+7 y^{21}+10 y^{20}+10 y^{19}+30 y^{18}+8 y^{17}+42 y^{16}+38 y^{15}+10 y^{14}+65 y^{13}+60 y^{12}+60 y^{10}+65 y^{9}+10 y^{8}+38 y^{7}+42 y^{6}+8 y^{5}+30 y^{4}+10 y^{3}+10 y^{2}+7 y +8) (y^{2}+y +1)^{2} (y -1)^{3} (y +1)^{3}.$

The cyclotomic roots are $\zeta_6^k$ for $k = 0, 1, \dots, 5$.

        \item $\mathrm{Resultant} (\mathit{L(x,y)} ,\mathit{L(-x^2,y^2)} ,x) $
         $=-8 y^{6} (y^{36}+9 y^{35}+35 y^{34}+39 y^{33}+93 y^{32}+159 y^{31}+216 y^{30}+263 y^{29}+387 y^{28}+458 y^{27}+425 y^{26}+564 y^{25}+675 y^{24}+515 y^{23}+652 y^{22}+879 y^{21}+636 y^{20}+698 y^{19}+992 y^{18}+698 y^{17}+636 y^{16}+879 y^{15}+652 y^{14}+515 y^{13}+675 y^{12}+564 y^{11}+425 y^{10}+458 y^{9}+387 y^{8}+263 y^{7}+216 y^{6}+159 y^{5}+93 y^{4}+39 y^{3}+35 y^{2}+9 y +1).$

It has no cyclotomic root.

        \item   $\mathrm{Resultant} (\mathit{L(x,y)} ,\mathit{L(x^2,-y^2)} ,x)$

        $=8 y^{6} (y^{2}+y +1) (y -1)^{3} (y^{2}-y +1)^{8} (y +1)^{15}.$

The cyclotomic roots are $\zeta_6^k$ for $k = 0, 1, \dots, 5$.

        \item  $\mathrm{Resultant} (\mathit{L(x,y)} ,\mathit{L(-x^2,-y^2)} ,x)$
        $=8 y^{7} (8 y^{34}-3 y^{33}-27 y^{32}+23 y^{31}+64 y^{30}-58 y^{29}-45 y^{28}+160 y^{27}+99 y^{26}-133 y^{25}-85 y^{24}+242 y^{23}+77 y^{22}-378 y^{21}-64 y^{20}+343 y^{19}-27 y^{18}-328 y^{17}-27 y^{16}+343 y^{15}-64 y^{14}-378 y^{13}+77 y^{12}+242 y^{11}-85 y^{10}-133 y^{9}+99 y^{8}+160 y^{7}-45 y^{6}-58 y^{5}+64 y^{4}+23 y^{3}-27 y^{2}-3 y +8)$.

It has no cyclotomic root.
      \end{enumerate}

In summary, the cyclotomic roots for $y$ are $\zeta_6^k$ for $k = 0, 1, \dots, 5$. The only solution  suitable for an $a^4b$-tiling is $y = e^{\frac{\mathrm{i}\pi}{3}}$ and $f=12$. Solving $L(x,e^{\frac{\mathrm{i}\pi}{3}})=0$ similarly, we obtain $x=\zeta_{3}^k$ for $k=1,2$. The only suitable solution is $x=e^{\frac{4\mathrm{i}\pi}{3}}$ and $\ddd =\frac{2\pi}{3}$. This implies that all angles are $\frac{2\pi}{3}$, contradicting $\bbb\neq \ccc$. Hence we put `None' for this case in the Appendix, saying that this case admits no rational $a^4b$-pentagon satisfying the basic constraints. 

\end{example}
\medskip
We now present an example of a polynomial with coefficients in $\mathbb{Q}^{ab}$.
%%ad2_bd2_c5
\begin{example}\label{example1.3}
Case $\{\aaa\ddd^2,\bbb\ddd^2,\ccc^5\}$ in the first row of Table \ref{two a2b}. 
	
	By solving the linear system of angle sums, we get
	\[
	\aaa=2\pi-2\ddd, \quad \bbb=2\pi-2\ddd, \quad \ccc=\frac{2\pi}{5}, \quad \eee=-\frac{7\pi}5+3\ddd+\frac{4\pi}{f}, 
	\]		
for $\ddd\in\left(\frac{7\pi}{15}-\frac{4\pi}{3f},\pi \right)$. The trigonometric equation \eqref{coolsaet_eq6}  reduces to
\begin{align*}
	Q&(\delta, f):=  -\sin\left(\frac{2 \pi}{5}+3 \delta+\frac{2 \pi}{f}\right)+2 \sin\left(\frac{2 \pi}{5}+5 \delta+\frac{6 \pi}{f}\right)\\
	-&\cos\left(\frac{3 \pi}{10}+5 \delta+\frac{6 \pi}{f}\right)+2 \cos\left(\frac{3 \pi}{10}+3 \delta+\frac{2 \pi}{f}\right)-\sin\left(5 \delta+\frac{6 \pi}{f}\right)\\
	+&\sin\left(\frac{\pi}{5}+3 \delta-\frac{2 \pi}{f}\right)-\cos\left(\frac{3 \pi}{10}+5 \delta+\frac{2 \pi}{f}\right)+2 \cos\left(\frac{3 \pi}{10}+\delta+\frac{2 \pi}{f}\right)\\
	-&2 \sin\left(\frac{\pi}{5}+\delta-\frac{2 \pi}{f}\right)+\sin\left(\frac{2 \pi}{5}+\delta+\frac{2 \pi}{f}\right)-\sin\left(\frac{\pi}{5}+5 \delta+\frac{2 \pi}{f}\right)\\
	-&\cos\left(\frac{\pi}{10}+\delta-\frac{2 \pi}{f}\right)+2 \sin\left(\frac{\pi}{5}+3 \delta+\frac{2 \pi}{f}\right) =0.
\end{align*}
By setting $x={\mathrm e}^{2\mathrm{i} \delta},y={\mathrm e}^{\frac{4\mathrm{i} \pi}f}$,
$Q(\delta, f)$  is transformed to a polynomial
\begin{align*}
 L(x,y) &:= (\zeta_{5}^{2}-2 \zeta_{5}+1) x^{5} y^{3}-(\zeta_{5}^{3}-\zeta_{5}^{2}) x^{5} y^{2}+(2 \zeta_{5}^{3}-2 \zeta_{5}^{2}+\zeta_{5}) x^{4} y^{2}\\
+\zeta_{5}^{3} x^{4} y&- (2 \zeta_{5}^{2}+\zeta_{5}) x^{3} y^{2}-(\zeta_{5}^{4}+2 \zeta_{5}^{3}) x^{3} y +(2 \zeta_{5}^{2}+\zeta_{5}) x^{2} y^{2}+(\zeta_{5}^{4}+2 \zeta_{5}^{3}) \\
x^{2} y-\zeta_{5}^{2}& x y^{2}-(\zeta_{5}^{4}-2 \zeta_{5}^{3}+2 \zeta_{5}^{2}) x y -(\zeta_{5}^{3}-\zeta_{5}^{2}) y +2 \zeta_{5}^{4}-\zeta_{5}^{3}-1.
\end{align*}
The norm $N(L)$ of $L$ is the product $\prod_{k=1}^4 L_k$, where $L_k$ is $L$ with $\zeta_5$ replaced by $\zeta_5^k$.  We then obtain
$
\mathit{N(L)} := 25 x^{20} y^{12}-25 x^{20} y^{11}+50 x^{19} y^{10}+50 x^{18} y^{11}+5 x^{20} y^{8}+5 x^{19} y^{9}-115 x^{18} y^{10}-50 x^{17} y^{11}-30 x^{19} y^{8}-70 x^{18} y^{9}+90 x^{17} y^{10}-10 x^{19} y^{7}+60 x^{18} y^{8}+175 x^{17} y^{9}+100 x^{16} y^{10}+75 x^{18} y^{7}+30 x^{17} y^{8}-160 x^{16} y^{9}-150 x^{15} y^{10}+10 x^{18} y^{6}-230 x^{17} y^{7}-229 x^{16} y^{8}-130 x^{15} y^{9}+100 x^{14} y^{10}-80 x^{17} y^{6}+274 x^{16} y^{7}+158 x^{15} y^{8}+345 x^{14} y^{9}-35 x^{13} y^{10}-5 x^{17} y^{5}+281 x^{16} y^{6}+97 x^{15} y^{7}+386 x^{14} y^{8}-275 x^{13} y^{9}+10 x^{12} y^{10}+39 x^{16} y^{5}-472 x^{15} y^{6}-496 x^{14} y^{7}-660 x^{13} y^{8}+25 x^{12} y^{9}+x^{16} y^{4}-143 x^{15} y^{5}+201 x^{14} y^{6}+5 x^{13} y^{7}+335 x^{12} y^{8}+85 x^{11} y^{9}-7 x^{15} y^{4}+264 x^{14} y^{5}+590 x^{13} y^{6}+805 x^{12} y^{7}+253 x^{11} y^{8}-70 x^{10} y^{9}+26 x^{14} y^{4}-150 x^{13} y^{5}-845 x^{12} y^{6}-563 x^{11} y^{7}-431 x^{10} y^{8}+25 x^{9} y^{9}-40 x^{13} y^{4}-400 x^{12} y^{5}-57 x^{11} y^{6}-504 x^{10} y^{7}+243 x^{9} y^{8}-5 x^{8} y^{9}-25 x^{12} y^{4}+892 x^{11} y^{5}+819 x^{10} y^{6}+892 x^{9} y^{7}-25 x^{8} y^{8}-5 x^{12} y^{3}+243 x^{11} y^{4}-504 x^{10} y^{5}-57 x^{9} y^{6}-400 x^{8} y^{7}-40 x^{7} y^{8}+25 x^{11} y^{3}-431 x^{10} y^{4}-563 x^{9} y^{5}-845 x^{8} y^{6}-150 x^{7} y^{7}+26 x^{6} y^{8}-70 x^{10} y^{3}+253 x^{9} y^{4}+805 x^{8} y^{5}+590 x^{7} y^{6}+264 x^{6} y^{7}-7 x^{5} y^{8}+85 x^{9} y^{3}+335 x^{8} y^{4}+5 x^{7} y^{5}+201 x^{6} y^{6}-143 x^{5} y^{7}+x^{4} y^{8}+25 x^{8} y^{3}-660 x^{7} y^{4}-496 x^{6} y^{5}-472 x^{5} y^{6}+39 x^{4} y^{7}+10 x^{8} y^{2}-275 x^{7} y^{3}+386 x^{6} y^{4}+97 x^{5} y^{5}+281 x^{4} y^{6}-5 x^{3} y^{7}-35 x^{7} y^{2}+345 x^{6} y^{3}+158 x^{5} y^{4}+274 x^{4} y^{5}-80 x^{3} y^{6}+100 x^{6} y^{2}-130 x^{5} y^{3}-229 x^{4} y^{4}-230 x^{3} y^{5}+10 x^{2} y^{6}-150 x^{5} y^{2}-160 x^{4} y^{3}+30 x^{3} y^{4}+75 x^{2} y^{5}+100 x^{4} y^{2}+175 x^{3} y^{3}+60 x^{2} y^{4}-10 x y^{5}+90 x^{3} y^{2}-70 x^{2} y^{3}-30 x y^{4}-50 x^{3} y -115 x^{2} y^{2}+5 x y^{3}+5 y^{4}+50 x^{2} y +50 x y^{2}-25 y +25$.

{$ \displaystyle  $}
Similar to Example \ref{example1.2}, we first find that the only suitable solution for $y$ is $ e^{\frac{\mathrm{i}\pi}{5}}$ and $f=20$. But then $N(L) (x, y=e^{\frac{\mathrm{i}\pi}{5}}) = 0$ has no cyclotomic root. 
\end{example}

\subsection{The cyclotomic roots of a $3$-variable polynomial}

There is only one case in this paper inducing a $3$-variable polynomial.
%%ade_2bc
\begin{example}

        \label{ex:first}
Case $\{\aaa\ddd\eee,\bbb^2\ccc\}$ in Table \ref{ade 2bc}.
By solving the linear system of angle sums, we get
\[
\aaa=2\pi-\ddd-\eee, \quad \bbb=\pi-\frac{4\pi}{f}, \quad \ccc=\frac{8\pi}{f}. 
\]	
The trigonometric equation \eqref{coolsaet_eq6}  reduces to
\begin{align*}
Q&(\delta, \eee, f):= -2-2 \cos\left(\frac{12 \pi}{f}\right)+\cos\left(\delta+\epsilon+\frac{12 \pi}{f}\right)+\cos\left(2 \delta-\frac{4 \pi}{f}\right)\\
+&\cos\left(\delta+\epsilon-\frac{12 \pi}{f}\right)-\cos\left(\delta+\epsilon+\frac{4 \pi}{f}\right)+\cos\left(\delta+\epsilon-\frac{8 \pi}{f}\right)\\
-&\cos\left(\delta+\epsilon-\frac{4 \pi}{f}\right)+\cos\left(2 \delta+\frac{4 \pi}{f}\right)+2 \cos\left(2 \epsilon\right)+\cos\left(\delta+\epsilon+\frac{8 \pi}{f}\right)\\
-&\cos\left(2 \epsilon+\frac{8 \pi}{f}\right)-2 \cos\left(\delta+\epsilon\right)-\cos\left(2 \epsilon-\frac{8 \pi}{f}\right)+2 \cos\left(2 \delta\right)=0
\end{align*}
By setting $x={\mathrm e}^{\frac{4\mathrm{i} \pi}f},y={\mathrm e}^{\mathrm{i} \delta},z={\mathrm e}^{\mathrm{i} \eee}$, $Q(\ddd,\eee,f)$  is transformed to a polynomial
\begin{align*}
P(x,y,z)&:=(x+1)^{2} (x^{4} y^{3} z^{3}-x^{3} y^{3} z^{3}-x^{3} y^{2} z^{4}-2 x^{4} y^{2} z^{2}+x^{2} y^{4} z^{2}+\\
2 x^{2} y^{2} z^{4}&+4 x^{3} y^{2} z^{2}-x y^{3} z^{3}-x y^{2} z^{4}+x^{4} y z-6 x^{2} y^{2} z^{2}+y^{3} z^{3}-x^{3} y^{2}\\
-x^{3} y z&+4 x y^{2} z^{2}+2 x^{2} y^{2}+x^{2} z^{2}-2 y^{2} z^{2}-x y^{2}-x y z+y z).
\end{align*}
Again we only need to consider the second factor. Setting $ u=yz, w=y/z $, it is transformed to a full polynomial $L(x,u,w)=
 x^{4} u^{2} w-2 x^{4} u w-x^{3} u^{2} w+x^{2} u^{2} w^{2}+x^{4} w-x^{3} u^{2}+4 x^{3} u w-x^{3} w^{2}-x^{3} w+2 x^{2} u^{2}-6 x^{2} u w+2 x^{2} w^{2}-x u^{2} w-x u^{2}+4 x u w-x w^{2}+u^{2} w+x^{2}-x w-2 u w+w.$

Following the algorithm in \cite{AS}, we computed $15$ resultants of $L(x,u,w)$ with each of the following:
$L(-x,u,w)$, $L(x,-u,w)$, $L(x,u,-w)$,
$L(-x,\\
-u,w)$, $L(-x,u,-w)$, $L(x,-u,-w)$,
$L(-x,-u,-w)$, $L(x^2,u^2,w^2)$, $L(-x^2,u^2,w^2)$,
$L(x^2,-u^2,w^2)$,
$L(x^2,u^2,-w^2)$, $L(-x^2,-u^2,w^2)$, $L(-x^2,\\
u^2,
-w^2)$, $L(x^2,-u^2,-w^2)$,
and $L(-x^2,-u^2,-w^2)$ in the Appendix. 

For example, 
$  \mathrm{Resultant} (L(-x,u,w) ,L(x,u,w) ,w) $
is
$
-4 x^{2} (x^{2}+1)^{2}$ $ 
(u -1)^{4} (x^{3} u^{2}+x^{4} u +x^{2} u^{2}-x u^{2}+x^{2} u +x^{3}+x^{2}+u -x ) (x^{3} u^{2}-x^{4} u -x^{2} u^{2}-x u^{2}-x^{2} u +x^{3}-x^{2}-u -x )
$.
Similar to Example \ref{example1.2}, the only suitable solution is $x = e^{\frac{\mathrm{i}\pi}{3}}$ and $f=12$. This implies $\bbb=\ccc=\frac{2\pi}{3}$, contradicting $\bbb\neq \ccc$. 

Solving all 15 resultants yields four rational angle sets (Table \ref{ade 2bc}):
\begin{align*}
&(2,4,2,5,3)\pi/5,  &(4,8,4,9,7)\pi/10, \\
&(16,8,4,1,3)\pi/10, &(10,12,6,5,15)\pi/15.
\end{align*}
The first two violate Lemma \ref{geometry1} ($\beta > \gamma$ and $\delta > \epsilon$). The third has $\aaa>\frac{3\pi}{2}$, contradicting \cite[Table 2]{slw}. These yield non-simple pentagons. The fourth gives a simple pentagon admitting:
\begin{itemize}
    \item A non-symmetric 3-layer earth map tiling: $T(20\alpha\delta\epsilon,\ 10\beta^2\gamma,\ 2\gamma^5)$
    \item Its unique flip modification: $T(20\alpha\delta\epsilon,\ 8\beta^2\gamma,\ 4\beta\gamma^3)$ 
\end{itemize}
(as per main theorem).
\end{example}

\section{Rational $a^4b$-pentagons suitable for spherical tilings}\label{cases}

We often express angles in units of $\pi$ radians for simplicity. So the sum of all angles at a vertex is $2$. When an almost equilateral pentagon is symmetric about an axis through
$\aaa$, i.e. $\bbb=\ccc,\ddd=\eee$ it can be divided into two congruent quadrilaterals along this axis, and the corresponding classification has already been completed by \cite{lw1, lw2, lw3}. 

In this section we only need to consider non-symmetric pentagons ($\beta \neq \gamma$, $\delta \neq \epsilon$) with $f \geq 12$ even and angles in $(0,2\pi)$ in $a^4b$-tilings. Per \cite{slw} (Tables 4-6,10-12,14,21,23,28,30) according to $a^2b$-vertex types, the 258 vertex combinations yield 43 rational pentagons satisfying the above constraints. The computation details are shown in the Appendix. Only 12 are simple (Figure \ref{easy pentagon}), with just one (first row, fourth pentagon) admitting two distinct $a^4b$-tilings.

\begin{figure}[htp]
 \centering
 \begin{overpic}[scale=0.20]{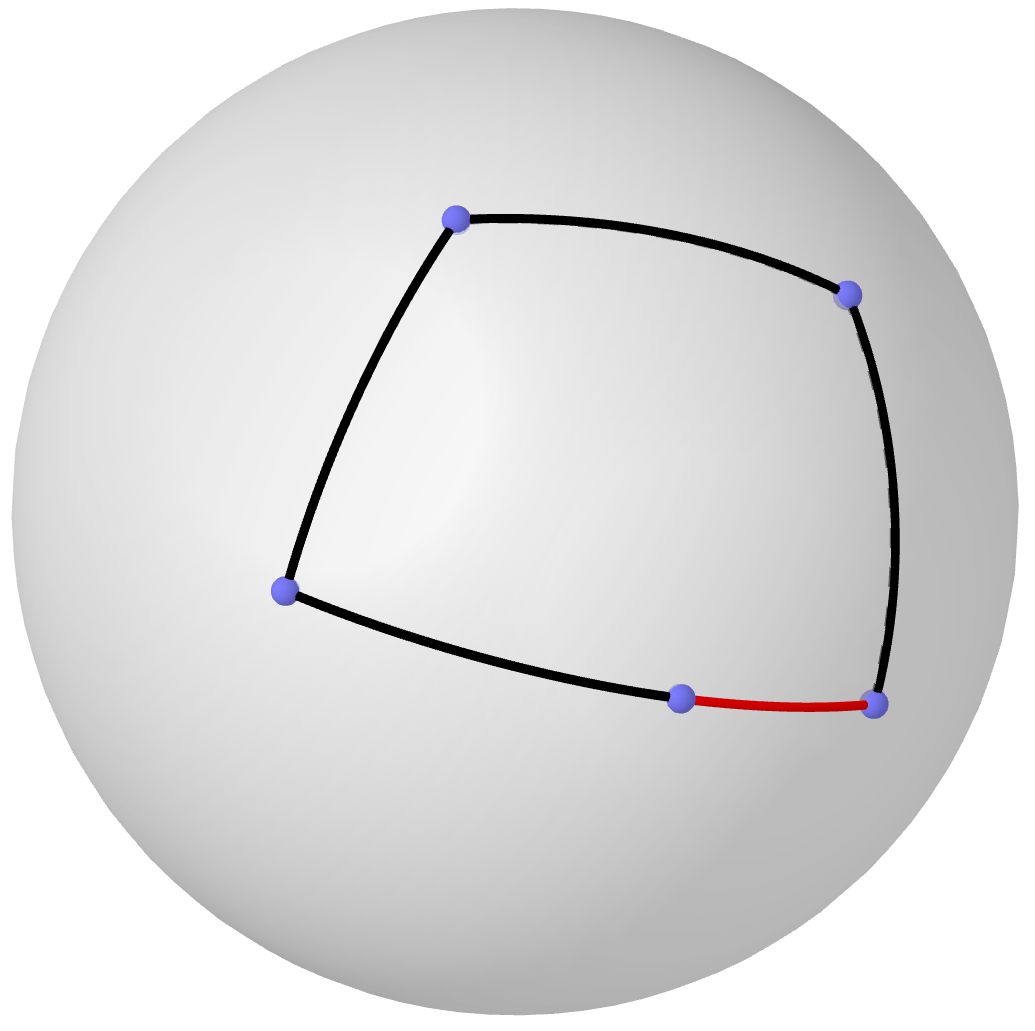}
   \put(0,-15.5){\small  $(4,3,4,6,3)/6$}
      \end{overpic}\hspace{10pt}
 \begin{overpic}[scale=0.19]{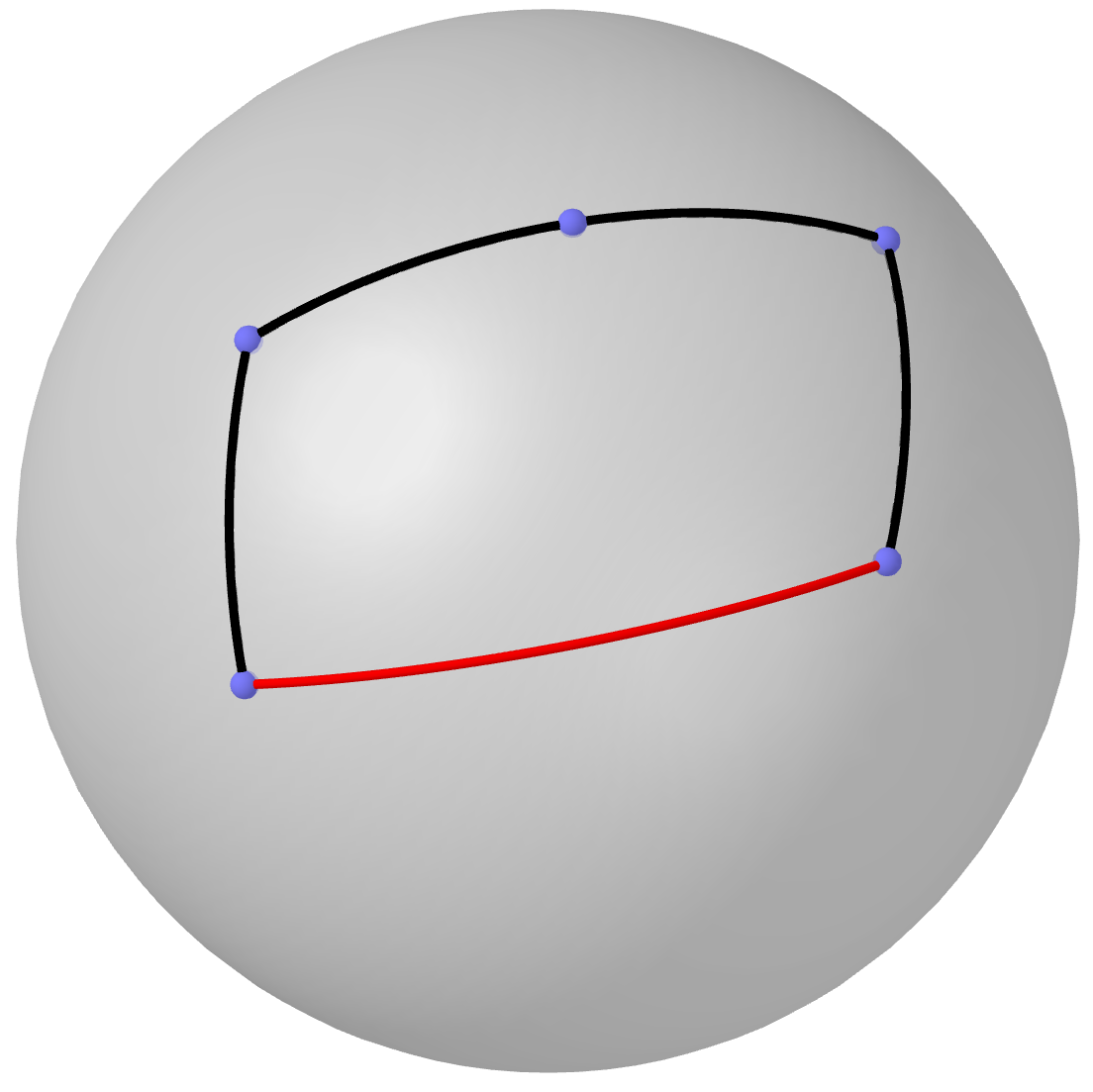}
   \put(0,-15.5){\small  $(6,4,3,3,4)/6$}
     \end{overpic}\hspace{10pt}
 \begin{overpic}[scale=0.18]{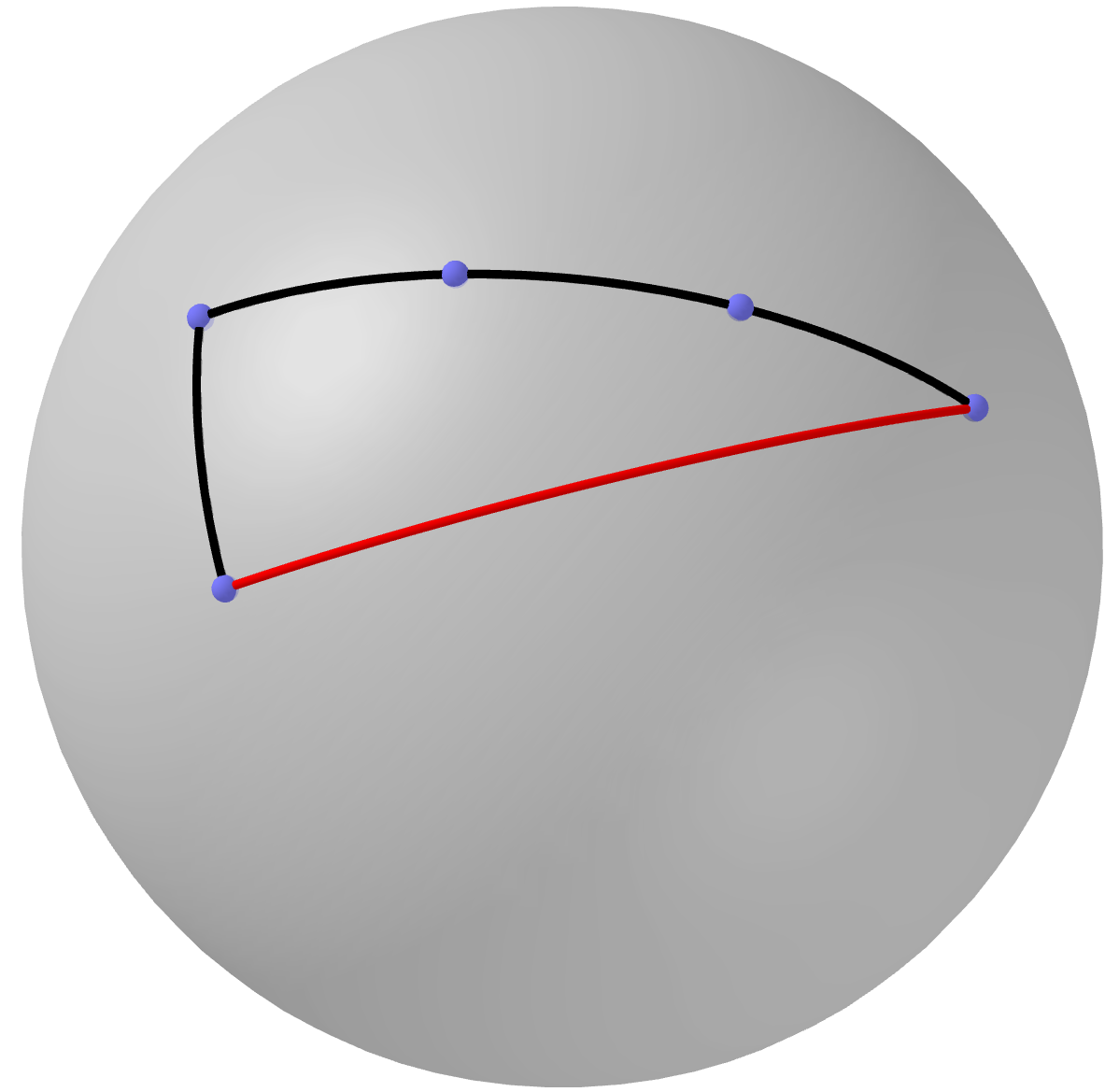}
    \put(5,-15.5){\small  $(6,3,6,3,1)/6$}
      \end{overpic}\hspace{10pt}
        \begin{overpic}[scale=0.033]{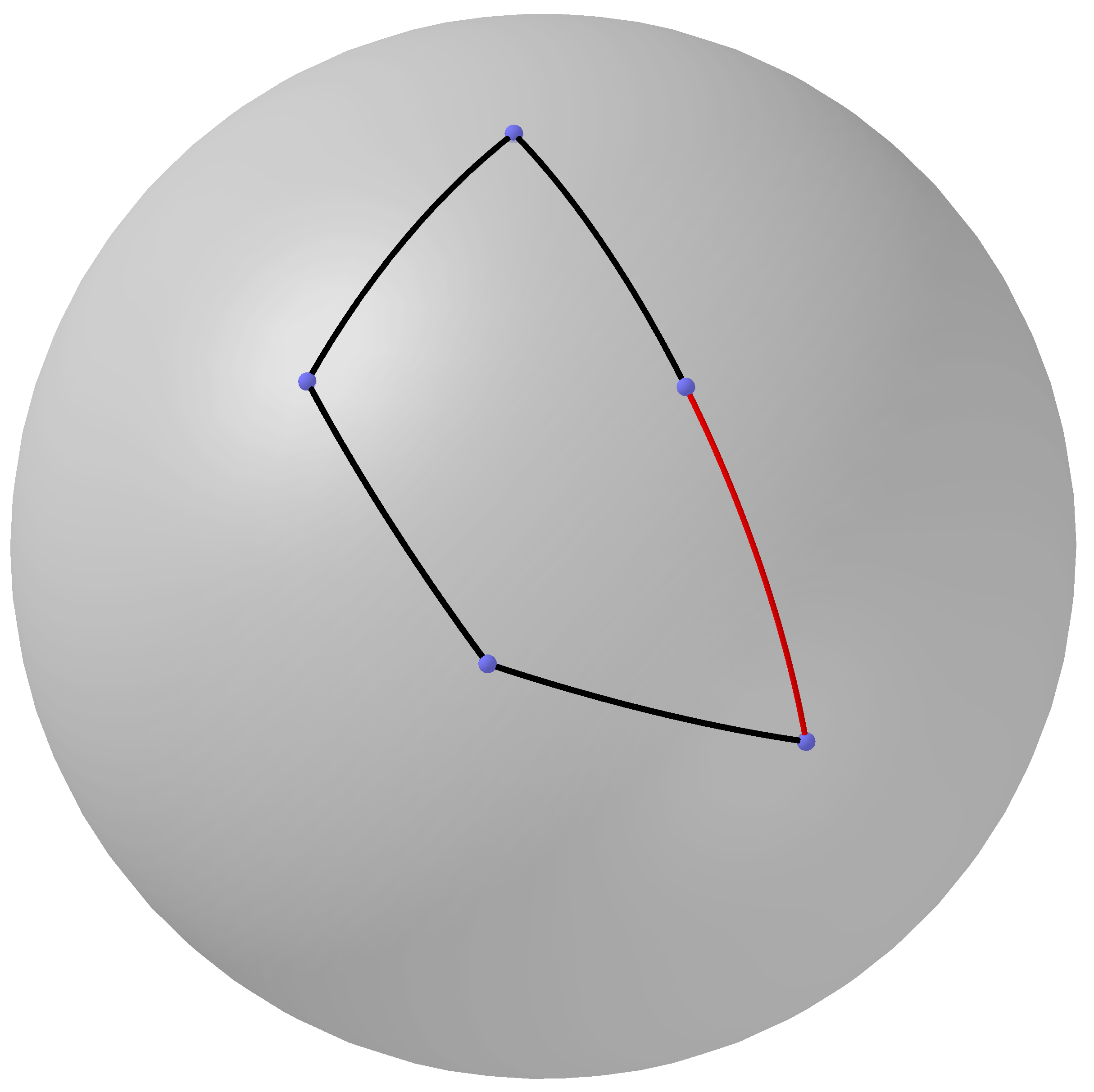}
       \put(-5,-15.5){\small  $(10,12,6,5,15)/15$}
          \end{overpic}\hspace{60pt}

        \begin{overpic}[scale=0.21]{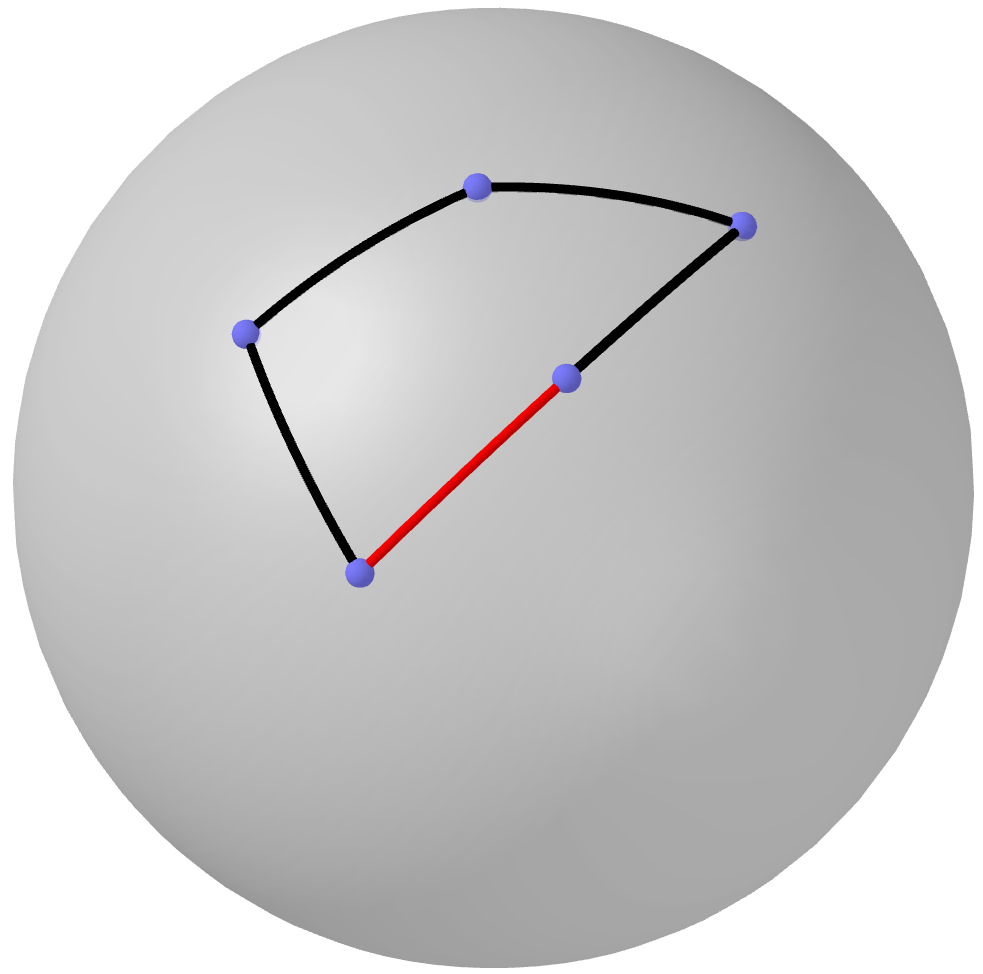}
    \put(0,-15.5){\small  $(6,4,2,3,7)/7$}
     \end{overpic}\hspace{10pt}
   \begin{overpic}[scale=0.20]{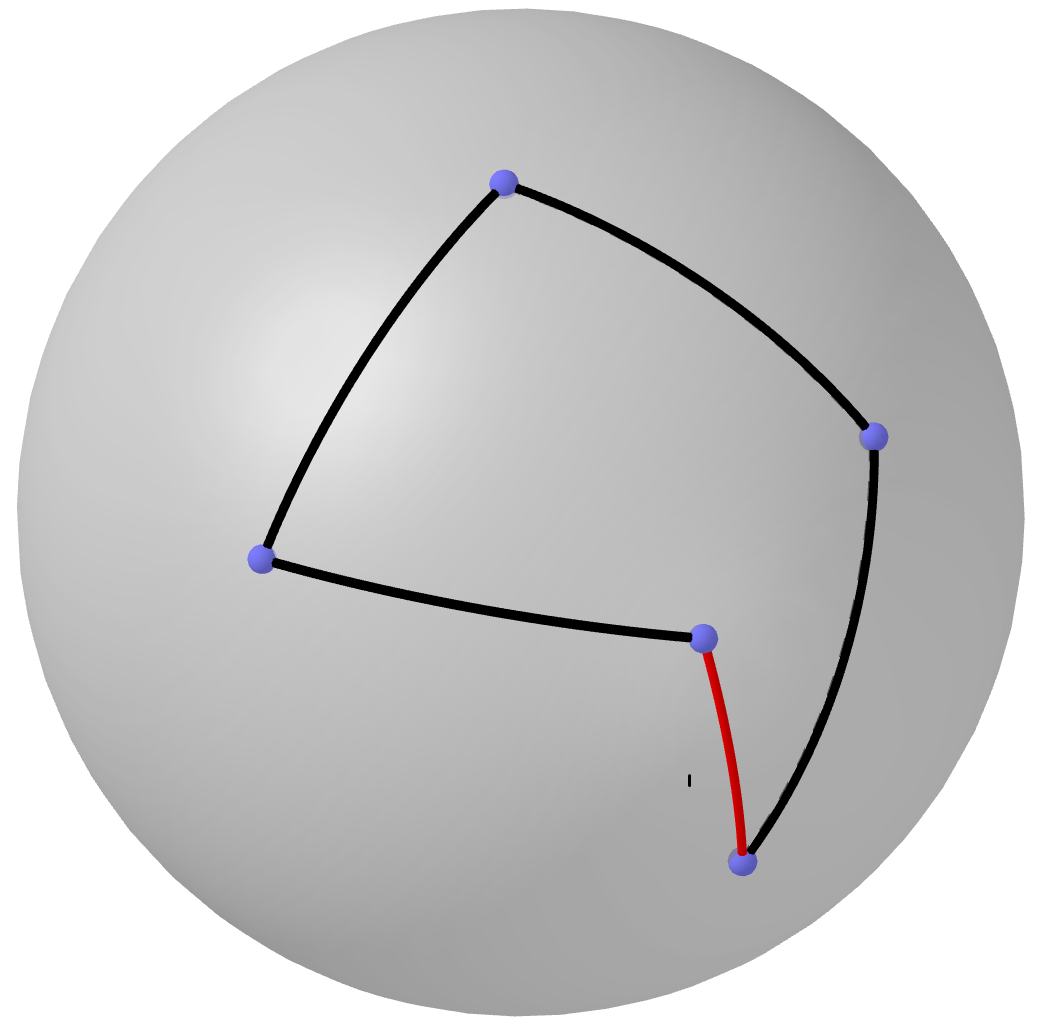}
    \put(0,-15.5){\small  $(8,6,10,19,3)/14$}
     \end{overpic}\hspace{10pt}    
   \begin{overpic}[scale=0.195]{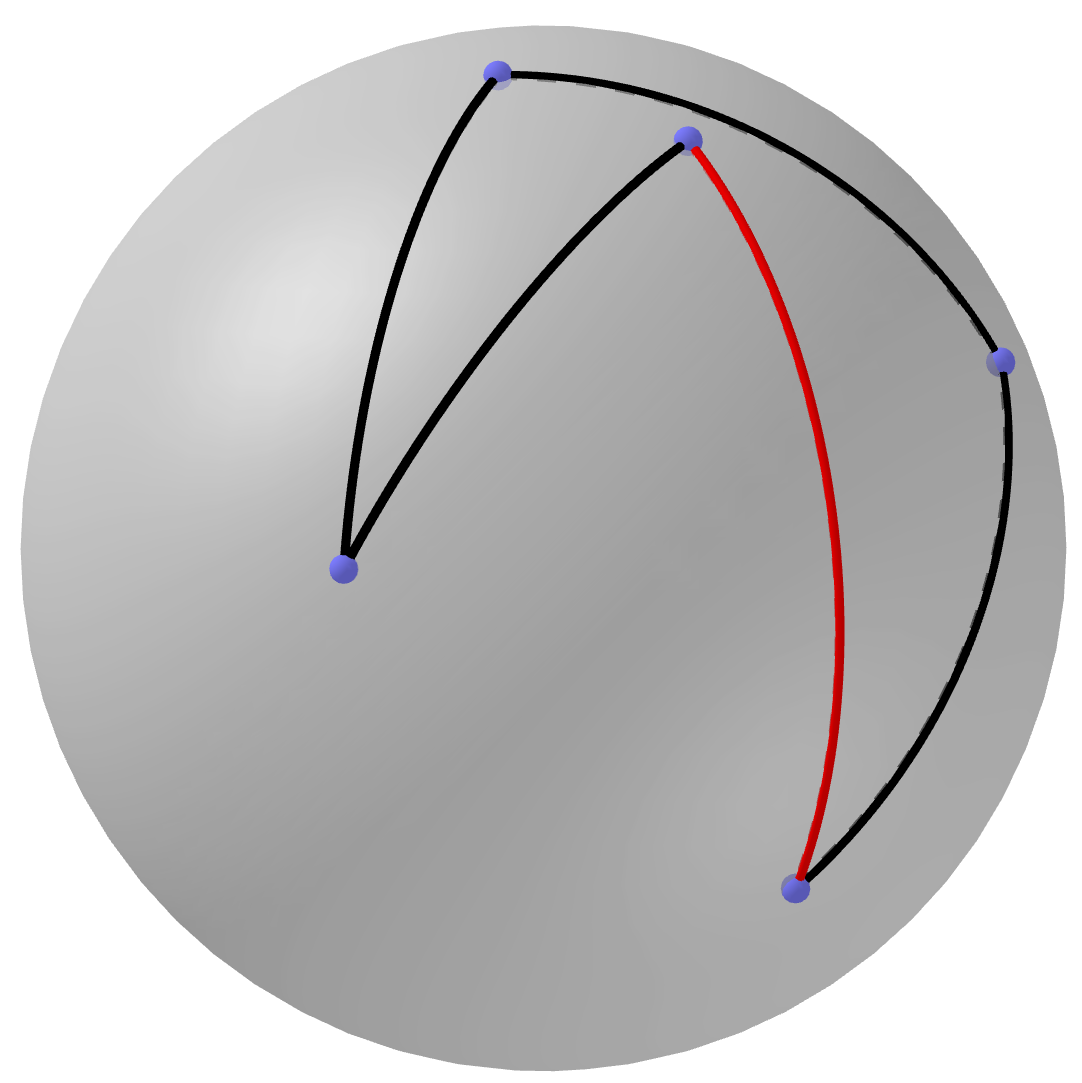}
       \put(8,-15.5){\small  $(8,2,10,23,3)/14$}
       \end{overpic}\hspace{10pt}
      \begin{overpic}[scale=0.22]{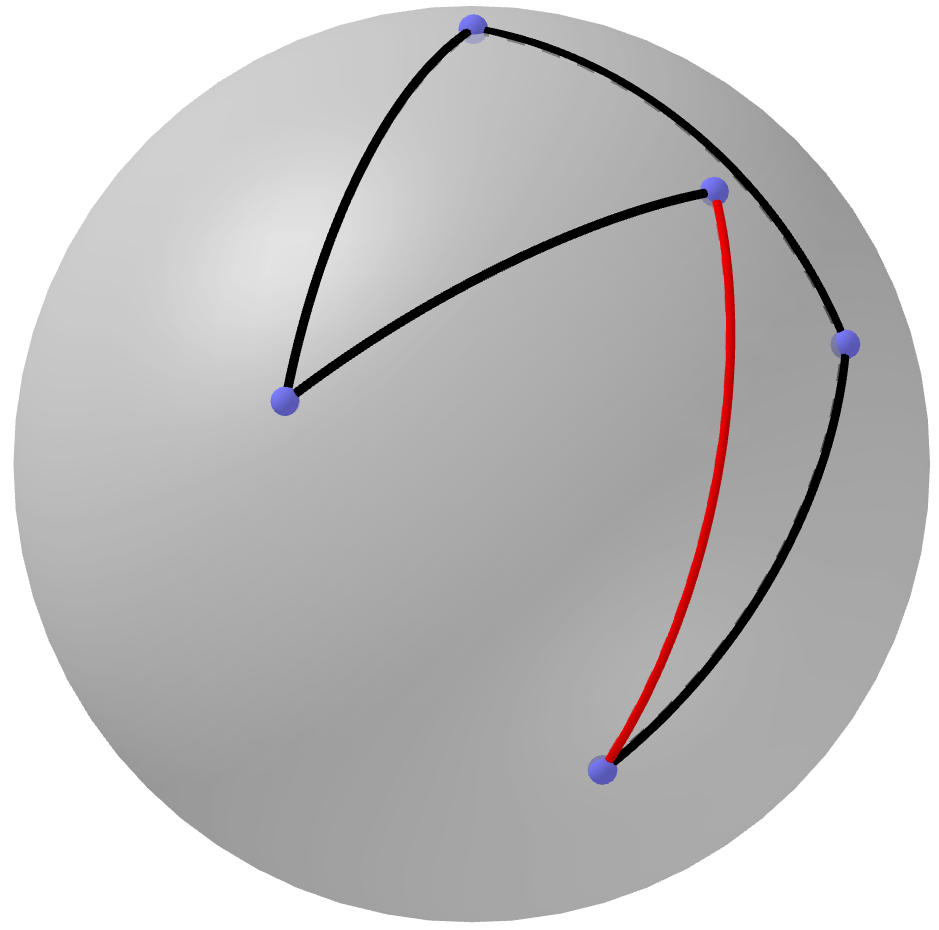}
     \put(10,-15.5){\small  $(4,2,6,13,1)/8$}
       \end{overpic}\hspace{60pt}

   \begin{overpic}[scale=0.21]{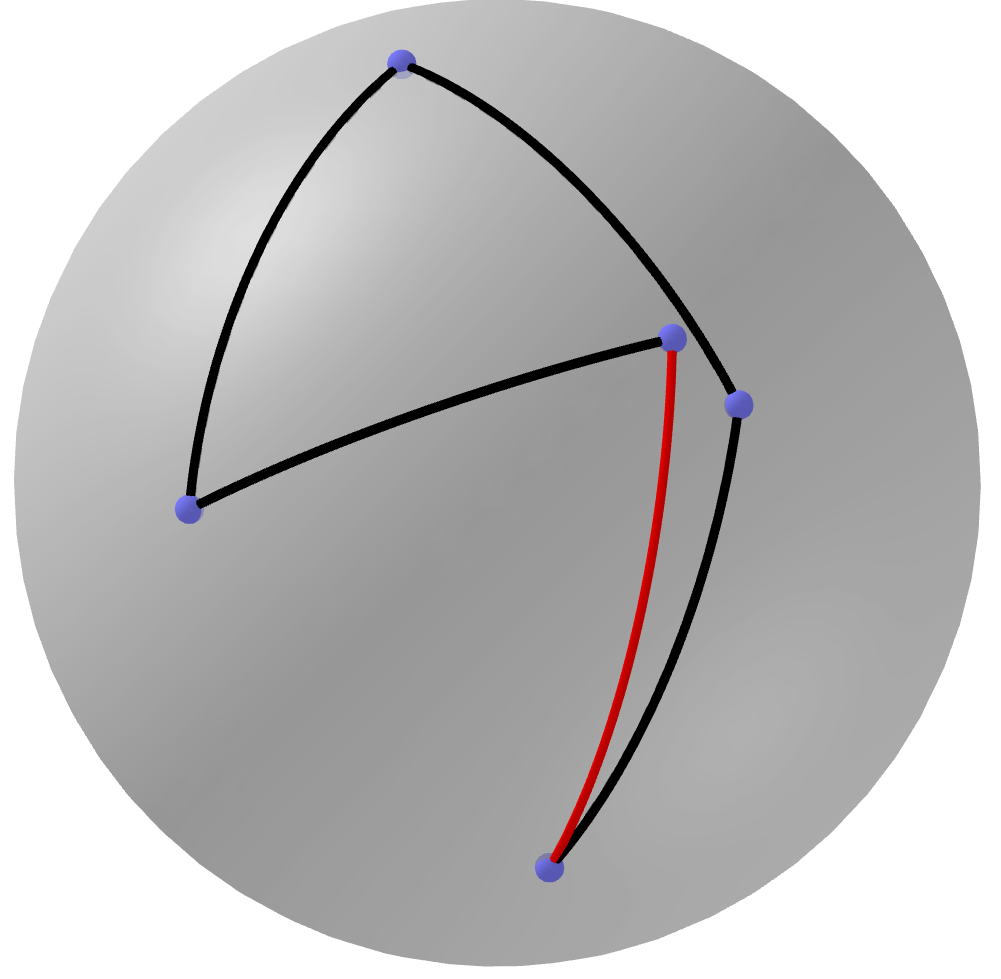}
     \put(0,-15.5){\small  $(8,6,14,29,1)/14$}
      \end{overpic}\hspace{10pt}
      \begin{overpic}[scale=0.175]{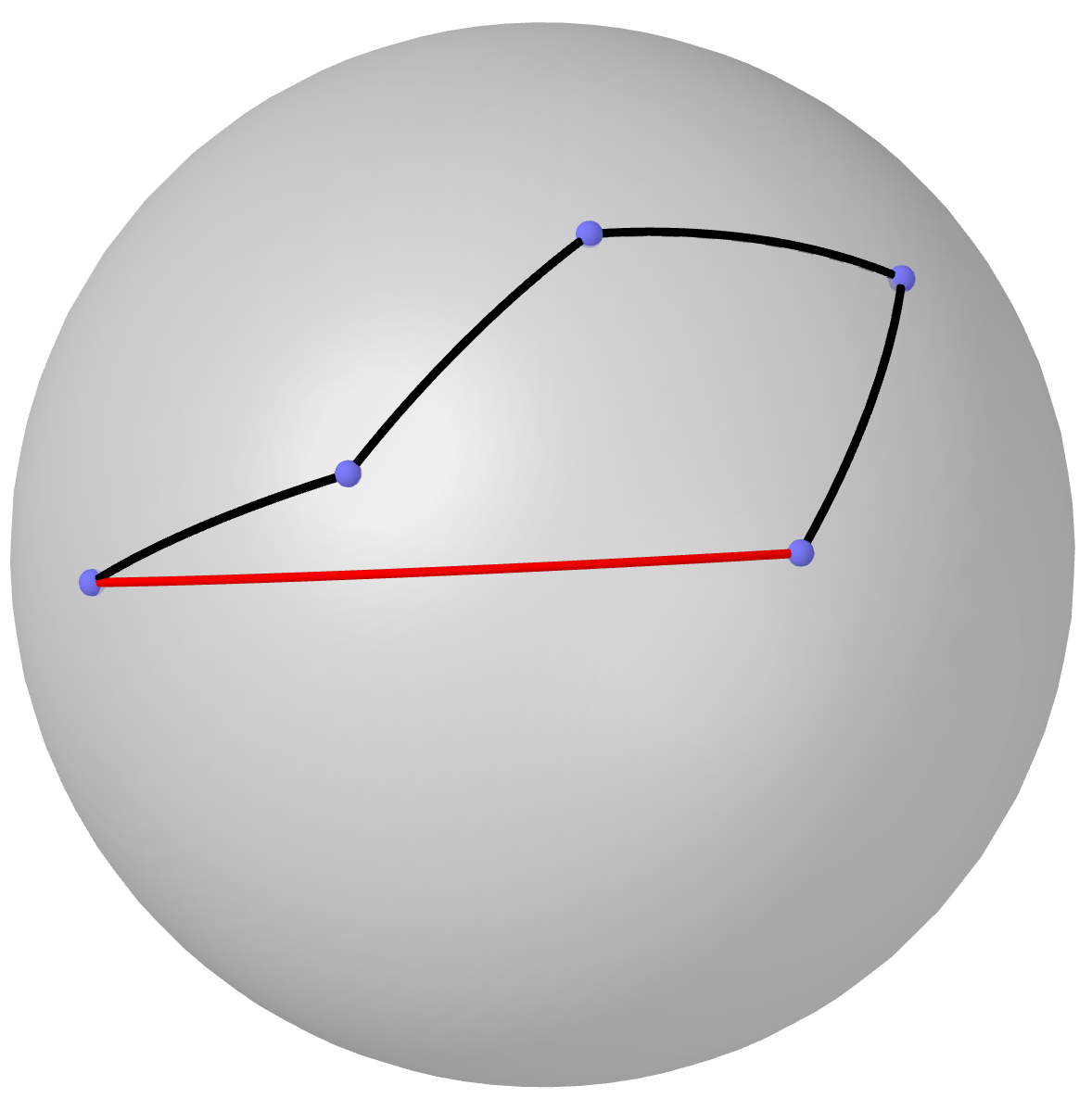}
    \put(-5,-15.5){\small  $(8,12,4,1,7)/10$}
      \end{overpic}\hspace{10pt}
      \begin{overpic}[scale=0.24]{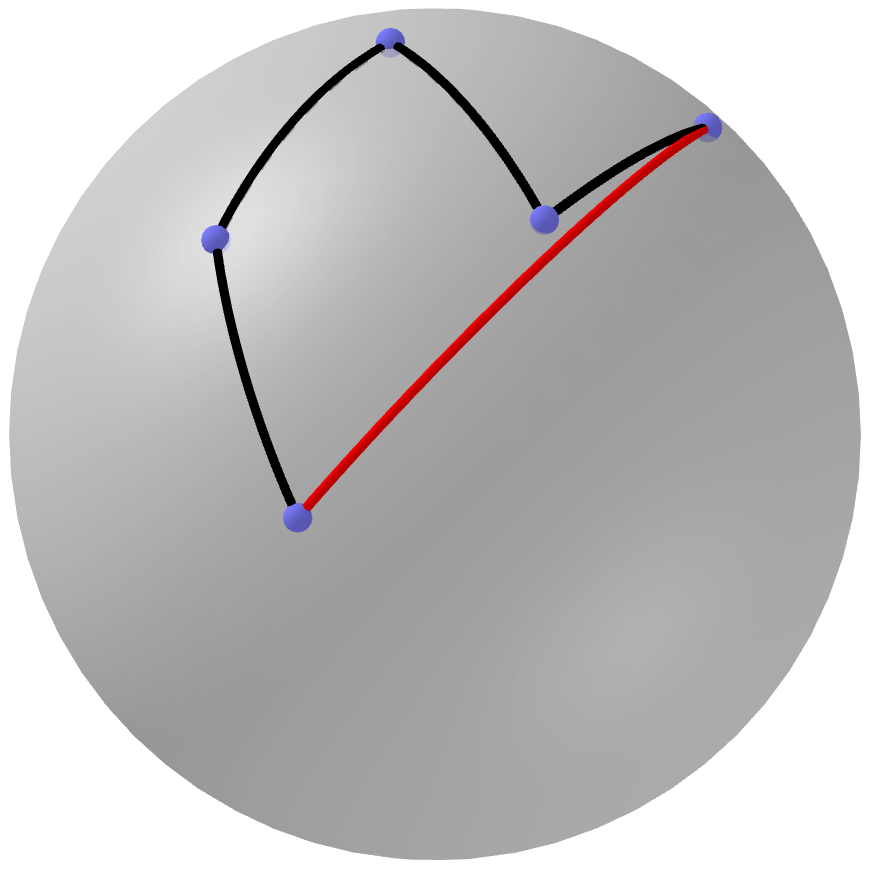}
     \put(-12,-15.5){\small  $(12,48,24,1,11)/30$}
       \end{overpic}\hspace{10pt}
      \begin{overpic}[scale=0.18]{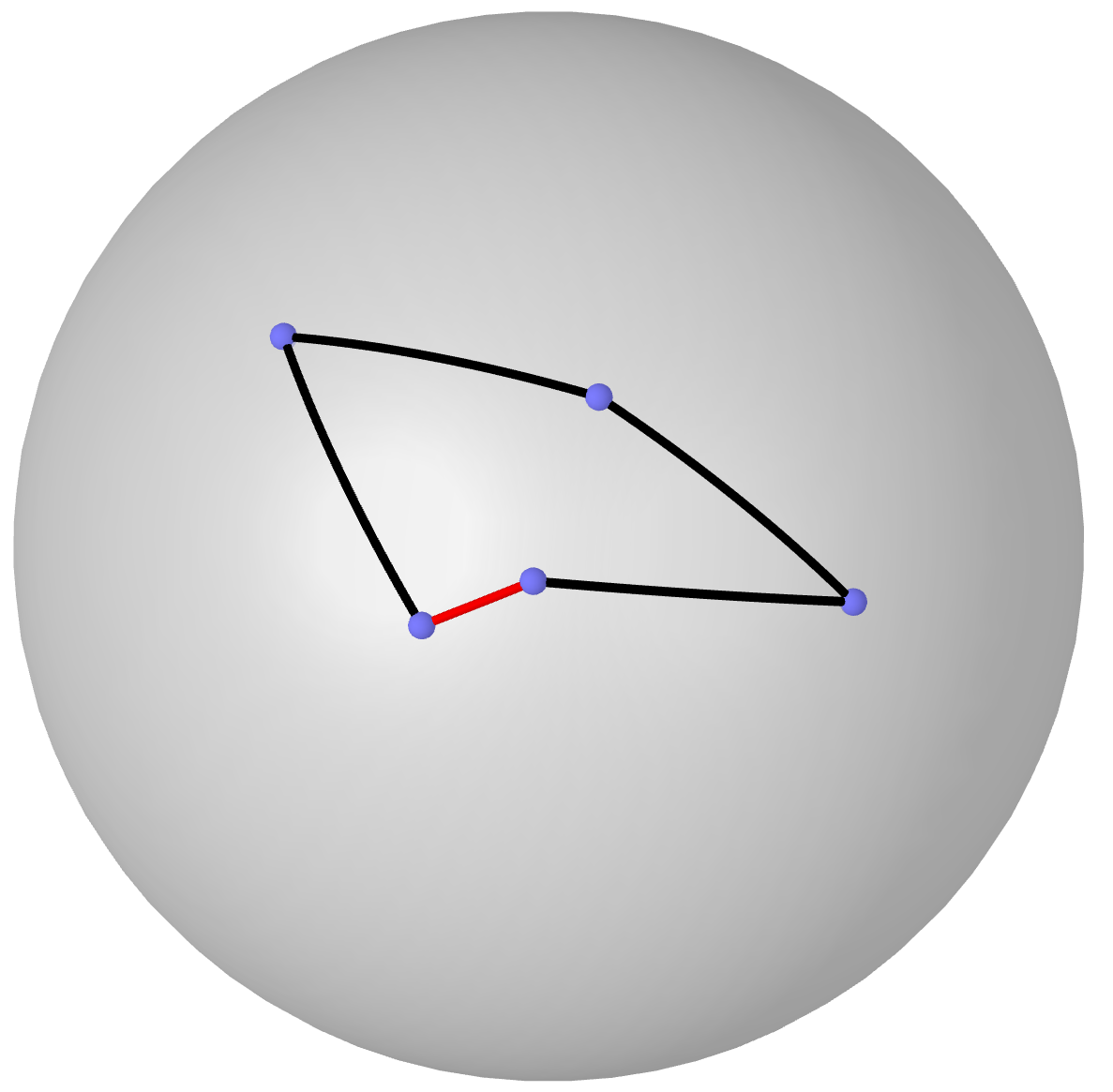}
   \put(-7,-15.5){\small  $(18,6,4,11,23)/20$}
  \end{overpic}\\
  \hspace{40pt}
 \caption{All simple non-symmetric rational $a^4b$-pentagons with area $\tfrac{4\pi}{f}$ for some even integer $f \ge 12$.} 
 \label{easy pentagon} 
\end{figure}

\subsection{Tilings with  three $a^2b$-vertex types}

By the proof of \cite[Lemma 14]{slw}, it is impossible to have more than three $a^2b$-vertex types in any non-symmetric $a^4b$-tiling, and there are only $6$ cases with exactly three $a^2b$-vertex types in Table \ref{3a2b1}. None of them admits any rational pentagon solution, with computation details given in the Appendix.

\begin{table}[htp]        
	\centering 
	\caption{Tilings with  three $a^2b$-vertex types}
	\label{3a2b1}
	\begin{tabular}{|c|c|}
		\hline  
		Vertex & Contradiction\\
		\hline\hline
		$\aaa\ddd\eee,\bbb\ddd^2,\ccc\eee^2$&\multirow{4}{*}{ No rational solution }\\
		\cline{1-1}
		$\aaa\ddd\eee,\bbb\ddd\eee,\ccc\eee^2$&\\
		\cline{1-1}
		$\aaa\ddd^2,\bbb\ddd^2,\ccc\eee^2$&\\
		\cline{1-1}
		$\aaa\ddd^2,\bbb\ddd\eee,\ccc\eee^2$&(or simply "No solution")\\
		\cline{1-1}
		$\aaa\eee^2,\bbb\ddd^2,\ccc\eee^2$&\\
		\cline{1-1}
		$\aaa\eee^2,\bbb\ddd\eee,\ccc\eee^2$&\\
		\hline
	\end{tabular}      	
\end{table}

\subsection{Tilings with two $a^2b$-vertex types}

By  \cite[Lemma 15]{slw},  we  derive $110$ distinct cases. These results are summarized in  Table \ref{two a2b}.

\begin{table}[htp]        
	\centering 
	\caption{Tilings with two $a^2b$-vertex types}\label{two a2b}
	\resizebox{\linewidth}{!}{\begin{tabular}{|c|c|c|c|}
			\hline  
			\multicolumn{2}{|c|}{Vertex}& $(\aaa,\bbb,\ccc,\ddd,\eee),f$ & Contradiction\\
			\hline\hline
			\multirow{4}{*}{$\aaa\ddd^2,\bbb\ddd^2$}& $\aaa^2\ccc,\aaa\bbb\ccc,\bbb^2\ccc,\ccc^4,\ccc^5$&& No solution \\
			\cline{2-4}
			&$\ccc^3$&$(4,4,6,7,9)/9,12$& \multirow{3}{*}{ $\bbb<\ccc$ and $\ddd<\eee$ }\\
			\cline{2-3}
			&$\aaa\ccc^2,\bbb\ccc^2$&$(4,4,f-2,f-2,f)/f$&\\
			\cline{2-3}
			&$\aaa\ccc^3,\bbb\ccc^3$&$(4,4,6,9,11)/11,44$&\\
			\hline
			\multirow{5}{*}{$\aaa\ddd^2,\bbb\ddd\eee$} & $\aaa\bbb\ccc,\bbb^2\ccc,\aaa\ccc^3,\bbb\ccc^3,\ccc^5$&& No solution \\
			\cline{2-4}			
			&$\aaa\ccc^2$&$(4,3,5,5,6)/7,14$&$\bbb<\ccc$ and $\ddd<\eee$\\
			\cline{2-4}
			&$\ccc^3$&$(10,8,12,13,15)/18,18$&$\bbb<\ccc$ and $\ddd<\eee$\\
			\cline{2-4}
			&$\bbb\ccc^2,\ccc^4$&$(4,6,3,4,2)/6,24$&$\bbb>\ccc$ and $\ddd>\eee$\\			     
			\cline{2-4}
			&$\aaa^2\ccc$&$(18,6,4,11,23)/20,40$&Balance
Lemma\\
			\hline
			\multirow{8}{*}{$\aaa\ddd^2,\bbb\eee^2$} & $\aaa^2\ccc,\bbb^2\ccc$&& No solution \\
			\cline{2-4}
			&$\aaa\ccc^2$&$(6,4,3,3,4)/6,12$&AAD\\
			\cline{2-4}
			&$\aaa\bbb\ccc,\ccc^3$&$(2,6,4,5,3)/6,12$&\multirow{6}{*}{ $\bbb>\ccc$ and $\ddd>\eee$ }\\
			\cline{2-3}
			&$\aaa\ccc^3$ &$(6,18,10,15,9)/18,18$&\\
			\cline{2-3}
			&$\bbb\ccc^2,\ccc^4$ &$(2,6,3,5,3)/6,24$&\\
			\cline{2-3}
			&$\bbb\ccc^2,\aaa\ccc^3$&$(2,6,4,6,4)/7,28$&\\
			\cline{2-3}
			&$\bbb\ccc^3$&$(6,42,2,21,3)/24,48$&\\
			\cline{2-3}
			&$\ccc^5$&$(10,30,12,25,15)/30,60$&\\
			\hline
			\multirow{2}{*}{$\aaa\ddd^2,\ccc\ddd^2$} & $\aaa^2\bbb,\aaa\bbb\ccc,\bbb\ccc^2,\bbb^3,\aaa\bbb^3,\bbb^3\ccc,\bbb^5$&& No solution \\
			\cline{2-4}
			&$\aaa\bbb^2,\bbb^2\ccc,\bbb^4$ &$(6,3,6,3,1)/6,24$&AAD\\
			\hline
			\multirow{3}{*}{$\aaa\ddd^2,\ccc\ddd\eee$} & $\aaa^2\bbb,\bbb^3,\bbb^2\ccc,\aaa\bbb^3,\bbb^4,\bbb^3\ccc,\bbb^5$&& No solution \\
			\cline{2-4}
			&$\aaa\bbb\ccc$ &$(2,6,4,5,3)/6,12$&\multirow{2}{*}{$\bbb>\ccc$ and $\ddd>\eee$}\\
			\cline{2-3}
			&$\aaa\bbb^2,\bbb\ccc^2$ &$(2,6,4,6,4)/7,28$&\\
			\hline
			\multirow{2}{*}{$\aaa\ddd^2,\ccc\eee^2$} & $\aaa^2\bbb,\aaa\bbb\ccc,\bbb^3,\bbb^2\ccc,\bbb\ccc^2,\aaa\bbb^3,\bbb^4,\bbb^3\ccc,\bbb^5$&& No solution \\
			\cline{2-4}
			&$\aaa\bbb^2$ &$(6,1,2,1,3)/4,16$&$\bbb<\ccc$ and $\ddd<\eee$\\
			\hline
			$\aaa\ddd\eee,\bbb\ddd^2$ & $\aaa^2\ccc,\aaa\bbb\ccc,\aaa\ccc^2,\bbb^2\ccc,\bbb\ccc^2,\ccc^3,\aaa\ccc^3,\bbb\ccc^3,\ccc^4,\ccc^5$&&\multirow{2}{*}{ No solution }  \\
			\cline{1-3}
			\multirow{3}{*}{$\aaa\ddd\eee,\bbb\ddd\eee$} & $\aaa^2\ccc,\aaa\bbb\ccc,\bbb^2\ccc,\ccc^3,\aaa\ccc^3,\bbb\ccc^3,\ccc^4,\ccc^5$&&  \\
			\cline{2-4}
			&\multirow{2}{*}{$\aaa\ccc^2,\bbb\ccc^2$} &$(2,2,4,3,5)/5,20$& \multirow{2}{*}{ $\bbb<\ccc$ and $\ddd<\eee$ }  \\
			\cline{3-3}
			& &$(4,4,8,7,9)/10,20$&\\
			\hline
			\multirow{2}{*}{$\aaa\ddd\eee,\bbb\eee^2$} &
			$\aaa^2\ccc,\aaa\bbb\ccc,\aaa\ccc^2,\bbb^2\ccc,\bbb\ccc^2,\ccc^3,\aaa\ccc^3,\ccc^4$&& No solution \\
			\cline{2-4}
			&$\bbb\ccc^3,\ccc^5$&$(2,4,2,5,3)/5,20$&$\bbb>\ccc$ and $\ddd>\eee$\\
			\hline
			$\bbb\ddd^2,\ccc\eee^2$ &
			$\aaa^3,\aaa^2\bbb,\aaa^2\ccc,\aaa\bbb\ccc,\aaa\bbb^2,\aaa\ccc^2,\aaa^4,\aaa^3\bbb,\aaa^3\ccc,\aaa^5$&& \multirow{2}{*}{ No solution }  \\
			\cline{1-3}
			$\bbb\ddd\eee,\ccc\eee^2$ &
			$\aaa^3,\aaa^2\bbb,\aaa^2\ccc,\aaa\bbb^2,\aaa\bbb\ccc,\aaa\ccc^2,\aaa^4,\aaa^3\bbb,\aaa^3\ccc,\aaa^5$&&  \\
			\hline
	\end{tabular} }
	
\end{table}

There are only three simple rational $a^4b$-pentagons to consider:

For Case $\{\aaa\ddd^2,\bbb\ddd\eee,\aaa^2\ccc,(18,6,4,11,23)/20,40\}$, by $\aaa\ddd^2$ and the parity lemma, we obtain that  $\eee^2\cdots$ is a vertex,  a contradiction.

For Case $\{\aaa\ddd^2,\bbb\eee^2,\aaa\ccc^2,(6,4,3,3,4)/6,12\}$, by the parity lemma, we have $\text{AVC}\sub\{\aaa\ccc^2,\aaa\ddd^2,\bbb^3,\bbb\eee^2,\ccc^4,\ccc^2\ddd^2,\ddd^4\}$, which admits a unique solution satisfying the balance Lemma: $\{6\aaa\ccc^2,6\aaa\ddd^2,6\bbb\eee^2,2\bbb^3\}$. This yields a unique AAD $\bbb\eee^2=\thin^{\aaa}\bbb^{\ddd}\thin^{\ccc}\eee^{\ddd}\thick^{\ddd}\eee^{\ccc}\thin$ inducing a vertex $\ccc\ddd\cdots$, which contradicts the AVC. 

For Case $\{\aaa\ddd^2,\ccc\ddd^2,\aaa\bbb^2(\text{or}\, \bbb^2\ccc, \bbb^4 ),(6,3,6,3,1)/6,24\}$, we obtain\begin{align*}
\text{AVC}\sub\{&\aaa\ddd^2,\ccc\ddd^2,\aaa\bbb^2,\bbb^2\ccc,\bbb^4,\bbb^2\ddd^2,
\ddd^4,\aaa\ddd\eee^3,\ccc\ddd\eee^3,
\bbb^2\ddd\eee^3,\ddd^3\eee^3,\aaa\eee^6,\ccc\eee^6,\\
&\bbb^2\eee^6, \ddd^2\eee^6,\ddd\eee^9,\eee^{12}\}.
\end{align*}
From this, we conclude that there are no vertices of the form $\aaa^2\cdots$, $\aaa\ccc^2\cdots$,or $\ccc^2\cdots$. This further implies (by AAD) that $\text{AVC}\sub\{\aaa\bbb^2,\aaa\ddd^2,\\
\bbb^2\ccc,\ccc\ddd^2,\bbb^4,\bbb^2\ddd^2,\ddd^4,\ddd^3\eee^3\}$. No solution satisfying the balance lemma exists in this case.
%\begin{table}[htp]        
%	\centering 
%	\caption{ $\cos\frac{1}{2}(\ddd+\eee-\aaa)=0$.}
%	\label{d+e-a}
%	\begin{tabular}{|c|c|c|}
%			\hline  
%			Vertex& $(\aaa,\bbb,\ccc,\ddd,\eee),f$ & Contradiction\\
%			\hline\hline
%			$\aaa\ddd^2,\bbb\ddd\eee,\ccc^3$&$(10,8,12,13,15)/18,18$&\multirow{1}{*}{$\bbb<\ccc$ and $\ddd<\eee$ } \\
%			\cline{1-3}
%			$\aaa\ddd^2,\bbb\eee^2,\aaa\ccc^3$&$(6,18,10,15,9)/18,18$& \multirow{4}{*}{$\bbb>\ccc$ and $\ddd>\eee$ }\\   
%			\cline{1-2}
%		    $\aaa\ddd^2,\bbb\eee^2,\bbb\ccc^2$&\multirow{3}{*}{$(2,6,3,5,3)/6,24$}& \\
%			\cline{1-1}
%			 $\aaa\ddd^2,\bbb\eee^2,\ccc^4$&&\\
%			 	  			
%			\cline{1-1}
%			$\aaa\ddd^2,\eee^4,\bbb\ccc^2$&&\\	  			
%		\cline{1-3}
%					
%
%			
%	\end{tabular}
%	
%\end{table}

\subsection{Tilings with a unique $a^2b$-vertex type  $\aaa\ddd^2$, $\bbb\ddd^2$ or $\ccc\ddd^2$}

By   \cite[Proposition 3]{slw}, we obtain  $22$ distinct cases.  These results are summarized in  Table \ref{ad2 bd2 cd2}.

\begin{table}[htp]        
	\centering 
	\caption{Tilings with a unique $a^2b$-vertex type $\aaa\ddd^2$, $\bbb\ddd^2$ or $\ccc\ddd^2$}\label{ad2 bd2 cd2}
	\begin{tabular}{|c|c|c|c|}
			\hline  
			\multicolumn{2}{|c|}{Vertex}& $(\aaa,\bbb,\ccc,\ddd,\eee),f$ & Contradiction\\
			\hline\hline
			$\aaa\ddd^2,\ddd\eee^3$& $\aaa\bbb^2,\aaa\bbb\ccc,\bbb^3,\bbb^2\ccc,\bbb\ccc^2$&& \multirow{2}{*}{ No solution } \\
			\cline{1-3}
			\multirow{2}{*}{$\aaa\ddd^2,\eee^4$}& $\aaa\bbb^2,\aaa\bbb\ccc,\bbb^3,\bbb^2\ccc$&&  \\
			\cline{2-4}
			&$\bbb\ccc^2$&$(2,6,3,5,3)/6,24$&$\bbb>\ccc$ and $\ddd>\eee$\\
			\hline
			$\bbb\ddd^2,\ddd\eee^3$& $\aaa^2\bbb,\aaa\bbb^2,\aaa\bbb\ccc$&& \multirow{4}{*}{ No solution } \\
		     \cline{1-3}
			$\bbb\ddd^2,\eee^4$& $\aaa^2\bbb,\aaa\bbb^2,\aaa\bbb\ccc$&&  \\
			\cline{1-3}
			$\ccc\ddd^2,\ddd\eee^3$& $\aaa^2\bbb,\aaa\bbb^2,\aaa\bbb\ccc$&& \\
			\cline{1-3}
			$\ccc\ddd^2,\eee^4$& $\aaa^2\bbb,\aaa\bbb^2,\aaa\bbb\ccc$&& \\
			\hline
	\end{tabular} 
	
\end{table}

There is  only one rational pentagon solution, but it is not simple.

\subsection{Tilings with a unique $a^2b$-vertex type  $\bbb\ddd\eee$}
\begin{table}[htp]        
	\centering 
	\caption{Tilings with a unique  $a^2b$-vertex type $\bbb\ddd\eee$}\label{bde}
	\resizebox{\linewidth}{!}{\begin{tabular}{|c|c|c|c|}
			\hline  
			\multicolumn{2}{|c|}{Vertex}& $(\aaa,\bbb,\ccc,\ddd,\eee),f$ & Contradiction\\
			\hline\hline
			\multirow{2}{*}{$\bbb\ddd\eee,\aaa^2\bbb$}& $\aaa\bbb\ccc,\bbb\ccc^2,\bbb^2\ccc,\ccc^3,\aaa\ccc^3,\bbb\ccc^3,\ccc^4,\ccc^5$&& No solution \\
			\cline{2-4}
			&$\aaa\ccc^2$&$(2,10,6,3,1)/7,28$&$\bbb>\ccc$ and $\ddd>\eee$\\
			\hline
			$\bbb\ddd\eee,\aaa^3$& $\aaa\bbb^2$&& \multirow{4}{*}{ No solution } \\
			\cline{1-3}
			\multirow{9}{*}{$\bbb\ddd\eee,\aaa^2\ccc$}& $ \aaa\bbb\ccc, \bbb^2\ccc, \bbb^3,  \aaa^2\bbb^2, \aaa\bbb^3, \aaa\bbb^2\ccc, \bbb\ccc\ddd^2,\ccc^2\ddd^2$&&\\
			&$\aaa^2\bbb^3, \aaa\bbb^4,  \aaa\bbb^2\ccc^2, \aaa\bbb^2\ddd^2, \aaa\bbb^2\eee^2, \aaa\bbb\ccc^3,\aaa\bbb\ccc\ddd^2$&&\\
			&$\aaa\bbb\ccc\eee^2,\aaa\ccc^2\ddd^2, \aaa\ccc^2\ddd\eee, \aaa\ccc^2\eee^2,  \bbb^2\ccc\ddd^2,  \ccc^3\ddd^2$&&\\
			\cline{2-4}
			&$\aaa\ccc^2,\ccc^3$&$(4,3,4,6,3)/6,12$& Balance
Lemma\\
			\cline{2-4}
			&\multirow{2}{*}{$\bbb\ccc^2,\aaa\ccc^3,\ccc^3\ddd\eee$}&$(8,12,4,7,1)/10,20$&$\bbb>\ccc$ and $\ddd>\eee$			 \\
			\cline{3-4}
			&&$(8,12,4,1,7)/10,20$&Balance
Lemma\\
			\cline{2-4}
			&$\aaa\bbb^2,\aaa\bbb\ccc^2,\bbb^2\ddd^2,\ccc^2\ddd\eee,\aaa\ccc^4,\bbb\ccc^2\ddd^2$&$(6,4,2,3,7)/7,28$&AAD\\
			\cline{2-4}
			&$\aaa\bbb^3\ccc,\bbb^3\ddd^2$&$(18,6,4,11,23)/20,40$&Balance
Lemma\\
			\cline{2-4}
			&$\aaa\bbb^3\ccc$&$(18,6,4,31,3)/20,40$&$\bbb>\ccc$ and $\ddd>\eee$\\
			\hline
			\multirow{5}{*}{$\bbb\ddd\eee,\aaa\bbb\ccc$}& $\aaa^3, \aaa^4,\aaa^3\bbb, \aaa^3\ccc, \aaa^2\ddd\eee, \aaa^2\ddd^2, \aaa\ccc\ddd^2, \ddd^3\eee,\aaa^2\bbb\ddd^2$&&\multirow{2}{*}{ No solution }\\
			&$\aaa^5,\aaa^4\bbb,\aaa^3\bbb^2,\aaa^3\ddd^2, \aaa^2\ccc\eee^2, \aaa\ccc^2\ddd^2,\ccc\ddd^3\eee,\aaa^2\bbb\eee^2$&&\\
			\cline{2-4}  			
			&$\aaa^4\ccc,\aaa^3\ddd\eee,\aaa^2\ccc\ddd^2,\aaa\ddd^3\eee,\ccc\ddd^4$&$(1,3,4,1,4)/4,16$&\multirow{3}{*}{$\bbb<\ccc$ and $\ddd<\eee$}\\
			\cline{2-3}  			
			&$\aaa^2\eee^2,\aaa^4\ccc,\aaa^3\ddd\eee$&$(1,3,4,2,3)/4,16$&\\
			\cline{2-3}
			&$\aaa^3\eee^2$&$(2,8,10,5,7)/10,20$&\\
			\hline
			\multirow{16}{*}{$\bbb\ddd\eee,\aaa\ccc^2$}& $ \aaa\ccc\ddd\eee, \aaa\ddd^2\eee^2, \aaa\ccc\eee^2, \ddd^2\eee^2,\aaa^3\bbb\ccc,  \aaa^2\bbb^2\ccc$&&\multirow{2}{*}{ No solution }\\
			&$\aaa\bbb\ccc\eee^2,\aaa^2\ccc\ddd^2, \aaa\ddd\eee^3, \ccc\ddd^3\eee,\ccc\ddd\eee^3,\aaa\ddd^3\eee$&&\\
			\cline{2-4}
		       &$\aaa^3$&$(4,3,4,6,3)/6,12$&Balance
Lemma\\
                       \cline{2-4}
			&$\aaa\ccc\ddd^2,\aaa^4\ccc,\aaa^2\ccc\ddd\eee,\ccc\ddd^2\eee^2$&$(2,10,6,3,1)/7,28$& $\bbb>\ccc$ and $\ddd>\eee$\\
			\cline{2-4}
			&\multirow{2}{*}{$\aaa\ccc\ddd^2$}&$(8,6,10,5,17)/14,14$			& \multirow{2}{*}{$\bbb<\ccc$ and $\ddd<\eee$}  \\
			\cline{3-3}
			&&$(8,32,38,19,33)/42,42$&\\
			\cline{2-4}
			&$\ddd\eee^3$&$(8,6,10,19,3)/14,14$	&  Balance
Lemma \\
			\cline{2-4}
			&$\aaa^2\bbb\ccc,\ddd\eee^3,\aaa^2\ccc\eee^2$&$(4,2,6,13,1)/8,16$&Balance
Lemma\\
			\cline{2-4}
			&\multirow{2}{*}{$\aaa^3\ccc,\aaa^2\ccc\ddd\eee$}&$(12,48,24,1,11)/30,20$&	Balance
Lemma	\\
			\cline{3-4}
			&&$(12,48,24,11,1)/30,20$&$\bbb>\ccc$ and $\ddd>\eee$\\
			\cline{2-4}
			&$\aaa^3\ccc$&$(2,2,4,3,5)/5,20$& \multirow{2}{*}{ $\bbb<\ccc$ and $\ddd<\eee$ }  \\
                          \cline{2-3}
                          &$\aaa^3\ccc,\aaa^2\bbb\ccc$&$(4,4,8,7,9)/10,20$&\\
                          \cline{2-4}
			&\multirow{3}{*}{$\aaa^2\bbb\ccc$}&$(8,2,10,23,3)/14,14$& Balance
Lemma\\
			\cline{3-4}
			&&$(8,6,14,29,1)/18,18$& Balance
Lemma\\
			\cline{3-4}
			&&$(16,2f-24,2f-8,f+4,f+20)/2f$&\multirow{2}{*}{ $\bbb<\ccc$ and $\ddd<\eee$ }  \\
		    \cline{2-3}
			&$\aaa^2\bbb\ccc,\aaa^2\ddd\eee,\aaa^4\ccc$&$(2,4,6,4,6)/7,28$&\\
			\cline{2-4}
			&$\aaa^2\bbb\ccc,\aaa^3\ddd\eee$&$(2,6,8,5,7)/9,36$&$\bbb<\ccc$ and $\ddd<\eee$\\
			\hline
	\end{tabular} }
\end{table}

According to  \cite{slw} Tables $12$,$14$,$21$, and $23$, we obtain $92$ distinct cases. These results are summarized in  Table \ref{bde}. 

There are only ten simple rational $a^4b$-pentagons to consider:

For Case $\{\bbb\ddd\eee,\aaa^2\ccc,\bbb\ccc^2(\text{or}\,\aaa\ccc^3,  \ccc^3\ddd\eee),(8,12,4,1,7)/10,20\}$, by the balance lemma and the parity lemma, we obtain $\text{AVC}\sub\{\bbb\ddd\eee,\aaa^2\ccc,\bbb\ccc^2,\\
\ccc^3\ddd\eee,\aaa\ccc^3,
\aaa\ccc\ddd\eee,\ccc\ddd^2\eee^2,\ccc^5\}$. There is no solution satisfying the balance lemma.

For Case $\{\bbb\ddd\eee,\aaa^2\ccc,\aaa\bbb^2(\text{or}\,  \aaa\bbb\ccc^2,\bbb^2\ddd^2,\ccc^2\ddd\eee,\aaa\ccc^4,\bbb\ccc^2\ddd^2)$,
$(6,4,2,3,7)\\
/7,28\}$, the balance lemma and the parity lemma imply \[\text{AVC}\sub\{\bbb\ddd\eee,\aaa^2\ccc,\aaa\bbb^2,  \aaa\bbb\ccc^2,\ccc^2\ddd\eee,\bbb^3\ccc,\aaa\ccc^4,\bbb^2\ccc^3,\bbb\ccc^5,\ccc^7\}.\] If $\ccc^7$ is a vertex, then the AAD at $\ccc^7$ must be either $\ccc^{\eee}\thin^{\aaa}\ccc\cdots$ or $\ccc^{\eee}\thin^{\eee}\ccc\cdots$. This results in a vertex $\aaa\eee\cdots$ or $\eee\thin\eee\cdots$, which is not in the AVC. Therefore, $\ccc^7$ is not a vertex. Similarly, $\aaa\ccc^4$ and $\bbb\ccc^5$ are also excluded as vertices. If $\bbb^3\ccc$ is a vertex,  the AAD at $\ccc^7$ must be one of the following: $\bbb^{\ddd}\thin^{\aaa}\bbb\cdots$, $\bbb^{\ddd}\thin^{\ddd}\bbb\cdots$, $\bbb^{\ddd}\thin^{\aaa}\ccc\cdots$ or $\bbb^{\ddd}\thin^{\eee}\ccc\cdots$. This leads to a vertex $\aaa\ddd\cdots$, $\ddd\thin\ddd\cdots$ or $\ddd\thin\eee\cdots$, none of which are in the AVC. Therefore, $\bbb^3\ccc$ is not  a vertex. Similarly, $\aaa\bbb\ccc^2$ and $\bbb^2\ccc^3$ are also excluded. Thus, the AVC reduces to
$\text{AVC}\sub\{\bbb\ddd\eee,\aaa^2\ccc,\aaa\bbb^2,\ccc^2\ddd\eee\}$.
 If $\ccc^2\ddd\eee$ is a vertex, the unique AAD configuration is $\ccc^2\ddd\eee=\thin^{\eee}\ccc^{\aaa}\thin^{\aaa}\ccc^{\eee}\thin^{\bbb}\ddd^{\eee}\thin^{\ddd}\eee^{\ccc}\thin$. This results in  a vertex $\aaa^{\ccc}\thin^{\ccc}\aaa\cdots=\thin^{\bbb}\aaa^{\ccc}\thin^{\ccc}\aaa^{\bbb}\thin^{\aaa}\ccc^{\eee}\thin$. The corresponding AAD is $\aaa^{\ccc}\thin^{\aaa}\bbb\cdots=\thin^{\bbb}\aaa^{\ccc}\thin^{\aaa}\bbb^{\ddd}\thin\bbb\thin$, which leads to  a vertex $\aaa\ddd\cdots$ or $\ddd\thin\ddd\cdots$,  neither of which are in the  AVC. Therefore, $\ccc^2\ddd\eee$  is not a vertex, and the AVC further reduces to $\text{AVC}\sub\{\bbb\ddd\eee,\aaa^2\ccc,\aaa\bbb^2\}$. No solution satisfying the balance lemma exists in this case.

For Case $\{\bbb\ddd\eee,\aaa^2\ccc,\aaa\bbb^3\ccc(\text{or}\,\bbb^3\ddd^2),(18,6,4,11,23)/20,40\}$,  the parity lemma implies that $\eee\cdots=\ddd\eee\cdots=\bbb\ddd\eee$. Consequently, we have $\#\bbb\ge\#\bbb\ddd\eee+3\#\aaa\bbb^3\ccc=f+3\#\aaa\bbb^3\ccc>f$ or $\#\ddd\ge\#\bbb\ddd\eee+2\#\bbb^3\ddd^2=f+2\#\bbb^3\ddd^2>f$, a contradiction. 

Similarly, the following cases are also dismissed for tilings:
	\begin{enumerate}
			 \item Case $\{\bbb\ddd\eee,\aaa^2\ccc,\aaa\ccc^2(\text{or}\,\ccc^3),(4,3,4,6,3)/6,12\}$.
		  \item Case $\{\bbb\ddd\eee,\aaa\ccc^2,\aaa^3,(4,3,4,6,3)/6,12\}$.
		  \item Case $\{\bbb\ddd\eee,\aaa\ccc^2,\aaa^3\ccc(\text{or}\,\aaa^2\ccc\ddd\eee),(12,48,24,1,11)/30,20\}$.
		  \item Case $\{\bbb\ddd\eee,\aaa\ccc^2,\ddd\eee^3, (8,6,10,19,3)/14,14\}$.
	      \item Case $\{\bbb\ddd\eee,\aaa\ccc^2,\aaa^2\bbb\ccc(\text{or}\,\ddd\eee^3,\aaa^2\ccc\eee^2),(4,2,6,13,1)/8,16\}$.
	      \item Case $\{\bbb\ddd\eee,\aaa\ccc^2,\aaa^2\bbb\ccc,(8,2,10,23,3)/14,14\}$.
	      \item Case $\{\bbb\ddd\eee,\aaa\ccc^2,\aaa^2\bbb\ccc,(8,6,14,29,1)/18,18\}$.

		\end{enumerate}

\subsection{Tilings with a unique $a^2b$-vertex type $\aaa\ddd\eee$}

According to \cite{slw} Tables $28$ and $30$, we obtain  $28$ distinct cases. These results are summarized in  Tables \ref{ade} and \ref{ade 2bc}.

\begin{table}[htp]        
	\centering 
	\caption{Tilings with a unique  $a^2b$-vertex type $\aaa\ddd\eee$}\label{ade}
	\resizebox{\linewidth}{!}{\begin{tabular}{|c|c|c|c|}
			\hline  
			\multicolumn{2}{|c|}{Vertex}& $(\aaa,\bbb,\ccc,\ddd,\eee),f$ & Contradiction\\
			\hline\hline
			\multirow{2}{*}{$\aaa\ddd\eee,\aaa^3$}& $\aaa^2\bbb,\aaa\bbb^2,\bbb^3,\aaa\bbb^3,\bbb^4$&& No solution \\
			\cline{2-4}
			&$\bbb^3\ccc,\bbb^5$&$(10,6,12,15,5)/15,20$&\\
			\hline
			$\aaa\ddd\eee,\aaa^2\ccc$& $\aaa\bbb^2,\aaa\bbb\ccc,\bbb^3,\aaa\bbb^3,\bbb^4,\bbb^3\ccc,\bbb^5$&&\multirow{2}{*}{ No solution }\\
			\cline{1-3}
			\multirow{3}{*}{$\aaa\ddd\eee,\aaa\ccc^2$}& $\aaa^2\bbb,\bbb^3,\aaa\bbb^3,\bbb^4$&& \\
			\cline{2-4}
			&\multirow{2}{*}{$\bbb^3\ccc,\bbb^5$}&$(2,2,4,3,5)/5,20$&\multirow{2}{*}{$\bbb<\ccc$ and $\ddd<\eee$}\\
			\cline{3-3}
			&&$(4,4,8,7,9)/10,20$&\\
			\hline
			$\aaa\ddd\eee,\ccc^3$& $\aaa^2\bbb,\aaa\bbb\ccc,\aaa\bbb^2,\aaa\bbb^3,\bbb^4,\bbb^3\ccc,\bbb^5$&& No solution \\
			\hline
	\end{tabular} }
	
\end{table}

\begin{table}[htp]        
 \centering 
 \caption{$\aaa\ddd\eee$, $\bbb^2\ccc$ appear as vertices}\label{ade 2bc}
 \begin{tabular}{|c|c|c|}
   \hline  
   Vertex& $(\aaa,\bbb,\ccc,\ddd,\eee),f$ & Contradiction\\
   \hline\hline
   \multirow{4}{*}{$\aaa\ddd\eee,\bbb^2\ccc$}&$(2,4,2,5,3)/5,20$&\multirow{2}{*}{$\bbb>\ccc$ and $\ddd>\eee$}  \\
   \cline{2-2}
   &$(4,8,4,9,7)/10,20$& \\   
   \cline{2-3}
   &$(16,8,4,1,3)/10,20$&self-intersecting \\
   \cline{2-3}
   &$(10,12,6,5,15)/15,20$&\\      
   \hline  
 \end{tabular}
 
\end{table}
See Example \ref{ex:first}  in Section \ref{region} for detailed discussion of four solutions in Table \ref{ade 2bc}. The second row of Table \ref{ade} and the last row of Table \ref{ade 2bc} give the same rational non-symmetric pentagon, which is the only one admitting tilings (see the third class in the main theorem).

\section{Conclusion}\label{conclusion}
Combining results from \cite{wy3, wy1, wy2, slw}, we complete the full classification of edge-to-edge tilings of the sphere by congruent pentagons. 
We summarize systematically the classification for all types: $a^2b^2c$, $a^3bc$, $a^3b^2$, $a^4b$, and $a^5$ (see  Figure \ref{pentagon}), which 
comprises four distinct classes:

\begin{figure}[htp]
\centering
\begin{tikzpicture}[>=latex,scale=0.9]
\begin{scope}[xshift=-2.5cm]

\draw
	(0,0.4) -- node[above=-2] {\small $a$} ++(1,0);

\draw[line width=1.5]
	(0,-0.1) -- node[above=-2] {\small $b$} ++(1,0);

\draw[dashed]
	(0,-0.6) -- node[above=-2] {\small $c$} ++(1,0);

\end{scope}

%% a^2b^2c

\draw
	(234:1) -- (162:1) -- (90:1);

\draw[line width=1.5]
	(-54:1) -- (18:1) -- (90:1);

\draw[densely dashed]
	(234:1) -- (-54:1);

\node at (90:0.75) {$\alpha$};
\node at (162:0.75) {$\beta$};
\node at (15:0.75) {$\gamma$};
\node at (234:0.75) {$\delta$};
\node at (-54:0.75) {$\epsilon$};

\node at (0,0) {\small $a^2b^2c$};
%% a^3bc

\begin{scope}[xshift=2.5cm]

\draw
	(162:1) -- (234:1) -- (-54:1) -- (18:1);

\draw[line width=1.5]
	(162:1) -- (90:1);

\draw[densely dashed]
	(18:1) -- (90:1);

\node at (90:0.75) {$\alpha$};
\node at (162:0.75) {$\beta$};
\node at (15:0.75) {$\gamma$};
\node at (234:0.75) {$\delta$};
\node at (-54:0.75) {$\epsilon$};
\node at (0,0) {\small $a^3bc$};
\end{scope}

%% a^3b^2

\begin{scope}[xshift=5cm]

\draw
	(162:1) -- (234:1) -- (-54:1) -- (18:1);

\draw[line width=1.5]
	(162:1) -- (90:1);

\draw[line width=1.5]
	(18:1) -- (90:1);

\node at (90:0.75) {$\alpha$};
\node at (162:0.75) {$\beta$};
\node at (15:0.75) {$\gamma$};
\node at (234:0.75) {$\delta$};
\node at (-54:0.75) {$\epsilon$};
\node at (0,0) {\small $a^3b^2$};
\end{scope}

%% a^4b

\begin{scope}[xshift=7.5cm]

\draw
	(162:1) -- (234:1) -- (-54:1) -- (18:1) -- (90:1);

\draw[line width=1.5]
	(234:1) -- (-54:1);
\draw
	(162:1) -- (90:1);

\node at (90:0.75) {$\alpha$};
\node at (162:0.75) {$\beta$};
\node at (15:0.75) {$\gamma$};
\node at (234:0.75) {$\delta$};
\node at (-54:0.75) {$\epsilon$};
\node at (0,0) {\small $a^4b$};
\end{scope}

%% a^5

\begin{scope}[xshift=10cm]

\draw
	(162:1) -- (234:1) -- (-54:1) -- (18:1) -- (90:1);

\draw
	(234:1) -- (-54:1);
\draw
	(162:1) -- (90:1);

\node at (90:0.75) {$\alpha$};
\node at (162:0.75) {$\beta$};
\node at (15:0.75) {$\gamma$};
\node at (234:0.75) {$\delta$};
\node at (-54:0.75) {$\epsilon$};
\node at (0,0) {\small $a^5$};
\end{scope}

\end{tikzpicture}
\caption{Edge combinations suitable for tiling, with $a,b,c$ distinct.}
\label{pentagon}
\end{figure}
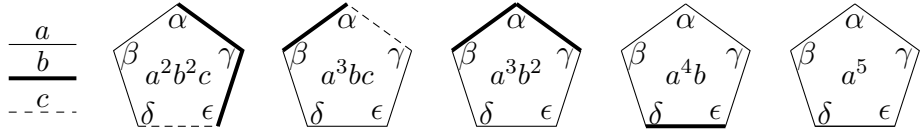

\begin{enumerate}
 \item  Three two-parameter families of pentagonal subdivisions of the Platonic solids, with $12$, $24$ and $60$ tiles as shown in Table \ref{Tab-summary-2} and Figure \ref{2a2bc}, with the moduli spaces \cite[Figure 22]{wy4} having the reductions from type $a^2b^2c$ to  $a^3b^2$ (see the first five pictures of  Figure \ref{3a2b}), $a^4b$ (see  \cite[Figure 4]{slw}), or $a^5$ (see the first three pictures of  Figure \ref{5a}), among which the following special $a^3b^2$-pentagons admit new flip modifications:

 \begin{itemize}
   \item The case $f=24$, $\ddd = \pi$ admits two modifications in Table \ref{Tab-summary-4} (see the sixth and seventh pictures of Figure \ref{3a2b}).
    \item The case $f=60$, $\ddd = \tfrac{6}{5}\pi$ admits four modifications in Table \ref{Tab-summary-4} (see the eighth to eleventh pictures of  Figure \ref{3a2b}).
     \item The case $f=60$, $\ddd = \tfrac{4}{5}\pi$ admits four modifications in Table \ref{Tab-summary-4} (see  the last four pictures of Figure \ref{3a2b}).
 \end{itemize}

 \begin{table*}[htp]                        
	\centering     
	\caption{All tilings of type $a^2b^2c$. }\label{Tab-summary-2}  
	~\\ 
	\begin{tabular}{c|c|c}	 
	\hline  	
	 $f$&$(\aaa,\bbb,\ccc,\ddd,\eee)$ &Tilings\\
	\hline	
	 $12$&$(\aaa,\frac{2}{3},\frac{2}{3},\ddd,2-\aaa-\ddd)$& $ T(12\aaa\ddd\eee,4\bbb^3,4\ccc^3)$\\
	\hline			
	$24$&$(\aaa,\frac{1}{2},\frac{2}{3},\ddd,2-\aaa-\ddd)$& $ T(24\aaa\ddd\eee,8\bbb^3,6\ccc^4)$\\    
	\hline
	$60$&$(\aaa,\frac{2}{5},\frac{2}{3},\ddd,2-\aaa-\ddd)$& $ T(60\aaa\ddd\eee,20\bbb^3,12\ccc^5)$\\
	\hline
	\end{tabular}
      
\end{table*}

\begin{figure}[htp]
 \centering
 \includegraphics[scale=0.0308]{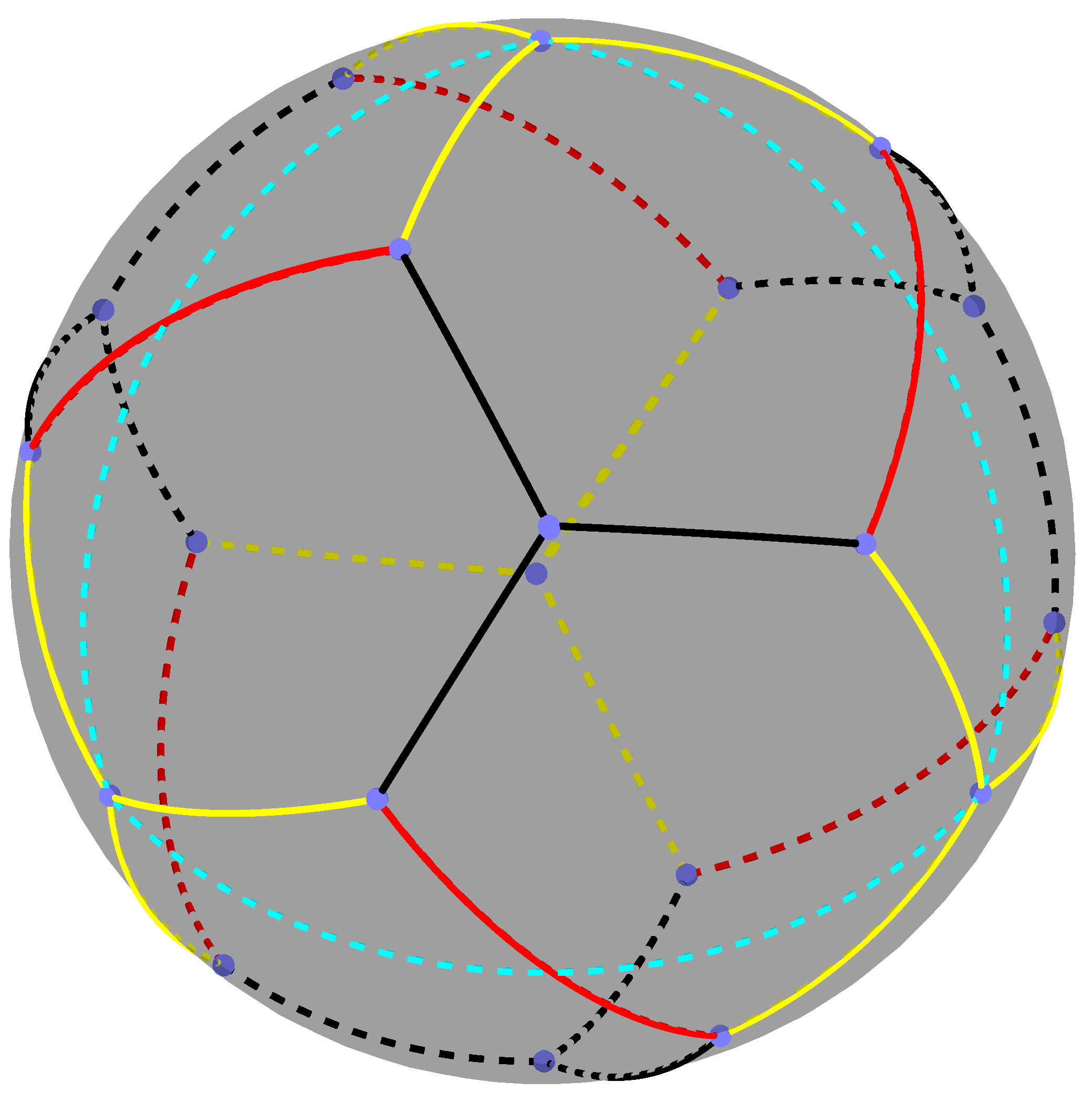}\hspace{10pt}
 \includegraphics[scale=0.0385]{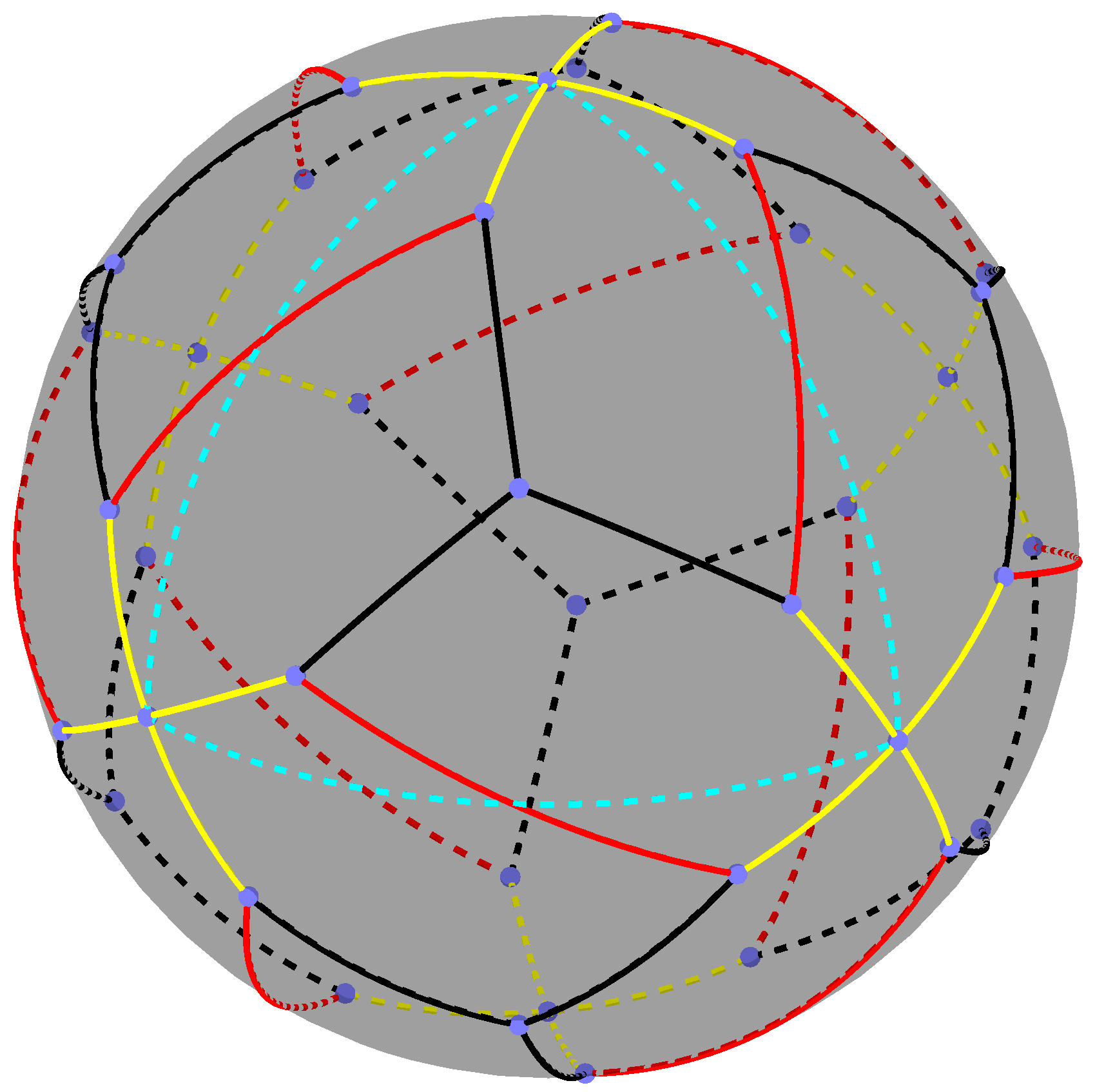}\hspace{10pt}
 \includegraphics[scale=0.0413]{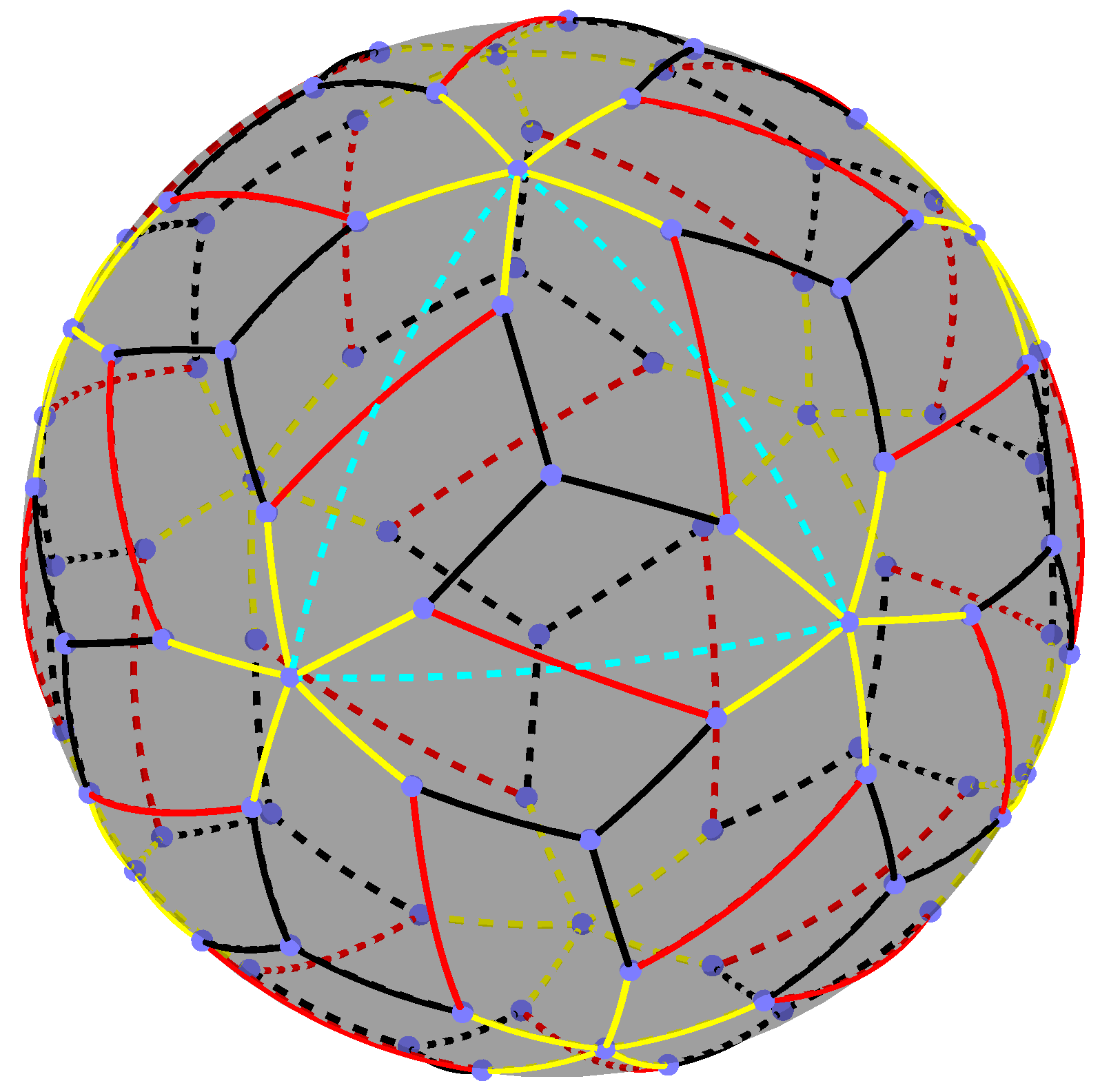} 
 \caption{Two-parameter families of pentagonal subdivisions for $a^2b^2c$.} 
 \label{2a2bc} 
\end{figure}

\begin{table}[htp]                  
	 \begin{center}   
	\caption{All tilings of type $a^3b^2$. }\label{Tab-summary-4}   
	\bgroup   
	   \resizebox{\textwidth}{40mm}{
	\begin{tabular}{c|c|c}	 
	\hline  	
	 $f$&$(\aaa,\bbb,\ccc,\ddd,\eee)$ &Tilings\\
	\hline
	 $12$&$(\frac{2}{3},\bbb,\ccc,2-\bbb-\ccc,\frac{2}{3})$& $ T(4\aaa^3,12\bbb\ccc\ddd,4\eee^3)$\\
	\hline
	$24$&$(\frac{2}{3},\bbb,\ccc,2-\bbb-\ccc,\frac{1}{2})$& $ T(8\aaa^3,24\bbb\ccc\ddd,6\eee^4)$\\
	\hline			
	$24$&$(\frac{1}{2},\bbb,\ccc,\frac{2}{3},2-\bbb-\ccc)$& $ T(6\aaa^4,24\bbb\ccc\eee,8\ddd^3)$\\    
	\hline
	$60$&$(\frac{2}{3},\bbb,\ccc,2-\bbb-\ccc,\frac{2}{5})$& $ T(20\aaa^3,60\bbb\ccc\ddd,12\eee^5)$\\
	\hline
	$60$&$(\frac{2}{5},\bbb,\ccc,\frac{2}{3},2-\bbb-\ccc)$& $ T(12\aaa^5,60\bbb\ccc\eee,20\ddd^3)$\\    
	\hline\hline
	\multirow{2}{*}{$24$}&\multirow{2}{*}{$(\frac{2}{3},\bbb,1-\bbb,1,\frac{1}{2}), \bbb \approx 0.3883$} & $ T(8\aaa^3,20\bbb\ccc\ddd,4\ddd\eee^2,4\bbb\ccc\eee^2,2\eee^4)$\\    
	\cline{3-3}
&& $ T(8\aaa^3,16\bbb\ccc\ddd,4\bbb^2\ccc^2,8\ddd\eee^2,2\eee^4)$\\    
	\hline\hline	
	 \multirow{4}{*}{$60$}& \multirow{4}{*}{$(\frac{2}{3},\bbb,\frac{4}{5}-\bbb,\frac{6}{5},\frac{2}{5}), \bbb \approx 0.2484$}& $ T(20\aaa^3,55\bbb\ccc\ddd,5\ddd\eee^2,5\bbb\ccc\eee^3,7\eee^5)$\\    
	\cline{3-3}
	& & $ T(20\aaa^3,50 \bbb\ccc\ddd, 10\ddd\eee^2,10\bbb\ccc\eee^3, 2\eee^5)$\\
      \cline{3-3}
	&  & $ T(20\aaa^3,50\bbb\ccc\ddd,10\ddd\eee^2,2\bbb^2\ccc^2\eee,6\bbb\ccc\eee^3,4\eee^5)$\\
	\cline{3-3} 
	 & & $ T(20\aaa^3,45\bbb\ccc\ddd,15\ddd\eee^2,6\bbb^2\ccc^2\eee,3\bbb\ccc\eee^3,3\eee^5)$\\
	\hline\hline
	 \multirow{4}{*}{$60$}& \multirow{4}{*}{$(\frac{2}{3},\bbb,\frac{6}{5}-\bbb,\frac{4}{5},\frac{2}{5}), \bbb \approx 0.4688$} & $ T(20\aaa^3,55\bbb\ccc\ddd,5\bbb\ccc\eee^2,5\ddd\eee^3,7\eee^5)$\\
	\cline{3-3}
	& & $ T(20\aaa^3,50\bbb\ccc\ddd,10\bbb\ccc\eee^2,10\ddd\eee^3,2\eee^5)$\\
	\cline{3-3}
	& & $ T(20\aaa^3,50\bbb\ccc\ddd,2\ddd^2\eee,10\bbb\ccc\eee^2,6\ddd\eee^3,4\eee^5)$\\
	\cline{3-3}
	&& $ T(20\aaa^3,45\bbb\ccc\ddd,6\ddd^2\eee,15\bbb\ccc\eee^2,3\ddd\eee^3,3\eee^5)$\\
	\hline
	\end{tabular}}
	  \egroup
  \end{center}
\end{table}

\begin{figure}[htp]
 \centering
 \includegraphics[scale=0.033]{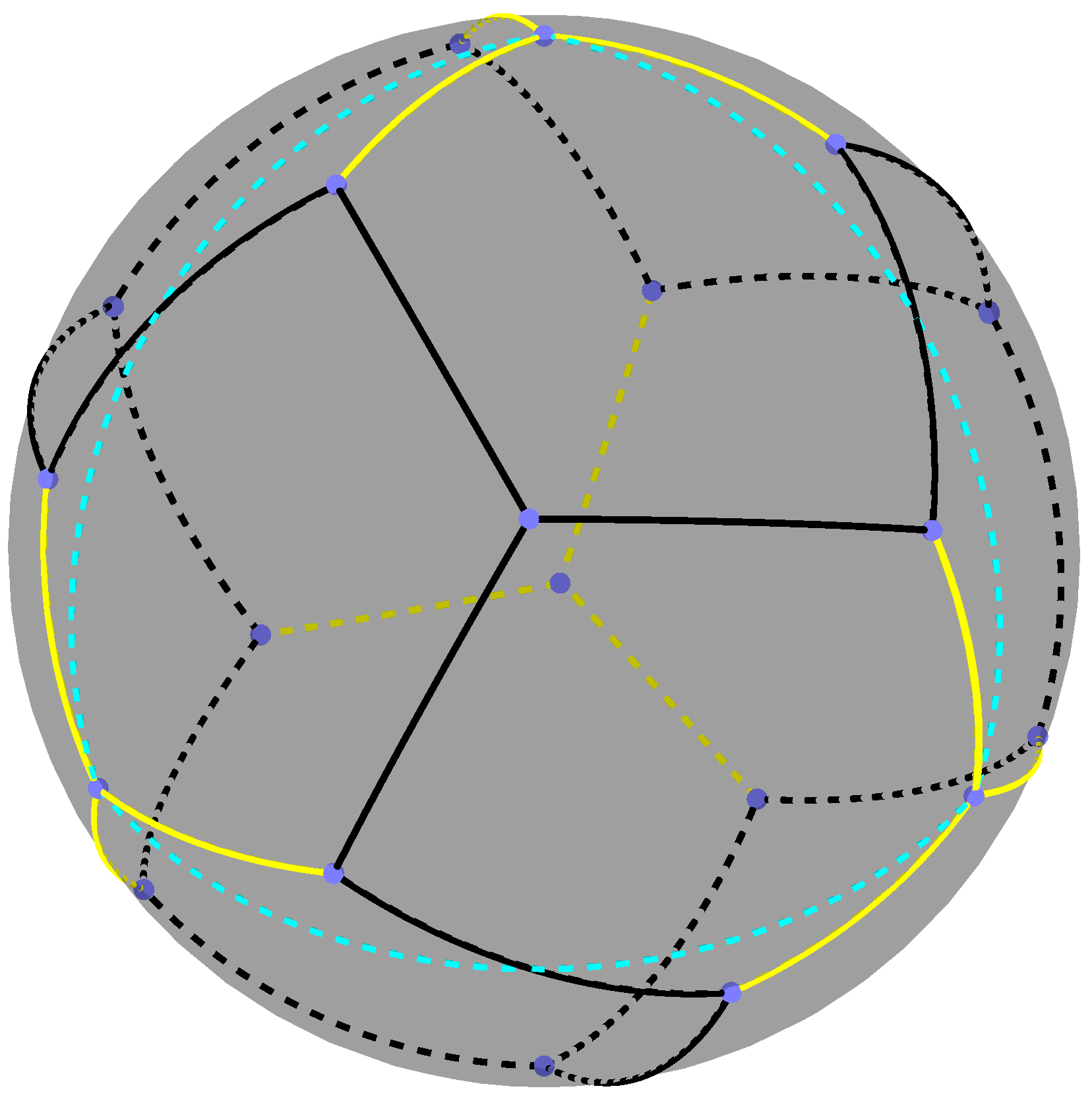}\hspace{10pt}
 \includegraphics[scale=0.033]{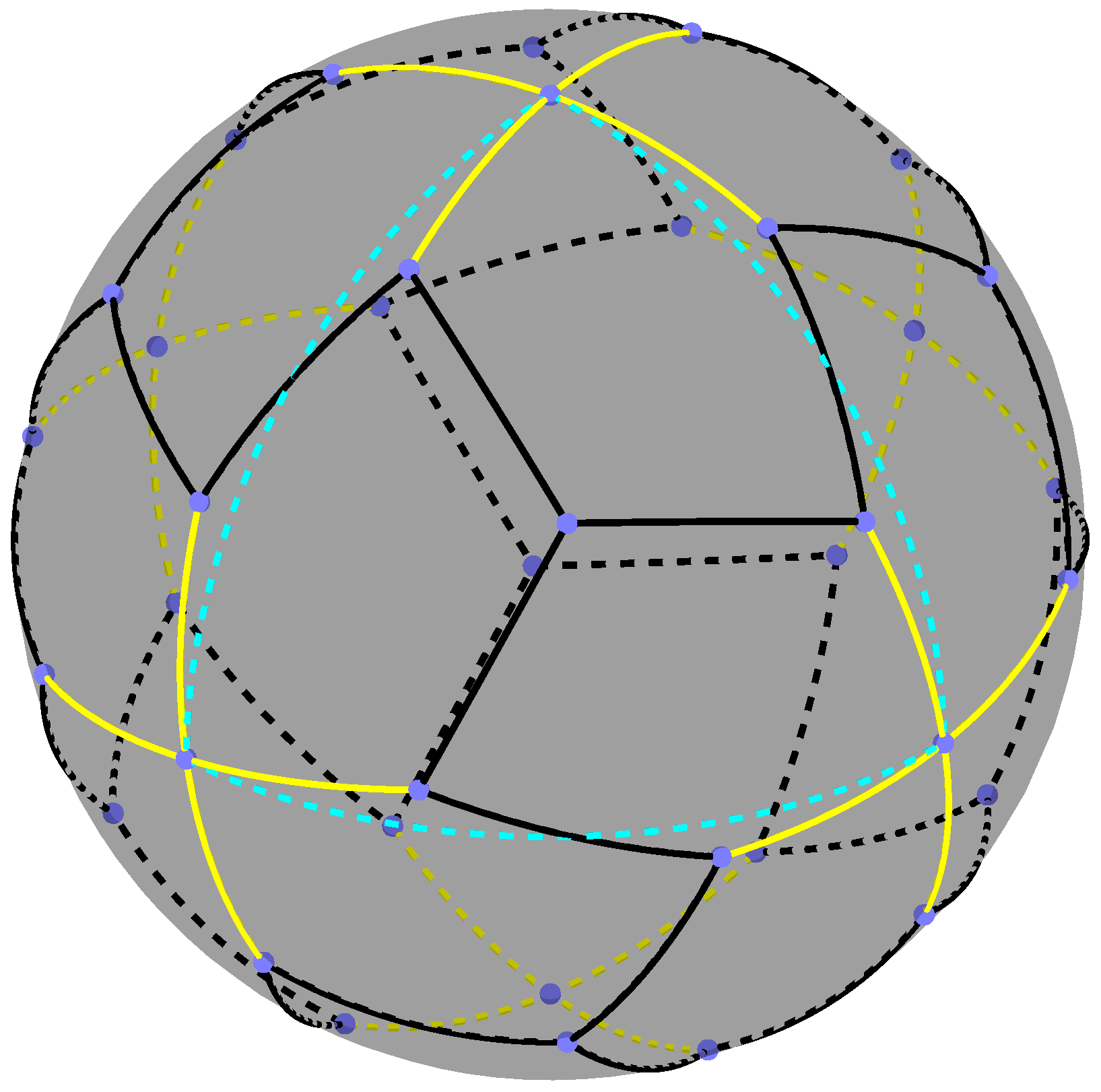}\hspace{14pt}
 \includegraphics[scale=0.033]{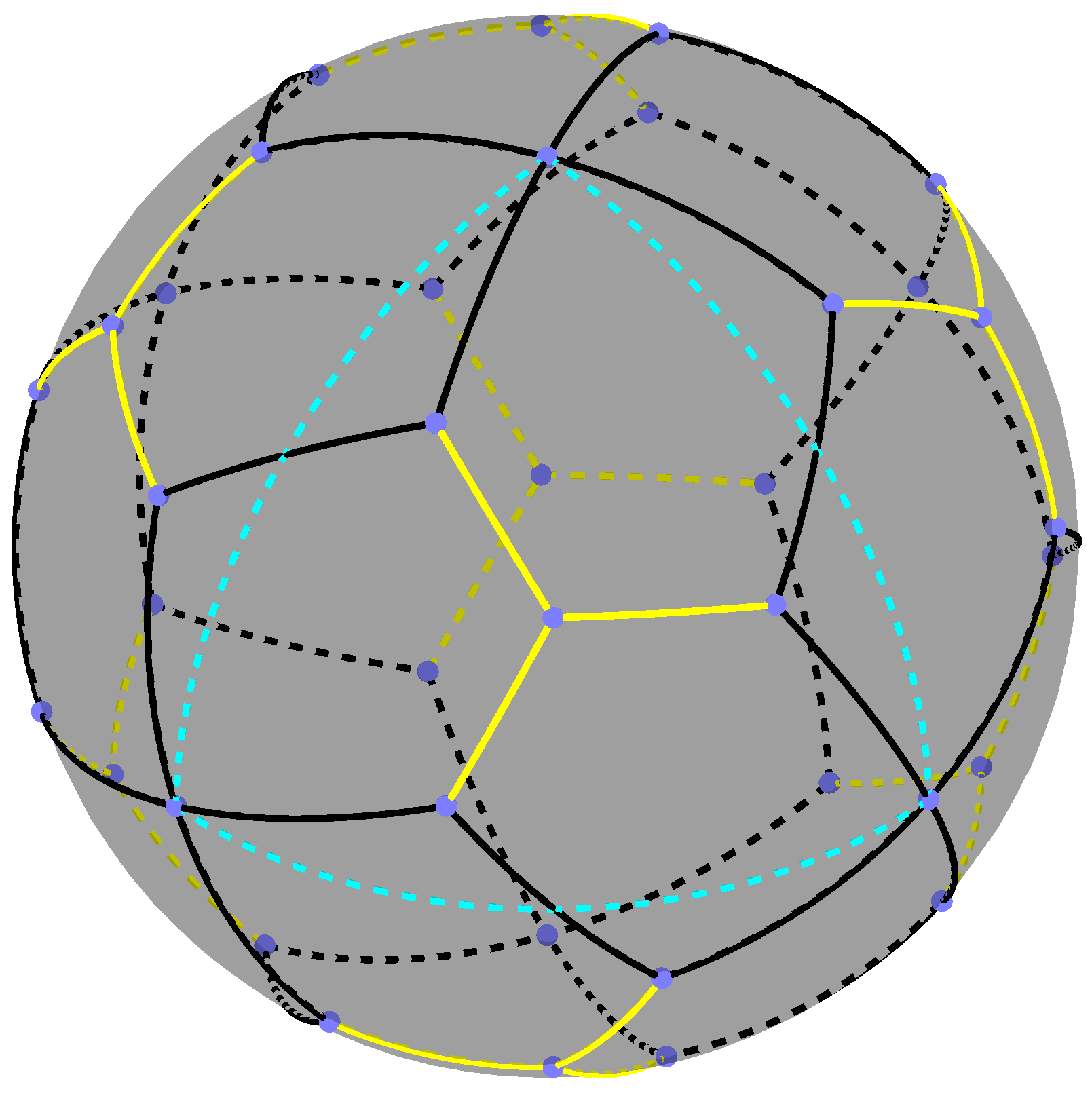}\hspace{10pt}
 \includegraphics[scale=0.033]{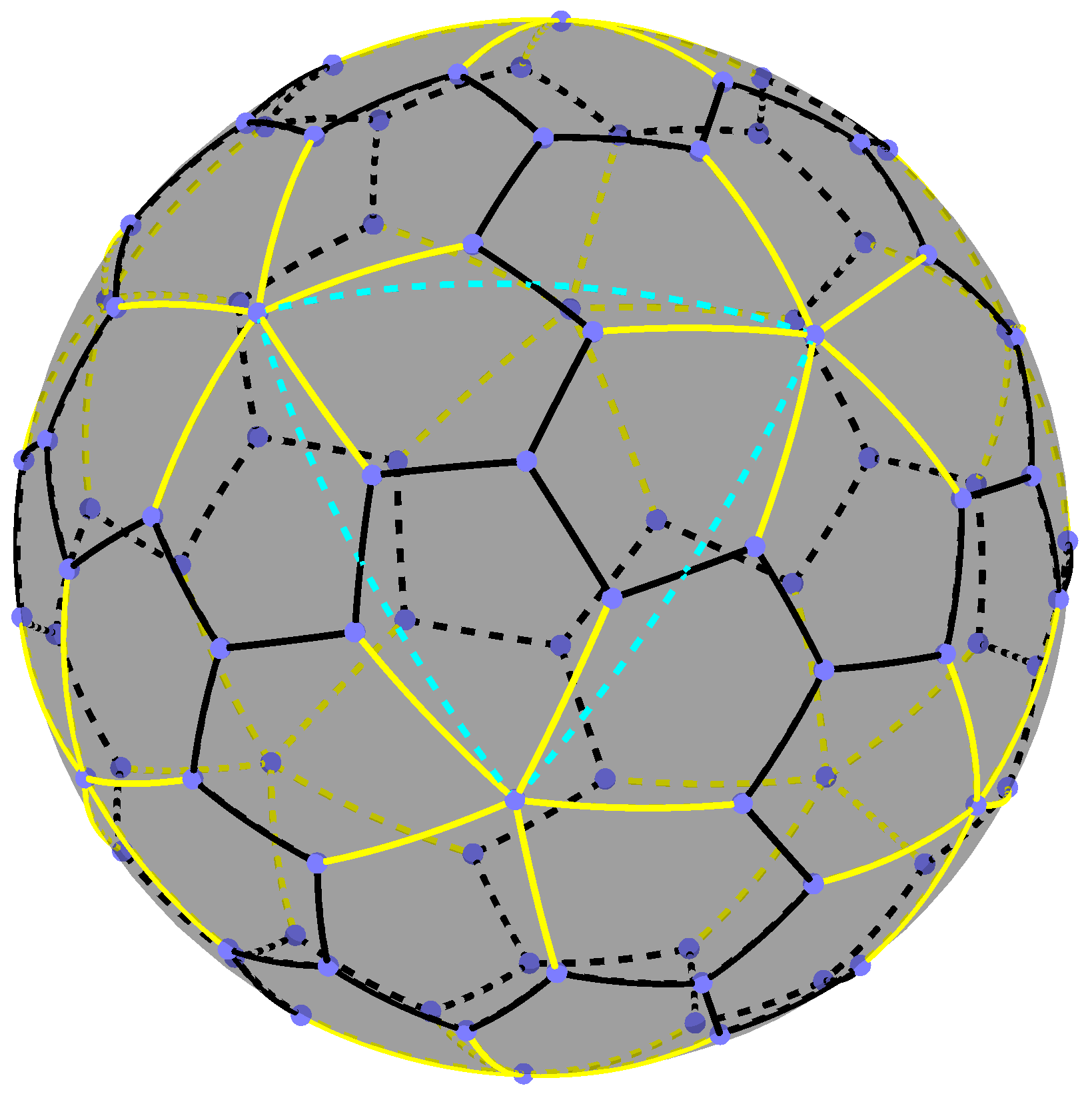}\hspace{10pt}
 \includegraphics[scale=0.033]{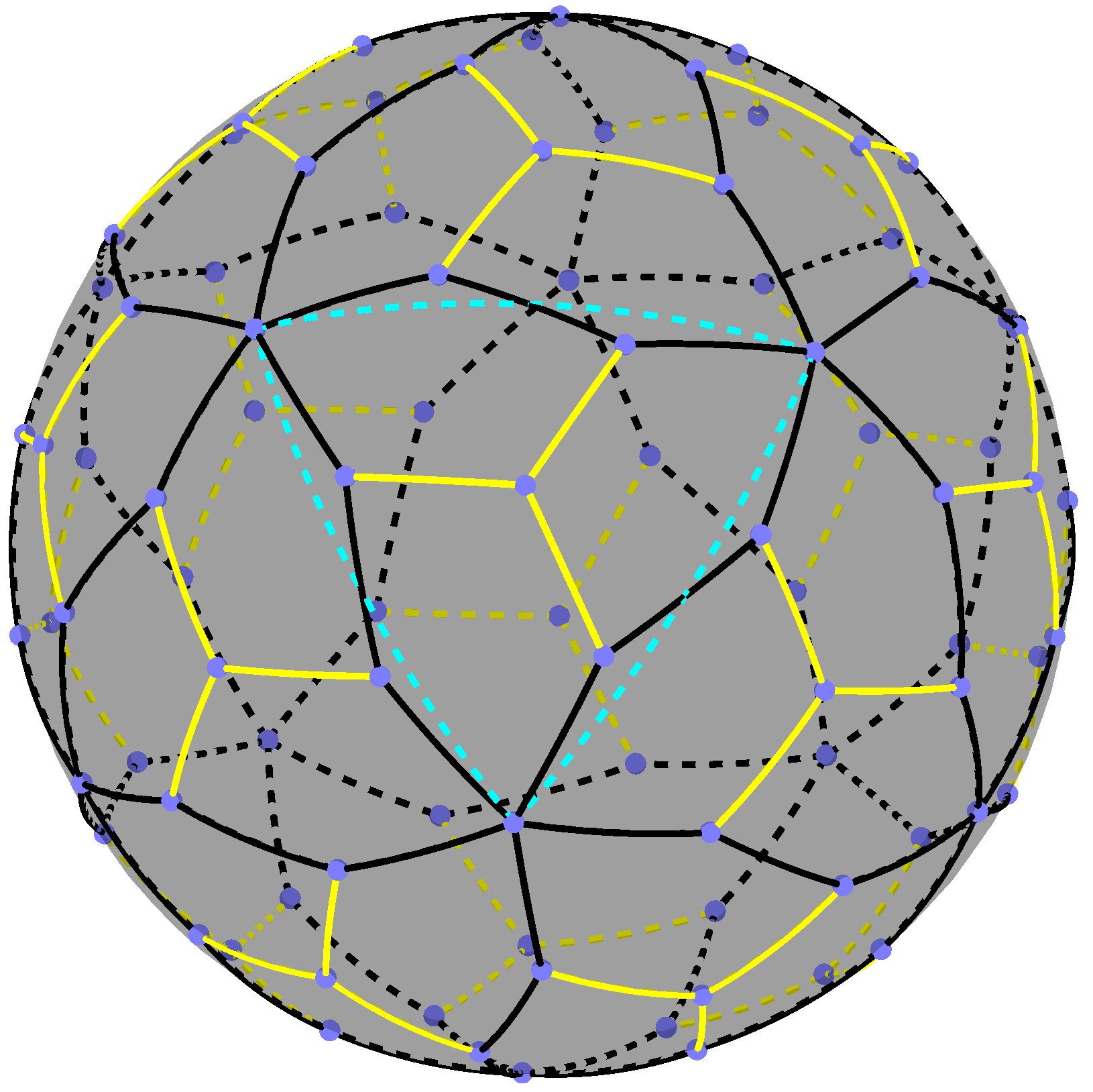}
 
 \includegraphics[scale=0.025]{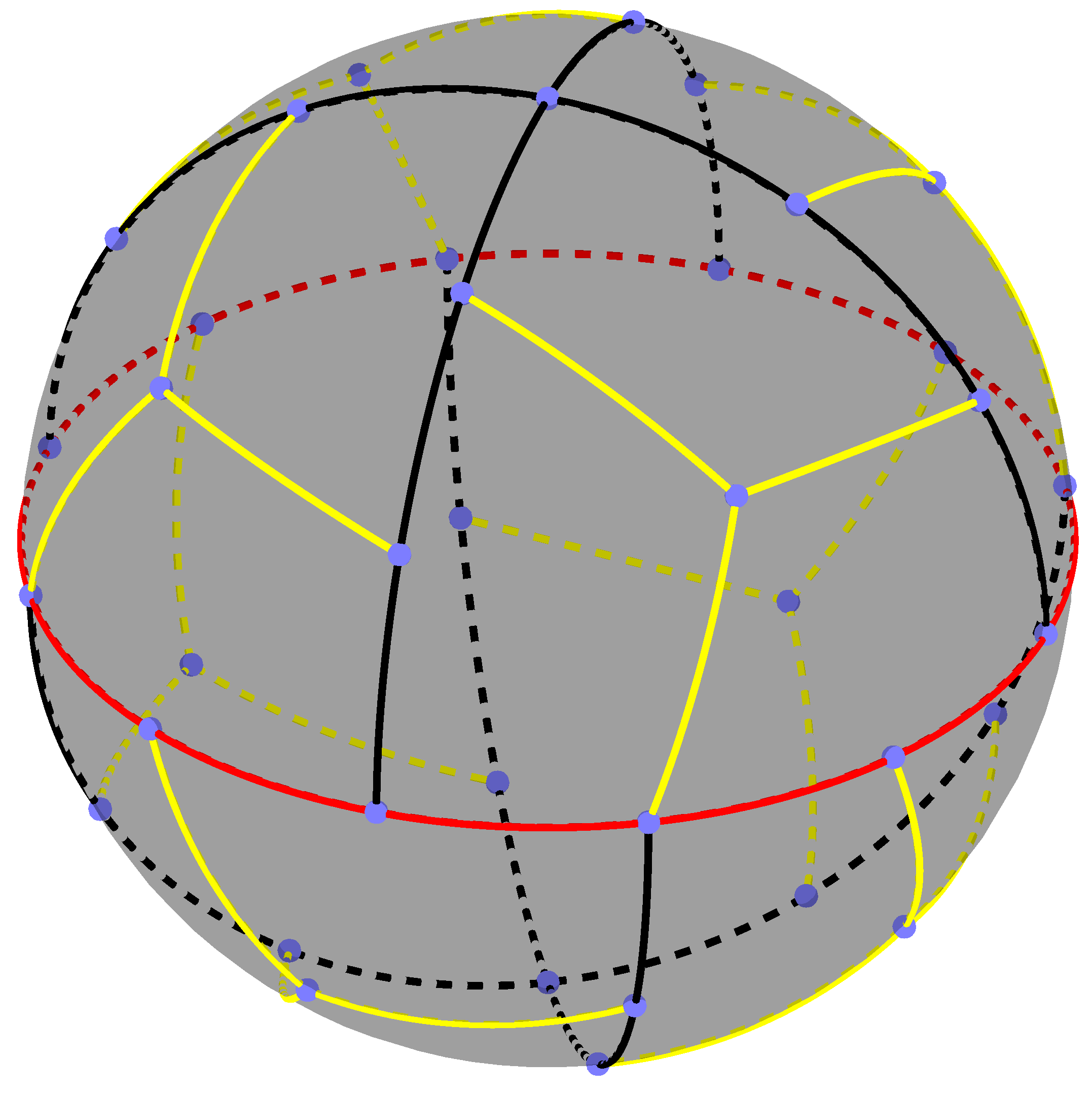}\hspace{10pt}
 \includegraphics[scale=0.040]{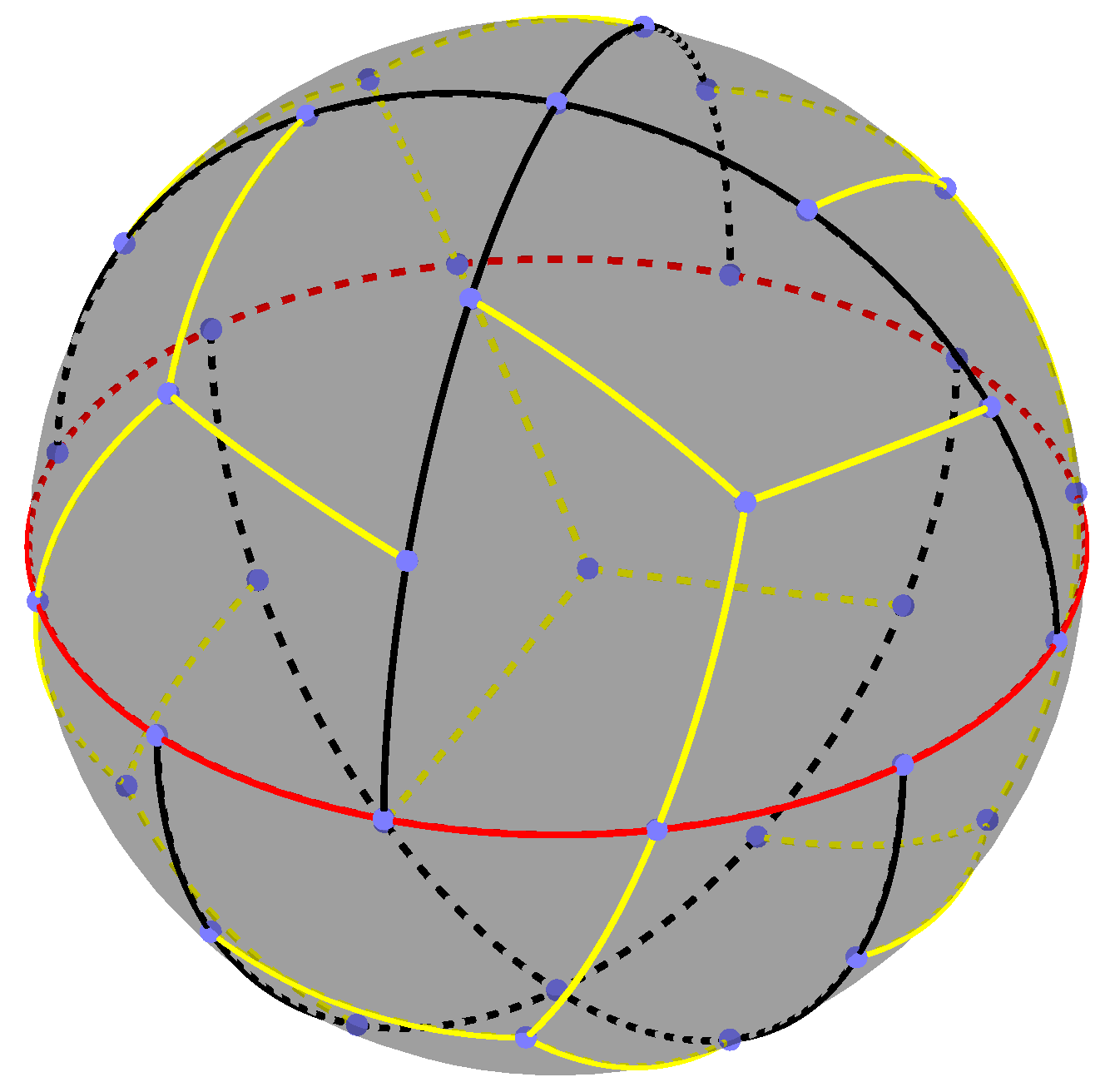}\hspace{17pt}
 \includegraphics[scale=0.033]{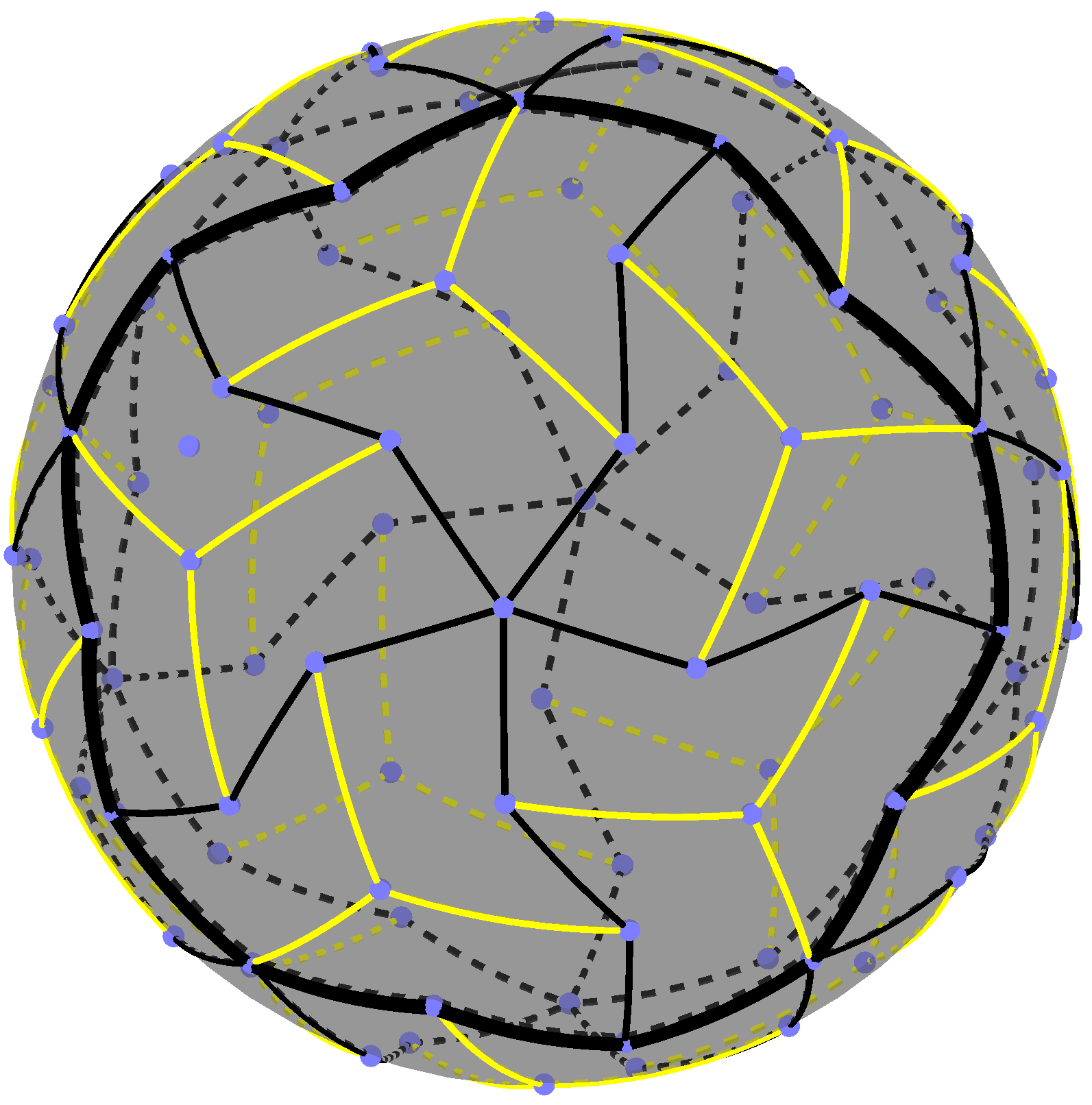}\hspace{10pt}
 \includegraphics[scale=0.0335]{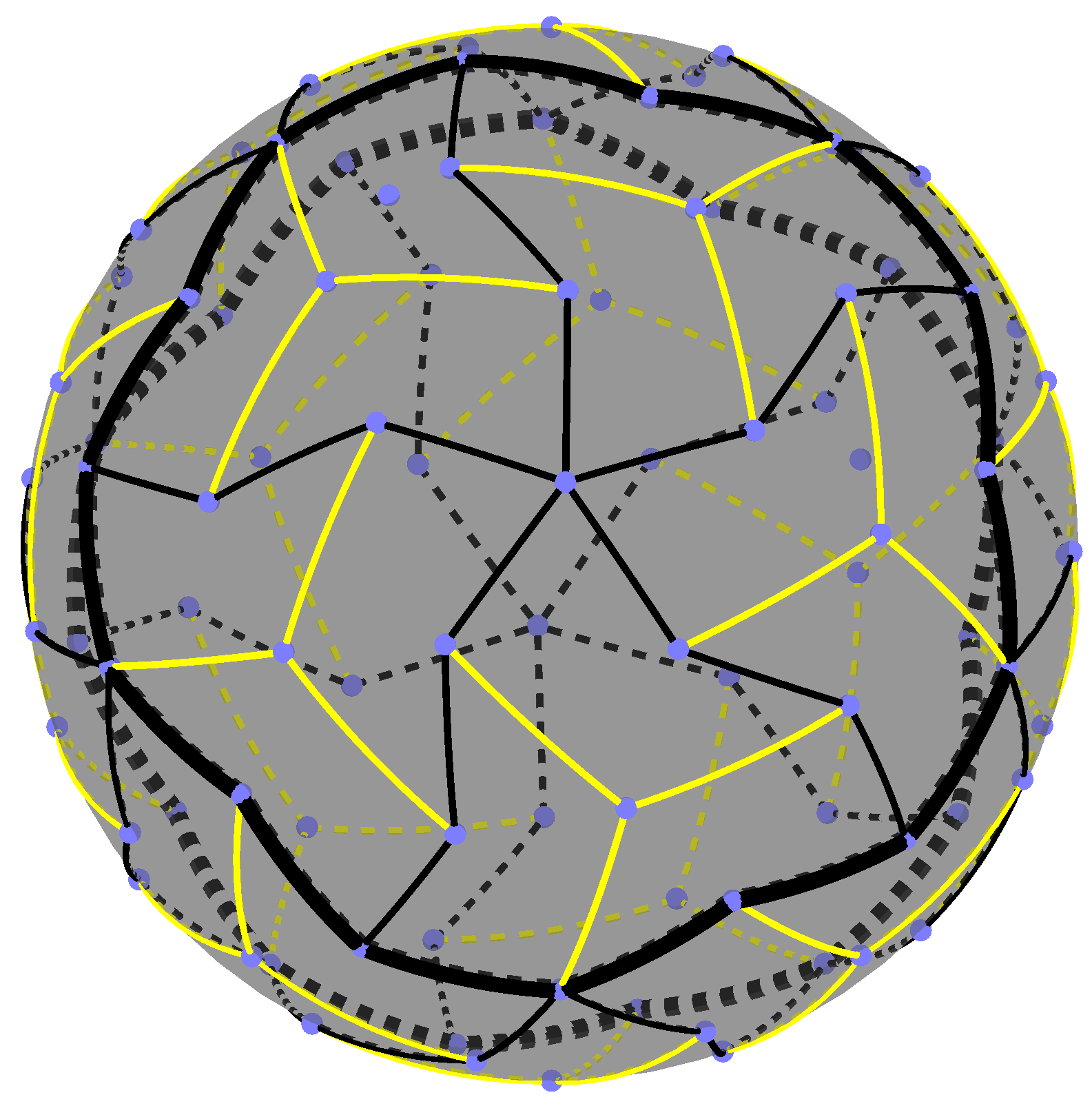}\hspace{10pt}
 \includegraphics[scale=0.033]{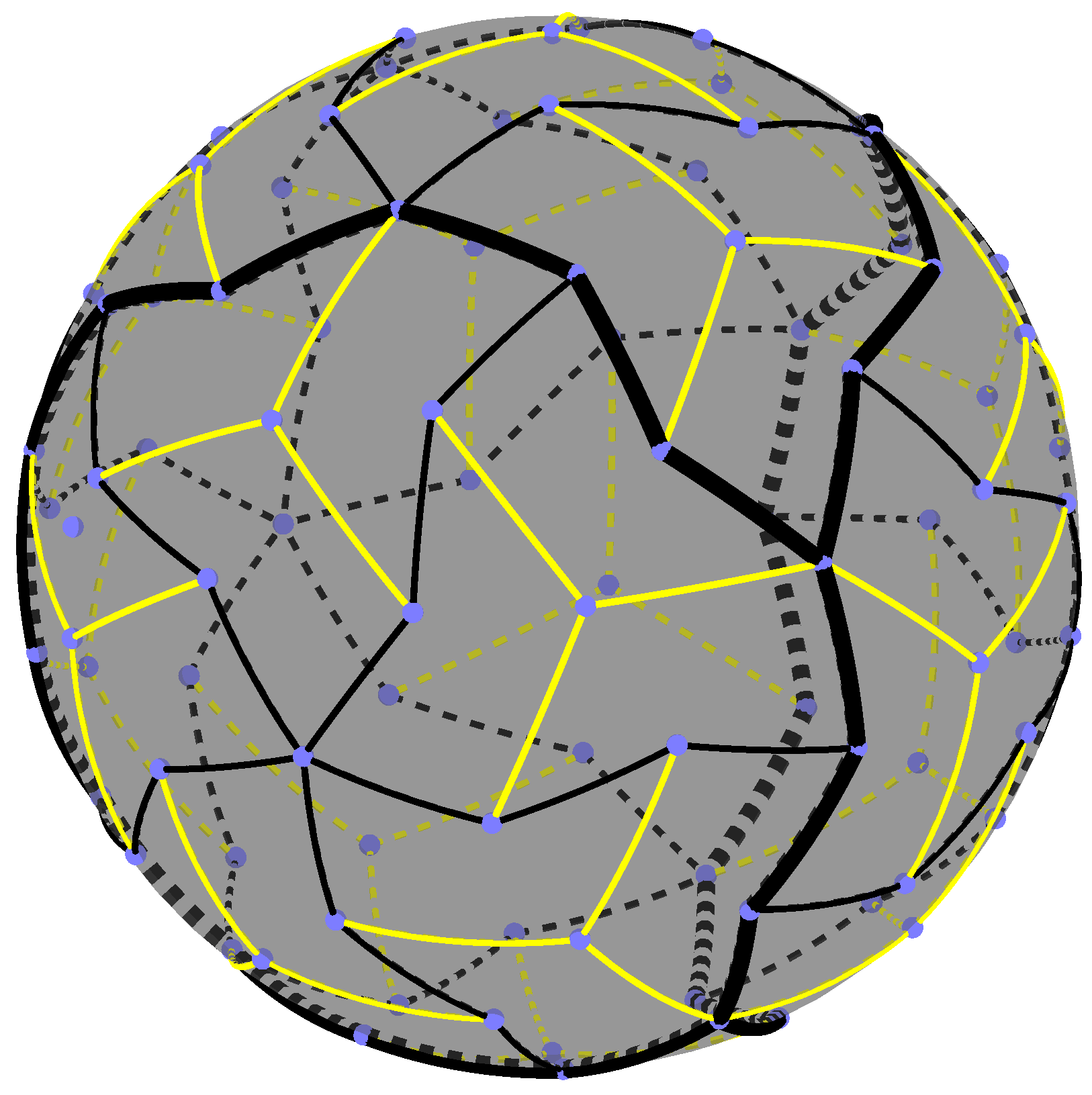}
 
 \includegraphics[scale=0.033]{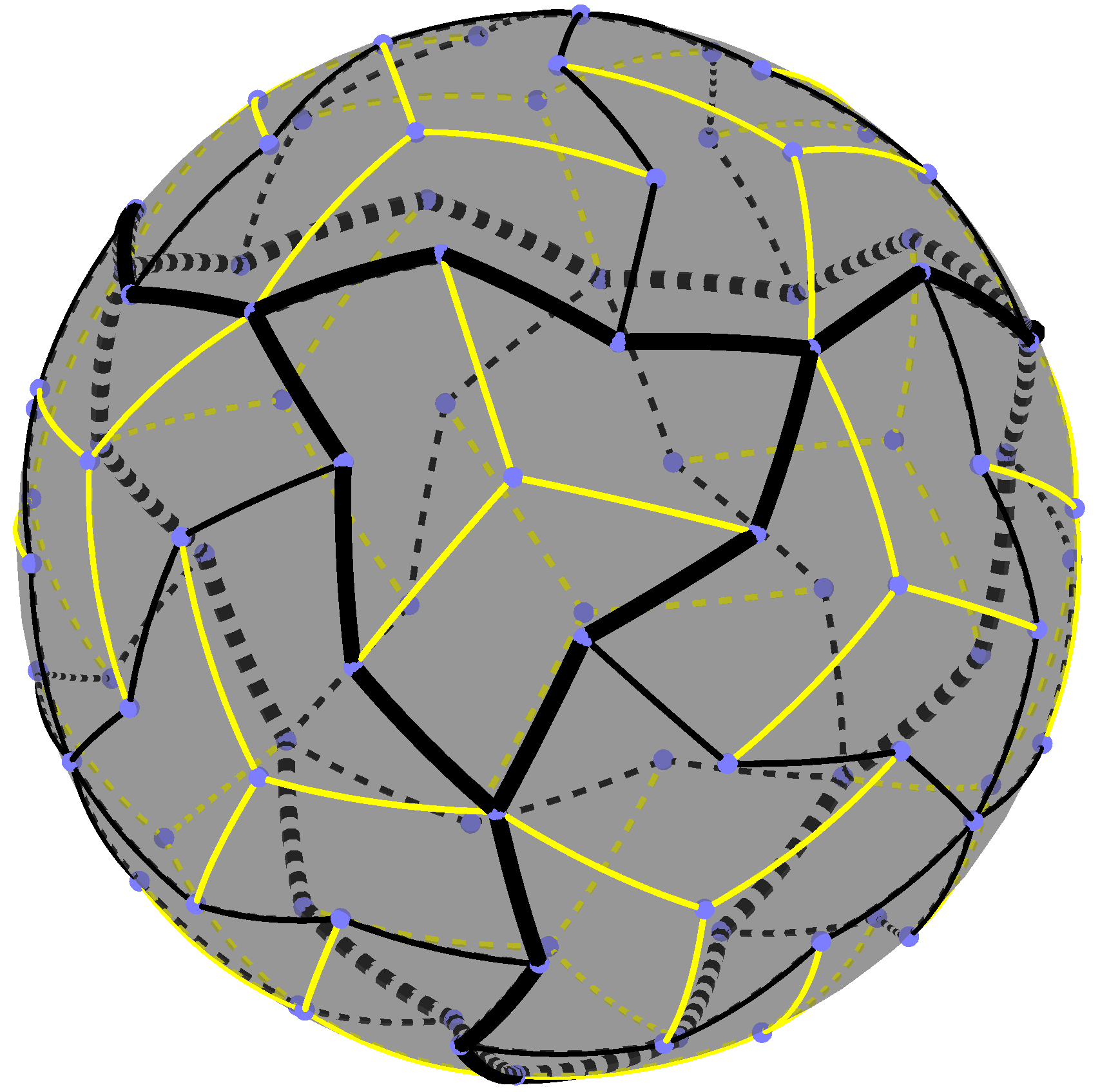}\hspace{10pt}
 \includegraphics[scale=0.032]{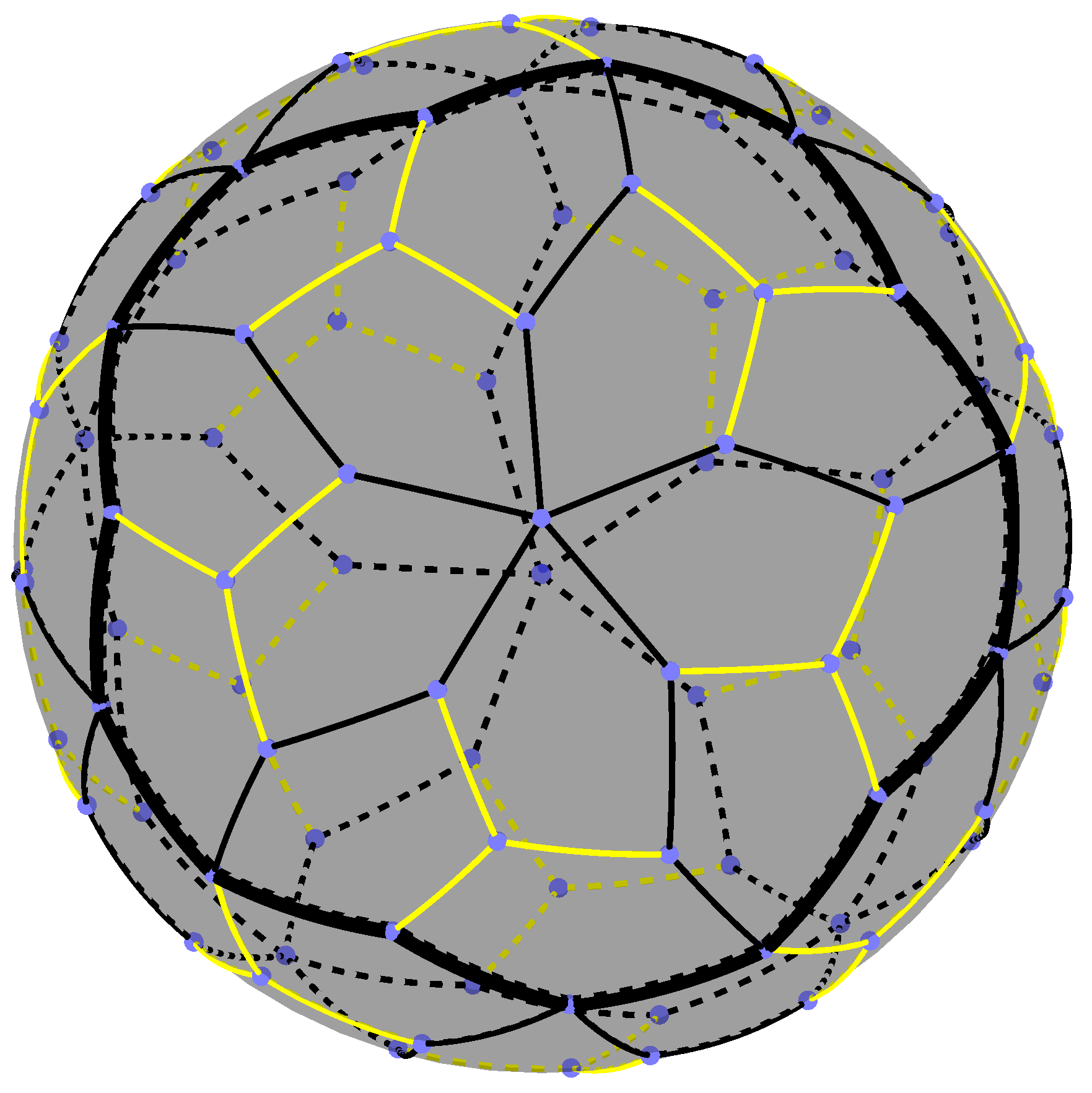}\hspace{10pt}
 \includegraphics[scale=0.0355]{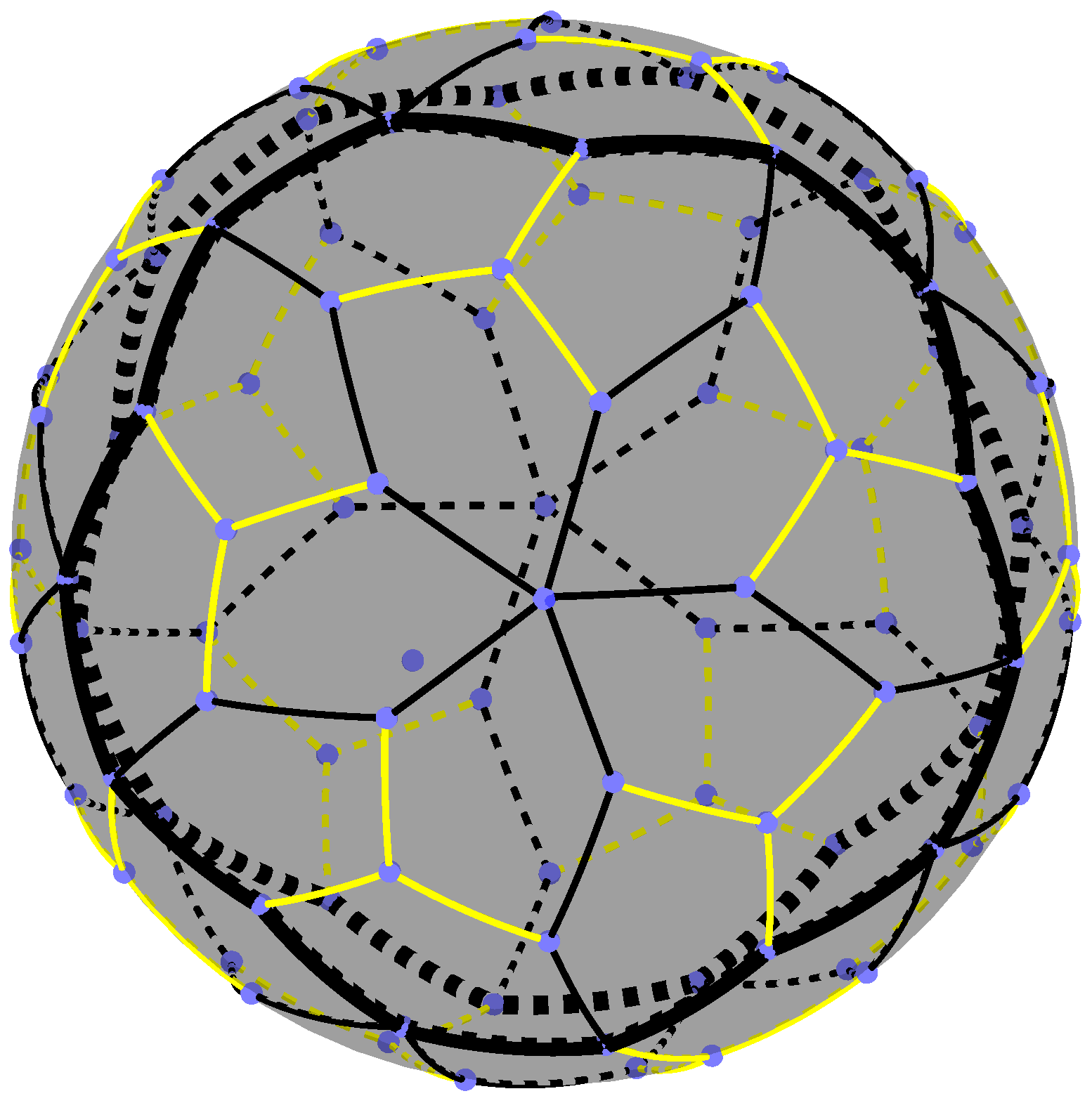}\hspace{10pt}
 \includegraphics[scale=0.033]{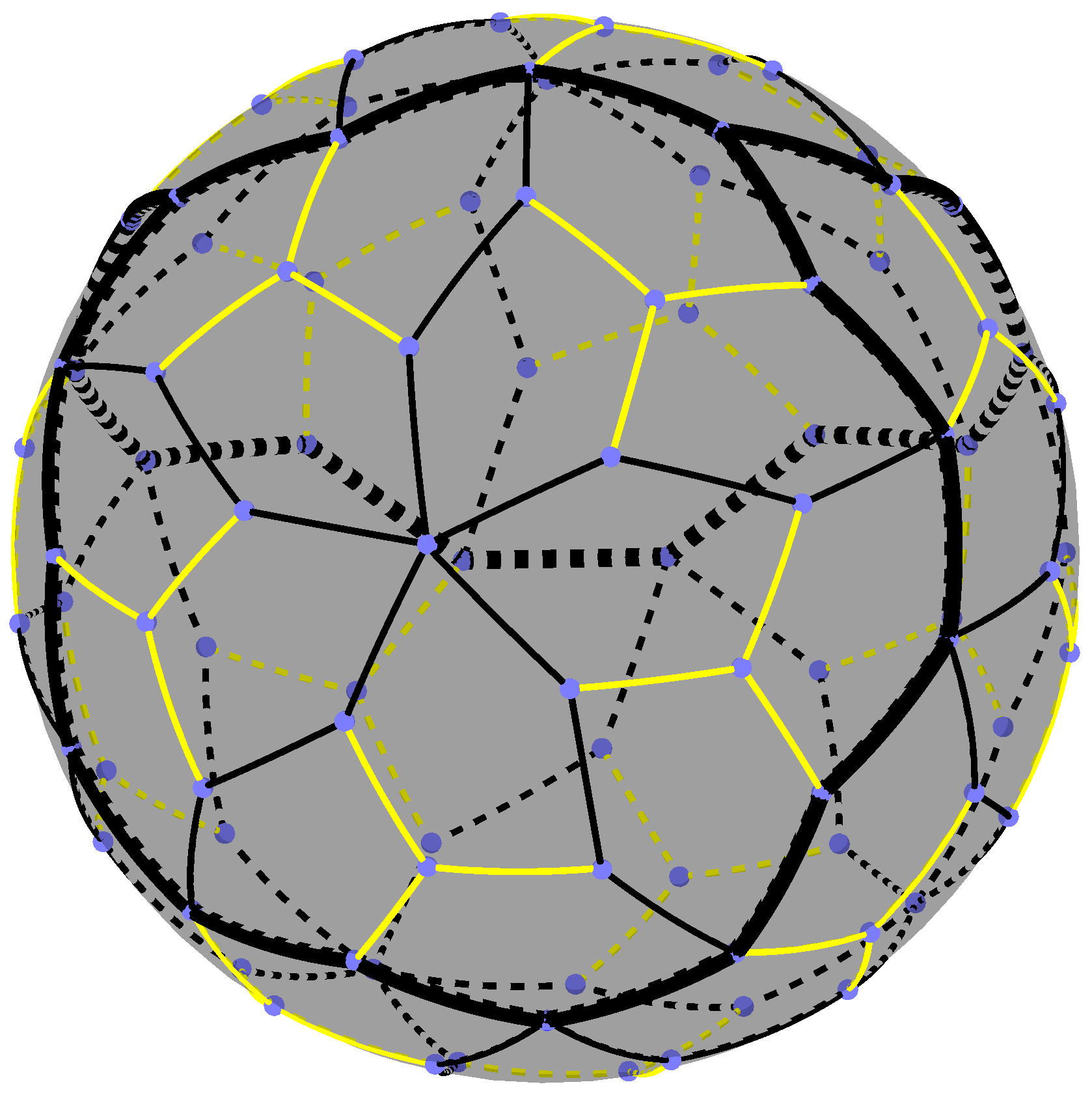}\hspace{10pt}
 \includegraphics[scale=0.034]{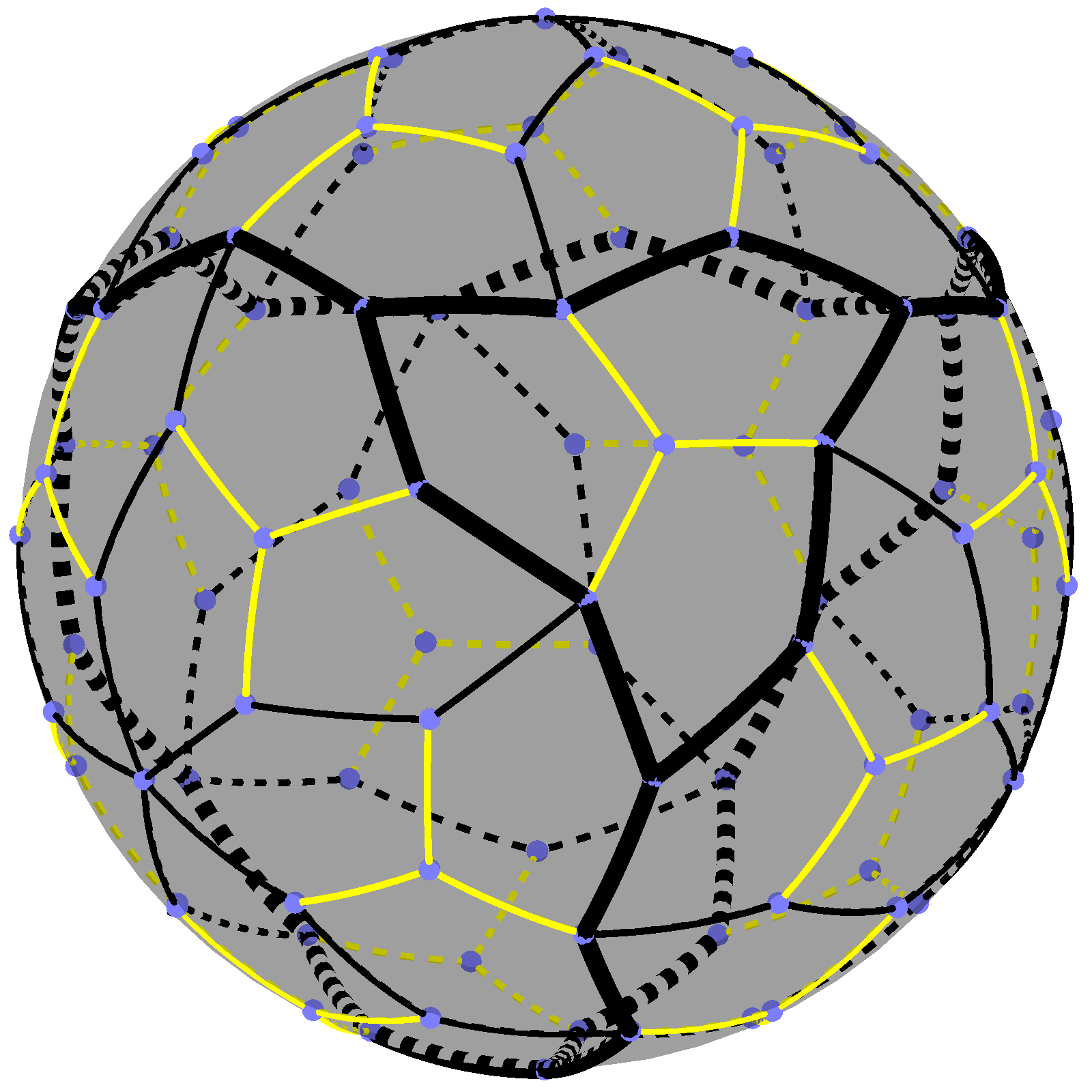}
 \caption{Five one-parameter families of pentagonal subdivision tilings, and ten flip modifications of three special cases of two pentagonal subdivision tilings for $a^3b^2$.} 
 \label{3a2b} 
\end{figure}

    \item A sequence of one-parameter families of $a^4b$-pentagons in \cite[Figure 5]{slw}, each admitting a non-symmetric $3$-layer earth map tiling with $4k$ tiles for any $k\ge4$, among which each odd $k=2n+1$ case always admits a standard flip modification in \cite[Figure 6]{slw}. Moreover, 
  \begin{itemize}
   \item The case $f=16$ and $\aaa=\tfrac{\pi}{2}$ admits a modification in Table \ref{Tab-summary-5} (see the second picture of \cite[Figure 7]{slw}).
    \item The case $\aaa=\bbb=1-\frac{4}{f}$ admits additional modifications in \cite[Table 3]{slw}  (see the second picture of \cite[Figure 9]{slw}).
    \item Special parameters induce five equilateral tilings in Table \ref{Tab-summary-6} (see the last five pictures of Figure \ref{5a}).    
    \end{itemize}

\begin{table*}[htp]                        
	\centering     
	\caption{All tilings of type $a^5$. }\label{Tab-summary-6}   
	~\\ 
	\begin{tabular}{c|c|c}	 
	\hline  	
	 $f$&$(\aaa,\bbb,\ccc,\ddd,\eee)$ &Tilings\\
	\hline	
	 $12$&$(\frac{2}{3},\frac{2}{3},\frac{2}{3},\frac{2}{3},\frac{2}{3})$& $ T(12\aaa\ddd\eee,4\ccc^3,4\bbb^3)$\\
	\hline			
	$24$&$(\alpha,\frac{1}{2},\frac{2}{3},\ddd,2-\aaa-\ddd),\aaa \approx 0.687,\ddd \approx 0.801$& $ T(24\aaa\ddd\eee,8\ccc^3,6\bbb^4)$\\    
	\hline
	$60$&$(\alpha,\frac{2}{5},\frac{2}{3},\ddd,2-\aaa-\ddd),\aaa \approx 0.684,\ddd \approx 0.905$& $ T(60\aaa\ddd\eee,20\ccc^3,12\bbb^5)$\\
	\hline\hline	
	 $16$&$(\alpha,\frac{3}{4},\frac{1}{2},\ddd,2-\aaa-\ddd),\aaa \approx 0.663,\ddd \approx 0.453$& $ T(16\aaa\ddd\eee,8\bbb^2\ccc,2\ccc^4)$\\
	\hline
	\multirow{2}{*}{$20$}&\multirow{2}{*}{$(\alpha,\frac{4}{5},\frac{2}{5},\ddd,2-\aaa-\ddd),\aaa \approx 0.628,\ddd \approx 0.309$}& $ T(20\aaa\ddd\eee,8\bbb^2\ccc,4\bbb\ccc^3)$\\
	\cline{3-3}
	&& $ T(20\aaa\ddd\eee,10\bbb^2\ccc,2\ccc^5)$\\
	\hline			
	\multirow{2}{*}{$24$}&$(\frac{2}{3}-\ddd,\frac{5}{6},\frac{1}{3},\ddd,\frac{4}{3}),\ddd \approx 0.144$& \multirow{2}{*}{$ T(24\aaa\ddd\eee,12\bbb^2\ccc,2\ccc^6)$}\\    
	\cline{2-2}
	&$(\frac{1}{2},\frac{5}{6},\frac{1}{3},\ddd,\frac{3}{2}-\ddd),\ddd \approx 0.119$& \\    
	\hline
	\end{tabular}
     
\end{table*}

\begin{figure}[htp]
 \centering
 \includegraphics[scale=0.043]{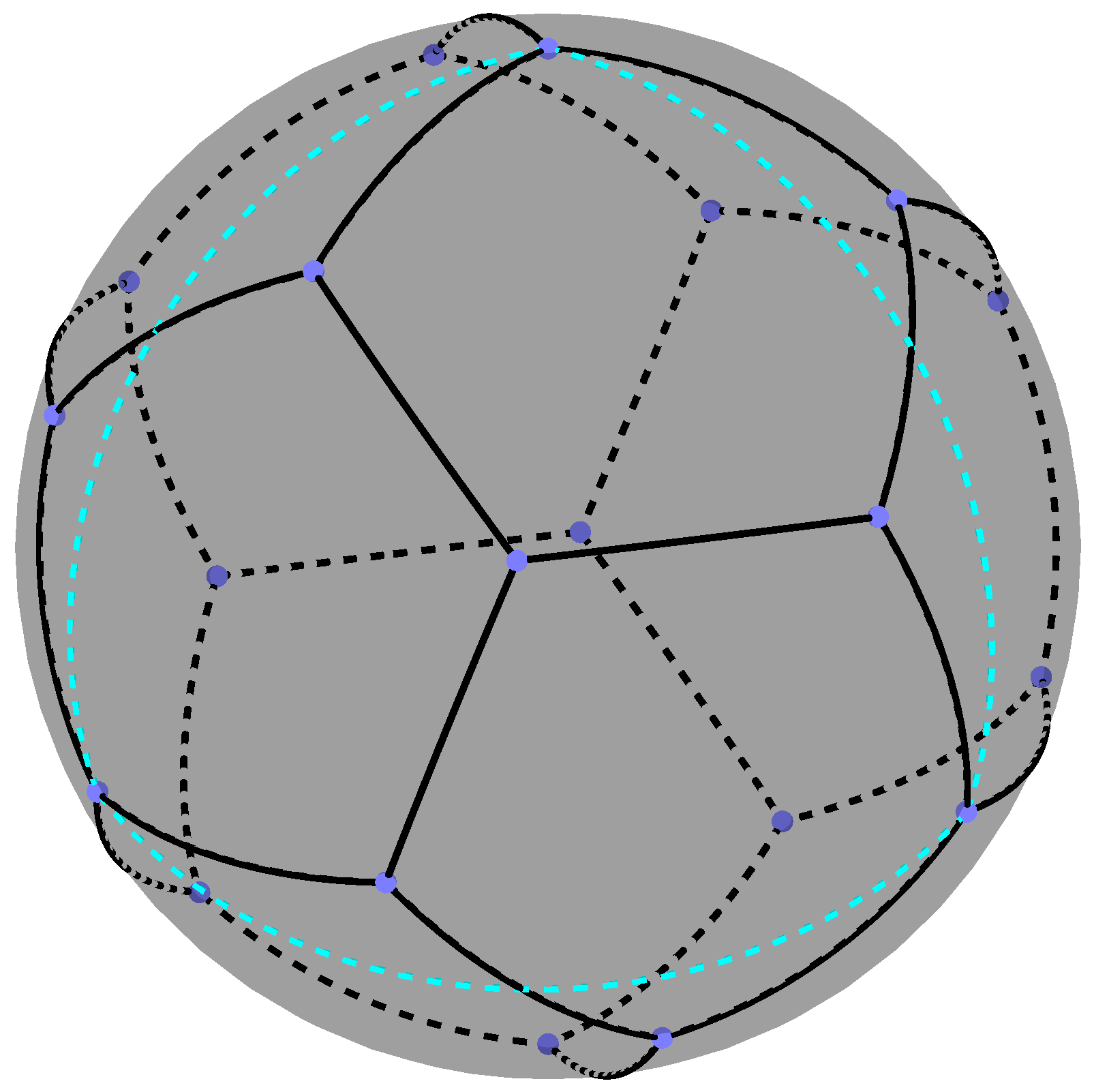}\hspace{20pt}
 \includegraphics[scale=0.044]{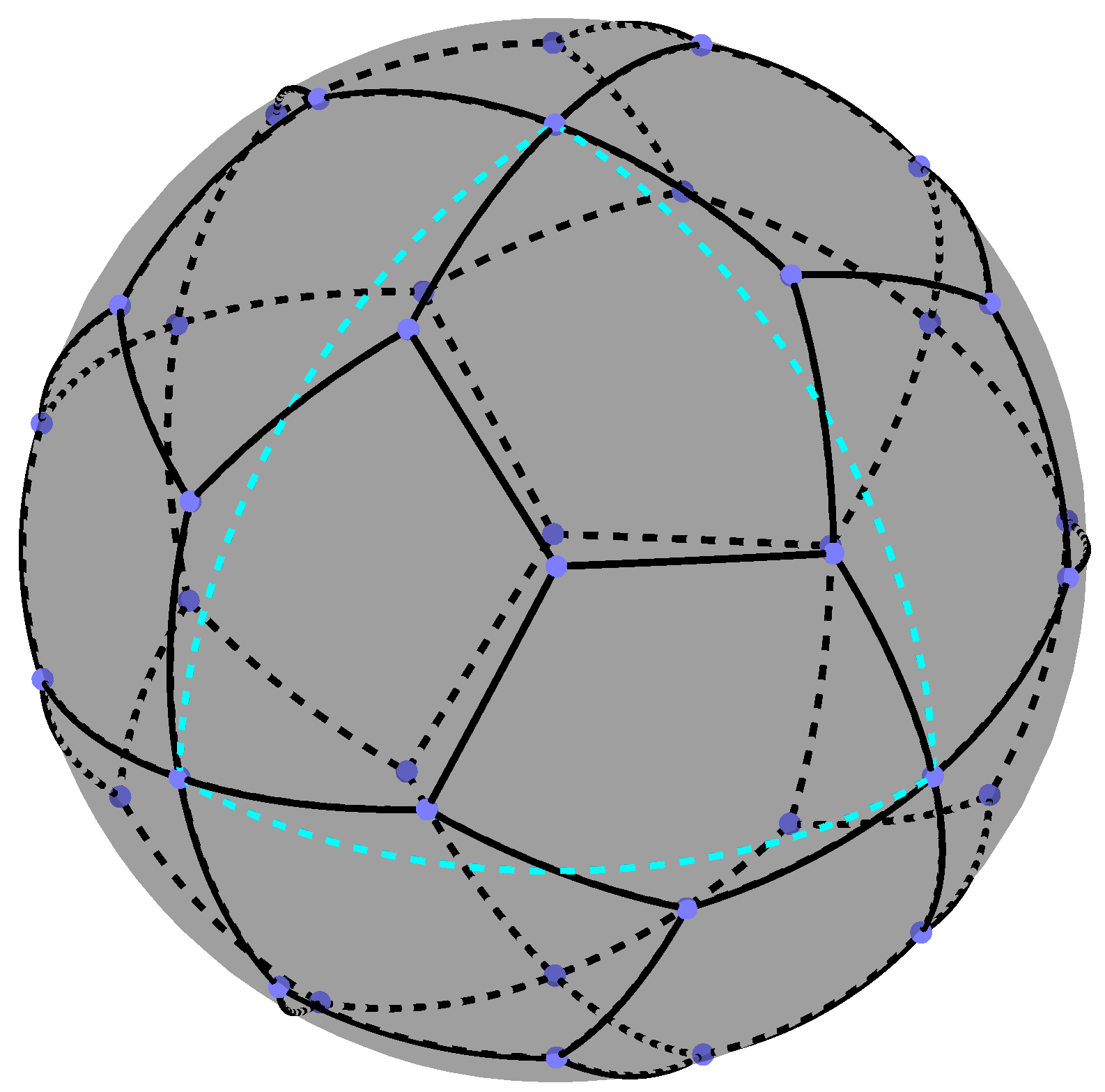}\hspace{20pt}
 \includegraphics[scale=0.042]{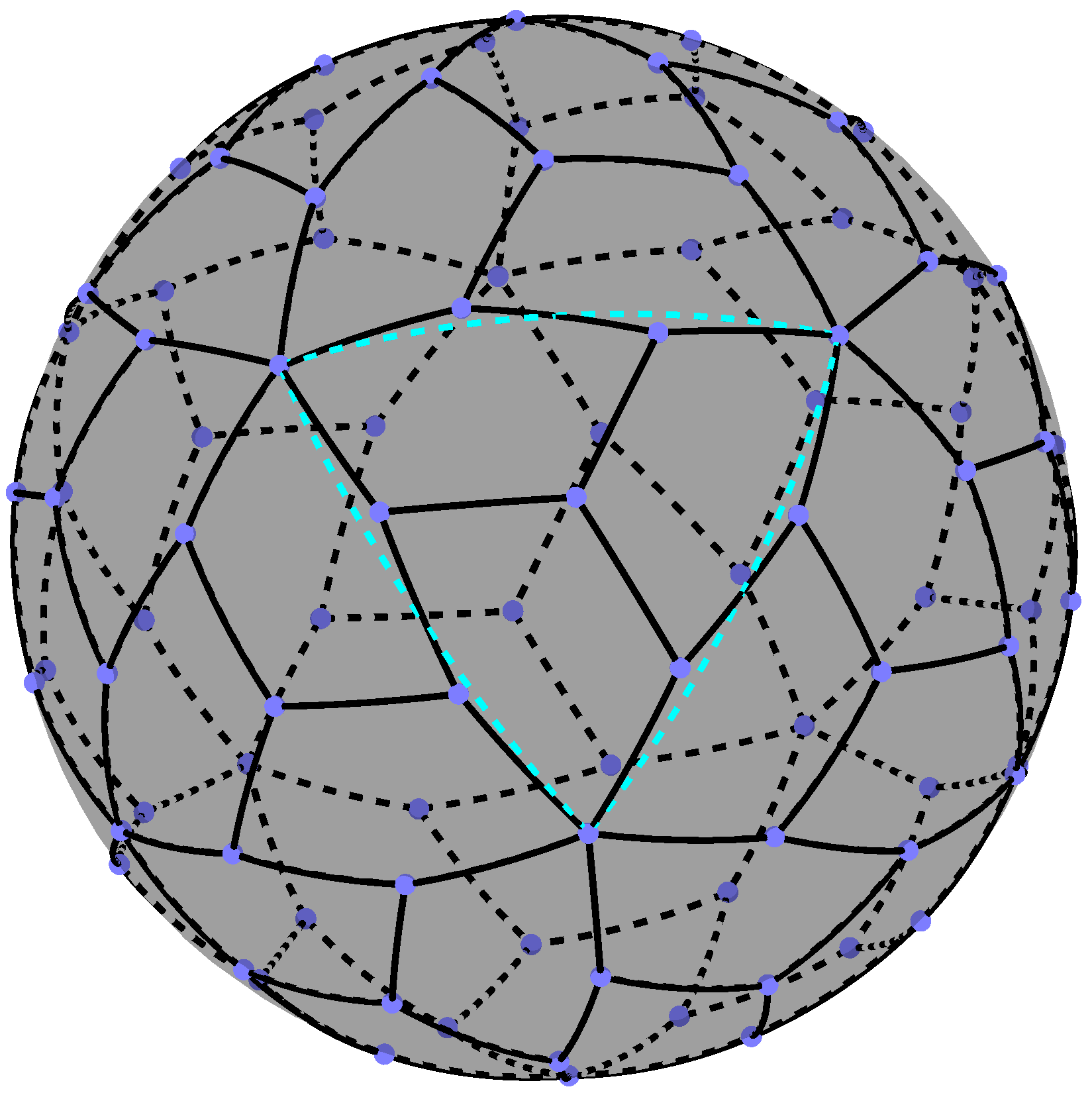}\hspace{20pt}
 \includegraphics[scale=0.043]{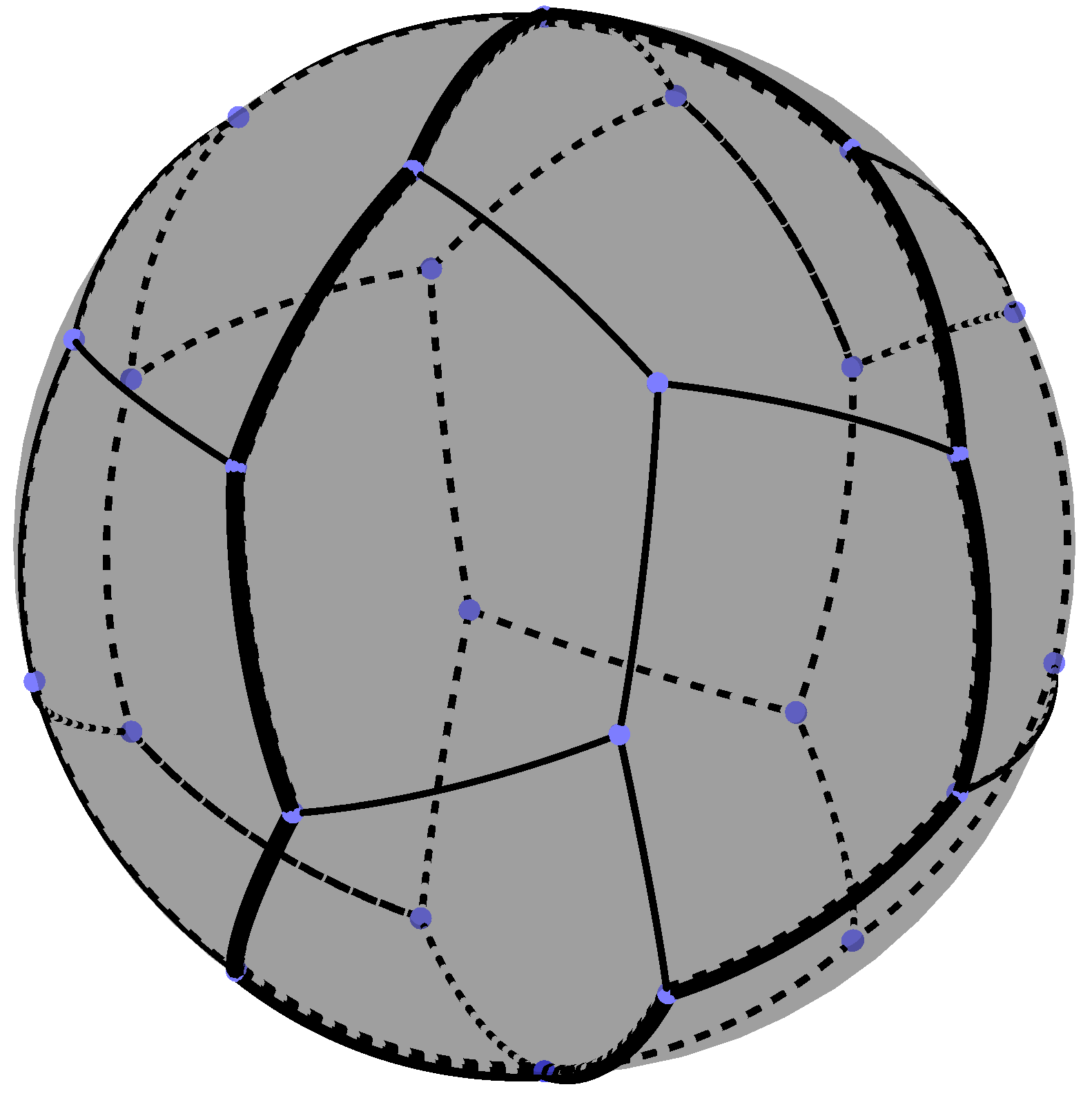}
 
\includegraphics[scale=0.043]{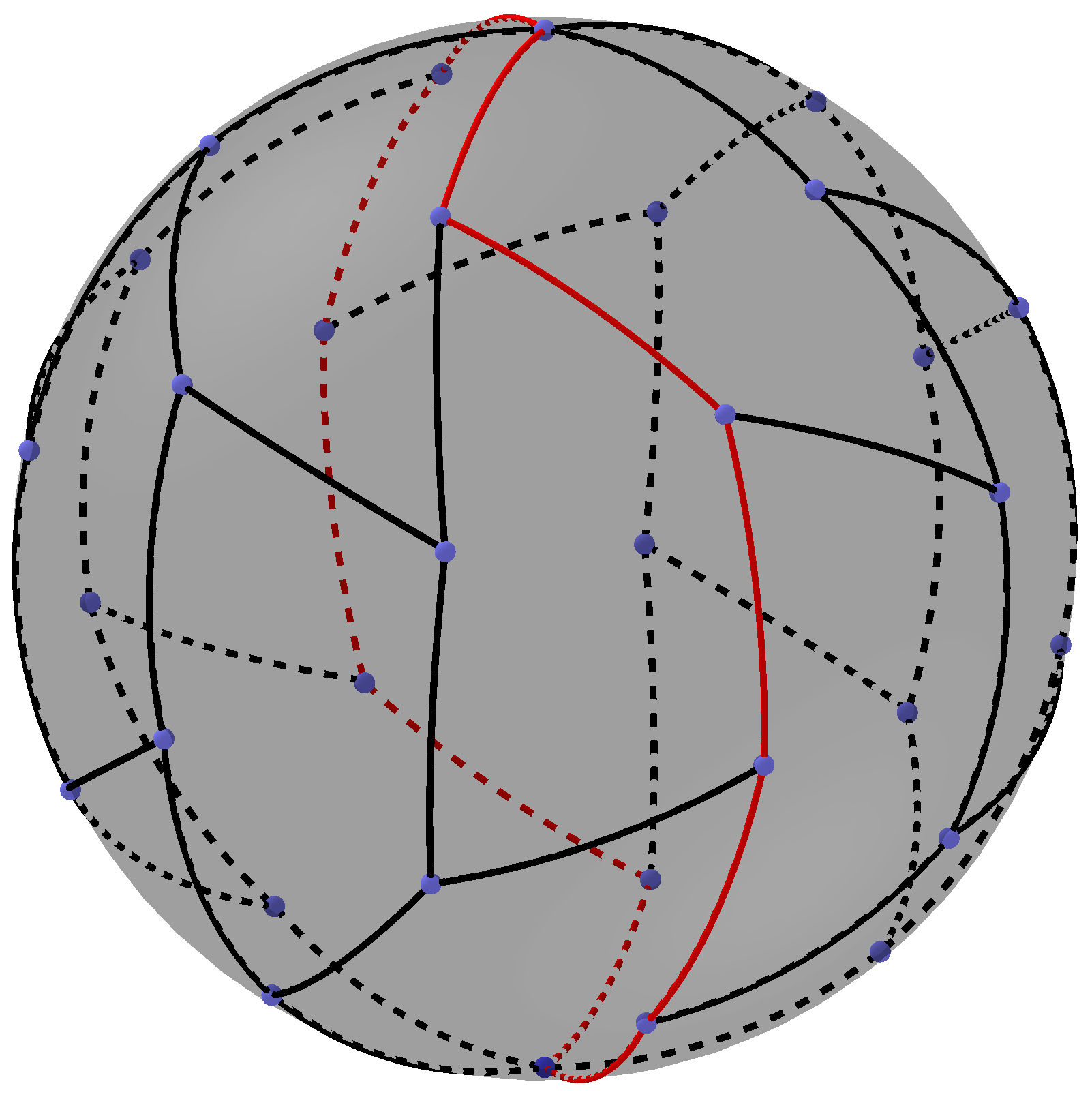}\hspace{19pt}
\includegraphics[scale=0.0504]{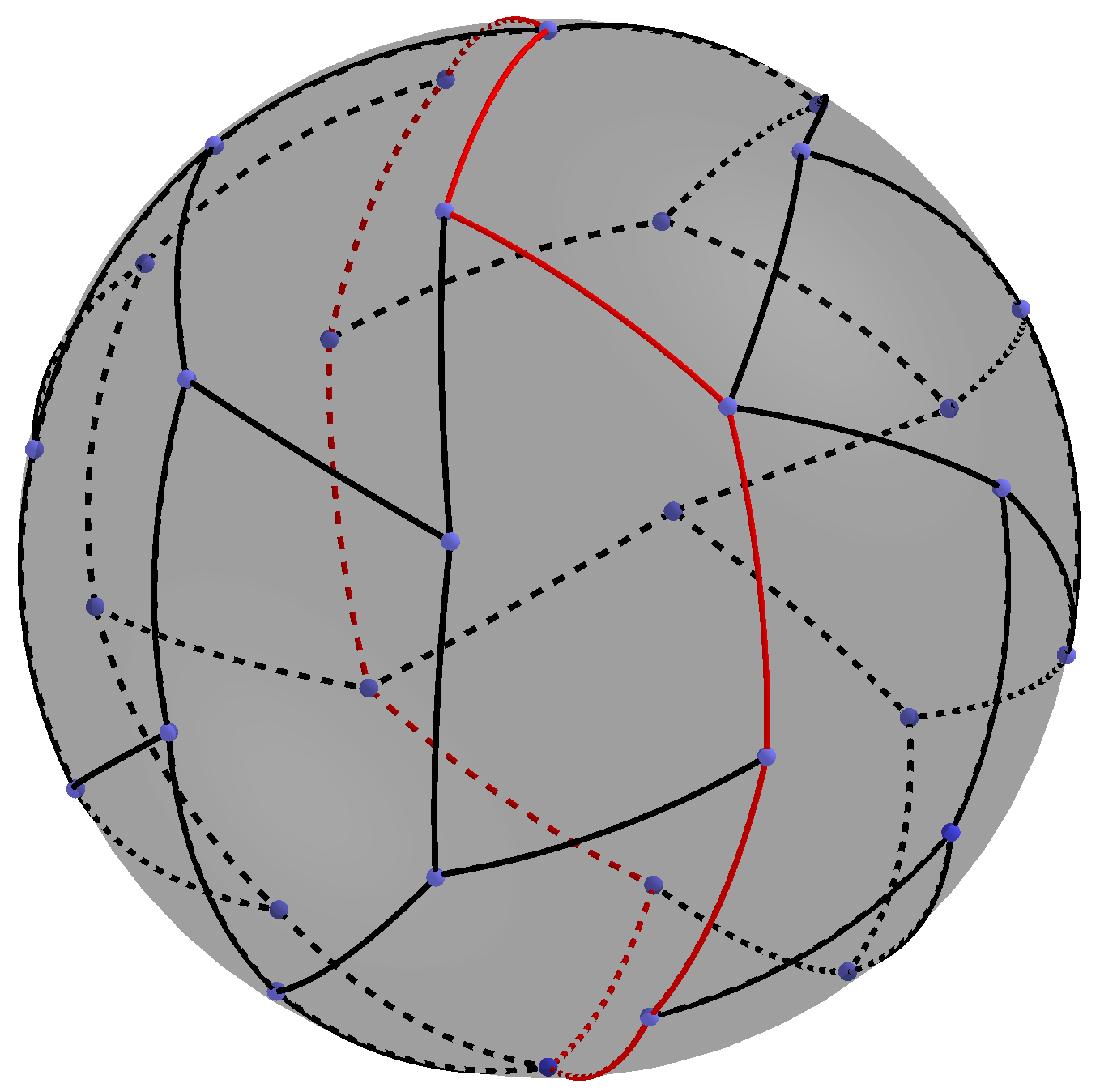}\hspace{20pt}
 \includegraphics[scale=0.044]{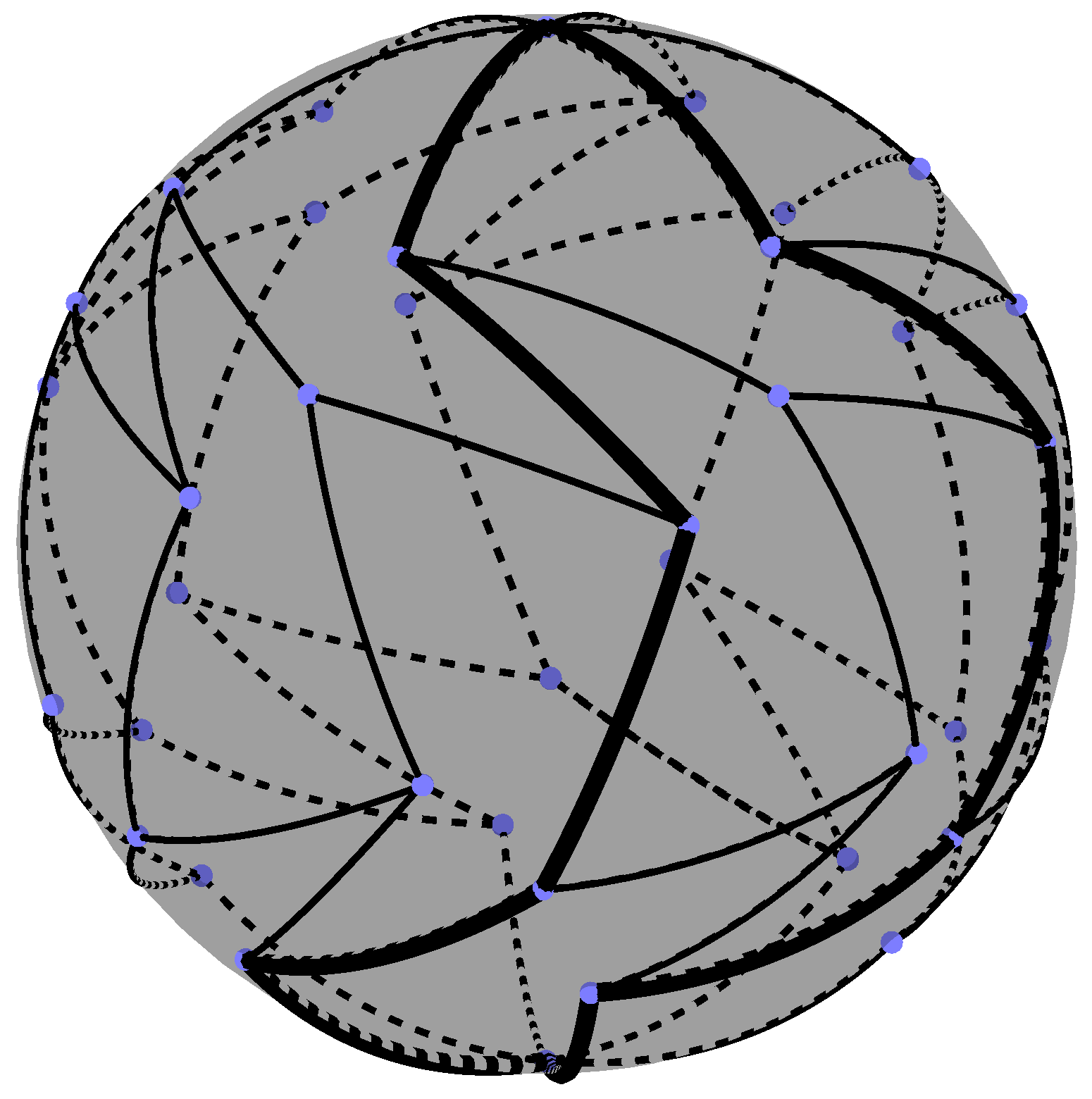}\hspace{20pt}
 \includegraphics[scale=0.042]{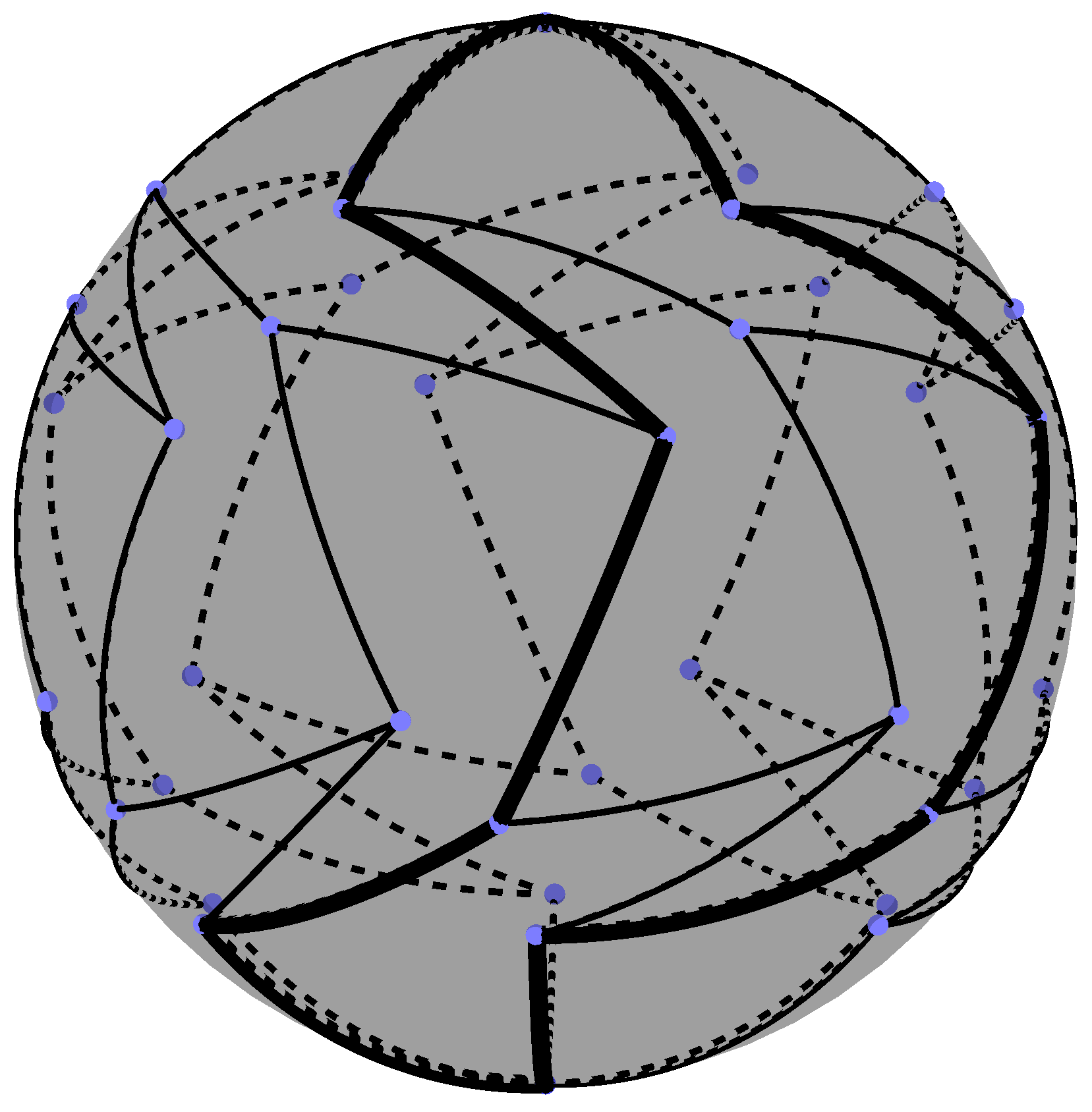}
 \caption{Three pentagonal subdivisions, four earth map tilings, and one flip modification of the earth map tiling for $a^5$.} 
 \label{5a} 
\end{figure}

  \item A sequence of unique symmetric $a^4b$-pentagons, 			each admitting a symmetric $3$-layer earth map tilings by $f=4k$ tiles for any $k\ge4$, among which each odd $k=2n+1$ case admits two standard flip modifications in Table \ref{Tab-summary-5} and \cite[Figure 2]{slw}.

\begin{table*}[htp]                        
	\centering     
	\caption{All tilings of type $a^4b$. }\label{Tab-summary-5}   
	\bgroup
   \resizebox{\linewidth}{!}{
   
	\begin{tabular}{c|c|c}	 
	\hline  	
	 $f$&$(\aaa,\bbb,\ccc,\ddd,\eee)$ &Tilings\\
	\hline	
	 $12$&$(\aaa,\frac{2}{3},\frac{2}{3},1-\frac{\aaa}{2},1-\frac{\aaa}{2})$& $ T(12\aaa\ddd\eee,4\bbb^3,4\ccc^3)$\\
	\hline			
	$24$&$(\alpha,\frac{1}{2},\frac{2}{3},\ddd,2-\aaa-\ddd)$& $ T(24\aaa\ddd\eee,8\ccc^3,6\bbb^4)$\\    
	\hline
	$60$&$(\alpha,\frac{2}{5},\frac{2}{3},\ddd,2-\aaa-\ddd)$& $ T(60\aaa\ddd\eee,20\ccc^3,12\bbb^5)$\\
	\hline\hline	
	 \multirow{3}{*}{$4k(k\geq 4)$}&\multirow{3}{*}{$(\frac{f}{8},1-\frac{4}{f},1-\frac{4}{f},\frac{1}{2}+\frac{2}{f},\frac{1}{2}+\frac{2}{f})$}& $ T(2k\,\aaa\bbb\ccc,2k\,\bbb\ddd\eee,2k\,\ccc\ddd\eee,2\aaa^{k})$\\
	\cline{3-3}
	&& $ T((2k-2)\aaa\bbb\ccc,2k\,\bbb\ddd\eee,2k\,\ccc\ddd\eee,2\aaa^{\frac{k+1}{2}}\bbb,2\aaa^{\frac{k+1}{2}}\ccc)$ ($k$ is odd)\\
	\cline{3-3}
	&& $ T(2k\,\aaa\bbb\ccc,(2k-1)\bbb\ddd\eee,(2k-1)\ccc\ddd\eee,\aaa^{\frac{k+1}{2}}\bbb,\aaa^{\frac{k+1}{2}}\ccc,2\aaa^{\frac{k-1}{2}}\ddd\eee)$ ($k$ is odd)\\
	\hline		
	\multirow{4}{*}{$4k(k\geq 4)$}&\multirow{4}{*}{$(\aaa,1-\frac{4}{f},\frac{8}{f},\ddd,2-\aaa-\ddd)$}&$T(4k\,\aaa\ddd\eee,2k\,\bbb^2\ccc,2\ccc^k)$\\     
	\cline{3-3}
	&& $T(4k\,\aaa\ddd\eee,(2k-2)\bbb^2\ccc,4\bbb\ccc^{\frac{k+1}{2}})$ ($k$ is odd)\\
	\cline{3-3}	
	&& $ T(8\aaa\ddd\eee,8\aaa\bbb^2,8\ccc\ddd\eee,2\ccc^4)$ when $k=4,\aaa=\frac{1}{2}$\\
	\cline{3-3}	
	&& Many more rearrangements when $k$ is odd and $\aaa=1-\frac{4}{f}$\\
	\hline	
	\end{tabular}}
	\egroup
     
\end{table*}

   \item Two unique double pentagonal subdivisions of the Platonic solids, with $48$ and $120$ tiles in Table \ref{Tab-summary-3} and Figure \ref{3abc}.

\begin{table*}[htp]                        
	\centering     
	\caption{All tilings of type $a^3bc$. }\label{Tab-summary-3}  
	~\\ 
	\begin{tabular}{c|c|c}	 
	\hline  	
	 $f$&$(\aaa,\bbb,\ccc,\ddd,\eee)$ &Tilings\\
	\hline	
	$48$&$(\frac{1}{2},\frac{3}{4},\frac{2}{3},\frac{2}{3},\frac{1}{2})$& $ T(12\aaa^4,24\bbb^2\eee,8\ddd^3,24\ccc^2\ddd,6\eee^4)$\\
	\hline			
	$120$&$(\frac{1}{2},\frac{4}{5},\frac{2}{3},\frac{2}{3},\frac{2}{5})$& $ T(30\aaa^4,60\bbb^2\eee,20\ddd^3,60\ccc^2\ddd,12\eee^4)$\\
	\hline
	\end{tabular}
      
\end{table*}

\begin{figure}[htp]
 \centering
 \includegraphics[scale=0.0436]{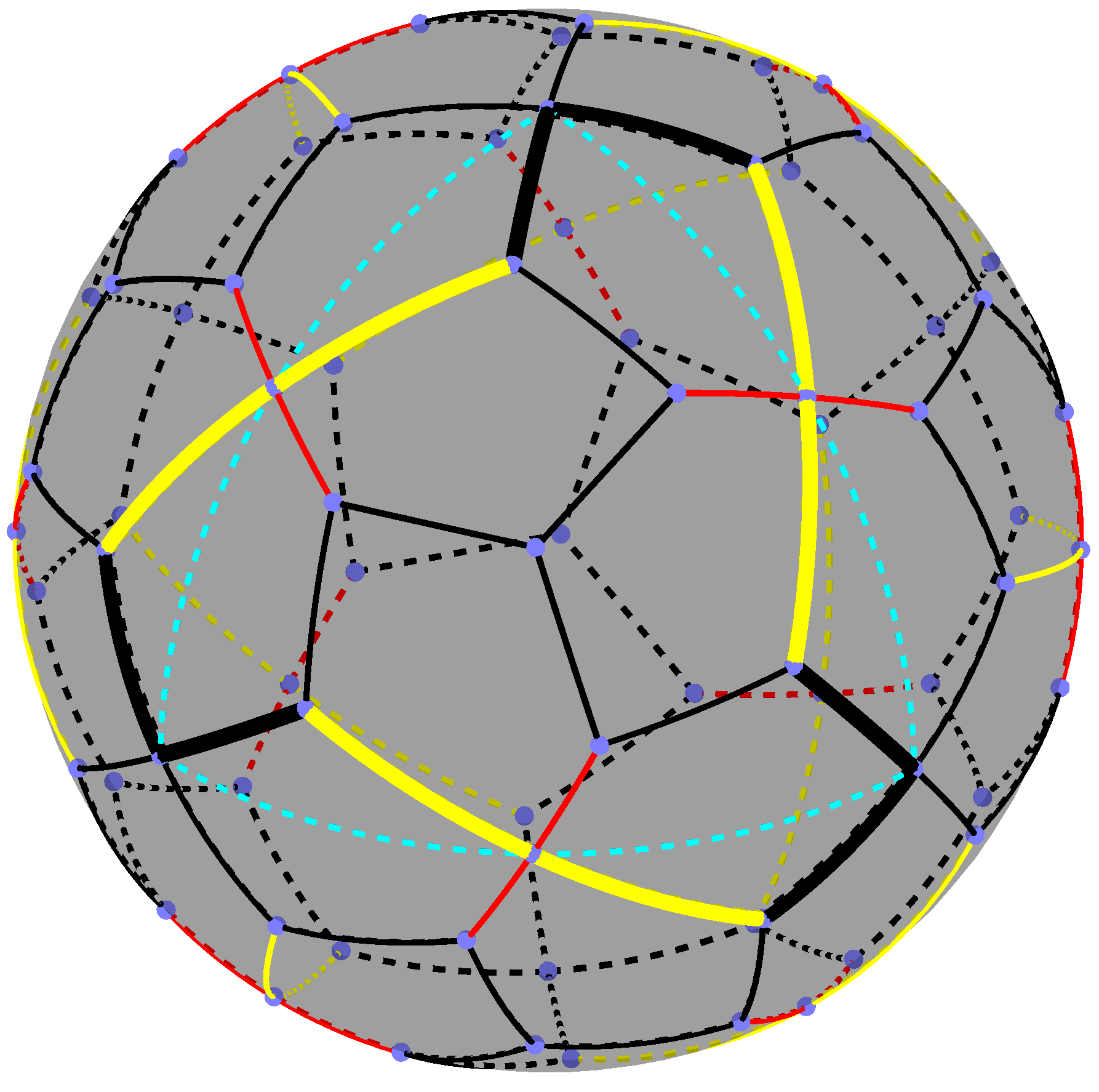}\hspace{20pt}
 \includegraphics[scale=0.0424]{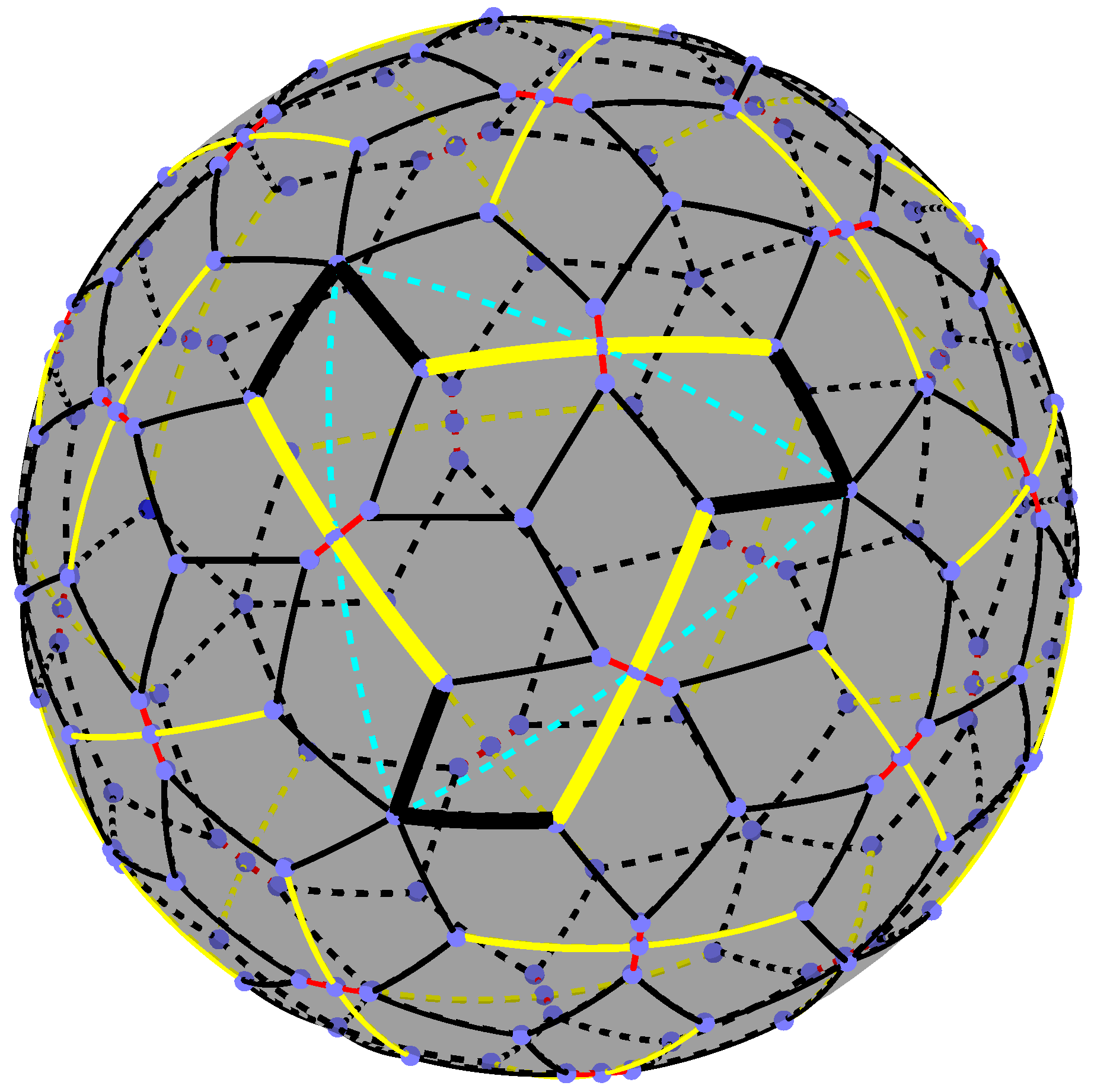}
 \caption{Two unique double pentagonal subdivisions of the Platonic solids  for $a^3bc$.} 
 \label{3abc} 
\end{figure}

\end{enumerate}

\medskip

Overall, the list of monohedral pentagonal tilings is shorter than the quadrilateral case, despite the notable complexity of geometric computations and the UFO tilings in \cite{slw}.

\subsection*{All induced non-edge-to-edge quadrilateral tilings}

\begin{figure}[htp]
 \centering
  \begin{overpic}[scale=0.038]{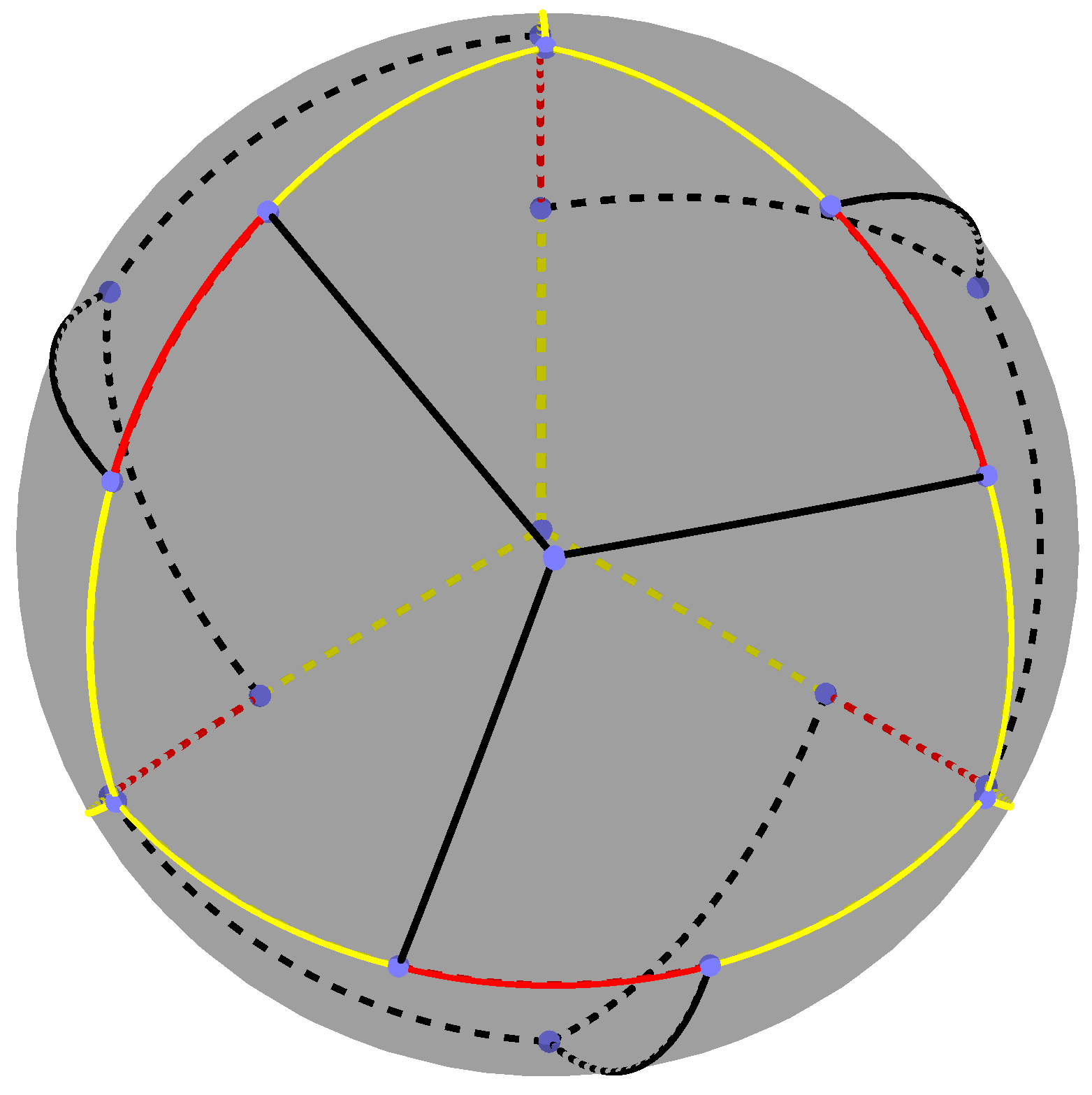}
 	\put(80,0){\small  $\eee=\pi$}
 \end{overpic}\hspace{20pt}
  \begin{overpic}[scale=0.083]{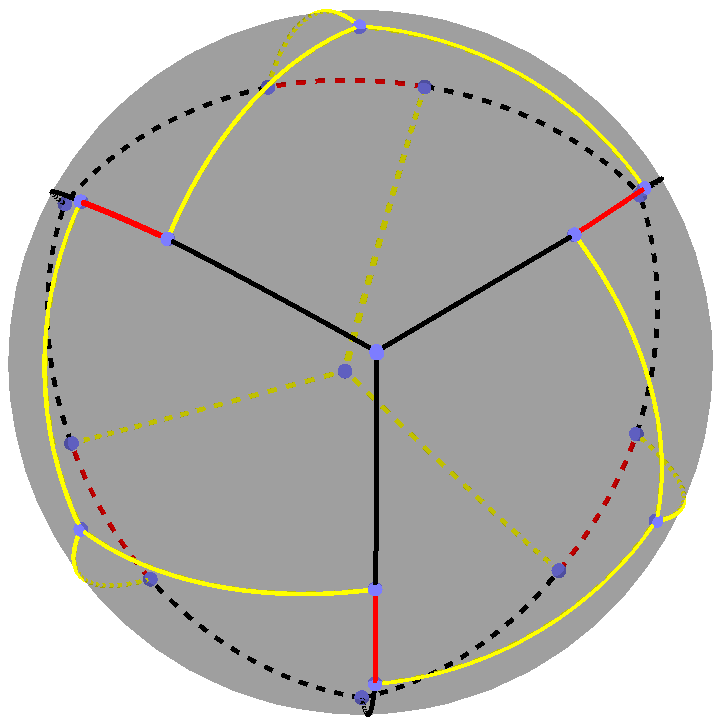}
	\put(80,0){\small  $\ddd=\pi$}
  \end{overpic}\hspace{20pt}
  \begin{overpic}[scale=0.0385]{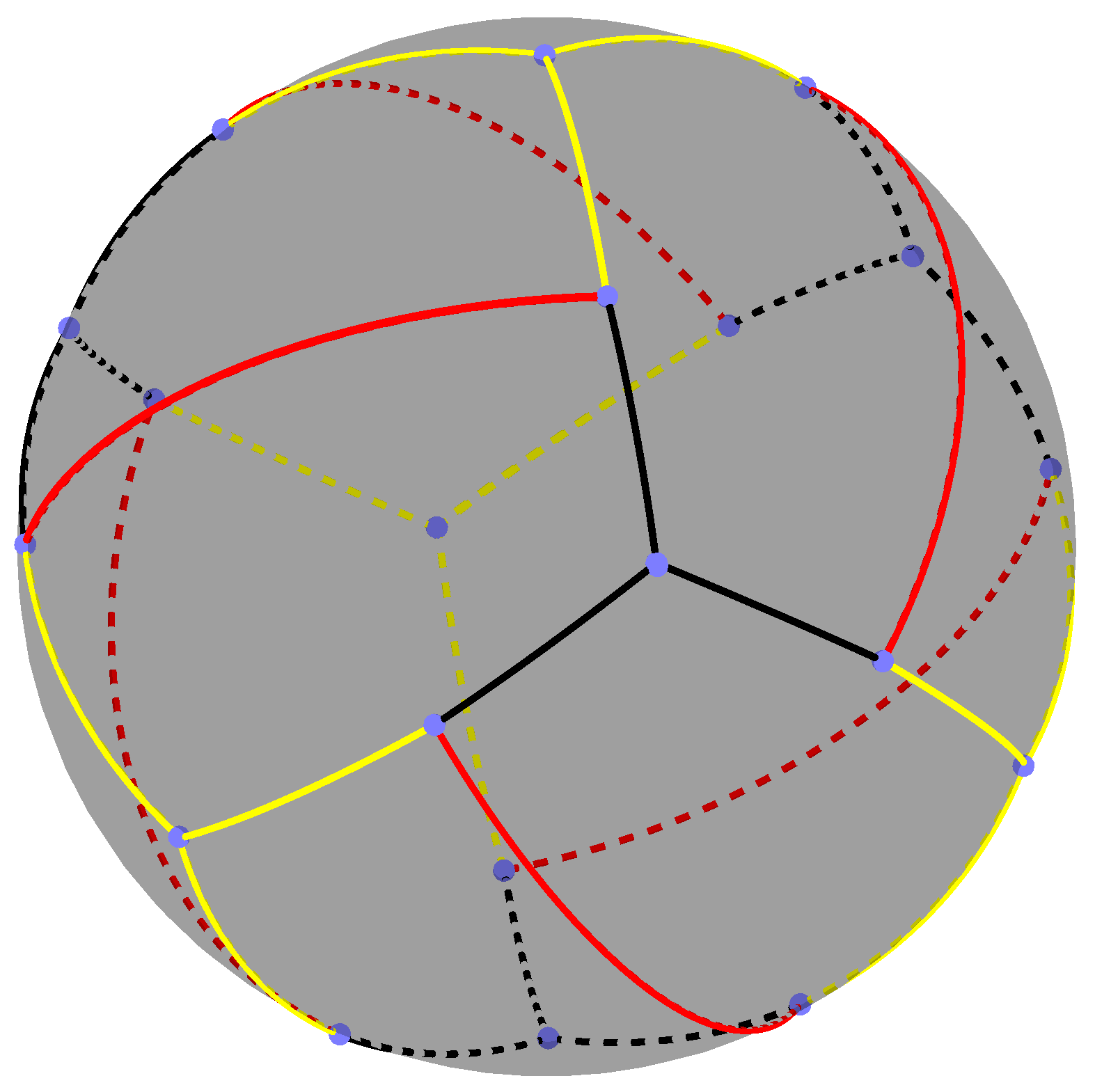}
	\put(80,0){\small  $\aaa=\pi$}
  \end{overpic}\hspace{20pt}

  \begin{overpic}[scale=0.036]{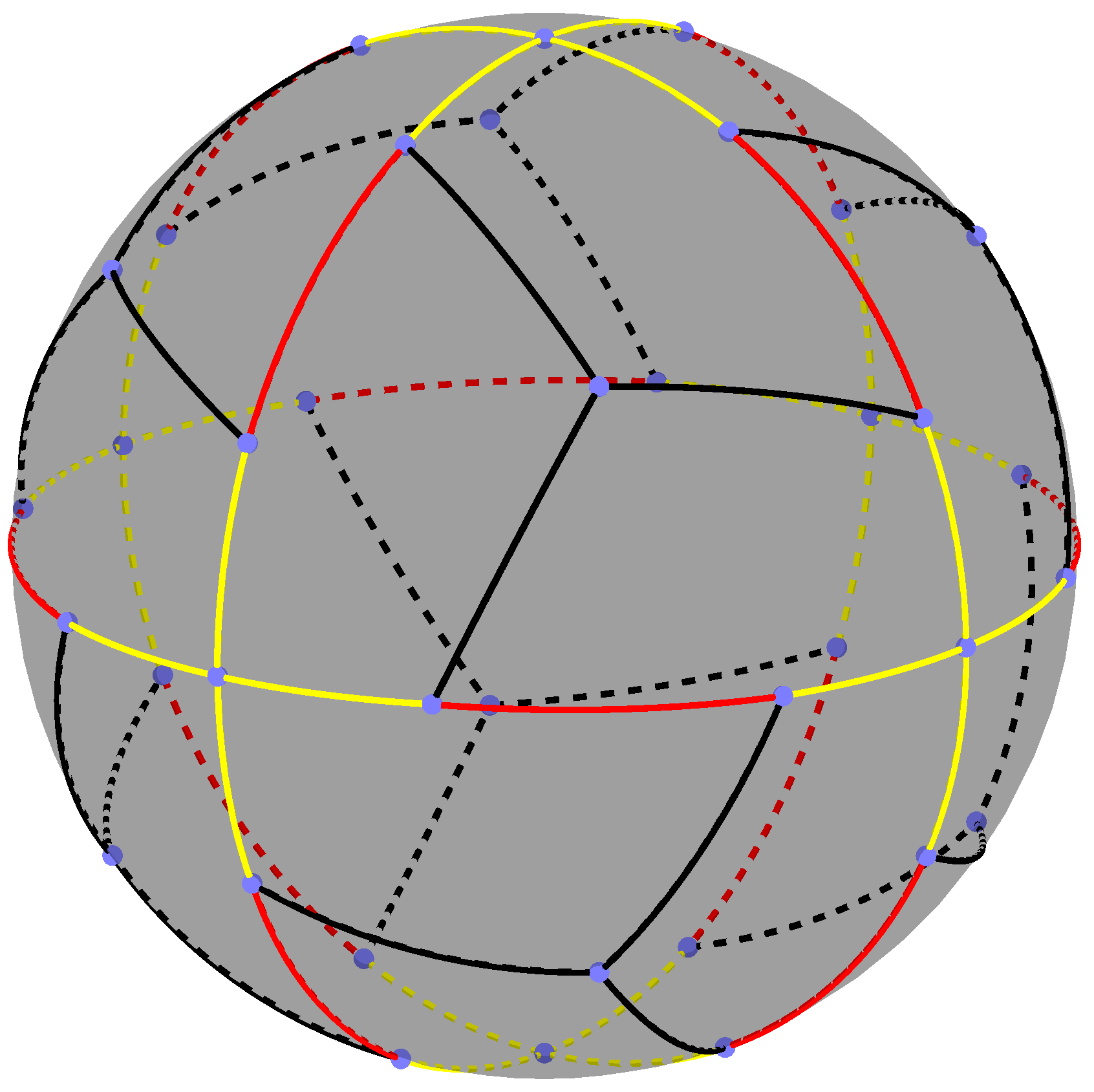}
  	\put(80,0){\small  $\eee=\pi$}
  \end{overpic}\hspace{20pt}
  \begin{overpic}[scale=0.036]{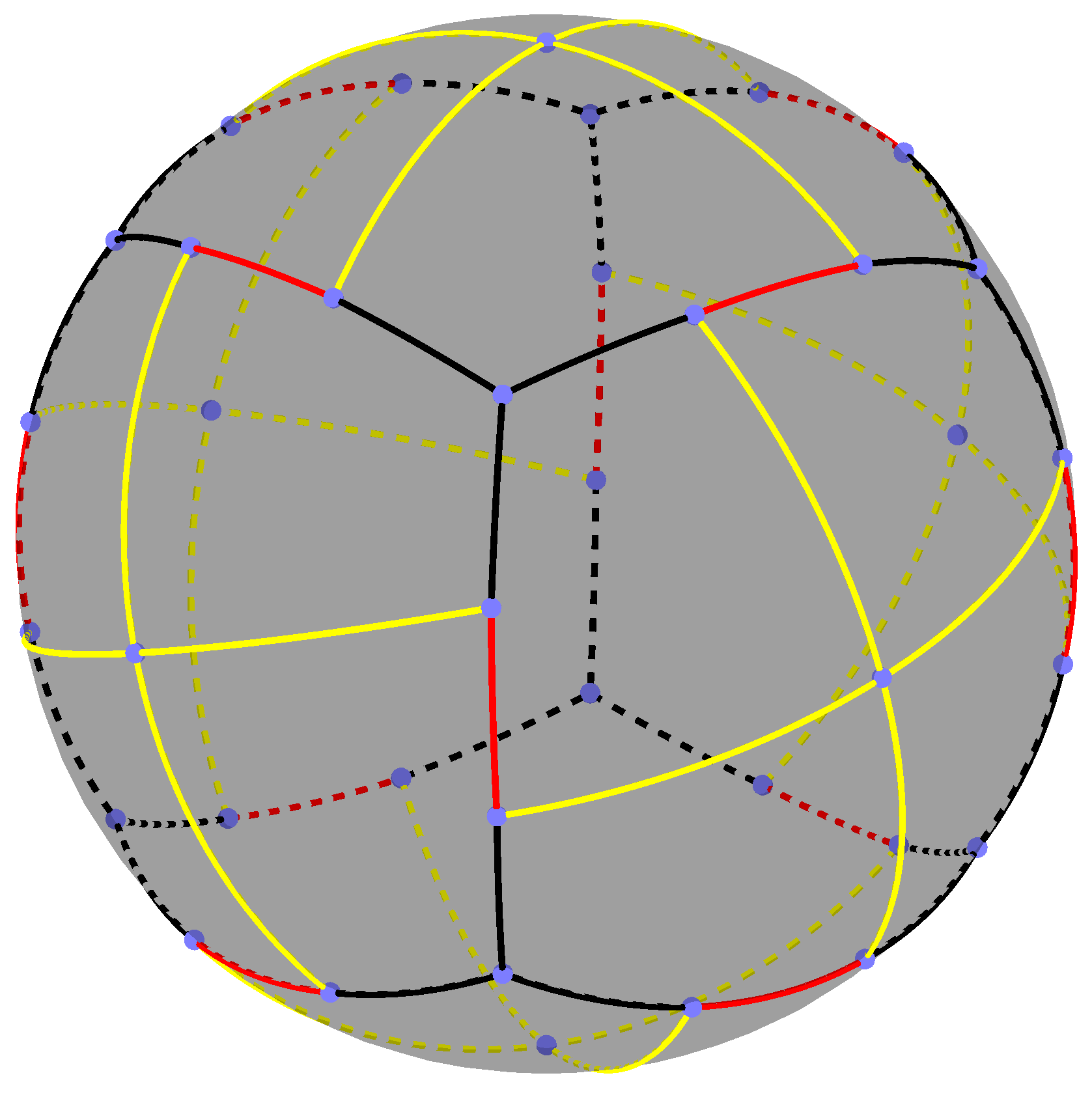}
  	\put(80,0){\small  $\ddd=\pi$}
  \end{overpic}\hspace{20pt}
  \begin{overpic}[scale=0.036]{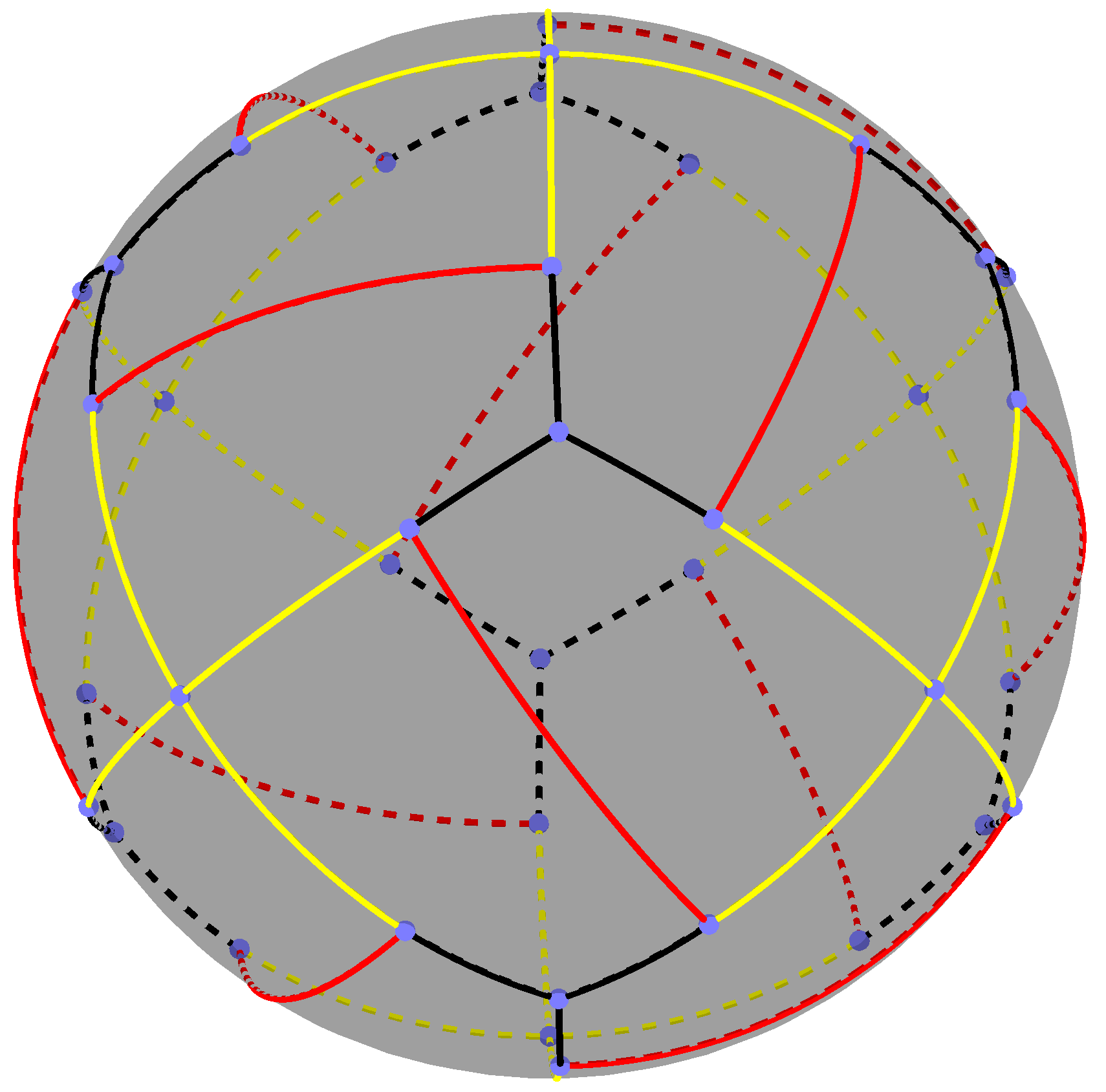} 
  	\put(80,0){\small  $\aaa=\pi$}
  \end{overpic}\hspace{20pt}

  \begin{overpic}[scale=0.041]{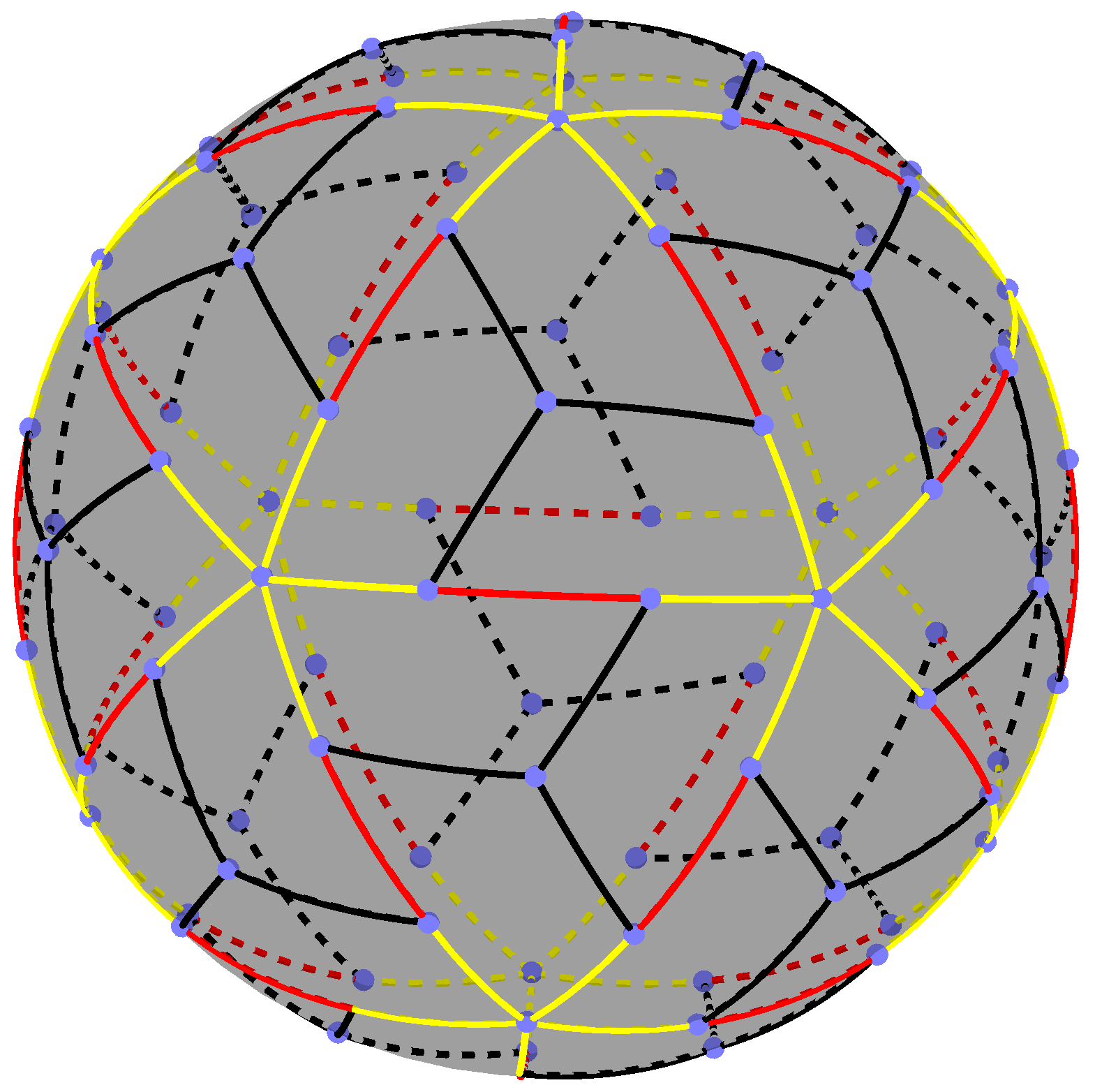}
  	\put(80,0){\small  $\eee=\pi$}
  \end{overpic}\hspace{20pt}
  \begin{overpic}[scale=0.040]{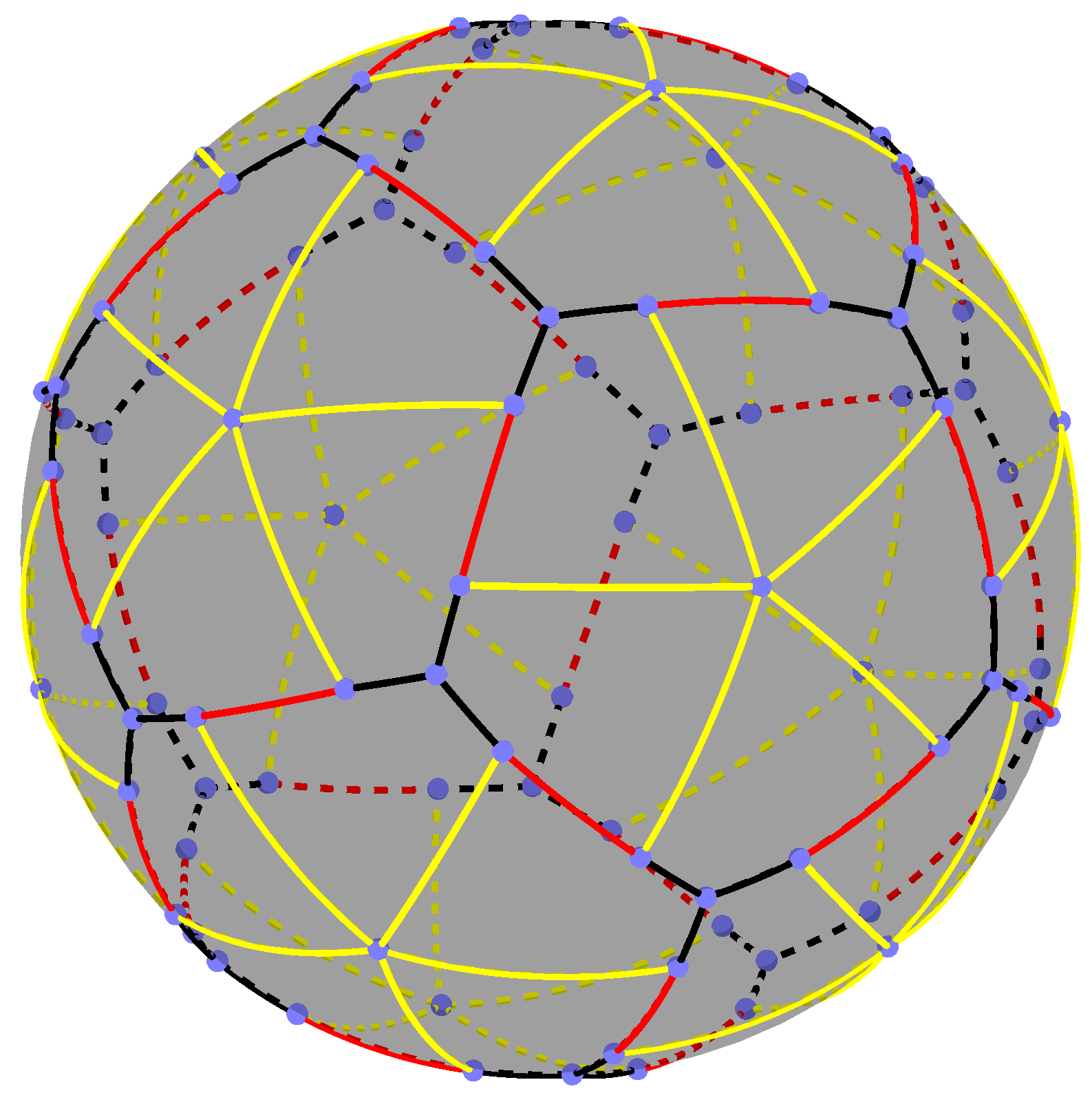}
  	\put(80,0){\small  $\ddd=\pi$}
  \end{overpic}\hspace{20pt}
  \begin{overpic}[scale=0.040]{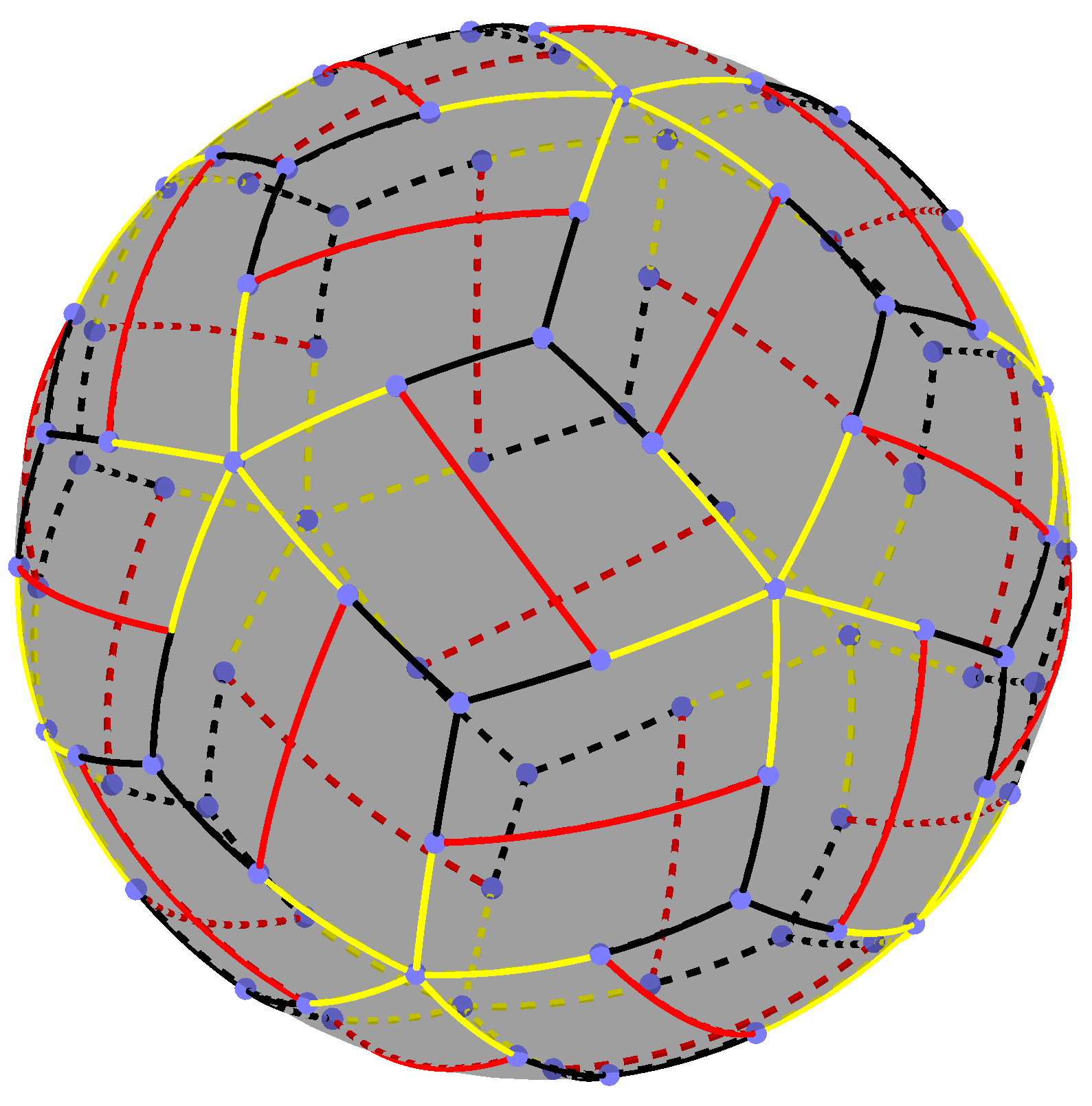} 
  	\put(80,0){\small  $\aaa=\pi$}
  \end{overpic}\hspace{20pt}

 \includegraphics[scale=0.043]{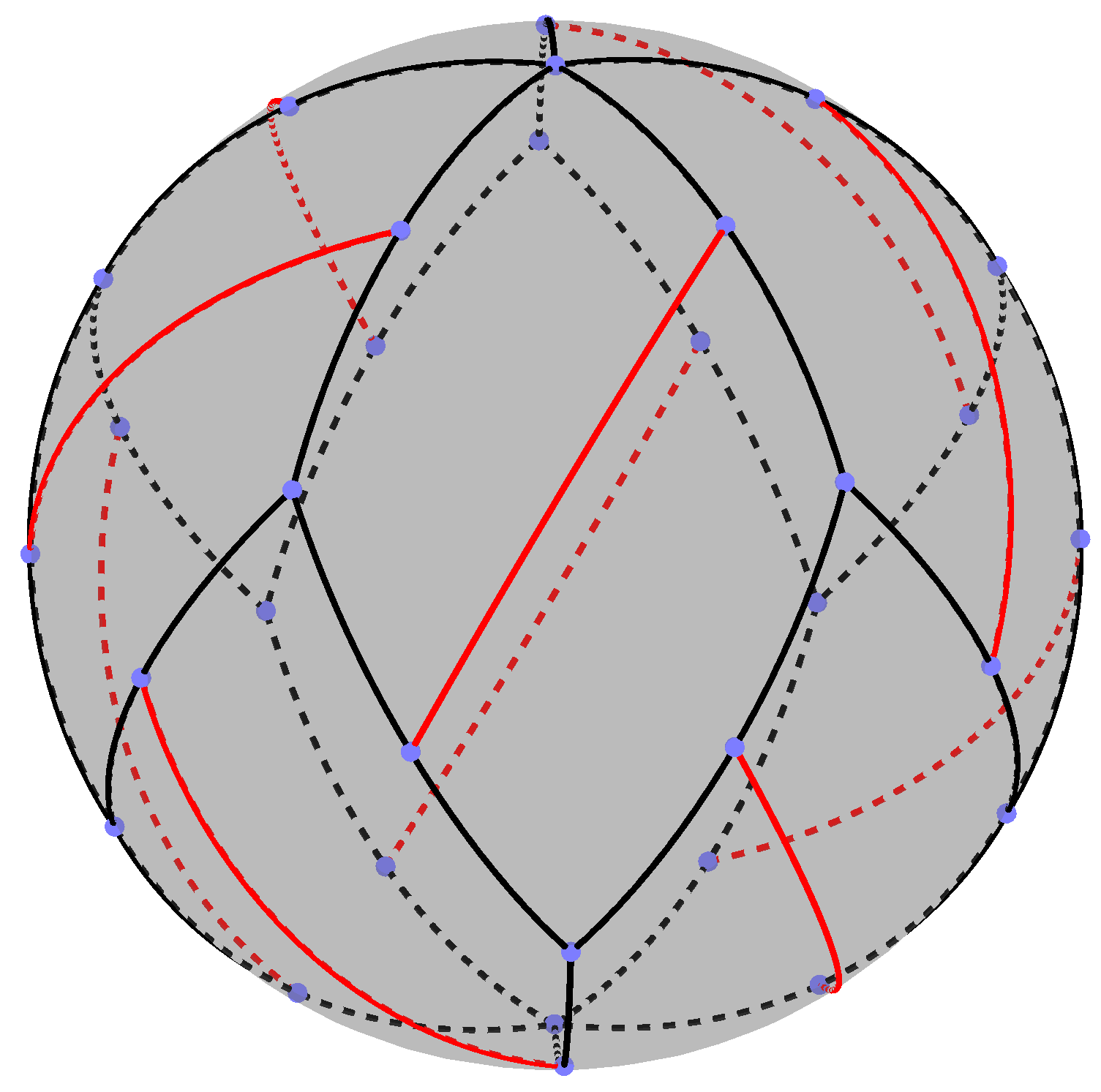}\hspace{10pt}
  \begin{overpic}[scale=0.028]{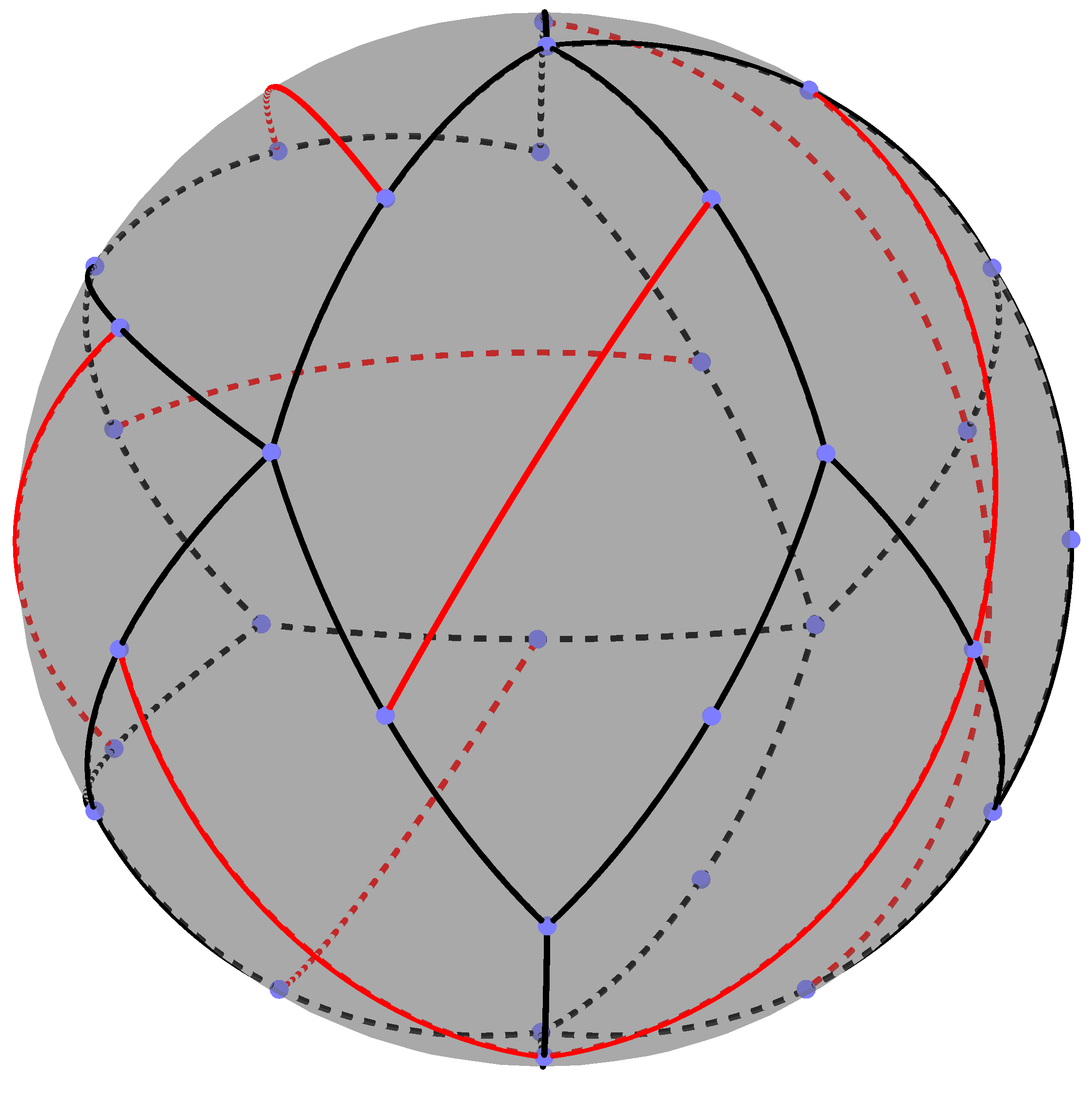}
   \put(-30,0){\small $\aaa=\pi$}
  \end{overpic} \hspace{10pt}
     \includegraphics[scale=0.038]{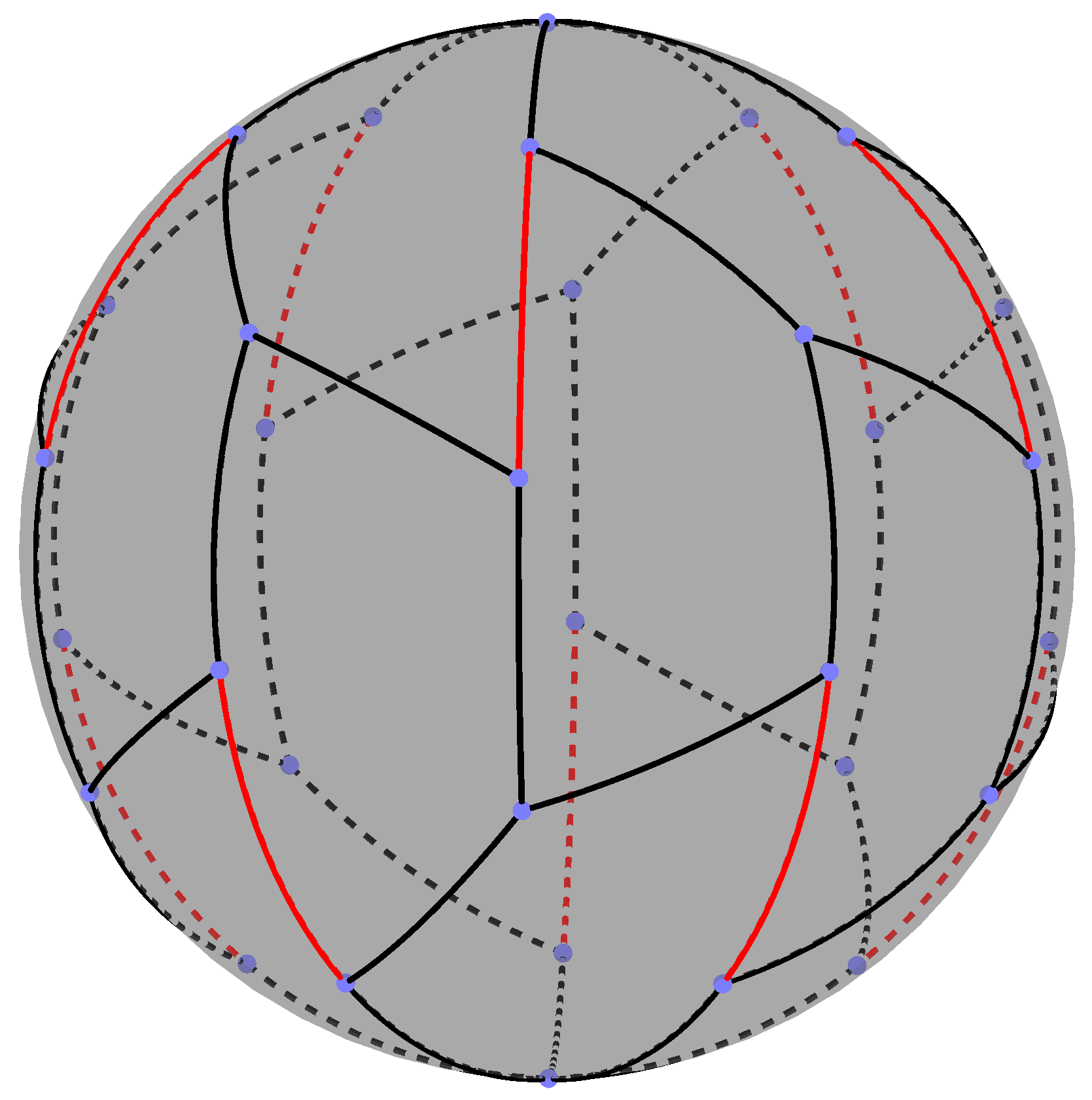}\hspace{10pt}
  \begin{overpic}[scale=0.040]{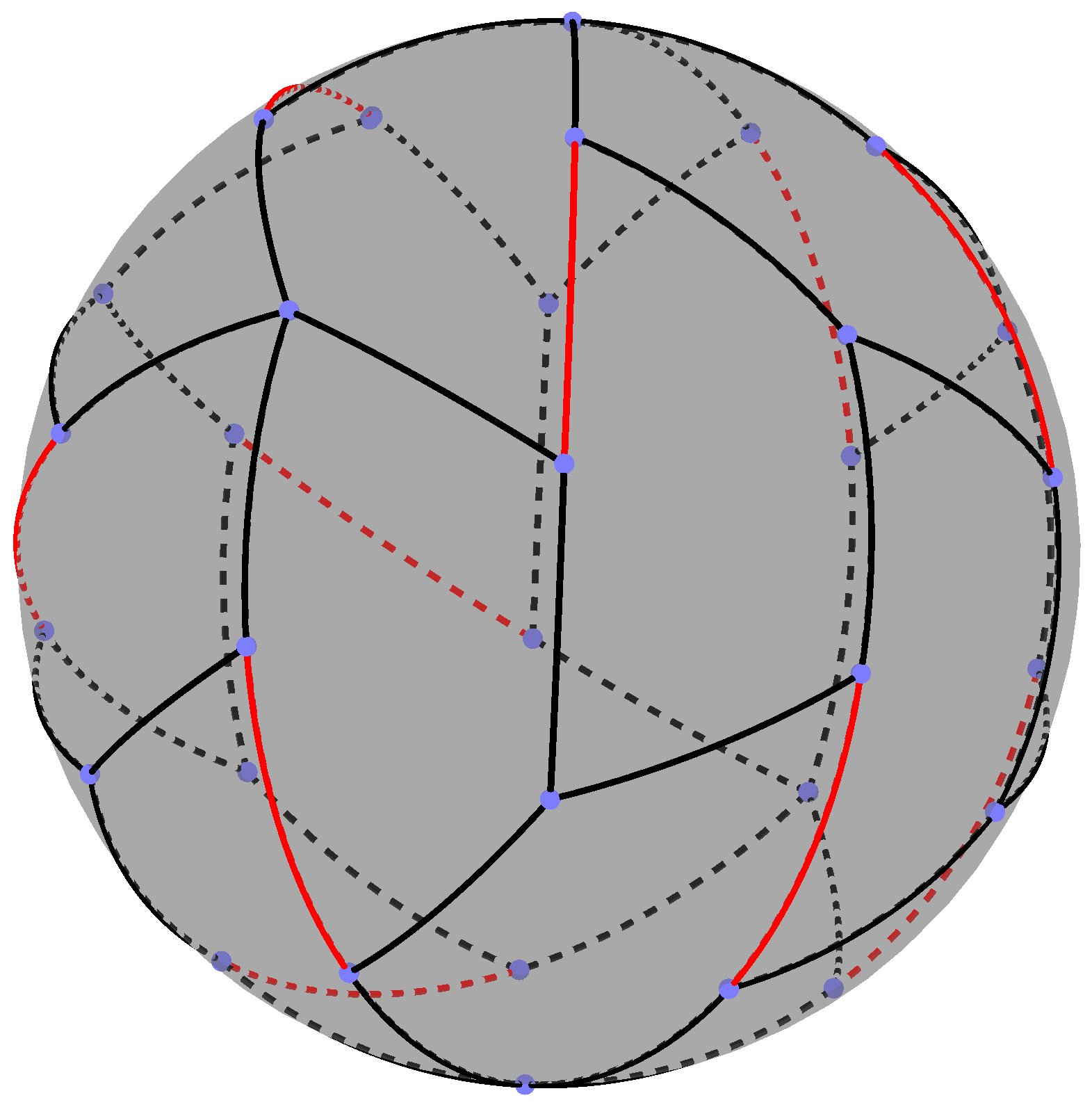}
   \put(-30,0){\small $\eee=\pi$}
  \end{overpic}

 \caption{Non-edge-to-edge quadrilateral tilings.} 
 \label{non} 
\end{figure}

When a pentagon has an angle $\pi$, it degenerates to a quadrilateral. This generates numerous new non-edge-to-edge quadrilateral tilings.
Figure \ref{non} shows nine different degenerations of three pentagonal subdivisions (still with a free parameter) and two different degenerations of the $a^4b$-pentagons in the non-symmetric 3-layer earth map tilings and their flips. For the pentagonal subdivisions of the tetrahedron ($f=12$), note that the first two degenerations in Figure \ref{non} are equivalent.  

For monohedral spherical tilings, there remain three open problems: non-edge-to-edge tilings of the sphere by congruent triangles, quadrilaterals, or pentagons. They prove to be more challenging after a few early explorations, and will be our next study goal. 

%%%%%%%%%%%%%%%%%%%%%%%%%%%%%%%%

\newpage
%%%%%%%%%%%%%%%%%%%%%%%%%%%%%%%%

\section*{Appendix:}

This Appendix finds all rational non-symmetric pentagons ($\bbb\neq\ccc$, $\ddd\neq\eee$) with  $f\ge12$ even and angles in $(0, 2\pi)$ in $a^4b$-tilings.

\subsection*{Rational angles solution for all cases with  three $a^2b$-vertex types: $6$ distinct cases}

\subsubsection*{Case $\{\aaa\ddd\eee,\bbb\ddd^2,\ccc\eee^2\}$, $(x={\mathrm e}^{2\mathrm{i} \delta},y={\mathrm e}^{\frac{2\mathrm{i} \pi}{f}})$}

$\zeta_{4} x^{4} y^{8}+x^{4} y^{7}+2 \zeta_{4} x^{4} y^{6}-x^{3} y^{7}-x^{4} y^{5}+2 \zeta_{4} x^{3} y^{6}+\zeta_{4} x^{4} y^{4}+2 x^{3} y^{5}-2 \zeta_{4} x^{3} y^{4}-2 x^{2} y^{5}+\zeta_{4}^{2} x^{3} y^{3}+8 \zeta_{4} x^{2} y^{4}+x y^{5}-2 \zeta_{4}^{2} x^{2} y^{3}-2 \zeta_{4} x y^{4}+2 \zeta_{4}^{2} x y^{3}+\zeta_{4} y^{4}+2 \zeta_{4} x y^{2}-\zeta_{4}^{2} y^{3}-\zeta_{4}^{2} x y+2 \zeta_{4} y^{2}+\zeta_{4}^{2} y+\zeta_{4}=0
$.

None.

\subsubsection*{Case $\{\aaa\ddd\eee,\bbb\ddd\eee,\ccc\eee^2\}$, $(x={\mathrm e}^{\frac{2\mathrm{i} \delta}{3}},y={\mathrm e}^{\frac{4\mathrm{i} \pi}{3f}})$}

$x^{6} y^{8}-x^{8} y^{3}-x^{7} y^{4}+2 x^{6} y^{5}-x^{5} y^{6}+x^{4} y^{7}-x^{6} y^{2}-2 x^{5} y^{3}-4 x^{4} y^{4}-2 x^{3} y^{5}-x^{2} y^{6}+x^{4} y-x^{3} y^{2}+2 x^{2} y^{3}-x y^{4}-y^{5}+x^{2}=0$.

None.

\subsubsection*{Case $\{\aaa\ddd^2,\bbb\ddd^2,\ccc\eee^2\}$, $(x={\mathrm e}^{2\mathrm{i} \delta},y={\mathrm e}^{\frac{4\mathrm{i} \pi}{f}})$}

$x^{10} y^{7}+x^{10} y^{6}-x^{9} y^{6}+x^{9} y^{5}+x^{8} y^{6}-2 x^{8} y^{5}-x^{8} y^{4}-2 x^{7} y^{5}+2 x^{7} y^{4}-3 x^{6} y^{4}-x^{6} y^{3}-2 x^{5} y^{4}+2 x^{5} y^{3}+x^{4} y^{4}+3 x^{4} y^{3}-2 x^{3} y^{3}+2 x^{3} y^{2}+x^{2} y^{3}+2 x^{2} y^{2}-x^{2} y-x y^{2}+x y-y-1=0$. 

None.

\subsubsection*{Case $\{\aaa\ddd^2,\bbb\ddd\eee,\ccc\eee^2\}$, $(x={\mathrm e}^{2\mathrm{i} \delta},y={\mathrm e}^{\frac{2\mathrm{i} \pi}{f}})$}

$x^{3} y^{7}+\zeta_{4} x^{2} y^{6}+x^{3} y^{5}-\zeta_{4} x^{3} y^{4}+x^{2} y^{5}+\zeta_{4} x^{2} y^{4}+x^{3} y^{3}+x y^{5}+2 \zeta_{4} x y^{4}+2 x^{2} y^{3}+\zeta_{4} x^{2} y^{2}+\zeta_{4} y^{4}+x y^{3}+\zeta_{4} x y^{2}-y^{3}+\zeta_{4} y^{2}+x y+\zeta_{4}=0$.

None.

\subsubsection*{Case $\{\aaa\eee^2,\bbb\ddd^2,\ccc\eee^2\}$, $(x={\mathrm e}^{\frac{2\mathrm{i} \delta}{3}},y={\mathrm e}^{\frac{4\mathrm{i} \pi}{3f}})$}

$x^{10} y^{7}+x^{9} y^{8}-x^{8} y^{9}+x^{10} y^{4}-x^{9} y^{5}-2 x^{8} y^{6}+2 x^{7} y^{7}-x^{6} y^{8}+x^{8} y^{3}-2 x^{7} y^{4}-3 x^{6} y^{5}+2 x^{5} y^{6}-2 x^{5} y^{3}+3 x^{4} y^{4}+2 x^{3} y^{5}-x^{2} y^{6}+x^{4} y-2 x^{3} y^{2}+2 x^{2} y^{3}+x y^{4}-y^{5}+x^{2}-x y-y^{2}=0$. 

None.

\subsubsection*{Case $\{\aaa\eee^2,\bbb\ddd\eee,\ccc\eee^2\}$, $(x={\mathrm e}^{\mathrm{i} \delta+\frac{\mathrm{i} \pi}{f}},y={\mathrm e}^{\mathrm{i} \delta-\frac{\mathrm{i} \pi}{f}})$}

$\zeta_{8}^{2} x^{3} y^{4}-\zeta_{8} x^{4} y^{2}+\zeta_{8}^{3} x^{3} y^{3}+\zeta_{8} x^{2} y^{4}+x^{5}+\zeta_{8}^{2} x^{4} y-x^{3} y^{2}+2 \zeta_{8}^{2} x^{2} y^{3}-x y^{4}-\zeta_{8}^{3} x^{4}+2 \zeta_{8} x^{3} y-\zeta_{8}^{3} x^{2} y^{2}+\zeta_{8} x y^{3}+\zeta_{8}^{3} y^{4}+\zeta_{8}^{2} x^{3}+x^{2} y-\zeta_{8}^{2} x y^{2}+\zeta_{8} x^{2}=0$. 

None.

\subsection*{Rational angles solution for all cases with two $a^2b$-vertex types: $110$ distinct cases}

\subsubsection*{Case $\{\aaa\ddd^2,\bbb\ddd^2,\aaa^2\ccc\,(\text{or}\,\aaa\bbb\ccc,\bbb^2\ccc)\}$, $(x={\mathrm e}^{2\mathrm{i} \delta},y={\mathrm e}^{\frac{4\mathrm{i} \pi}{f}})$}

$x^{7} y^{2}+x^{7} y-2 x^{6} y^{2}-x^{5} y^{3}-2 x^{6} y-x^{5} y^{2}+x^{6}+3 x^{5} y+2 x^{4} y^{2}+2 x^{3} y^{3}-2 x^{4} y+2 x^{3} y^{2}-2 x^{4}-2 x^{3} y-3 x^{2} y^{2}-x y^{3}+x^{2} y+2 x y^{2}+x^{2}+2 x y-y^{2}-y=0
$. 

None.

\subsubsection*{Case $\{\aaa\ddd^2,\bbb\ddd^2,\ccc^4\}$, $(x={\mathrm e}^{2\mathrm{i} \delta},y={\mathrm e}^{\frac{4\mathrm{i} \pi}{f}})$}

$2 x^{5} y^{3}-(\zeta_{4}+1) x^{5} y^{2}+(2 \zeta_{4}+1) x^{4} y^{2}+x^{4} y+(2 \zeta_{4}+1) x^{3} y^{2}-(\zeta_{4}+2) x^{3} y-(2 \zeta_{4}+1) x^{2} y^{2}+(\zeta_{4}+2) x^{2} y+\zeta_{4} x y^{2}+(\zeta_{4}+2) x y-(\zeta_{4}+1) y+2 \zeta_{4}=0$. 

None.

\subsubsection*{Case $\{\aaa\ddd^2,\bbb\ddd^2,\ccc^5\}$, $(x={\mathrm e}^{2\mathrm{i} \delta},y={\mathrm e}^{\frac{4\mathrm{i} \pi}{f}})$}

$(\zeta_{5}^{2}-2 \zeta_{5}+1) x^{5} y^{3}-(\zeta_{5}^{3}-\zeta_{5}^{2}) x^{5} y^{2}+(2 \zeta_{5}^{3}-2 \zeta_{5}^{2}+\zeta_{5}) x^{4} y^{2}+\zeta_{5}^{3} x^{4} y-(2 \zeta_{5}^{2}+\zeta_{5}) x^{3} y^{2}-(\zeta_{5}^{4}+2 \zeta_{5}^{3}) x^{3} y+(2 \zeta_{5}^{2}+\zeta_{5}) x^{2} y^{2}+(\zeta_{5}^{4}+2 \zeta_{5}^{3}) x^{2} y-\zeta_{5}^{2} x y^{2}-(\zeta_{5}^{4}-2 \zeta_{5}^{3}+2 \zeta_{5}^{2}) x y-(\zeta_{5}^{3}-\zeta_{5}^{2}) y+2 \zeta_{5}^{4}-\zeta_{5}^{3}-1=0$. 

None.

\subsubsection*{Case $\{\aaa\ddd^2,\bbb\ddd^2,\ccc^3\}$, $(x={\mathrm e}^{2\mathrm{i} \delta},y={\mathrm e}^{\frac{4\mathrm{i} \pi}{f}})$}

$(\zeta_{3}^{2}+\zeta_{3}-2) x^{5} y^{3}-(\zeta_{3}^{2}-\zeta_{3}) x^{5} y^{2}+(2 \zeta_{3}^{2}-2 \zeta_{3}+1) x^{4} y^{2}+\zeta_{3}^{2} x^{4} y-(2 \zeta_{3}+1) x^{3} y^{2}-(2 \zeta_{3}^{2}+1) x^{3} y+(2 \zeta_{3}+1) x^{2} y^{2}+(2 \zeta_{3}^{2}+1) x^{2} y-\zeta_{3} x y^{2}+(2 \zeta_{3}^{2}-2 \zeta_{3}-1) x y-(\zeta_{3}^{2}-\zeta_{3}) y+2-\zeta_{3}^{2}-\zeta_{3}=0$. 

$f=12,(\aaa,\bbb,\ccc,\ddd,\eee)=(4,4,6,7,9)/9$.

\subsubsection*{Case $\{\aaa\ddd^2,\bbb\ddd^2,\aaa\ccc^2\,(\text{or}\,\bbb\ccc^2)\}$, $(x={\mathrm e}^{\mathrm{i} \delta},y={\mathrm e}^{\frac{4\mathrm{i} \pi}{f}})$}

$(x-1) (x^{9} y^{2}-x^{8} y^{3}+x^{8} y^{2}-x^{7} y^{2}-x^{7} y-x^{6} y^{2}-x^{6} y-2 x^{5} y^{2}+x^{5} y-x^{4} y^{2}+2 x^{4} y+x^{3} y^{2}+x^{3} y+x^{2} y^{2}+x^{2} y-x y+x-y)=0
$.

$(\aaa,\bbb,\ccc,\ddd,\eee)=(4,4,f-2,f-2,f)/f$.
%for each $f\in\mathbb{Q}$.

\subsubsection*{Case $\{\aaa\ddd^2,\bbb\ddd^2,\aaa\ccc^3\,(\text{or}\,\bbb\ccc^3)\}$, $(x={\mathrm e}^{\frac{2\mathrm{i} \delta}{3}},y={\mathrm e}^{\frac{4\mathrm{i} \pi}{f}})$}

$(x-1) (x^{14} y^{2}-x^{13} y^{3}+x^{13} y^{2}+x^{12} y^{2}-x^{11} y^{2}-x^{11} y-x^{10} y^{2}-2 x^{10} y-2 x^{9} y^{2}-2 x^{9} y-3 x^{8} y^{2}-2 x^{7} y^{2}+2 x^{7} y+3 x^{6} y+2 x^{5} y^{2}+2 x^{5} y+2 x^{4} y^{2}+x^{4} y+x^{3} y^{2}+x^{3} y-x^{2} y-x y+x-y)=0$. 

$f=44,(\aaa,\bbb,\ccc,\ddd,\eee)=(4,4,6,9,11)/11$.

\subsubsection*{Case $\{\aaa\ddd^2,\bbb\ddd\eee,\aaa\bbb\ccc\}$, $(x={\mathrm e}^{2\mathrm{i} \delta},y={\mathrm e}^{\frac{4\mathrm{i} \pi}{f}})$}

$(y+1) (x^{3} y^{4}+x^{3} y^{3}+x^{3} y+x^{2} y^{2}+2 x y^{3}-2 x^{2} y-x y^{2}-y^{3}-y-1)=0
$. 

None.

\subsubsection*{Case $\{\aaa\ddd^2,\bbb\ddd\eee,\bbb^2\ccc\}$, $(x={\mathrm e}^{\mathrm{i} \ddd+\frac{2\mathrm{i} \pi}{f}},y={\mathrm e}^{\mathrm{i} \ddd-\frac{2\mathrm{i} \pi}{f}})$}

$x^{8} y^{2}+2 \zeta_{4} x^{7} y^{2}+\zeta_{4} x^{6} y^{3}+x^{6} y^{2}-2 x^{5} y^{3}-\zeta_{4} x^{6} y+3 \zeta_{4} x^{5} y^{2}-\zeta_{4} x^{4} y^{3}+3 x^{5} y+x^{4} y^{2}+\zeta_{4} x^{4} y+3 \zeta_{4} x^{3} y^{2}-x^{4}+3 x^{3} y-x^{2} y^{2}-2 \zeta_{4} x^{3}+\zeta_{4} x^{2} y+x^{2}+2 x y+\zeta_{4} y=0$. 

None.

\subsubsection*{Case $\{\aaa\ddd^2,\bbb\ddd\eee,\aaa\ccc^3\}$, $(x={\mathrm e}^{\mathrm{i} \eee+\frac{\mathrm{i} \pi}{f}},y={\mathrm e}^{\mathrm{i} \eee-\frac{\mathrm{i} \pi}{f}})$}

$\zeta_{8} x^{8} y^{4}-\zeta_{8} x^{6} y^{6}+2 \zeta_{8}^{3} x^{5} y^{7}+\zeta_{8} x^{2} y^{10}-\zeta_{8}^{3} x y^{11}+2 \zeta_{8}^{2} x^{9} y^{2}-x^{6} y^{5}-\zeta_{8}^{2} x^{5} y^{6}-x^{4} y^{7}+2 \zeta_{8}^{2} x^{3} y^{8}+2 x^{2} y^{9}+y^{11}+\zeta_{8}^{3} x^{10}+2 \zeta_{8}^{3} x^{8} y^{2}+2 \zeta_{8} x^{7} y^{3}-\zeta_{8}^{3} x^{6} y^{4}-\zeta_{8} x^{5} y^{5}-\zeta_{8}^{3} x^{4} y^{6}+2 \zeta_{8} x y^{9}-x^{9}+\zeta_{8}^{2} x^{8} y+2 x^{5} y^{4}-\zeta_{8}^{2} x^{4} y^{5}+\zeta_{8}^{2} x^{2} y^{7}=0
$. 

None.

\subsubsection*{Case $\{\aaa\ddd^2,\bbb\ddd\eee,\bbb\ccc^3\}$, $(x={\mathrm e}^{2\mathrm{i} \ddd},y={\mathrm e}^{\frac{4\mathrm{i} \pi}{f}})$}

$x^{12} y^{11}+x^{11} y^{10}+2 x^{10} y^{10}-x^{9} y^{10}+x^{9} y^{9}+2 x^{9} y^{8}+x^{9} y^{7}-x^{8} y^{8}+2 x^{8} y^{7}+x^{7} y^{7}+x^{7} y^{6}-2 y^{7} x^{6}+2 x^{6} y^{4}-x^{5} y^{5}-x^{5} y^{4}-2 x^{4} y^{4}+x^{4} y^{3}-x^{3} y^{4}-2 x^{3} y^{3}-x^{3} y^{2}+x^{3} y-2 x^{2} y-x y-1=0$. 

None.

\subsubsection*{Case $\{\aaa\ddd^2,\bbb\ddd\eee,\ccc^5\}$, $(x={\mathrm e}^{\mathrm{i} \eee+\frac{2\mathrm{i} \pi}{f}},y={\mathrm e}^{\mathrm{i} \eee-\frac{2\mathrm{i} \pi}{f}})$}

$\zeta_{20}^{4} x^{3} y^{2}+(2 \zeta_{20}^{8}-\zeta_{20}^{4}) x^{2} y^{3}-(\zeta_{20}^{6}-\zeta_{20}^{2}) x y^{4}+2 \zeta_{20}^{7} x^{3} y-(\zeta_{20}^{5}+\zeta_{20}) x^{2} y^{2}-(2 \zeta_{20}^{9}-2 \zeta_{20}^{5}+\zeta_{20}) x y^{3}-\zeta_{20}^{9} y^{4}-x^{3}-(\zeta_{20}^{8}-2 \zeta_{20}^{4}+2) x^{2} y-(\zeta_{20}^{8}+\zeta_{20}^{4}) x y^{2}+2 \zeta_{20}^{2} y^{3}+(\zeta_{20}^{7}-\zeta_{20}^{3}) x^{2}-(\zeta_{20}^{5}-2 \zeta_{20}) x y+\zeta_{20}^{5} y^{2}=0$. 

None.

\subsubsection*{Case $\{\aaa\ddd^2,\bbb\ddd\eee,\aaa\ccc^2\}$, $(x={\mathrm e}^{\mathrm{i} \eee},y={\mathrm e}^{\frac{4\mathrm{i} \pi}{f}})$}

$(y+1) (x^{3} y^{4}+2 x^{2} y^{5}+x y^{6}-2 x^{3} y^{3}-2 x^{2} y^{4}+x y^{5}+y^{6}+x^{2} y^{3}+x y^{4}-x^{2} y^{2}-x y^{3}-x^{3}-x^{2} y+2 x y^{2}+2 y^{3}-x^{2}-2 x y-y^{2})=0
$. 

$f=14,(\aaa,\bbb,\ccc,\ddd,\eee)=(4,3,5,5,6)/7$.

\subsubsection*{Case $\{\aaa\ddd^2,\bbb\ddd\eee,\ccc^3\}$, $(x={\mathrm e}^{\mathrm{i} \eee+\frac{2\mathrm{i} \pi}{f}},y={\mathrm e}^{\mathrm{i} \eee-\frac{2\mathrm{i} \pi}{f}})$}

$\zeta_{12} x^{3} y^{2}+(2 \zeta_{12}^{5}-\zeta_{12}) x^{2} y^{3}+(\zeta_{12}^{3}+\zeta_{12}) x y^{4}+2 \zeta_{12}^{2} x^{3} y-(\zeta_{12}^{4}+1) x^{2} y^{2}+(2 \zeta_{12}^{4}+2 \zeta_{12}^{2}-1) x y^{3}+\zeta_{12}^{2} y^{4}+\zeta_{12}^{3} x^{3}-(\zeta_{12}^{5}-2 \zeta_{12}^{3}-2 \zeta_{12}) x^{2} y-(\zeta_{12}^{5}+\zeta_{12}) x y^{2}+2 \zeta_{12}^{3} y^{3}+(\zeta_{12}^{4}+\zeta_{12}^{2}) x^{2}-(\zeta_{12}^{4}-2) x y+\zeta_{12}^{4} y^{2}=0
$.

$f=18,(\aaa,\bbb,\ccc,\ddd,\eee)=(10,8,12,13,15)/18$.

\subsubsection*{Case $\{\aaa\ddd^2,\bbb\ddd\eee,\bbb\ccc^2\}$, $(x={\mathrm e}^{2\mathrm{i} \ddd},y={\mathrm e}^{\frac{4\mathrm{i} \pi}{f}})$}

$(x y-1) (x y+1) (x^{6} y^{5}+2 x^{4} y^{4}-x^{4} y^{3}-x^{3} y^{4}-x^{4} y^{2}+x^{2} y^{3}+x^{3} y+x^{2} y^{2}-2 x^{2} y-1)=0$.

$f=24,(\aaa,\bbb,\ccc,\ddd,\eee)=(4,6,3,4,2)/6$.

\subsubsection*{Case $\{\aaa\ddd^2,\bbb\ddd\eee,\ccc^4\}$, $(x={\mathrm e}^{\mathrm{i} \eee+\frac{2\mathrm{i} \pi}{f}},y={\mathrm e}^{\mathrm{i} \eee-\frac{2\mathrm{i} \pi}{f}})$}

$\zeta_{8} x^{3} y^{2}+(2 \zeta_{8}^{3}-\zeta_{8}) x^{2} y^{3}-(\zeta_{8}^{3}-\zeta_{8}) x y^{4}+2 \zeta_{8}^{2} x^{3} y-(\zeta_{8}^{2}+1) x^{2} y^{2}+(2 \zeta_{8}^{2}+1) x y^{3}+y^{4}+\zeta_{8}^{3} x^{3}+(\zeta_{8}^{3}+2 \zeta_{8}) x^{2} y-(\zeta_{8}^{3}+\zeta_{8}) x y^{2}+2 \zeta_{8} y^{3}+(\zeta_{8}^{2}-1) x^{2}-(\zeta_{8}^{2}-2) x y+\zeta_{8}^{2} y^{2}=0
$.

$f=24,(\aaa,\bbb,\ccc,\ddd,\eee)=(4,6,3,4,2)/6$.

\subsubsection*{Case $\{\aaa\ddd^2,\bbb\ddd\eee,\aaa^2\ccc\}$, $(x={\mathrm e}^{\mathrm{i} \eee+\frac{2\mathrm{i} \pi}{f}},y={\mathrm e}^{\mathrm{i} \eee-\frac{2\mathrm{i} \pi}{f}})$}

$(x+y) (\zeta_{4} x^{7}-\zeta_{4} x^{6} y+2 \zeta_{4} x^{4} y^{3}-\zeta_{4} x^{3} y^{4}-2 x^{6}+2 x^{4} y^{2}-3 x^{3} y^{3}-3 \zeta_{4} x^{4} y+2 \zeta_{4} x^{3} y^{2}-2 \zeta_{4} x y^{4}-x^{4}+2 x^{3} y-x y^{3}+y^{4})=0$.

$f=40,(\aaa,\bbb,\ccc,\ddd,\eee)=(18,6,4,11,23)/20$.

\subsubsection*{Case $\{\aaa\ddd^2,\bbb\eee^2,\aaa^2\ccc\}$, $(x={\mathrm e}^{2\mathrm{i} \ddd},y={\mathrm e}^{\frac{4\mathrm{i} \pi}{f}})$}

$(x-1) (x^{3}-y^{2}) (x^{9}-x^{7}-2 x^{6} y-x^{5} y^{2}-x^{6}+x^{3} y^{3}+x^{4} y+2 x^{3} y^{2}+x^{2} y^{3}-y^{3}) =0$.

None.

\subsubsection*{Case $\{\aaa\ddd^2,\bbb\eee^2,\bbb^2\ccc\}$, $(x={\mathrm e}^{2\mathrm{i} \eee},y={\mathrm e}^{\frac{4\mathrm{i} \pi}{f}})$}

$(x-1) (x^{9} y+x^{8} y-x^{6} y^{3}-x^{8}-2 x^{5} y^{3}-x^{4} y^{4}+x^{6} y+2 x^{5} y^{2}-2 x^{4} y^{3}-x^{3} y^{4}+x^{5} y+2 x^{4} y^{2}+x y^{5}+x^{3} y^{2}-x y^{4}-y^{4})=0$.

None.

\subsubsection*{Case $\{\aaa\ddd^2,\bbb\eee^2,\aaa\ccc^2\}$, $(x={\mathrm e}^{\mathrm{i} \ddd},y={\mathrm e}^{\frac{4\mathrm{i} \pi}{f}})$}

$(y-1) (y+1) (x^{5} y^{3}-x^{5} y+x^{4} y^{2}+x^{5}-2 x^{3} y^{2}-x^{2} y^{3}-x^{4}+2 x^{3} y+2 x^{2} y^{2}-x y^{3}-x^{3}-2 x^{2} y+y^{3}+x y-y^{2}+1)=0$.

$f=12,(\aaa,\bbb,\ccc,\ddd,\eee)=(6,4,3,3,4)/6$.

\subsubsection*{Case $\{\aaa\ddd^2,\bbb\eee^2,\aaa\bbb\ccc\}$, $(x={\mathrm e}^{2\mathrm{i} \eee},y={\mathrm e}^{\frac{4\mathrm{i} \pi}{f}})$}

$(x-1) (x^{2} y^{7}-x^{3} y^{5}-x^{3} y^{4}-2 x^{2} y^{5}-x y^{6}+2 x^{2} y^{4}+x y^{5}+y^{6}+x^{4} y+x^{3} y^{2}+2 x^{2} y^{3}-x^{3} y-2 x^{2} y^{2}-x y^{3}-x y^{2}+x^{2})=0$.

$f=12,(\aaa,\bbb,\ccc,\ddd,\eee)=(2,6,4,5,3)/6$.

\subsubsection*{Case $\{\aaa\ddd^2,\bbb\eee^2,\ccc^3\}$, $(x={\mathrm e}^{2\mathrm{i} \ddd},y={\mathrm e}^{\frac{4\mathrm{i} \pi}{f}})$}

$\zeta_{3}^{2} x^{5} y^{5}-\zeta_{3} x^{4} y^{5}-(\zeta_{3}^{2}-1) x^{4} y^{4}-(\zeta_{3}^{2}-\zeta_{3}) x^{3} y^{5}-2 \zeta_{3} x^{4} y^{3}-(\zeta_{3}-2) x^{3} y^{4}+2 \zeta_{3} x^{3} y^{3}-x^{2} y^{4}-(2 \zeta_{3}^{2}-\zeta_{3}-1) x^{3} y^{2}-(\zeta_{3}^{2}-2 \zeta_{3}+1) x^{2} y^{3}+x^{3} y-2 \zeta_{3}^{2} x^{2} y^{2}+(\zeta_{3}^{2}-2) x^{2} y+2 \zeta_{3}^{2} x y^{2}-(\zeta_{3}^{2}-\zeta_{3}) x^{2}+(\zeta_{3}-1) x y+\zeta_{3}^{2} x-\zeta_{3}=0$.

$f=12,(\aaa,\bbb,\ccc,\ddd,\eee)=(2,6,4,5,3)/6$.

\subsubsection*{Case $\{\aaa\ddd^2,\bbb\eee^2,\aaa\ccc^3\}$, $(x={\mathrm e}^{\frac{2\mathrm{i} \ddd}{3}},y={\mathrm e}^{\frac{4\mathrm{i} \pi}{f}})$}

$(x-1)(xy^2-1) (x^8y^3 + x^7y^3 + x^6y^3 - x^7*y + x^6y^2 + x^5y^3 - x^6y + x^6 - x^5y - x^4y^2 + x^4y + x^3y^2 - x^2y^3 + x^2y^2 - x^3 - x^2y + xy^2 - x^2 - x - 1)=0$.

$f=18,(\aaa,\bbb,\ccc,\ddd,\eee)=(6,18,10,15,9)/18$.

$f=28,(\aaa,\bbb,\ccc,\ddd,\eee)=(2,6,4,6,4)/7$.

\subsubsection*{Case $\{\aaa\ddd^2,\bbb\eee^2,\bbb\ccc^2\}$, $(x={\mathrm e}^{\mathrm{i} \eee},y={\mathrm e}^{\frac{4\mathrm{i} \pi}{f}})$}

$(x+1) (x-1) (x^{3} y^{5}-x^{4} y^{3}-x^{3} y^{4}+x^{2} y^{5}-x^{3} y^{3}+x^{4} y-2 x^{2} y^{3}+2 x^{2} y^{2}-y^{4}+x y^{2}-x^{2}+x y+y^{2}-x)=0$.

$f=24,(\aaa,\bbb,\ccc,\ddd,\eee)=(2,6,3,5,3)/6$.

$f=28,(\aaa,\bbb,\ccc,\ddd,\eee)=(2,6,4,6,4)/7$.

\subsubsection*{Case $\{\aaa\ddd^2,\bbb\eee^2,\ccc^4\}$, $(x={\mathrm e}^{2\mathrm{i} \eee},y={\mathrm e}^{\frac{4\mathrm{i} \pi}{f}})$}

$\zeta_{4} x^{5} y^{5}-2 \zeta_{4} x^{4} y^{5}+(\zeta_{4}-1) x^{4} y^{4}+\zeta_{4} x^{3} y^{5}-x^{4} y^{3}+2 x^{3} y^{4}-2 \zeta_{4} x^{3} y^{3}-(\zeta_{4}+1) x^{2} y^{4}+(\zeta_{4}+2) x^{3} y^{2}+(2 \zeta_{4}+1) x^{2} y^{3}-(\zeta_{4}+1) x^{3} y-2 x^{2} y^{2}+2 \zeta_{4} x^{2} y-\zeta_{4} x y^{2}+x^{2}-(\zeta_{4}-1) x y-2 x+1=0$.

$f=24,(\aaa,\bbb,\ccc,\ddd,\eee)=(2,6,3,5,3)/6$.
\subsubsection*{Case $\{\aaa\ddd^2,\bbb\eee^2,\bbb\ccc^3\}$, $(x={\mathrm e}^{\frac{2\mathrm{i} \eee}{3}},y={\mathrm e}^{\frac{4\mathrm{i} \pi}{f}})$}

$x^{10} y^{5}-x^{10} y^{3}-x^{9} y^{4}+x^{8} y^{4}-2 x^{7} y^{5}+x^{9} y+x^{6} y^{4}-x^{8} y+2 x^{7} y^{2}+2 x^{6} y^{3}-2 x^{5} y^{4}+x^{4} y^{5}-x^{6} y^{2}+x^{4} y^{4}-x^{6} y+x^{4} y^{3}-x^{6}+2 x^{5} y-2 x^{4} y^{2}-2 x^{3} y^{3}+x^{2} y^{4}-x^{4} y-x y^{4}+2 x^{3}-x^{2} y+x y+y^{2}-1=0$.

$f=48,(\aaa,\bbb,\ccc,\ddd,\eee)=(6,42,2,21,3)/24$.

\subsubsection*{Case $\{\aaa\ddd^2,\bbb\eee^2,\ccc^5\}$, $(x={\mathrm e}^{2\mathrm{i} \eee},y={\mathrm e}^{\frac{4\mathrm{i} \pi}{f}})$}

$(x-1) (x^{4} y^{5}-x^{3} y^{5}-(\zeta_{5}^{2}-\zeta_{5}) x^{3} y^{4}-\zeta_{5}^{3} x^{3} y^{3}-(\zeta_{5}-1) x^{2} y^{4}-(\zeta_{5}^{3}-2 \zeta_{5}^{2}) x^{2} y^{3}+(2 \zeta_{5}^{3}-\zeta_{5}^{2}) x^{2} y^{2}-(\zeta_{5}^{4}-1) x^{2} y-\zeta_{5}^{2} x y^{2}+(\zeta_{5}^{4}-\zeta_{5}^{3}) x y-x+1)=0
$.

$f=60,(\aaa,\bbb,\ccc,\ddd,\eee)=(10,30,12,25,15)/30$.

\subsubsection*{Case $\{\aaa\ddd^2,\ccc\ddd^2,\aaa^2\bbb\,(\text{or}\,\aaa\bbb\ccc,\bbb\ccc^2)\}$, $(x={\mathrm e}^{2\mathrm{i} \eee},y={\mathrm e}^{\frac{4\mathrm{i} \pi}{f}})$}

$y^{15}-2 x y^{13}-x^{3} y^{10}-x^{2} y^{11}-x y^{12}+2 x^{4} y^{8}+3 x^{3} y^{9}+3 x^{2} y^{10}-x^{5} y^{6}+x^{4} y^{7}-x^{3} y^{8}+x^{2} y^{9}-3 x^{5} y^{5}-3 x^{4} y^{6}-2 x^{3} y^{7}+x^{6} y^{3}+x^{5} y^{4}+x^{4} y^{5}+2 x^{6} y^{2}-x^{7}=0$.

None.

\subsubsection*{Case $\{\aaa\ddd^2,\ccc\ddd^2,\bbb^3\}$, $(x={\mathrm e}^{2\mathrm{i} \ddd},y={\mathrm e}^{\frac{4\mathrm{i} \pi}{f}})$}

$(x-1) (x^{6} y^{3}-x^{5} y^{3}-(\zeta_{3}^{2}-\zeta_{3}) x^{5} y^{2}+(\zeta_{3}^{2}-\zeta_{3}+1) x^{4} y^{2}+(\zeta_{3}^{2}-2 \zeta_{3}) x^{3} y^{2}-(2 \zeta_{3}^{2}-\zeta_{3}) x^{3} y-(\zeta_{3}^{2}-\zeta_{3}-1) x^{2} y+(\zeta_{3}^{2}-\zeta_{3}) x y-x+1)=0$

None.

\subsubsection*{Case $\{\aaa\ddd^2,\ccc\ddd^2,\aaa\bbb^3\,(\text{or}\,\bbb^3\ccc)\}$, $(x={\mathrm e}^{\frac{2\mathrm{i} \ddd}{3}},y={\mathrm e}^{\frac{4\mathrm{i} \pi}{f}})$}

$(x-1) (x^{2}+x+1) (x^{14} y^{3}+x^{13} y^{3}-x^{13} y^{2}+x^{12} y^{3}+x^{10} y^{2}+x^{8} y^{2}+2 x^{7} y^{2}-2 x^{7} y-x^{6} y-x^{4} y-x^{2}+x y-x-1)=0
$.

None.

\subsubsection*{Case $\{\aaa\ddd^2,\ccc\ddd^2,\bbb^5\}$, $(x={\mathrm e}^{2\mathrm{i} \ddd},y={\mathrm e}^{\frac{4\mathrm{i} \pi}{f}})$}

$(x-1)(\zeta_{5}^{3} x^{6} y^{3}-\zeta_{5}^{3} x^{5} y^{3}+(\zeta_{5}^{4}-1) x^{5} y^{2}-(\zeta_{5}^{4}-\zeta_{5}^{3}-1) x^{4} y^{2}-(2 \zeta_{5}^{4}-1) x^{3} y^{2}+(\zeta_{5}^{4}-2) x^{3} y+(\zeta_{5}^{4}+\zeta_{5}-1) x^{2} y-(\zeta_{5}^{4}-1) x y-\zeta_{5} x+\zeta_{5})=0$.

None.

\subsubsection*{Case $\{\aaa\ddd^2,\ccc\ddd^2,\aaa\bbb^2\,(\text{or}\,\bbb^2\ccc)\}$, $(x={\mathrm e}^{\frac{\mathrm{i} \eee}{2}+\frac{2\mathrm{i} \pi}{f}},y={\mathrm e}^{\frac{\mathrm{i} \eee}{2}-\frac{2\mathrm{i} \pi}{f}})$}

$(y^{2}+1) (\zeta_{4} x^{3} y^{7}-\zeta_{4} x^{2} y^{8}+x^{2} y^{7}+\zeta_{4} x^{3} y^{5}-\zeta_{4} x^{2} y^{6}+x^{2} y^{5}+2 \zeta_{4} x^{2} y^{4}-2 \zeta_{4} x y^{5}-2 x^{2} y^{3}+2 x y^{4}+\zeta_{4} x y^{3}-x y^{2}+y^{3}+\zeta_{4} x y-x+y)=0$.

$f=24,(\aaa,\bbb,\ccc,\ddd,\eee)=(6,3,6,3,1)/6$.

\subsubsection*{Case $\{\aaa\ddd^2,\ccc\ddd^2,\bbb^4\}$, $(x={\mathrm e}^{\frac{2\mathrm{i} \eee}{3}},y={\mathrm e}^{\frac{4\mathrm{i} \pi}{3f}})$}

$(y^2 + x) (\zeta_{4} x^{5} y^{3}+\zeta_{4} x^{6} y+(\zeta_{4}-2) x^{3} y^{4}+x^{4} y^{2}+(\zeta_{4}-1) x y^{5}-(\zeta_{4}-1) x^{5}+\zeta_{4} x^{2} y^{3}-(2 \zeta_{4}-1) x^{3} y+y^{4}+x y^{2})=0$.

$f=24,(\aaa,\bbb,\ccc,\ddd,\eee)=(6,3,6,3,1)/6$.

\subsubsection*{Case $\{\aaa\ddd^2,\ccc\ddd\eee,\aaa^2\bbb\}$, $(x={\mathrm e}^{\mathrm{i} \eee},y={\mathrm e}^{\frac{2\mathrm{i} \pi}{f}})$}

$\zeta_{4} x^{4} y^{13}+\zeta_{4} x^{4} y^{11}-2 x^{3} y^{12}+\zeta_{4} x^{2} y^{13}+x^{5} y^{8}-\zeta_{4} x^{4} y^{9}-x^{3} y^{10}+2 \zeta_{4} x^{4} y^{7}+3 x^{3} y^{8}-3 \zeta_{4} x^{2} y^{9}-2 x^{3} y^{6}-2 \zeta_{4} x^{2} y^{7}-3 x^{3} y^{4}+3 \zeta_{4} x^{2} y^{5}+2 x y^{6}-\zeta_{4} x^{2} y^{3}-x y^{4}+\zeta_{4} y^{5}+x^{3}-2 \zeta_{4} x^{2} y+x y^{2}+x=0$.

None.

\subsubsection*{Case $\{\aaa\ddd^2,\ccc\ddd\eee,\bbb^3\}$, $(x={\mathrm e}^{\mathrm{i} \eee+\frac{2\mathrm{i} \pi}{f}},y={\mathrm e}^{\mathrm{i} \eee-\frac{2\mathrm{i} \pi}{f}})$}

$\zeta_{12}^{4} x^{2} y^{4}+(2 \zeta_{12}^{5}-\zeta_{12}^{3}) x^{2} y^{3}+(\zeta_{12}^{3}+\zeta_{12}) x y^{4}+(\zeta_{12}^{4}+2 \zeta_{12}^{2}+1) x^{3} y+(-\zeta_{12}^{4}-2) x^{2} y^{2}+(2 \zeta_{12}^{4}+2 \zeta_{12}^{2}-1) x y^{3}-(\zeta_{12}^{5}-2 \zeta_{12}^{3}-2 \zeta_{12}) x^{2} y-(2 \zeta_{12}^{5}+\zeta_{12}) x y^{2}+(\zeta_{12}^{5}+2 \zeta_{12}^{3}+\zeta_{12}) y^{3}+(\zeta_{12}^{4}+\zeta_{12}^{2}) x^{2}-(\zeta_{12}^{2}-2) x y+\zeta_{12} x=0$.

None.

\subsubsection*{Case $\{\aaa\ddd^2,\ccc\ddd\eee,\bbb^2\ccc\}$, $(x={\mathrm e}^{2\mathrm{i} \ddd},y={\mathrm e}^{\frac{4\mathrm{i} \pi}{f}})$}

$(y-1) (x y+1) (x^{7} y^{7}-x^{6} y^{6}-x^{5} y^{5}+2 x^{4} y^{4}+x^{4} y^{3}-x^{3} y^{4}-2 x^{3} y^{3}+x^{2} y^{2}+x y-1)=0$.

None.

\subsubsection*{Case $\{\aaa\ddd^2,\ccc\ddd\eee,\aaa\bbb^3\}$, $(x={\mathrm e}^{\mathrm{i} \eee+\frac{\mathrm{i} \pi}{f}},y={\mathrm e}^{\mathrm{i} \eee-\frac{\mathrm{i} \pi}{f}})$}

$\zeta_{8}^{2} x^{4} y^{9}+2 \zeta_{8}^{3} x^{5} y^{7}-\zeta_{8} x^{4} y^{8}+\zeta_{8} x^{2} y^{10}-\zeta_{8}^{3} x y^{11}-x^{10} y+2 \zeta_{8}^{2} x^{9} y^{2}+x^{8} y^{3}-2 x^{6} y^{5}-\zeta_{8}^{2} x^{5} y^{6}-x^{4} y^{7}+2 \zeta_{8}^{2} x^{3} y^{8}+2 x^{2} y^{9}+2 \zeta_{8}^{3} x^{8} y^{2}+2 \zeta_{8} x^{7} y^{3}-\zeta_{8}^{3} x^{6} y^{4}-\zeta_{8} x^{5} y^{5}-2 \zeta_{8}^{3} x^{4} y^{6}+\zeta_{8}^{3} x^{2} y^{8}+2 \zeta_{8} x y^{9}-\zeta_{8}^{3} y^{10}-x^{9}+\zeta_{8}^{2} x^{8} y-\zeta_{8}^{2} x^{6} y^{3}+2 x^{5} y^{4}+\zeta_{8} x^{6} y^{2}=0$.

None.

\subsubsection*{Case $\{\aaa\ddd^2,\ccc\ddd\eee,\bbb^4\}$, $(x={\mathrm e}^{\mathrm{i} \eee+\frac{2\mathrm{i} \pi}{f}},y={\mathrm e}^{\mathrm{i} \eee-\frac{2\mathrm{i} \pi}{f}})$}

$\zeta_{8}^{2} x^{2} y^{4}+(2 \zeta_{8}^{3}-\zeta_{8}) x^{2} y^{3}-(\zeta_{8}^{3}-\zeta_{8}) x y^{4}+2 \zeta_{8}^{2} x^{3} y-(\zeta_{8}^{2}+2) x^{2} y^{2}+(2 \zeta_{8}^{2}+1) x y^{3}+(\zeta_{8}^{3}+2 \zeta_{8}) x^{2} y-(2 \zeta_{8}^{3}+\zeta_{8}) x y^{2}+2 \zeta_{8} y^{3}+(\zeta_{8}^{2}-1) x^{2}-(\zeta_{8}^{2}-2) x y+\zeta_{8} x=0$.

none

\subsubsection*{Case $\{\aaa\ddd^2,\ccc\ddd\eee,\bbb^3\ccc\}$, $(x={\mathrm e}^{2\mathrm{i} \ddd},y={\mathrm e}^{\frac{4\mathrm{i} \pi}{f}})$}

$x^{14} y^{15}+x^{13} y^{13}-x^{12} y^{13}+x^{12} y^{12}+2 x^{11} y^{12}+2 x^{10} y^{10}-x^{9} y^{10}+2 x^{9} y^{9}-x^{8} y^{10}+2 x^{8} y^{9}+x^{8} y^{8}-2 x^{7} y^{9}+x^{8} y^{7}-x^{6} y^{8}+2 x^{7} y^{6}-x^{6} y^{7}-2 x^{6} y^{6}+x^{6} y^{5}-2 x^{5} y^{6}+x^{5} y^{5}-2 x^{4} y^{5}-2 x^{3} y^{3}-x^{2} y^{3}+x^{2} y^{2}-x y^{2}-1=0$.

None.

\subsubsection*{Case $\{\aaa\ddd^2,\ccc\ddd\eee,\bbb^5\}$, $(x={\mathrm e}^{\mathrm{i} \eee+\frac{2\mathrm{i} \pi}{f}},y={\mathrm e}^{\mathrm{i} \eee-\frac{2\mathrm{i} \pi}{f}})$}

$\zeta_{20}^{5} x^{2} y^{4}+(2 \zeta_{20}^{8}-\zeta_{20}^{2}) x^{2} y^{3}-(\zeta_{20}^{6}-\zeta_{20}^{2}) x y^{4}+(2 \zeta_{20}^{7}-\zeta_{20}^{3}+\zeta_{20}) x^{3} y-(\zeta_{20}^{5}+2 \zeta_{20}) x^{2} y^{2}-(2 \zeta_{20}^{9}-2 \zeta_{20}^{5}+\zeta_{20}) x y^{3}-(\zeta_{20}^{8}-2 \zeta_{20}^{4}+2) x^{2} y-(2 \zeta_{20}^{8}+\zeta_{20}^{4}) x y^{2}+(\zeta_{20}^{8}-\zeta_{20}^{6}+2 \zeta_{20}^{2}) y^{3}+(\zeta_{20}^{7}-\zeta_{20}^{3}) x^{2}-(\zeta_{20}^{7}-2 \zeta_{20}) x y+\zeta_{20}^{4} x=0$.

None.

\subsubsection*{Case $\{\aaa\ddd^2,\ccc\ddd\eee,\aaa\bbb\ccc\}$, $(x={\mathrm e}^{2\mathrm{i} \eee},y={\mathrm e}^{\frac{4\mathrm{i} \pi}{f}})$}

$x y^{11}-x^{2} y^{9}+2 x y^{10}+y^{11}+x^{3} y^{7}+x^{3} y^{6}+x^{2} y^{7}+2 x y^{8}+2 x^{3} y^{5}-2 x^{2} y^{6}-2 x^{4} y^{3}-x^{3} y^{4}-x^{2} y^{5}-x^{2} y^{4}-x^{5}-2 x^{4} y+x^{3} y^{2}-x^{4}=0$.

$f=12,(\aaa,\bbb,\ccc,\ddd,\eee)=(2,6,4,5,3)/6$.

\subsubsection*{Case $\{\aaa\ddd^2,\ccc\ddd\eee,\aaa\bbb^2\}$, $(x={\mathrm e}^{\mathrm{i} \eee},y={\mathrm e}^{\frac{4\mathrm{i} \pi}{f}})$}

$x^{3} y^{7}+2 x^{3} y^{6}+x^{3} y^{5}+2 x^{2} y^{6}+x y^{7}-x^{5} y^{2}-2 x^{4} y^{3}-2 x^{3} y^{4}+2 x^{2} y^{5}+x y^{6}+x^{4} y^{2}-x y^{5}-x^{4} y-2 x^{3} y^{2}+2 x^{2} y^{3}+2 x y^{4}+y^{5}-x^{4}-2 x^{3} y-x^{2} y^{2}-2 x^{2} y-x^{2}=0$.

$f=28,(\aaa,\bbb,\ccc,\ddd,\eee)=(2,6,4,6,4)/7$.

\subsubsection*{Case $\{\aaa\ddd^2,\ccc\ddd\eee,\bbb\ccc^2\}$, $(x={\mathrm e}^{\mathrm{i} \eee+\frac{2\mathrm{i} \pi}{f}},y={\mathrm e}^{\mathrm{i} \eee-\frac{2\mathrm{i} \pi}{f}})$}

$x^{8} y^{2}-x^{7} y^{3}+2 \zeta_{4} x^{7} y^{2}-x^{7} y+x^{6} y^{2}-2 x^{5} y^{3}-2 \zeta_{4} x^{6} y+2 \zeta_{4} x^{5} y^{2}+3 x^{5} y+2 x^{4} y^{2}-x^{3} y^{3}-\zeta_{4} x^{5}+2 \zeta_{4} x^{4} y+3 \zeta_{4} x^{3} y^{2}+2 x^{3} y-2 x^{2} y^{2}-2 \zeta_{4} x^{3}+\zeta_{4} x^{2} y-\zeta_{4} x y^{2}+2 x y-\zeta_{4} x+\zeta_{4} y=0$.

$f=28,(\aaa,\bbb,\ccc,\ddd,\eee)=(2,6,4,6,4)/7$.

\subsubsection*{Case $\{\aaa\ddd^2,\ccc\eee^2,\aaa^2\bbb\}$, $(x={\mathrm e}^{\frac{2\mathrm{i} \eee}{3}},y={\mathrm e}^{\frac{4\mathrm{i} \pi}{3f}})$}

$\zeta_{3}^{2} x^{13} y^{12}-\zeta_{3} x^{12} y^{10}+x^{11} y^{11}-\zeta_{3}^{2} x^{10} y^{12}-\zeta_{3} x^{9} y^{13}-\zeta_{3} x^{12} y^{7}-x^{11} y^{8}-\zeta_{3}^{2} x^{10} y^{9}+2 \zeta_{3} x^{9} y^{7}+3 x^{8} y^{8}+3 \zeta_{3}^{2} x^{7} y^{9}+x^{8} y^{5}+3 \zeta_{3}^{2} x^{7} y^{6}-3 \zeta_{3} x^{6} y^{7}-x^{5} y^{8}-3 \zeta_{3} x^{6} y^{4}-3 x^{5} y^{5}-2 \zeta_{3}^{2} x^{4} y^{6}+\zeta_{3} x^{3} y^{4}+x^{2} y^{5}+\zeta_{3}^{2} x y^{6}+\zeta_{3}^{2} x^{4}+\zeta_{3} x^{3} y-x^{2} y^{2}+\zeta_{3}^{2} x y^{3}-\zeta_{3} y=0$.

None.

\subsubsection*{Case $\{\aaa\ddd^2,\ccc\eee^2,\aaa\bbb\ccc\}$, $(x={\mathrm e}^{2\mathrm{i} \eee},y={\mathrm e}^{\frac{4\mathrm{i} \pi}{f}})$}

$(y+1) (-x^{2} y^{7}+x^{3} y^{5}+x y^{7}+y^{8}+x^{3} y^{4}+x^{2} y^{5}-x y^{6}-y^{7}-2 x^{3} y^{3}-4 x^{2} y^{4}-x y^{5}-y^{6}+x^{3} y^{2}+x^{2} y^{3}+4 x y^{4}+2 y^{5}+x^{3} y+x^{2} y^{2}-x y^{3}-y^{4}-x^{3}-x^{2} y-y^{3}+x y)=0$.

None.

\subsubsection*{Case $\{\aaa\ddd^2,\ccc\eee^2,\bbb^3\}$, $(x={\mathrm e}^{2\mathrm{i} \eee},y={\mathrm e}^{\frac{4\mathrm{i} \pi}{f}})$}

$(2 \zeta_{3}^{2}-\zeta_{3}-1) x^{3} y^{5}+(\zeta_{3}-1) x^{3} y^{4}+\zeta_{3}^{2} x^{3} y^{3}-(\zeta_{3}^{2}+\zeta_{3}-2) x^{2} y^{4}+\zeta_{3}^{2} x^{3} y^{2}-3 \zeta_{3} x^{2} y^{3}+(\zeta_{3}^{2}-1) x \,y^{4}-(2 \zeta_{3}^{2}+1) x^{2} y^{2}+(2 \zeta_{3}+1) x \,y^{3}-(\zeta_{3}-1) x^{2} y+3 \zeta_{3}^{2} x \,y^{2}-\zeta_{3} y^{3}+(\zeta_{3}^{2}+\zeta_{3}-2) x y-\zeta_{3} y^{2}-(\zeta_{3}^{2}-1) y+\zeta_{3}^{2}-2 \zeta_{3}+1=0$.

None.

\subsubsection*{Case $\{\aaa\ddd^2,\ccc\eee^2,\bbb^2\ccc\}$, $(x={\mathrm e}^{\mathrm{i} \eee},y={\mathrm e}^{\frac{4\mathrm{i} \pi}{f}})$}

$(x-1) (x^{5} y^{3}+x^{4} y^{4}-x^{3} y^{5}+x^{5} y^{2}+2 x^{4} y^{3}+x^{3} y^{4}-x^{5} y+x^{4} y^{2}+3 x^{3} y^{3}-x^{4} y-x^{3} y^{2}+x^{2} y^{3}+x y^{4}-3 x^{2} y^{2}-x y^{3}+y^{4}-x^{2} y-2 x y^{2}-y^{3}+x^{2}-x y-y^{2})=0$.

None.

\subsubsection*{Case $\{\aaa\ddd^2,\ccc\eee^2,\bbb\ccc^2\}$, $(x={\mathrm e}^{2\mathrm{i} \eee},y={\mathrm e}^{\frac{4\mathrm{i} \pi}{f}})$}

$x^{11} y-x^{11}-x^{10} y+x^{9} y^{2}-x^{8} y^{3}-x^{9} y-x^{8} y^{2}-x^{6} y^{4}+2 x^{9}+2 x^{8} y+2 x^{7} y^{2}+x^{5} y^{4}+x^{4} y^{5}-x^{7} y-3 x^{6} y^{2}+3 x^{5} y^{3}+x^{4} y^{4}-x^{7}-x^{6} y-2 x^{4} y^{3}-2 x^{3} y^{4}-2 x^{2} y^{5}+x^{5} y+x^{3} y^{3}+x^{2} y^{4}+x^{3} y^{2}-x^{2} y^{3}+x y^{4}+y^{5}-y^{4}=0$.

None.

\subsubsection*{Case $\{\aaa\ddd^2,\ccc\eee^2,\aaa\bbb^3\}$, $(x={\mathrm e}^{2\mathrm{i} \eee},y={\mathrm e}^{\frac{4\mathrm{i} \pi}{f}})$}

$x^{10} y^{19}-2 x^{9} y^{17}+x^{8} y^{16}+x^{8} y^{15}-x^{7} y^{16}+x^{8} y^{14}-2 x^{7} y^{14}+x^{6} y^{14}-x^{7} y^{12}-x^{6} y^{13}+x^{6} y^{12}+3 x^{6} y^{11}-2 x^{5} y^{11}-x^{6} y^{9}+x^{4} y^{11}-x^{6} y^{8}+x^{4} y^{10}+2 x^{5} y^{8}-3 x^{4} y^{8}-x^{4} y^{7}+x^{4} y^{6}+x^{3} y^{7}-x^{4} y^{5}+2 x^{3} y^{5}-x^{2} y^{5}+x^{3} y^{3}-x^{2} y^{4}-x^{2} y^{3}+2 x y^{2}-1=0$.

None.

\subsubsection*{Case $\{\aaa\ddd^2,\ccc\eee^2,\bbb^4\}$, $(x={\mathrm e}^{2\mathrm{i} \eee},y={\mathrm e}^{\frac{4\mathrm{i} \pi}{f}})$}

$2 \zeta_{4} x^{3} y^{5}-(\zeta_{4}-1) x^{3} y^{4}+x^{3} y^{3}-2 x^{2} y^{4}+x^{3} y^{2}+3 \zeta_{4} x^{2} y^{3}+(\zeta_{4}+1) x y^{4}-(\zeta_{4}+2) x^{2} y^{2}-(2 \zeta_{4}+1) x y^{3}+(\zeta_{4}+1) x^{2} y+3 x y^{2}+\zeta_{4} y^{3}-2 \zeta_{4} x y+\zeta_{4} y^{2}+(\zeta_{4}-1) y+2=0$.

None.

\subsubsection*{Case $\{\aaa\ddd^2,\ccc\eee^2,\bbb^3\ccc\}$, $(x={\mathrm e}^{\frac{2\mathrm{i} \eee}{3}},y={\mathrm e}^{\frac{4\mathrm{i} \pi}{f}})$}

$(x-1) (x^{8} y^{3}+x^{7} y^{4}-x^{6} y^{5}+x^{8} y^{2}+2 x^{7} y^{3}+x^{6} y^{4}+2 x^{7} y^{2}+3 x^{6} y^{3}+x^{5} y^{4}-x^{7} y+2 x^{6} y^{2}+4 x^{5} y^{3}-x^{6} y+2 x^{4} y^{3}+x^{3} y^{4}-x^{5} y-2 x^{4} y^{2}+x^{2} y^{4}-4 x^{3} y^{2}-2 x^{2} y^{3}+x y^{4}-x^{3} y-3 x^{2} y^{2}-2 x y^{3}-x^{2} y-2 x y^{2}-y^{3}+x^{2}-x y-y^{2})=0
$.

None.

\subsubsection*{Case $\{\aaa\ddd^2,\ccc\eee^2,\bbb^5\}$, $(x={\mathrm e}^{2\mathrm{i} \eee},y={\mathrm e}^{\frac{4\mathrm{i} \pi}{f}})$}

$(\zeta_{5}^{4}+\zeta_{5}-2) x^{3} y^{5}-(\zeta_{5}^{2}-\zeta_{5}) x^{3} y^{4}-\zeta_{5}^{3} x^{3} y^{3}+(\zeta_{5}^{2}-2 \zeta_{5}+1) x^{2} y^{4}-\zeta_{5}^{3} x^{3} y^{2}+3 \zeta_{5}^{2} x^{2} y^{3}+(\zeta_{5}-1) x \,y^{4}+(\zeta_{5}^{4}+2 \zeta_{5}^{3}) x^{2} y^{2}-(2 \zeta_{5}^{2}+\zeta_{5}) x y^{3}-(\zeta_{5}^{4}-1) x^{2} y-3 \zeta_{5}^{3} x \,y^{2}+\zeta_{5}^{2} y^{3}+(2 \zeta_{5}^{4}-\zeta_{5}^{3}-1) x y+\zeta_{5}^{2} y^{2}-(\zeta_{5}^{4}-\zeta_{5}^{3}) y-\zeta_{5}^{4}-\zeta_{5}+2=0$.

None.

\subsubsection*{Case $\{\aaa\ddd^2,\ccc\eee^2,\aaa\bbb^2\}$, $(x={\mathrm e}^{\mathrm{i} \ddd},y={\mathrm e}^{\frac{4\mathrm{i} \pi}{f}})$}

$x^{7} y^{4}-2 x^{6} y^{4}+x^{5} y^{5}-x^{3} y^{7}+x^{6} y^{3}+x^{5} y^{4}-x^{4} y^{5}-x^{3} y^{6}-2 x^{5} y^{3}+3 x^{4} y^{4}+2 x^{3} y^{5}-x^{2} y^{6}-x^{6} y-x^{5} y^{2}+x^{4} y^{3}-x^{3} y^{4}+x^{2} y^{5}+x y^{6}+x^{5} y-2 x^{4} y^{2}-3 x^{3} y^{3}+2 x^{2} y^{4}+x^{4} y+x^{3} y^{2}-x^{2} y^{3}-x y^{4}+x^{4}-x^{2} y^{2}+2 x y^{3}-y^{3}=0$.

$f=16,(\aaa,\bbb,\ccc,\ddd,\eee)=(6,1,2,1,3)/4$.

\subsubsection*{Case $\{\aaa\ddd\eee,\bbb\ddd^2,\aaa^2\ccc\}$, $(x={\mathrm e}^{\mathrm{i} \ddd},y={\mathrm e}^{\frac{2\mathrm{i} \pi}{f}})$}

$x^{10} y^{7}+2 \zeta_{4} x^{9} y^{6}-x^{8} y^{7}+\zeta_{4} x^{7} y^{8}+\zeta_{4} x^{9} y^{4}-\zeta_{4} x^{7} y^{4}-x^{6} y^{5}+2 \zeta_{4} x^{5} y^{6}+\zeta_{4}^{2} x^{6} y^{3}+6 \zeta_{4} x^{5} y^{4}+x^{4} y^{5}+2 \zeta_{4} x^{5} y^{2}-\zeta_{4}^{2} x^{4} y^{3}-\zeta_{4} x^{3} y^{4}+\zeta_{4} x y^{4}+\zeta_{4} x^{3}-\zeta_{4}^{2} x^{2} y+2 \zeta_{4} x y^{2}+\zeta_{4}^{2} y=0$.

None.

\subsubsection*{Case $\{\aaa\ddd\eee,\bbb\ddd^2,\aaa\bbb\ccc\}$, $(x={\mathrm e}^{2\mathrm{i} \ddd},y={\mathrm e}^{\frac{4\mathrm{i} \pi}{f}})$}

$x^{4} y^{5}+2 x^{4} y^{4}-x^{3} y^{5}-y^{6} x^{2}+2 x^{4} y^{3}+x^{3} y^{4}-x^{4} y^{2}-3 x^{3} y^{3}-2 x^{2} y^{4}-2 x y^{5}+x^{3} y^{2}+6 x^{2} y^{3}+x y^{4}-2 x^{3} y-2 x^{2} y^{2}-3 x y^{3}-y^{4}+x y^{2}+2 y^{3}-x^{2}-x y+2 y^{2}+y=0$.

None.

\subsubsection*{Case $\{\aaa\ddd\eee,\bbb\ddd^2,\aaa\ccc^2\}$, $(x={\mathrm e}^{2\mathrm{i} \ddd},y={\mathrm e}^{\frac{4\mathrm{i} \pi}{f}})$}

$x^{8} y^{10}+x^{8} y^{8}+2 x^{7} y^{9}-x^{7} y^{8}+x^{7} y^{7}+x^{6} y^{8}-2 x^{6} y^{7}-2 x^{6} y^{6}+x^{5} y^{7}-x^{6} y^{5}-x^{5} y^{6}+2 x^{5} y^{5}+x^{4} y^{6}-6 x^{4} y^{5}+x^{4} y^{4}+2 x^{3} y^{5}-x^{3} y^{4}-x^{2} y^{5}+x^{3} y^{3}-2 x^{2} y^{4}-2 x^{2} y^{3}+x^{2} y^{2}+x y^{3}-x y^{2}+2 x y+y^{2}+1=0$.

None.

\subsubsection*{Case $\{\aaa\ddd\eee,\bbb\ddd^2,\bbb^2\ccc\}$, $(x={\mathrm e}^{\mathrm{i} \eee},y={\mathrm e}^{\frac{2\mathrm{i} \pi}{f}})$}

$(y^{2}+1) (x^{3} y^{10}+\zeta_{4} x^{4} y^{7}-x^{3} y^{8}+2 \zeta_{4} x^{2} y^{9}-2 \zeta_{4} x^{4} y^{5}-3 \zeta_{4} x^{2} y^{7}+\zeta_{4}^{2} x y^{8}+\zeta_{4} x^{4} y^{3}-x^{3} y^{4}+6 \zeta_{4} x^{2} y^{5}-\zeta_{4}^{2} x y^{6}+\zeta_{4} y^{7}+x^{3} y^{2}-3 \zeta_{4} x^{2} y^{3}-2 \zeta_{4} y^{5}+2 \zeta_{4} x^{2} y-\zeta_{4}^{2} x y^{2}+\zeta_{4} y^{3}+\zeta_{4}^{2} x)=0$.

None.

\subsubsection*{Case $\{\aaa\ddd\eee,\bbb\ddd^2,\bbb\ccc^2\}$, $(x={\mathrm e}^{\mathrm{i} \eee},y={\mathrm e}^{\frac{4\mathrm{i} \pi}{f}})$}

$(y+1) (x^{2} y^{6}-x^{4} y^{3}+x^{3} y^{4}+x y^{6}-x^{3} y^{3}-3 x^{2} y^{4}-x y^{5}+6 x^{2} y^{3}-x^{3} y-3 x^{2} y^{2}-x y^{3}+x^{3}+x y^{2}-y^{3}+x^{2})=0$.

None.

\subsubsection*{Case $\{\aaa\ddd\eee,\bbb\ddd^2,\ccc^3\}$, $(x={\mathrm e}^{\mathrm{i} \eee+\frac{2\mathrm{i} \pi}{f}},y={\mathrm e}^{\mathrm{i} \eee-\frac{2\mathrm{i} \pi}{f}})$}

$(\zeta_{12}^{5}+2 \zeta_{12}^{3}+\zeta_{12}) x^{3} y^{3}-(\zeta_{12}^{4}-1) x^{3} y^{2}-(\zeta_{12}^{6}-\zeta_{12}^{4}-2 \zeta_{12}^{2}) x^{2} y^{3}+(\zeta_{12}^{6}+\zeta_{12}^{4}) x y^{4}+\zeta_{12} x^{4}+(2 \zeta_{12}^{5}+2 \zeta_{12}) x^{3} y-6 \zeta_{12}^{3} x^{2} y^{2}+(2 \zeta_{12}^{5}+2 \zeta_{12}) x y^{3}+\zeta_{12}^{5} y^{4}+(\zeta_{12}^{2}+1) x^{3}+(2 \zeta_{12}^{4}+\zeta_{12}^{2}-1) x^{2} y+(\zeta_{12}^{6}-\zeta_{12}^{2}) x y^{2}+(\zeta_{12}^{5}+2 \zeta_{12}^{3}+\zeta_{12}) x y=0$.

None.

\subsubsection*{Case $\{\aaa\ddd\eee,\bbb\ddd^2,\aaa\ccc^3\}$, $(x={\mathrm e}^{2\mathrm{i} \ddd},y={\mathrm e}^{\frac{4\mathrm{i} \pi}{f}})$}

$x^{12} y^{14}+2 x^{11} y^{13}-x^{11} y^{11}+x^{10} y^{12}+x^{10} y^{11}-x^{10} y^{10}+2 x^{9} y^{10}-x^{8} y^{10}+x^{8} y^{9}-2 x^{8} y^{8}-x^{7} y^{9}-x^{8} y^{7}+2 x^{7} y^{7}-6 x^{6} y^{7}+2 x^{5} y^{7}-x^{4} y^{7}-x^{5} y^{5}-2 x^{4} y^{6}+x^{4} y^{5}-x^{4} y^{4}+2 x^{3} y^{4}-x^{2} y^{4}+x^{2} y^{3}+x^{2} y^{2}-x y^{3}+2 x y+1=0$.

None.

\subsubsection*{Case $\{\aaa\ddd\eee,\bbb\ddd^2,\bbb\ccc^3\}$, $(x={\mathrm e}^{\mathrm{i} \eee+\frac{\mathrm{i} \pi}{f}},y={\mathrm e}^{\mathrm{i} \eee-\frac{\mathrm{i} \pi}{f}})$}

$x^{8} y^{6}-2 \zeta_{8}^{2} x^{7} y^{7}+\zeta_{8}^{4} x^{6} y^{8}-\zeta_{8} x^{9} y^{4}+\zeta_{8}^{3} x^{8} y^{5}-\zeta_{8}^{3} x^{6} y^{7}-2 \zeta_{8} x^{5} y^{8}+\zeta_{8}^{3} x^{4} y^{9}-\zeta_{8}^{3} x^{2} y^{11}-\zeta_{8} x y^{12}-\zeta_{8}^{2} x^{12}-2 \zeta_{8}^{2} x^{10} y^{2}+2 x^{9} y^{3}+6 \zeta_{8}^{2} x^{6} y^{6}+2 \zeta_{8}^{4} x^{3} y^{9}-2 \zeta_{8}^{2} x^{2} y^{10}-\zeta_{8}^{2} y^{12}-\zeta_{8}^{3} x^{11}-\zeta_{8} x^{10} y+\zeta_{8} x^{8} y^{3}-2 \zeta_{8}^{3} x^{7} y^{4}-\zeta_{8} x^{6} y^{5}+\zeta_{8} x^{4} y^{7}-\zeta_{8}^{3} x^{3} y^{8}+x^{6} y^{4}-2 \zeta_{8}^{2} x^{5} y^{5}+\zeta_{8}^{4} x^{4} y^{6}=0$.

None.

\subsubsection*{Case $\{\aaa\ddd\eee,\bbb\ddd^2,\ccc^4\}$, $(x={\mathrm e}^{\mathrm{i} \eee+\frac{2\mathrm{i} \pi}{f}},y={\mathrm e}^{\mathrm{i} \eee-\frac{2\mathrm{i} \pi}{f}})$}

$2 \zeta_{8}^{2} x^{3} y^{3}-(\zeta_{8}^{3}-\zeta_{8}) x^{3} y^{2}+2 \zeta_{8} x^{2} y^{3}+(\zeta_{8}^{3}+\zeta_{8}) x y^{4}+\zeta_{8}^{2} x^{4}+(2 \zeta_{8}^{2}-2) x^{3} y-6 \zeta_{8}^{2} x^{2} y^{2}-(2 \zeta_{8}^{4}-2 \zeta_{8}^{2}) x y^{3}+\zeta_{8}^{2} y^{4}+(\zeta_{8}^{3}+\zeta_{8}) x^{3}+2 \zeta_{8}^{3} x^{2} y+(\zeta_{8}^{3}-\zeta_{8}) x y^{2}+2 \zeta_{8}^{2} x y=0$.

None.

\subsubsection*{Case $\{\aaa\ddd\eee,\bbb\ddd^2,\ccc^5\}$, $(x={\mathrm e}^{\mathrm{i} \eee+\frac{2\mathrm{i} \pi}{f}},y={\mathrm e}^{\mathrm{i} \eee-\frac{2\mathrm{i} \pi}{f}})$}

$(\zeta_{20}^{9}-2 \zeta_{20}^{5}+\zeta_{20}) x^{3} y^{3}+(\zeta_{20}^{8}-\zeta_{20}^{4}) x^{3} y^{2}-(\zeta_{20}^{8}-\zeta_{20}^{6}+2 \zeta_{20}^{2}) x^{2} y^{3}+(\zeta_{20}^{10}-\zeta_{20}^{6}) x y^{4}-\zeta_{20}^{7} x^{4}-(2 \zeta_{20}^{7}-2 \zeta_{20}) x^{3} y+6 \zeta_{20}^{5} x^{2} y^{2}+(2 \zeta_{20}^{9}-2 \zeta_{20}^{3}) x y^{3}-\zeta_{20}^{3} y^{4}-(\zeta_{20}^{4}-1) x^{3}-(2 \zeta_{20}^{8}-\zeta_{20}^{4}+\zeta_{20}^{2}) x^{2} y-(\zeta_{20}^{6}-\zeta_{20}^{2}) x y^{2}+(\zeta_{20}^{9}-2 \zeta_{20}^{5}+\zeta_{20}) x y=0$.

None.

\subsubsection*{Case $\{\aaa\ddd\eee,\bbb\ddd\eee,\aaa^2\ccc\,(\text{or}\,\aaa\bbb\ccc,\bbb^2\ccc)\}$, $(x={\mathrm e}^{2\mathrm{i} \ddd},y={\mathrm e}^{\frac{4\mathrm{i} \pi}{f}})$}

$(y+1) (x y^{6}+x y^{5}+y^{6}-x^{2} y^{3}-3 x y^{4}-2 y^{5}+x^{2} y^{2}+4 x y^{3}+y^{4}-2 x^{2} y-3 x y^{2}-y^{3}+x^{2}+x y+x)=0$.

None.

\subsubsection*{Case $\{\aaa\ddd\eee,\bbb\ddd\eee,\ccc^3\}$, $(x={\mathrm e}^{2\mathrm{i} \ddd},y={\mathrm e}^{\frac{4\mathrm{i} \pi}{f}})$}

$(\zeta_{3}^{2}+\zeta_{3}-2) x^{2} y^{4}-\zeta_{3} x^{2} y^{3}+(\zeta_{3}^{2}-1) x y^{4}-2 \zeta_{3}^{2} x^{2} y^{2}+(\zeta_{3}^{2}-2 \zeta_{3}+3) x y^{3}-x^{2} y+\zeta_{3}^{2} x y^{2}-\zeta_{3} y^{3}+(\zeta_{3}^{2}+3 \zeta_{3}-2) x y-2 \zeta_{3}^{2} y^{2}+(\zeta_{3}^{2}-\zeta_{3}) x-y+\zeta_{3}^{2}-2 \zeta_{3}+1=0$.

None.

\subsubsection*{Case $\{\aaa\ddd\eee,\bbb\ddd\eee,\aaa\ccc^3\,(\text{or}\,\bbb\ccc^3)\}$, $(x={\mathrm e}^{2\mathrm{i} \ddd},y={\mathrm e}^{\frac{2\mathrm{i} \pi}{f}})$}

$x^{2} y^{14}-2 \zeta_{4} x^{2} y^{13}+x y^{14}-x^{2} y^{12}-\zeta_{4} x y^{13}+x^{2} y^{10}+3 \zeta_{4} x y^{11}+2 x y^{10}-\zeta_{4} x y^{9}-\zeta_{4}^{2} y^{10}+2 \zeta_{4} x^{2} y^{7}-x y^{8}-2 \zeta_{4} x y^{7}-\zeta_{4}^{2} x y^{6}+2 \zeta_{4} y^{7}-x^{2} y^{4}-\zeta_{4} x y^{5}+2 \zeta_{4}^{2} x y^{4}+3 \zeta_{4} x y^{3}+\zeta_{4}^{2} y^{4}-\zeta_{4} x y-\zeta_{4}^{2} y^{2}+\zeta_{4}^{2} x-2 \zeta_{4} y+\zeta_{4}^{2}=0$.

None.

\subsubsection*{Case $\{\aaa\ddd\eee,\bbb\ddd\eee,\ccc^4\}$, $(x={\mathrm e}^{2\mathrm{i} \ddd},y={\mathrm e}^{\frac{4\mathrm{i} \pi}{f}})$}

$2 \zeta_{4} x^{2} y^{4}-x^{2} y^{3}+(\zeta_{4}-1) x y^{4}-2 \zeta_{4} x^{2} y^{2}-(2 \zeta_{4}+2) x y^{3}+x^{2} y+2 \zeta_{4} x y^{2}+\zeta_{4}^{2} y^{3}-(2 \zeta_{4}^{2}+2 \zeta_{4}) x y-2 \zeta_{4} y^{2}-(\zeta_{4}^{2}-\zeta_{4}) x-\zeta_{4}^{2} y+2 \zeta_{4}=0$.

None.

\subsubsection*{Case $\{\aaa\ddd\eee,\bbb\ddd\eee,\ccc^5\}$, $(x={\mathrm e}^{2\mathrm{i} \ddd},y={\mathrm e}^{\frac{4\mathrm{i} \pi}{f}})$}

$(2 \zeta_{5}^{4}-\zeta_{5}^{3}-1) x^{2} y^{4}+x^{2} y^{3}+(\zeta_{5}^{4}-\zeta_{5}^{3}) x y^{4}+2 \zeta_{5} x^{2} y^{2}-(3 \zeta_{5}^{4}+\zeta_{5}-2) x y^{3}+\zeta_{5}^{2} x^{2} y-(\zeta_{5}^{2}+2 \zeta_{5}+1) x y^{2}+y^{3}-(3 \zeta_{5}^{3}-2 \zeta_{5}^{2}+\zeta_{5}) x y+2 \zeta_{5} y^{2}-(\zeta_{5}^{4}-\zeta_{5}^{3}) x+\zeta_{5}^{2} y-\zeta_{5}^{4}+2 \zeta_{5}^{3}-\zeta_{5}^{2}=0$.

None.

\subsubsection*{Case $\{\aaa\ddd\eee,\bbb\ddd\eee,\aaa\ccc^2\,(\text{or}\,\bbb\ccc^2)\}$, $(x={\mathrm e}^{2\mathrm{i} \ddd},y={\mathrm e}^{\frac{4\mathrm{i} \pi}{f}})$}

$(y+1) (x^{2} y^{8}+xy^{8}-xy^{7}-x^{2} y^{5}-2 x y^{6}+2 x^{2} y^{4}+3 x y^{5}-x^{2} y^{3}-4 x y^{4}-y^{5}+3 x y^{3}+2 y^{4}-2 x y^{2}-y^{3}-xy+x+1)=0$.

$f=20,(\aaa,\bbb,\ccc,\ddd,\eee)=(2,2,4,3,5)/5$.

$f=20,(\aaa,\bbb,\ccc,\ddd,\eee)=(4,4,8,7,9)/10$.

\subsubsection*{Case $\{\aaa\ddd\eee,\bbb\eee^2,\aaa^2\ccc\}$, $(x={\mathrm e}^{\mathrm{i} \eee+\frac{2\mathrm{i} \pi}{f}},y={\mathrm e}^{\mathrm{i} \eee-\frac{2\mathrm{i} \pi}{f}})$}

$x^{8} y^{3}+\zeta_{4} x^{7} y^{3}-x^{7} y^{2}-\zeta_{4} x^{6} y^{2}-2 \zeta_{4} x^{5} y^{3}-\zeta_{4} x^{4} y^{4}-x^{5} y^{2}+\zeta_{4}^{2} x^{4} y^{3}+\zeta_{4} x^{5} y-4 \zeta_{4} x^{4} y^{2}+\zeta_{4} x^{3} y^{3}+x^{4} y-\zeta_{4}^{2} x^{3} y^{2}-\zeta_{4} x^{4}-2 \zeta_{4} x^{3} y-\zeta_{4} x^{2} y^{2}-\zeta_{4}^{2} x y^{2}+\zeta_{4} x y+\zeta_{4}^{2} y=0$.

None.

\subsubsection*{Case $\{\aaa\ddd\eee,\bbb\eee^2,\aaa\bbb\ccc\}$, $(x={\mathrm e}^{2\mathrm{i} \eee},y={\mathrm e}^{\frac{4\mathrm{i} \pi}{f}})$}

$x^{3} y^{4}+x^{3} y^{3}+x^{3} y^{2}+2 x^{2} y^{3}+x y^{4} -2 x^{2} y^{2}-x y^{3} -y^{4}+x^{3}+x^{2} y+2 x y^{2}-x^{2}-2 x y-y^{2}-y-1=0$.

None.

\subsubsection*{Case $\{\aaa\ddd\eee,\bbb\eee^2,\aaa\ccc^2\}$, $(x={\mathrm e}^{2\mathrm{i} \eee},y={\mathrm e}^{\frac{4\mathrm{i} \pi}{f}})$}

$(x-1) (x^{7} y^{8}-x^{7} y^{7}-x^{6} y^{8}-x^{6} y^{6}+x^{5} y^{6}+x^{5} y^{5}+2 x^{4} y^{5}-2 x^{4} y^{4}+2 x^{3} y^{4}-2 x^{3} y^{3}-x^{2} y^{3}-x^{2} y^{2}+x y^{2}+x+y-1)=0$.

None.

\subsubsection*{Case $\{\aaa\ddd\eee,\bbb\eee^2,\bbb^2\ccc\}$, $(x={\mathrm e}^{\mathrm{i} \ddd},y={\mathrm e}^{\frac{2\mathrm{i} \pi}{f}})$}

$(y^{2}+1)(x^{3} y^{10}-x^{3} y^{8}+\zeta_{4} x^{2} y^{9}-\zeta_{4} x^{4} y^{5}-2 \zeta_{4} x^{2} y^{7}+\zeta_{4}^{2} x y^{8}-x^{3} y^{4}+4 \zeta_{4} x^{2} y^{5}-\zeta_{4}^{2} x y^{6}+x^{3} y^{2}-2 \zeta_{4} x^{2} y^{3}-\zeta_{4} y^{5}+\zeta_{4} x^{2} y-\zeta_{4}^{2} x y^{2}+\zeta_{4}^{2} x)=0$.

None.

\subsubsection*{Case $\{\aaa\ddd\eee,\bbb\eee^2,\bbb\ccc^2\}$, $(x={\mathrm e}^{\mathrm{i} \ddd},y={\mathrm e}^{\frac{4\mathrm{i} \pi}{f}})$}

$(y-1)(y+1) (x^{4} y^{2}+x^{3} y^{2}+x^{2} y^{3}+x y^{4}+x^{3} y-2 x^{2} y^{2}+x y^{3}+x^{3}+x^{2} y+x y^{2}+y^{2})=0$.

None.

\subsubsection*{Case $\{\aaa\ddd\eee,\bbb\eee^2,\ccc^3\}$, $(x={\mathrm e}^{\mathrm{i} \ddd},y={\mathrm e}^{\frac{2\mathrm{i} \pi}{f}})$}

$\zeta_{12}^{5} x^{4} y^{5}-2 \zeta_{12}^{3} x^{4} y^{3}+(\zeta_{12}^{4}-1) x^{3} y^{4}-(2 \zeta_{12}^{5}+\zeta_{12}^{3}+\zeta_{12}) x^{2} y^{5}-(\zeta_{12}^{2}+1) x y^{6}+\zeta_{12} x^{4} y+(\zeta_{12}^{6}-\zeta_{12}^{4}-2 \zeta_{12}^{2}) x^{3} y^{2}+4 \zeta_{12}^{3} x^{2} y^{3}-(2 \zeta_{12}^{4}+\zeta_{12}^{2}-1) x y^{4}+\zeta_{12}^{5} y^{5}-(\zeta_{12}^{6}+\zeta_{12}^{4}) x^{3}-(\zeta_{12}^{5}+\zeta_{12}^{3}+2 \zeta_{12}) x^{2} y-(\zeta_{12}^{6}-\zeta_{12}^{2}) x y^{2}-2 \zeta_{12}^{3} y^{3}+\zeta_{12} y=0$.

None.

\subsubsection*{Case $\{\aaa\ddd\eee,\bbb\eee^2,\aaa\ccc^3\}$, $(x={\mathrm e}^{2\mathrm{i} \eee},y={\mathrm e}^{\frac{4\mathrm{i} \pi}{f}})$}

$(x-1) (x^{11} y^{12}-x^{10} y^{12}+x^{10} y^{10}+x^{9} y^{9}-x^{8} y^{9}-x^{7} y^{8}+x^{7} y^{7}+x^{6} y^{7}-2 x^{6} y^{6}+2 x^{5} y^{6}-x^{5} y^{5}-x^{4} y^{5}+x^{4} y^{4}+x^{3} y^{3}-x^{2} y^{3}-x y^{2}+x-1)=0$.

None.

\subsubsection*{Case $\{\aaa\ddd\eee,\bbb\eee^2,\ccc^4\}$, $(x={\mathrm e}^{\mathrm{i} \ddd},y={\mathrm e}^{\frac{2\mathrm{i} \pi}{f}})$}

$x^{4} y^{5}+2 \zeta_{8}^{2} x^{4} y^{3}-(\zeta_{8}^{3}-\zeta_{8}) x^{3} y^{4}+(2 \zeta_{8}^{2}-2) x^{2} y^{5}+(\zeta_{8}^{3}+\zeta_{8}) x y^{6}+\zeta_{8}^{4} x^{4} y+2 \zeta_{8} x^{3} y^{2}-4 \zeta_{8}^{2} x^{2} y^{3}+2 \zeta_{8}^{3} x y^{4}+y^{5}+(\zeta_{8}^{3}+\zeta_{8}) x^{3}-(2 \zeta_{8}^{4}-2 \zeta_{8}^{2}) x^{2} y+(\zeta_{8}^{3}-\zeta_{8}) x y^{2}+2 \zeta_{8}^{2} y^{3}+\zeta_{8}^{4} y=0$.

None.

\subsubsection*{Case $\{\aaa\ddd\eee,\bbb\eee^2,\bbb\ccc^3\}$, $(x={\mathrm e}^{\mathrm{i} \ddd},y={\mathrm e}^{\frac{\mathrm{i} \pi}{f}})$}

$\zeta_{8}^{3} x y^{22}+x^{4} y^{17}+\zeta_{8}^{2} x^{2} y^{19}+\zeta_{8} x y^{20}+\zeta_{8} x^{3} y^{16}-2 x^{2} y^{17}-\zeta_{8}^{3} x^{3} y^{14}+\zeta_{8}^{2} x^{2} y^{15}-\zeta_{8} x y^{16}+y^{17}+2 \zeta_{8}^{2} x^{4} y^{11}+2 \zeta_{8}^{3} x y^{14}+\zeta_{8}^{3} x^{3} y^{10}-4 \zeta_{8}^{2} x^{2} y^{11}+\zeta_{8} x y^{12}+2 \zeta_{8} x^{3} y^{8}+2 \zeta_{8}^{2} y^{11}+\zeta_{8}^{4} x^{4} y^{5}-\zeta_{8}^{3} x^{3} y^{6}+\zeta_{8}^{2} x^{2} y^{7}-\zeta_{8} x y^{8}-2 \zeta_{8}^{4} x^{2} y^{5}+\zeta_{8}^{3} x y^{6}+\zeta_{8}^{3} x^{3} y^{2}+\zeta_{8}^{2} x^{2} y^{3}+\zeta_{8}^{4} y^{5}+\zeta_{8} x^{3}=0$.

$f=20,(\aaa,\bbb,\ccc,\ddd,\eee)=(2,4,2,5,3)/5$.

\subsubsection*{Case $\{\aaa\ddd\eee,\bbb\eee^2,\ccc^5\}$, $(x={\mathrm e}^{\mathrm{i} \ddd},y={\mathrm e}^{\frac{2\mathrm{i} \pi}{f}})$}

$\zeta_{20} x^{4} y^{5}+2 \zeta_{20}^{5} x^{4} y^{3}-(\zeta_{20}^{8}-\zeta_{20}^{4}) x^{3} y^{4}+(\zeta_{20}^{7}+\zeta_{20}^{5}-2 \zeta_{20}) x^{2} y^{5}+(\zeta_{20}^{4}-1) x y^{6}+\zeta_{20}^{9} x^{4} y+(\zeta_{20}^{8}-\zeta_{20}^{6}+2 \zeta_{20}^{2}) x^{3} y^{2}-4 \zeta_{20}^{5} x^{2} y^{3}+(2 \zeta_{20}^{8}-\zeta_{20}^{4}+\zeta_{20}^{2}) x y^{4}+\zeta_{20} y^{5}-(\zeta_{20}^{10}-\zeta_{20}^{6}) x^{3}-(2 \zeta_{20}^{9}-\zeta_{20}^{5}-\zeta_{20}^{3}) x^{2} y+(\zeta_{20}^{6}-\zeta_{20}^{2}) x y^{2}+2 \zeta_{20}^{5} y^{3}+\zeta_{20}^{9} y=0$.

$f=20,(\aaa,\bbb,\ccc,\ddd,\eee)=(2,4,2,5,3)/5$.

\subsubsection*{Case $\{\bbb\ddd^2,\ccc\eee^2,\aaa^3\}$, $(x={\mathrm e}^{2\mathrm{i} \ddd},y={\mathrm e}^{\frac{4\mathrm{i} \pi}{f}})$}

$x^{4} y^{7}-(\zeta_{3}-1) x^{4} y^{6}-(\zeta_{3}^{2}-\zeta_{3}) x^{4} y^{5}+\zeta_{3} x^{3} y^{6}-\zeta_{3}^{2} x^{4} y^{4}-(2 \zeta_{3}^{2}+\zeta_{3}) x^{3} y^{5}+(2 \zeta_{3}^{2}+1) x^{3} y^{4}-x^{3} y^{3}-(4 \zeta_{3}^{2}+2) x^{2} y^{4}+(4 \zeta_{3}+2) x^{2} y^{3}+x y^{4}-(2 \zeta_{3}+1) x y^{3}+(\zeta_{3}^{2}+2 \zeta_{3}) x y^{2}+\zeta_{3} y^{3}-\zeta_{3}^{2} x y-(\zeta_{3}^{2}-\zeta_{3}) y^{2}+(\zeta_{3}^{2}-1) y-1=0$.

None.

\subsubsection*{Case $\{\bbb\ddd^2,\ccc\eee^2,\aaa^2\bbb\,(\text{or}\,\aaa^2\ccc)\}$, $(x={\mathrm e}^{\mathrm{i} \ddd},y={\mathrm e}^{\frac{4\mathrm{i} \pi}{f}})$}

$x^{6} y^{5}+x^{5} y^{6}-x^{4} y^{7}+x^{6} y^{4}-x^{5} y^{5}-x^{4} y^{6}-x^{5} y^{4}+2 x^{4} y^{5}-x^{3} y^{6}+x^{5} y^{3}-2 x^{4} y^{4}+x^{3} y^{5}-4 x^{4} y^{3}+2 x^{3} y^{4}-2 x^{3} y^{3}+4 x^{2} y^{4}-x^{3} y^{2}+2 x^{2} y^{3}-x y^{4}+x^{3} y-2 x^{2} y^{2}+x y^{3}+x^{2} y+x y^{2}-y^{3}+x^{2}-x y-y^{2}=0$.

None.

\subsubsection*{Case $\{\bbb\ddd^2,\ccc\eee^2,\aaa\bbb\ccc\}$, $(x={\mathrm e}^{2\mathrm{i} \ddd},y={\mathrm e}^{\frac{4\mathrm{i} \pi}{f}})$}

$x^{4} y^{3}-x^{3} y^{4}-4 x^{2} y^{5}-x y^{6}+y^{7}+2 y^{2} x^{4}+3 x^{3} y^{3}+2 x^{2} y^{4}+3 x y^{5}+2 y^{6}-2 x^{4} y-3 x^{3} y^{2}-2 x^{2} y^{3}-3 x y^{4}-2 y^{5}-x^{4}+x^{3} y+4 x^{2} y^{2}+x y^{3}-y^{4}=0$.

None.

\subsubsection*{Case $\{\bbb\ddd^2,\ccc\eee^2,\aaa\bbb^2\,(\text{or}\,\aaa\ccc^2)\}$, $(x={\mathrm e}^{2\mathrm{i} \ddd},y={\mathrm e}^{\frac{4\mathrm{i} \pi}{f}})$}

$x^{8} y-x^{7} y^{2}-4 x^{6} y^{3}-x^{5} y^{4}+x^{4} y^{5}+x^{8}+x^{7} y+x^{6} y^{2}+3 x^{5} y^{3}+3 x^{4} y^{4}+2 x^{3} y^{5}+x^{2} y^{6}-x^{6} y-2 x^{5} y^{2}-3 x^{4} y^{3}-3 x^{3} y^{4}-x^{2} y^{5}-x y^{6}-y^{7}-x^{4} y^{2}+x^{3} y^{3}+4 x^{2} y^{4}+x y^{5}-y^{6}=0$.

None.

\subsubsection*{Case $\{\bbb\ddd^2,\ccc\eee^2,\aaa^4\}$, $(x={\mathrm e}^{2\mathrm{i} \ddd},y={\mathrm e}^{\frac{4\mathrm{i} \pi}{f}})$}

$\zeta_{4} x^{4} y^{7}+(\zeta_{4}+1) x^{4} y^{6}+(\zeta_{4}-1) x^{4} y^{5}-x^{3} y^{6}+\zeta_{4} x^{4} y^{4}+(2 \zeta_{4}+1) x^{3} y^{5}-(2 \zeta_{4}-1) x^{3} y^{4}+\zeta_{4}^{2} x^{3} y^{3}+(4 \zeta_{4}-2) x^{2} y^{4}-(2 \zeta_{4}^{2}-4 \zeta_{4}) x^{2} y^{3}+x y^{4}+(\zeta_{4}^{2}-2 \zeta_{4}) x y^{3}+(\zeta_{4}^{2}+2 \zeta_{4}) x y^{2}+\zeta_{4} y^{3}-\zeta_{4}^{2} x y-(\zeta_{4}^{2}-\zeta_{4}) y^{2}+(\zeta_{4}^{2}+\zeta_{4}) y+\zeta_{4}=0$.

None.

\subsubsection*{Case $\{\bbb\ddd^2,\ccc\eee^2,\aaa^3\bbb\,(\text{or}\,\aaa^3\ccc)\}$, $(x={\mathrm e}^{\frac{2\mathrm{i} \ddd}{3}},y={\mathrm e}^{\frac{4\mathrm{i} \pi}{f}})$}

$x^{10} y^{5}+x^{9} y^{6}-x^{8} y^{7}+x^{10} y^{4}-x^{9} y^{5}-x^{8} y^{6}-x^{8} y^{4}+2 x^{7} y^{5}-x^{6} y^{6}+x^{8} y^{3}-2 x^{7} y^{4}+x^{6} y^{5}-4 x^{6} y^{3}+2 x^{5} y^{4}-2 x^{5} y^{3}+4 x^{4} y^{4}-x^{4} y^{2}+2 x^{3} y^{3}-x^{2} y^{4}+x^{4} y-2 x^{3} y^{2}+x^{2} y^{3}+x^{2} y+x y^{2}-y^{3}+x^{2}-x y-y^{2}=0$.

None.

\subsubsection*{Case $\{\bbb\ddd^2,\ccc\eee^2,\aaa^5\}$, $(x={\mathrm e}^{2\mathrm{i} \ddd},y={\mathrm e}^{\frac{4\mathrm{i} \pi}{f}})$}

$\zeta_{5}^{2} x^{4} y^{7}-(\zeta_{5}^{3}-\zeta_{5}^{2}) x^{4} y^{6}-(\zeta_{5}^{4}-\zeta_{5}^{3}) x^{4} y^{5}+\zeta_{5}^{3} x^{3} y^{6}-\zeta_{5}^{4} x^{4} y^{4}-(2 \zeta_{5}^{4}+\zeta_{5}^{3}) x^{3} y^{5}+(2 \zeta_{5}^{4}+1) x^{3} y^{4}-x^{3} y^{3}-(4 \zeta_{5}^{4}+2) x^{2} y^{4}+(4 \zeta_{5}+2) x^{2} y^{3}+x y^{4}-(2 \zeta_{5}+1) x y^{3}+(\zeta_{5}^{2}+2 \zeta_{5}) x y^{2}+\zeta_{5} y^{3}-\zeta_{5}^{2} x y-(\zeta_{5}^{2}-\zeta_{5}) y^{2}-(\zeta_{5}^{3}-\zeta_{5}^{2}) y-\zeta_{5}^{3}=0$.

None.

\subsubsection*{Case $\{\bbb\ddd\eee,\ccc\eee^2,\aaa^3\}$, $(x={\mathrm e}^{\mathrm{i} \ddd+\frac{2\mathrm{i} \pi}{f}},y={\mathrm e}^{\mathrm{i} \ddd-\frac{2\mathrm{i} \pi}{f}})$}

$ (\zeta_{12} x+1) (\zeta_{12} x y^{3}+(\zeta_{12}^{2}-1) x y^{2}-y^{3}-\zeta_{12}^{3} x^{2}-(\zeta_{12}^{3}-\zeta_{12}) x y+\zeta_{12}^{2} x)=0$.

None.
\subsubsection*{Case $\{\bbb\ddd\eee,\ccc\eee^2,\aaa^2\bbb\}$, $(x={\mathrm e}^{\frac{2\mathrm{i} \ddd}{3}},y={\mathrm e}^{\frac{4\mathrm{i} \pi}{3f}})$}

$\zeta_{3}^{2} x^{3} y^{13}-\zeta_{3} x^{7} y^{6}-\zeta_{3}^{2} x^{6} y^{7}+x^{5} y^{8}-\zeta_{3} x^{4} y^{9}+\zeta_{3}^{2} x^{3} y^{10}-x^{2} y^{11}+\zeta_{3}^{2} x^{6} y^{4}+2 x^{5} y^{5}-2 \zeta_{3} x^{4} y^{6}+2 \zeta_{3}^{2} x^{3} y^{7}-2 x^{2} y^{8}-\zeta_{3} x y^{9}+x^{5} y^{2}-\zeta_{3} x^{4} y^{3}+\zeta_{3}^{2} x^{3} y^{4}-x^{2} y^{5}+\zeta_{3} x y^{6}+\zeta_{3}^{2} y^{7}-\zeta_{3} x^{4}=0$.

None.

\subsubsection*{Case $\{\bbb\ddd\eee,\ccc\eee^2,\aaa^2\ccc\}$, $(x={\mathrm e}^{\mathrm{i} \ddd},y={\mathrm e}^{\frac{4\mathrm{i} \pi}{f}})$}

$(x y+1) (x y^{4}-x^{3} y-x^{2} y^{2}+x y^{3}+x^{2} y-x y^{2}-y^{3}+x^{2})=0$.

None.

\subsubsection*{Case $\{\bbb\ddd\eee,\ccc\eee^2,\aaa\bbb^2\}$, $(x={\mathrm e}^{\mathrm{i} \eee+\frac{2\mathrm{i} \pi}{f}},y={\mathrm e}^{\mathrm{i} \eee-\frac{2\mathrm{i} \pi}{f}})$}

$(x+y) (\zeta_{4} x^{2} y^{3}+2 x^{2} y^{2}-\zeta_{4} x^{2} y-y^{2}+2 \zeta_{4} y+1)=0
$.

None.

\subsubsection*{Case $\{\bbb\ddd\eee,\ccc\eee^2,\aaa\bbb\ccc\}$, $(x={\mathrm e}^{\mathrm{i} \ddd},y={\mathrm e}^{\frac{4\mathrm{i} \pi}{f}})$}

$(y-1) (y+1) (x-y) (x+y) (x^{2}+y)=0$.

None.

\subsubsection*{Case $\{\bbb\ddd\eee,\ccc\eee^2,\aaa\ccc^2\}$, $(x={\mathrm e}^{\mathrm{i} \ddd+\frac{2\mathrm{i} \pi}{f}},y={\mathrm e}^{\mathrm{i} \ddd-\frac{2\mathrm{i} \pi}{f}})$}

$\zeta_{4} x^{2} y^{4}-x^{3} y^{2}+x^{2} y^{3}-x y^{4}+2 \zeta_{4} x^{3} y+2 \zeta_{4} x^{2} y^{2}+3 \zeta_{4} x y^{3}+\zeta_{4} y^{4}+x^{3}+3 x^{2} y+2 x y^{2}+2 y^{3}-\zeta_{4} x^{2}+\zeta_{4} x y-\zeta_{4} y^{2}+x=0$.

None.

\subsubsection*{Case $\{\bbb\ddd\eee,\ccc\eee^2,\aaa^4\}$, $(x={\mathrm e}^{\mathrm{i} \ddd+\frac{2\mathrm{i} \pi}{f}},y={\mathrm e}^{\mathrm{i} \ddd-\frac{2\mathrm{i} \pi}{f}})$}

$\zeta_{8}^{2} x^{3} y^{3}+(\zeta_{8}^{3}-\zeta_{8}) x^{3} y^{2}+\zeta_{8} x^{2} y^{3}+x^{4}+(\zeta_{8}^{2}+1) x^{3} y+(2 \zeta_{8}^{2}-2) x^{2} y^{2}-x y^{3}-\zeta_{8}^{3} x^{3}-(2 \zeta_{8}^{3}-2 \zeta_{8}) x^{2} y+(\zeta_{8}^{3}+\zeta_{8}) x y^{2}+\zeta_{8}^{3} y^{3}+\zeta_{8}^{2} x^{2}-(\zeta_{8}^{2}-1) x y+\zeta_{8} x=0$.

None.

\subsubsection*{Case $\{\bbb\ddd\eee,\ccc\eee^2,\aaa^3\bbb\}$, $(x={\mathrm e}^{\frac{2\mathrm{i} \ddd}{5}},y={\mathrm e}^{\frac{4\mathrm{i} \pi}{5f}})$}

$\zeta_{5} x^{6} y^{20}-\zeta_{5}^{4} x^{12} y^{9}+\zeta_{5}^{3} x^{10} y^{11}-x^{9} y^{12}-\zeta_{5}^{4} x^{7} y^{14}+\zeta_{5} x^{6} y^{15}+x^{4} y^{17}-\zeta_{5}^{3} x^{10} y^{6}+2 \zeta_{5}^{2} x^{8} y^{8}-2 \zeta_{5}^{4} x^{7} y^{9}-2 \zeta_{5}^{3} x^{5} y^{11}+2 x^{4} y^{12}-\zeta_{5}^{4} x^{2} y^{14}+\zeta_{5}^{2} x^{8} y^{3}+\zeta_{5} x^{6} y^{5}-\zeta_{5}^{3} x^{5} y^{6}-\zeta_{5}^{2} x^{3} y^{8}+\zeta_{5}^{4} x^{2} y^{9}-\zeta_{5}^{3} y^{11}+\zeta_{5} x^{6}=0$.

None.

\subsubsection*{Case $\{\bbb\ddd\eee,\ccc\eee^2,\aaa^3\ccc\}$, $(x={\mathrm e}^{\mathrm{i} \ddd+\frac{\mathrm{i} \pi}{f}},y={\mathrm e}^{\mathrm{i} \ddd-\frac{\mathrm{i} \pi}{f}})$}

$\zeta_{8}^{2} x^{6} y^{6}+\zeta_{8}^{3} x^{7} y^{4}-\zeta_{8} x^{6} y^{5}+\zeta_{8} x^{4} y^{7}+x^{10}+x^{8} y^{2}+\zeta_{8}^{2} x^{7} y^{3}+2 \zeta_{8}^{2} x^{5} y^{5}-2 x^{4} y^{6}-x^{2} y^{8}-\zeta_{8}^{3} x^{8} y-2 \zeta_{8}^{3} x^{6} y^{3}+2 \zeta_{8} x^{5} y^{4}+\zeta_{8} x^{3} y^{6}+\zeta_{8}^{3} x^{2} y^{7}+\zeta_{8}^{3} y^{9}+\zeta_{8}^{2} x^{6} y^{2}-\zeta_{8}^{2} x^{4} y^{4}+x^{3} y^{5}+\zeta_{8} x^{4} y^{3}=0$.

None.

\subsubsection*{Case $\{\bbb\ddd\eee,\ccc\eee^2,\aaa^5\}$, $(x={\mathrm e}^{\mathrm{i} \ddd+\frac{2\mathrm{i} \pi}{f}},y={\mathrm e}^{\mathrm{i} \ddd-\frac{2\mathrm{i} \pi}{f}})$
}

$\zeta_{20}^{6} x^{3} y^{3}+(\zeta_{20}^{9}-\zeta_{20}^{3}) x^{3} y^{2}+\zeta_{20}^{3} x^{2} y^{3}+\zeta_{20}^{2} x^{4}+(\zeta_{20}^{6}+\zeta_{20}^{2}) x^{3} y+(2 \zeta_{20}^{6}-2) x^{2} y^{2}-x y^{3}-\zeta_{20}^{9} x^{3}-(2 \zeta_{20}^{9}-2 \zeta_{20}^{3}) x^{2} y+(\zeta_{20}^{7}+\zeta_{20}^{3}) x y^{2}+\zeta_{20}^{7} y^{3}+\zeta_{20}^{6} x^{2}-(\zeta_{20}^{6}-1) x y+\zeta_{20}^{3} x=0$.

None.

\subsection*{Rational angles solution for all cases with a unique $a^2b$-vertex type  $\aaa\ddd^2$, $\bbb\ddd^2$ or $\ccc\ddd^2$: $22$ distinct cases}

\subsubsection*{Case $\{\aaa\ddd^2,\ddd\eee^3,\aaa\bbb^2\}$, $(x={\mathrm e}^{\frac{\mathrm{i} \ddd}{3}},y={\mathrm e}^{\frac{4\mathrm{i} \pi}{f}})$}

$x^{22} y-2 x^{19} y-\zeta_{3} x^{18} y+x^{16} y^{2}-\zeta_{3}^{2} x^{17}+x^{16} y+2 \zeta_{3} x^{15} y+\zeta_{3} x^{15}-x^{13} y^{2}+\zeta_{3}^{2} x^{14}+2 x^{13} y-2 \zeta_{3} x^{12} y^{2}-\zeta_{3} x^{12} y+\zeta_{3}^{2} x^{11} y^{2}-4 \zeta_{3}^{2} x^{11} y+\zeta_{3}^{2} x^{11}-x^{10} y-2 x^{10}+2 \zeta_{3} x^{9} y+\zeta_{3}^{2} x^{8} y^{2}-\zeta_{3} x^{9}+x^{7} y^{2}+2 x^{7} y+\zeta_{3} x^{6} y-\zeta_{3}^{2} x^{5} y^{2}+\zeta_{3} x^{6}-x^{4} y-2 \zeta_{3} x^{3} y+\zeta_{3} y=0$.

None.

\subsubsection*{Case $\{\aaa\ddd^2,\ddd\eee^3,\aaa\bbb\ccc\}$, $(x={\mathrm e}^{\mathrm{i} \ccc},y={\mathrm e}^{\frac{4\mathrm{i} \pi}{f}})$}

$(y+1)(x y^{12} -2 x y^{11}+2 x^{2} y^{9}+2 x y^{10}+x^{3} y^{7}-4 x^{2} y^{8}-2 x y^{9}+y^{10}-x^{3} y^{6}+4 x^{2} y^{7}+2 x y^{8}-y^{9}+x^{4} y^{4}-x^{3} y^{5}-2 x^{2} y^{6}-x y^{7}+y^{8}-x^{4} y^{3}+2 x^{3} y^{4}+4 x^{2} y^{5}-x y^{6}+x^{4} y^{2}-2 x^{3} y^{3}-4 x^{2} y^{4}+x y^{5}+2 x^{3} y^{2}+2 x^{2} y^{3}-2 x^{3} y+x^{3})=0$.

None.

\subsubsection*{Case $\{\aaa\ddd^2,\ddd\eee^3,\bbb^3\}$, $(x={\mathrm e}^{\frac{\mathrm{i} \ccc}{2}+\frac{2\mathrm{i} \pi}{f}},y={\mathrm e}^{\frac{\mathrm{i} \ccc}{2}-\frac{2\mathrm{i} \pi}{f}})$}

$(\zeta_{12}^{3}+\zeta_{12}) x^{2} y^{7}-(2 \zeta_{12}^{5}+\zeta_{12}^{3}-\zeta_{12}) x y^{8}+\zeta_{12}^{4} x^{2} y^{6}-2 \zeta_{12}^{3} x^{2} y^{5}+(\zeta_{12}^{5}-\zeta_{12}^{3}-2 \zeta_{12}) x y^{6}-x^{2} y^{4}+2 \zeta_{12}^{2} x y^{5}+(\zeta_{12}^{5}-\zeta_{12}) x^{2} y^{3}+4 \zeta_{12}^{3} x y^{4}-(\zeta_{12}^{5}-\zeta_{12}) y^{5}+2 \zeta_{12}^{4} x y^{3}-\zeta_{12}^{6} y^{4}-(2 \zeta_{12}^{5}+\zeta_{12}^{3}-\zeta_{12}) x y^{2}-2 \zeta_{12}^{3} y^{3}+\zeta_{12}^{2} y^{2}+(\zeta_{12}^{5}-\zeta_{12}^{3}-2 \zeta_{12}) x+(\zeta_{12}^{5}+\zeta_{12}^{3}) y=0$.

None.

\subsubsection*{Case $\{\aaa\ddd^2,\ddd\eee^3,\bbb^2\ccc\}$, $(x={\mathrm e}^{\frac{\mathrm{i} \ccc}{4}+\frac{2\mathrm{i} \pi}{f}},y={\mathrm e}^{\frac{\mathrm{i} \ccc}{4}-\frac{2\mathrm{i} \pi}{f}})$}

$(y+1) (x y+1) (x^{4} y^{6}-x^{4} y^{5}+x^{4} y^{4}-x^{3} y^{5}+x^{2} y^{6}-2 x^{2} y^{5}-x^{3} y^{3}+3 x^{2} y^{4}-x^{2} y^{3}+3 x^{2} y^{2}-x y^{3}-2 x^{2} y+x^{2}-x y+y^{2}-y+1)=0$.

None.

\subsubsection*{Case $\{\aaa\ddd^2,\ddd\eee^3,\bbb\ccc^2\}$, $(x={\mathrm e}^{\frac{\mathrm{i} \ccc}{2}+\frac{2\mathrm{i} \pi}{f}},y={\mathrm e}^{\frac{\mathrm{i} \ccc}{2}-\frac{2\mathrm{i} \pi}{f}})$}

$(x y-1) (\zeta_{4} x^{10} y^{3}-\zeta_{4} x^{9} y^{4}+x^{8} y^{4}+2 \zeta_{4} x^{9} y^{2}-\zeta_{4} x^{8} y^{3}+2 x^{7} y^{3}+\zeta_{4}^{2} x^{6} y^{4}+\zeta_{4} x^{8} y-\zeta_{4} x^{5} y^{4}+3 x^{6} y^{2}+2 \zeta_{4}^{2} x^{5} y^{3}-\zeta_{4} x^{6} y+2 \zeta_{4} x^{5} y^{2}-\zeta_{4} x^{4} y^{3}+2 x^{5} y+3 \zeta_{4}^{2} x^{4} y^{2}-\zeta_{4} x^{5}+\zeta_{4} x^{2} y^{3}+x^{4}+2 \zeta_{4}^{2} x^{3} y-\zeta_{4} x^{2} y+2 \zeta_{4} x y^{2}+\zeta_{4}^{2} x^{2}-\zeta_{4} x+\zeta_{4} y)=0$.

None.

\subsubsection*{Case $\{\aaa\ddd^2,\eee^4,\aaa\bbb^2\}$, $(x={\mathrm e}^{\mathrm{i} \ddd},y={\mathrm e}^{\frac{4\mathrm{i} \pi}{f}})$}

$x^{7} y-(\zeta_{4}+2) x^{6} y+x^{5} y^{2}+x^{6}+(2 \zeta_{4}+1) x^{5} y-(\zeta_{4}+1) x^{4} y^{2}+(\zeta_{4}-1) x^{5}-(\zeta_{4}-4) x^{4} y-(\zeta_{4}-1) x^{3} y^{2}+(\zeta_{4}-1) x^{4}+(4 \zeta_{4}-1) x^{3} y-(\zeta_{4}-1) x^{2} y^{2}-(\zeta_{4}+1) x^{3}+(\zeta_{4}+2) x^{2} y+\zeta_{4} x y^{2}+\zeta_{4} x^{2}-(2 \zeta_{4}+1) x y+\zeta_{4} y=0$.

None.

\subsubsection*{Case $\{\aaa\ddd^2,\eee^4,\aaa\bbb\ccc\}$, $(x={\mathrm e}^{\mathrm{i} \ccc},y={\mathrm e}^{\frac{4\mathrm{i} \pi}{f}})$}

$(y+1) (x y^{8}+x^{3} y^{5}+2 x^{2} y^{6}-2 x y^{7}+x^{4} y^{3}-x^{3} y^{4}-4 x^{2} y^{5}+2 x y^{6}+y^{7}-x^{3} y^{3}+8 x^{2} y^{4}-x y^{5}+x^{4} y+2 x^{3} y^{2}-4 x^{2} y^{3}-x y^{4}+y^{5}-2 x^{3} y+2 x^{2} y^{2}+x y^{3}+x^{3})=0$.

None.

\subsubsection*{Case $\{\aaa\ddd^2,\eee^4,\bbb^3\}$, $(x={\mathrm e}^{\mathrm{i} \ccc},y={\mathrm e}^{\frac{4\mathrm{i} \pi}{f}})$}

$(\zeta_{3}^{2}-1) x^{5} y-(\zeta_{3}^{2}+\zeta_{3}) x^{4} y^{2}+(\zeta_{3}+1) x^{3} y^{3}-(\zeta_{3}^{2}-1) x^{2} y^{4}-(\zeta_{3}^{2}-2 \zeta_{3}+1) x^{5}-(\zeta_{3}^{2}+\zeta_{3}-2) x^{4} y+4 \zeta_{3}^{2} x^{3} y^{2}-4 \zeta_{3} x^{2} y^{3}+(\zeta_{3}^{2}+\zeta_{3}-2) x y^{4}-(2 \zeta_{3}^{2}-\zeta_{3}-1) y^{5}+(\zeta_{3}-1) x^{3} y-(\zeta_{3}^{2}+1) x^{2} y^{2}+(\zeta_{3}^{2}+\zeta_{3}) x y^{3}-(\zeta_{3}-1) y^{4}=0$.

None.

\subsubsection*{Case $\{\aaa\ddd^2,\eee^4,\bbb^2\ccc\}$, $(x={\mathrm e}^{\frac{\mathrm{i} \ccc}{2}},y={\mathrm e}^{\frac{4\mathrm{i} \pi}{f}})$}

$(x+1) (x+y) (x^{6} y+x^{4} y^{3}-x^{5} y-2 x^{4} y^{2}-x^{3} y^{3}+x^{5}-x^{4} y+6 x^{3} y^{2}-x^{2} y^{3}+x y^{4}-x^{3} y-2 x^{2} y^{2}-x y^{3}+x^{2} y+y^{3})=0$.

None.

\subsubsection*{Case $\{\aaa\ddd^2,\eee^4,\bbb\ccc^2\}$, $(x={\mathrm e}^{\mathrm{i} \ccc},y={\mathrm e}^{\frac{4\mathrm{i} \pi}{f}})$}

$(x-1) (x^{7} y^{5}-x^{7} y^{4}+2 x^{6} y^{5}-x^{6} y^{4}+x^{5} y^{5}+x^{6} y^{3}+2 x^{5} y^{3}-x^{4} y^{4}-x^{5} y^{2}+4 x^{4} y^{3}-x^{3} y^{4}-2 x^{4} y^{2}+2 x^{3} y^{3}+x^{4} y-4 x^{3} y^{2}+x^{2} y^{3}+x^{3} y-2 x^{2} y^{2}-x y^{2}-x^{2}+x y-2 x+y-1)=0$.

$f=24,(\aaa,\bbb,\ccc,\ddd,\eee)=(2,6,3,5,3)/6$.

\subsubsection*{Case $\{\bbb\ddd^2,\ddd\eee^3,\aaa^2\bbb\}$, $(x={\mathrm e}^{\mathrm{i} \ccc},y={\mathrm e}^{\frac{4\mathrm{i} \pi}{f}})$}

$(x+y) (x^{18} y-x^{17} y^{2}+x^{15} y^{4}-2 x^{14} y^{5}+x^{13} y^{6}-x^{12} y^{7}+x^{11} y^{8}-x^{9} y^{10}+x^{8} y^{11}-x^{7} y^{12}+x^{6} y^{13}+x^{18}-x^{17} y+x^{16} y^{2}-x^{14} y^{4}+3 x^{13} y^{5}-5 x^{12} y^{6}+2 x^{11} y^{7}-2 x^{10} y^{8}+2 x^{8} y^{10}-2 x^{7} y^{11}+5 x^{6} y^{12}-3 x^{5} y^{13}+x^{4} y^{14}-x^{2} y^{16}+x y^{17}-y^{18}-x^{12} y^{5}+x^{11} y^{6}-x^{10} y^{7}+x^{9} y^{8}-x^{7} y^{10}+x^{6} y^{11}-x^{5} y^{12}+2 x^{4} y^{13}-x^{3} y^{14}+x y^{16}-y^{17})=0$.

None.

\subsubsection*{Case $\{\bbb\ddd^2,\ddd\eee^3,\aaa\bbb^2\}$, $(x={\mathrm e}^{\frac{2\mathrm{i} \ddd}{3}},y={\mathrm e}^{\frac{4\mathrm{i} \pi}{f}})$}

$x^{17}-\zeta_{3} x^{15}+x^{14} y+2 \zeta_{3}^{2} x^{13} y-\zeta_{3}^{2} x^{13}-3 \zeta_{3} x^{12} y-x^{11} y^{2}+\zeta_{3} x^{12}-3 x^{11} y+\zeta_{3}^{2} x^{10} y^{2}+\zeta_{3} x^{9} y^{2}-\zeta_{3}^{2} x^{10}+3 \zeta_{3} x^{9} y+3 x^{8} y-\zeta_{3}^{2} x^{7} y^{2}+x^{8}+\zeta_{3}^{2} x^{7}-3 \zeta_{3} x^{6} y+x^{5} y^{2}-\zeta_{3} x^{6}-3 x^{5} y-\zeta_{3}^{2} x^{4} y^{2}+2 \zeta_{3}^{2} x^{4} y+\zeta_{3} x^{3} y-x^{2} y^{2}+\zeta_{3} y^{2}=0$.

None.

\subsubsection*{Case $\{\bbb\ddd^2,\ddd\eee^3,\aaa\bbb\ccc\}$, $(x={\mathrm e}^{\mathrm{i} \ccc},y={\mathrm e}^{\frac{4\mathrm{i} \pi}{f}})$}

$(y+1) (x^{2} y^{9}+x^{2} y^{8}+x^{3} y^{6}-2 x^{2} y^{7}-x y^{8}+x^{2} y^{6}+3 x y^{7}+2 x^{2} y^{5}-y^{7}+x^{3} y^{3}-2 x^{2} y^{4}-2 x y^{5}+y^{6}-x^{3} y^{2}+2 x y^{4}+3 x^{2} y^{2}+x y^{3}-x^{2} y-2 x y^{2}+y^{3}+x y+x)=0$.

None.

\subsubsection*{Case $\{\bbb\ddd^2,\eee^4,\aaa^2\bbb\}$, $(x={\mathrm e}^{\mathrm{i} \ddd},y={\mathrm e}^{\frac{4\mathrm{i} \pi}{f}})$}

$x^{6} y^{2}+x^{6} y+(\zeta_{4}-1) x^{5} y^{2}-(\zeta_{4}+1) x^{5} y-(\zeta_{4}-1) x^{4} y^{2}+(5 \zeta_{4}-2) x^{4} y-(\zeta_{4}-1) x^{3} y^{2}-x^{4}+(2 \zeta_{4}+2) x^{3} y-\zeta_{4} x^{2} y^{2}+(\zeta_{4}-1) x^{3}-(2 \zeta_{4}-5) x^{2} y+(\zeta_{4}-1) x^{2}-(\zeta_{4}+1) x y-(\zeta_{4}-1) x+\zeta_{4} y+\zeta_{4}=0$.

None.

\subsubsection*{Case $\{\bbb\ddd^2,\eee^4,\aaa\bbb^2\}$, $(x={\mathrm e}^{\mathrm{i} \ddd},y={\mathrm e}^{\frac{4\mathrm{i} \pi}{f}})$}

$x^{11}-\zeta_{4} x^{10}+x^{9} y+x^{9}-5 \zeta_{4} x^{8} y-x^{7} y^{2}+\zeta_{4} x^{8}-3 x^{7} y-\zeta_{4} x^{6} y^{2}+x^{7}+3 \zeta_{4} x^{6} y-x^{5} y^{2}-\zeta_{4} x^{6}+3 x^{5} y+\zeta_{4} x^{4} y^{2}-x^{5}-3 \zeta_{4} x^{4} y+x^{3} y^{2}-\zeta_{4} x^{4}-5 x^{3} y+\zeta_{4} x^{2} y^{2}+\zeta_{4} x^{2} y-x y^{2}+\zeta_{4} y^{2}=0$.

None.

\subsubsection*{Case $\{\bbb\ddd^2,\eee^4,\aaa\bbb\ccc\}$, $(x={\mathrm e}^{\mathrm{i} \ccc},y={\mathrm e}^{\frac{4\mathrm{i} \pi}{f}})$}

$x^{2} y^{7}+x^{3} y^{5}+2 x^{2} y^{6}+x^{3} y^{4}-2 x^{2} y^{5}-x y^{6}+x^{3} y^{3}+3 x^{2} y^{4}+5 x y^{5}-x^{3} y^{2}-2 x^{2} y^{3}-2 x y^{4}-y^{5}+5 x^{2} y^{2}+3 x y^{3}+y^{4}-x^{2} y-2 x y^{2}+y^{3}+2 x y+y^{2}+x=0$.

None.

\subsubsection*{Case $\{\ccc\ddd^2,\ddd\eee^3,\aaa^2\bbb\}$, $(x={\mathrm e}^{\frac{2\mathrm{i} \ddd}{3}},y={\mathrm e}^{\frac{4\mathrm{i} \pi}{f}})$}

$(x-1) (x^{2}+x+1) (\zeta_{3} x^{14} y^{5}+\zeta_{3}^{2} x^{12} y^{4}-\zeta_{3} x^{11} y^{4}-x^{10} y^{3}+\zeta_{3}^{2} x^{9} y^{4}-2 \zeta_{3}^{2} x^{9} y^{3}-\zeta_{3}^{2} x^{9} y^{2}+\zeta_{3} x^{8} y^{3}+x^{7} y^{3}+x^{7} y^{2}+\zeta_{3}^{2} x^{6} y^{2}-\zeta_{3} x^{5} y^{3}-2 \zeta_{3} x^{5} y^{2}+\zeta_{3} x^{5} y-x^{4} y^{2}-\zeta_{3}^{2} x^{3} y+\zeta_{3} x^{2} y+\zeta_{3}^{2})=0$.

None.

\subsubsection*{Case $\{\ccc\ddd^2,\ddd\eee^3,\aaa\bbb^2\}$, $(x={\mathrm e}^{\frac{2\mathrm{i} \ddd}{3}},y={\mathrm e}^{\frac{2\mathrm{i} \pi}{f}})$}

$\zeta_{12}^{4} x^{13} y^{6}+\zeta_{12}^{4} x^{10} y^{8}+\zeta_{12}^{2} x^{11} y^{6}-2 \zeta_{12}^{4} x^{10} y^{6}+2 x^{9} y^{6}-\zeta_{12}^{4} x^{7} y^{8}+\zeta_{12}^{4} x^{10} y^{4}-2 \zeta_{12}^{2} x^{8} y^{6}-2 x^{9} y^{4}+\zeta_{12}^{4} x^{7} y^{6}-2 \zeta_{12} x^{7} y^{5}-2 x^{6} y^{6}+\zeta_{12}^{4} x^{7} y^{4}+\zeta_{12}^{2} x^{5} y^{6}+\zeta_{12}^{2} x^{8} y^{2}+x^{6} y^{4}-2 \zeta_{12}^{4} x^{7} y^{2}-2 \zeta_{12}^{3} x^{6} y^{3}+x^{6} y^{2}-2 \zeta_{12}^{4} x^{4} y^{4}-2 \zeta_{12}^{2} x^{5} y^{2}+x^{3} y^{4}-x^{6}+2 \zeta_{12}^{4} x^{4} y^{2}-2 x^{3} y^{2}+\zeta_{12}^{2} x^{2} y^{2}+x^{3}+y^{2}=0$.

None.

\subsubsection*{Case $\{\ccc\ddd^2,\ddd\eee^3,\aaa\bbb\ccc\}$, $(x={\mathrm e}^{\mathrm{i} \bbb},y={\mathrm e}^{\frac{4\mathrm{i} \pi}{f}})$}

$(y+1) (y^{2}-y+1) (x^{2} y^{8}+x y^{9}+x^{2} y^{7}+x^{3} y^{5}-x y^{7}-2 x^{2} y^{5}+x y^{6}+2 x^{2} y^{4}+2 x y^{5}+x^{2} y^{3}-2 x y^{4}-x^{2} y^{2}+y^{4}+x y^{2}+x^{2}+x y)=0$.

None.

\subsubsection*{Case $\{\ccc\ddd^2,\eee^4,\aaa^2\bbb\}$, $(x={\mathrm e}^{\mathrm{i} \ddd},y={\mathrm e}^{\frac{4\mathrm{i} \pi}{f}})$}

$(x-1) (x+1) (x^{8} y^{5}+\zeta_{4} x^{7} y^{4}-x^{6} y^{4}+x^{6} y^{3}+\zeta_{4} x^{5} y^{4}-3 \zeta_{4} x^{5} y^{3}-\zeta_{4} x^{5} y^{2}-x^{4} y^{3}+\zeta_{4}^{2} x^{4} y^{2}+\zeta_{4} x^{3} y^{3}+3 \zeta_{4} x^{3} y^{2}-\zeta_{4} x^{3} y-\zeta_{4}^{2} x^{2} y^{2}+\zeta_{4}^{2} x^{2} y-\zeta_{4} x y-\zeta_{4}^{2})=0$.

None.

\subsubsection*{Case $\{\ccc\ddd^2,\eee^4,\aaa\bbb^2\}$, $(x={\mathrm e}^{\mathrm{i} \ddd},y={\mathrm e}^{\frac{4\mathrm{i} \pi}{f}})$}

$(x-1) (x+1) (x^{5} y^{3}-\zeta_{4} x^{4} y^{3}+x^{3} y^{4}-\zeta_{4} x^{4} y^{2}-3 x^{3} y^{3}+2 x^{3} y^{2}+\zeta_{4} x^{2} y^{3}-x^{3} y-2 \zeta_{4} x^{2} y^{2}+3 \zeta_{4} x^{2} y+x y^{2}-\zeta_{4} x^{2}+x y-\zeta_{4} y)=0$.

None.

\subsubsection*{Case $\{\ccc\ddd^2,\eee^4,\aaa\bbb\ccc\}$, $(x={\mathrm e}^{\mathrm{i} \bbb},y={\mathrm e}^{\frac{4\mathrm{i} \pi}{f}})$}

$(y^{2}+1) (x^{2} y^{5}+x^{3} y^{3}+x^{2} y^{4}+x y^{5}-2 x^{2} y^{3}-x y^{4}+3 x^{2} y^{2}+3 x y^{3}-x^{2} y-2 x y^{2}+x^{2}+x y+y^{2}+x)=0$.

None.

\subsection*{Rational angles solution for all case with a unique $a^2b$-vertex type  $\bbb\ddd\eee$: $92$ distinct cases}

\subsubsection*{Case $\{\bbb\ddd\eee,\aaa^2\bbb,\aaa\bbb\ccc\,(\text{or}\,\bbb\ccc^2)\}$, $(x={\mathrm e}^{2\mathrm{i} \ddd},y={\mathrm e}^{\frac{2\mathrm{i} \pi}{f}})$}

$\zeta_{4} x y^{9}-x y^{8}-2 y^{8}-\zeta_{4} x^{2} y^{5}-x y^{6}-2 \zeta_{4} y^{7}-\zeta_{4} x y^{5}-y^{6}-\zeta_{4} x^{2} y^{3}-x y^{4}-2 x^{2} y^{2}-\zeta_{4} x y^{3}-y^{4}-2 \zeta_{4} x^{2} y-\zeta_{4} x y+x=0$.

None.

\subsubsection*{Case $\{\bbb\ddd\eee,\aaa^2\bbb,\bbb^2\ccc\}$, $(x={\mathrm e}^{2\mathrm{i} \ddd},y={\mathrm e}^{\frac{4\mathrm{i} \pi}{5f}})$}

$x y^{15}+\zeta_{5} x y^{14}-\zeta_{5}^{2} x y^{13}-2 \zeta_{5}^{3} x y^{12}+\zeta_{5}^{4} x y^{11}+\zeta_{5} y^{14}-\zeta_{5}^{4} y^{11}-\zeta_{5} x^{2} y^{9}+\zeta_{5}^{2} x^{2} y^{8}-2 \zeta_{5}^{4} x^{2} y^{6}-x y^{10}-\zeta_{5} x y^{9}+2 \zeta_{5}^{2} x y^{8}-2 \zeta_{5}^{3} x y^{7}+\zeta_{5}^{4} x y^{6}-2 y^{10}+2 \zeta_{5} y^{9}-\zeta_{5}^{3} y^{7}+\zeta_{5}^{4} y^{6}+2 x^{2} y^{5}+\zeta_{5} x^{2} y^{4}-\zeta_{5}^{4} x^{2} y+x y^{5}-\zeta_{5} x y^{4}+2 \zeta_{5}^{2} x y^{3}+\zeta_{5}^{3} x y^{2}-\zeta_{5}^{4} x y-x=0$.

None.

\subsubsection*{Case $\{\bbb\ddd\eee,\aaa^2\bbb,\ccc^3\}$, $(x={\mathrm e}^{2\mathrm{i} \ddd},y={\mathrm e}^{\frac{4\mathrm{i} \pi}{f}})$}

$\zeta_{3}^{2} x y^{7}-(\zeta_{3}^{2}-1) x y^{6}-\zeta_{3} y^{7}+\zeta_{3} x^{2} y^{4}-(\zeta_{3}+1) x y^{5}-(2 \zeta_{3}^{2}-\zeta_{3}-1) y^{6}-(2 \zeta_{3}^{2}-2 \zeta_{3}+1) x y^{4}+2 y^{5}-2 x^{2} y^{2}-(2 \zeta_{3}^{2}-2 \zeta_{3}-1) x y^{3}-(\zeta_{3}^{2}-2 \zeta_{3}+1) x^{2} y+(\zeta_{3}^{2}+1) x y^{2}-\zeta_{3}^{2} y^{3}+\zeta_{3}^{2} x^{2}+(\zeta_{3}-1) x y-\zeta_{3} x=0$.

None.

\subsubsection*{Case $\{\bbb\ddd\eee,\aaa^2\bbb,\aaa\ccc^3\}$, $(x={\mathrm e}^{2\mathrm{i} \ddd},y={\mathrm e}^{\frac{2\mathrm{i} \pi}{f}})$}

$\zeta_{4} x y^{20}-y^{21}-\zeta_{4} x y^{18}+y^{19}-x y^{17}-2 \zeta_{4} y^{18}-y^{17}+x y^{15}+\zeta_{4} x y^{14}-2 y^{15}-\zeta_{4} x^{2} y^{12}+x y^{13}-2 \zeta_{4} x y^{12}-2 x y^{11}-2 \zeta_{4} x y^{10}-2 x y^{9}+\zeta_{4} x y^{8}-y^{9}-2 \zeta_{4} x^{2} y^{6}+x y^{7}+\zeta_{4} x y^{6}-\zeta_{4} x^{2} y^{4}-2 x^{2} y^{3}-\zeta_{4} x y^{4}+\zeta_{4} x^{2} y^{2}-x y^{3}-\zeta_{4} x^{2}+x y=0$.

None.

\subsubsection*{Case $\{\bbb\ddd\eee,\aaa^2\bbb,\bbb\ccc^3\}$, $(x={\mathrm e}^{2\mathrm{i} \ddd},y={\mathrm e}^{\frac{4\mathrm{i} \pi}{5f}})$}

$x y^{25}+\zeta_{5} x y^{22}-\zeta_{5}^{4} y^{23}+\zeta_{5} y^{22}-x y^{20}-\zeta_{5}^{2} x y^{19}-2 y^{20}-\zeta_{5} x y^{17}+\zeta_{5}^{4} y^{18}-2 \zeta_{5}^{3} x y^{16}+2 \zeta_{5} y^{17}+\zeta_{5}^{2} x^{2} y^{14}+2 \zeta_{5}^{2} x y^{14}+\zeta_{5}^{4} x y^{13}-\zeta_{5} x y^{12}-2 \zeta_{5}^{3} x y^{11}-\zeta_{5}^{3} y^{11}-2 \zeta_{5}^{4} x^{2} y^{8}+2 \zeta_{5}^{2} x y^{9}-\zeta_{5} x^{2} y^{7}+\zeta_{5}^{4} x y^{8}+2 x^{2} y^{5}+\zeta_{5}^{3} x y^{6}+x y^{5}-\zeta_{5}^{4} x^{2} y^{3}+\zeta_{5} x^{2} y^{2}-\zeta_{5}^{4} x y^{3}-x=0$.

None.

\subsubsection*{Case $\{\bbb\ddd\eee,\aaa^2\bbb,\ccc^4\}$, $(x={\mathrm e}^{2\mathrm{i} \ddd},y={\mathrm e}^{\frac{4\mathrm{i} \pi}{f}})$}

$\zeta_{4} x y^{7}-(\zeta_{4}+1) x y^{6}-y^{7}-\zeta_{4} x^{2} y^{4}+(\zeta_{4}+1) x y^{5}-2 \zeta_{4} y^{6}-(2 \zeta_{4}+1) x y^{4}-2 y^{5}-2 \zeta_{4} x^{2} y^{2}-(\zeta_{4}+2) x y^{3}-2 x^{2} y+(\zeta_{4}+1) x y^{2}-y^{3}-\zeta_{4} x^{2}-(\zeta_{4}+1) x y+x=0$.

None.

\subsubsection*{Case $\{\bbb\ddd\eee,\aaa^2\bbb,\ccc^5\}$, $(x={\mathrm e}^{2\mathrm{i} \ddd},y={\mathrm e}^{\frac{4\mathrm{i} \pi}{f}})$}

$x y^{7}+(\zeta_{5}-1) x y^{6}-\zeta_{5}^{4} y^{7}+\zeta_{5}^{2} x^{2} y^{4}-(\zeta_{5}^{2}+\zeta_{5}) x y^{5}+(\zeta_{5}^{4}+\zeta_{5}-2) y^{6}-(2 \zeta_{5}^{3}-2 \zeta_{5}^{2}+\zeta_{5}) x y^{4}+2 \zeta_{5} y^{5}-2 \zeta_{5}^{4} x^{2} y^{2}+(\zeta_{5}^{4}-2 \zeta_{5}^{3}+2 \zeta_{5}^{2}) x y^{3}-(\zeta_{5}^{4}+\zeta_{5}-2) x^{2} y+(\zeta_{5}^{4}+\zeta_{5}^{3}) x y^{2}-\zeta_{5}^{3} y^{3}+\zeta_{5} x^{2}-(\zeta_{5}^{4}-1) x y-x=0$.

None.

\subsubsection*{Case $\{\bbb\ddd\eee,\aaa^2\bbb,\aaa\ccc^2\}$, $(x={\mathrm e}^{2\mathrm{i} \ddd},y={\mathrm e}^{\frac{4\mathrm{i} \pi}{f}})$}

$(y+1) (x y^{11}+y^{12}-3 x y^{10}-3 y^{11}+4 x y^{9}+3 y^{10}-4 x y^{8}-4 y^{9}+x^{2} y^{6}+4 x y^{7}+3 y^{8}-2 x^{2} y^{5}-2 x y^{6}-2 y^{7}+3 x^{2} y^{4}+4 x y^{5}+y^{6}-4 x^{2} y^{3}-4 x y^{4}+3 x^{2} y^{2}+4 x y^{3}-3 x^{2} y-3 x y^{2}+x^{2}+x y)=0$.

$f=28,(\aaa,\bbb,\ccc,\ddd,\eee)=(2,10,6,3,1)/7$.

\subsubsection*{Case $\{\bbb\ddd\eee,\aaa^3,\aaa\bbb^2\}$, $(x={\mathrm e}^{2\mathrm{i} \ddd},y={\mathrm e}^{\frac{4\mathrm{i} \pi}{f}})$}

$x^{2} y^{2}+(2 \zeta_{3}^{2}+\zeta_{3}-1) x^{2} y+(\zeta_{3}-3) x y^{2}+\zeta_{3}^{2} x^{2}-\zeta_{3}^{2} x y+\zeta_{3}^{2} y^{2}-(3 \zeta_{3}-1) x+(2 \zeta_{3}^{2}-\zeta_{3}+1) y+\zeta_{3}=0$.

None.

\subsubsection*{Case $\{\bbb\ddd\eee,\aaa^2\ccc,\aaa\bbb^2\}$, $(x={\mathrm e}^{2\mathrm{i} \ddd},y={\mathrm e}^{\frac{2\mathrm{i} \pi}{f}})$}

$\zeta_{4} x^{2} y^{9}-2 \zeta_{4} x y^{9}-\zeta_{4} x^{2} y^{7}+2 x^{2} y^{6}-\zeta_{4} x y^{7}+y^{8}-\zeta_{4} x^{2} y^{5}+\zeta_{4} x y^{5}+x y^{4}-y^{4}+\zeta_{4} x^{2} y-x y^{2}+2 \zeta_{4} y^{3}-y^{2}-2 x+1=0$.

None.

\subsubsection*{Case $\{\bbb\ddd\eee,\aaa^2\ccc,\aaa\bbb\ccc\,(\text{or}\,\bbb^2\ccc)\}$, $(x={\mathrm e}^{2\mathrm{i} \ddd},y={\mathrm e}^{\frac{4\mathrm{i} \pi}{f}})$}

$(y+1) (x y^{6}+x y^{5}+y^{6}-x^{2} y^{3}-3 x y^{4}-2 y^{5}+x^{2} y^{2}+4 x y^{3}+y^{4}-2 x^{2} y-3 x y^{2}-y^{3}+x^{2}+x y+x)=0$.

None.

\subsubsection*{Case $\{\bbb\ddd\eee,\aaa^2\ccc,\bbb^3\}$, $(x={\mathrm e}^{2\mathrm{i} \ddd},y={\mathrm e}^{\frac{4\mathrm{i} \pi}{f}})$}

$x^{2} y^{4}+(2 \zeta_{3}-1) x y^{5}+(2 \zeta_{3}^{2}-\zeta_{3}-1) x^{2} y^{3}-(\zeta_{3}^{2}-\zeta_{3}+1) x y^{4}-y^{5}-2 x^{2} y^{2}-(\zeta_{3}^{2}+2 \zeta_{3}-1) x y^{3}+(2 \zeta_{3}^{2}+\zeta_{3}-1) x y^{2}+2 y^{3}+x^{2}-(\zeta_{3}^{2}-\zeta_{3}-1) x y+(\zeta_{3}^{2}-2 \zeta_{3}+1) y^{2}-(2 \zeta_{3}^{2}-1) x-y=0$.

None.

\subsubsection*{Case $\{\bbb\ddd\eee,\aaa^2\ccc,\aaa^2\bbb^2\}$, $(x={\mathrm e}^{2\mathrm{i} \ddd},y={\mathrm e}^{\frac{4\mathrm{i} \pi}{f}})$}

$x^{2} y^{6}-x^{2} y^{5}+x y^{5}-x^{2} y^{3}-2x y^{4}+x^{2} y^{2}-y^{4}+2 x y^{2}+y^{3}-x y+y-1=0$.

None.

\subsubsection*{Case $\{\bbb\ddd\eee,\aaa^2\ccc,\aaa\bbb^3\}$, $(x={\mathrm e}^{2\mathrm{i} \ddd},y={\mathrm e}^{\frac{4\mathrm{i} \pi}{3f}})$}

$\zeta_{3}^{2} x^{2} y^{13}+2 \zeta_{3} x y^{14}+\zeta_{3}^{2} x y^{13}+x y^{12}+\zeta_{3}^{2} y^{13}-\zeta_{3}^{2} x^{2} y^{10}+\zeta_{3} x y^{11}-2 x^{2} y^{9}+\zeta_{3}^{2} x y^{10}-\zeta_{3} x^{2} y^{8}+x y^{9}-2 \zeta_{3}^{2} x^{2} y^{7}-2 \zeta_{3} x y^{8}-2 \zeta_{3}^{2} x y^{7}-2 x y^{6}-2 \zeta_{3}^{2} y^{7}+\zeta_{3} x y^{5}-y^{6}+\zeta_{3}^{2} x y^{4}-2 \zeta_{3} y^{5}+x y^{3}-\zeta_{3}^{2} y^{4}+\zeta_{3}^{2} x^{2} y+\zeta_{3} x y^{2}+\zeta_{3}^{2} x y+2 x+\zeta_{3}^{2} y=0$.

None.

\subsubsection*{Case $\{\bbb\ddd\eee,\aaa^2\ccc,\aaa\bbb^2\ccc\}$, $(x={\mathrm e}^{2\mathrm{i} \ddd},y={\mathrm e}^{\frac{2\mathrm{i} \pi}{f}})$}

$x y^{13}-2 \zeta_{4} x y^{12}+y^{13}+x y^{11}-\zeta_{4} x y^{10}+2 \zeta_{4} x^{2} y^{8}-2 x y^{9}+2 x^{2} y^{7}+2 \zeta_{4} x y^{8}-2 y^{9}-\zeta_{4} x^{2} y^{6}-x y^{7}-\zeta_{4} x y^{6}-y^{7}-2 \zeta_{4} x^{2} y^{4}+2 x y^{5}+2 \zeta_{4} y^{6}-2 \zeta_{4} x y^{4}+2 y^{5}-x y^{3}+\zeta_{4} x y^{2}+\zeta_{4} x^{2}-2 x y+\zeta_{4} x=0$.

None.

\subsubsection*{Case $\{\bbb\ddd\eee,\aaa^2\ccc,\bbb\ccc\ddd^2\}$, $(x={\mathrm e}^{\mathrm{i} \ddd},y={\mathrm e}^{\frac{4\mathrm{i} \pi}{f}})$}

$x^{3} y^{11}+x^{3} y^{10}+x^{2} y^{11}-x^{3} y^{9}-4 x^{2} y^{9}+x^{3} y^{7}-2 x^{2} y^{8}+3 x^{2} y^{7}-x^{2} y^{6}-3 x^{2} y^{5}+3 x y^{6}+x y^{5}-3 x y^{4}+2 x y^{3}-y^{4}+4  x y^{2}+y^{2}-x-y-1=0$.

None.

\subsubsection*{Case $\{\bbb\ddd\eee,\aaa^2\ccc,\ccc^2\ddd^2\}$, $(x={\mathrm e}^{\mathrm{i} \eee},y={\mathrm e}^{\frac{4\mathrm{i} \pi}{f}})$}

$(y+1) (y-1) (x^{3} y^{5}+x^{3} y^{4}-x y^{6}-2 x^{2} y^{4}-x y^{5}-x^{3} y^{2}-x y^{4}-x^{2} y^{2}-y^{4}-x^{2} y-2 x y^{2}-x^{2}+y^{2}+y)=0$.

None.

\subsubsection*{Case $\{\bbb\ddd\eee,\aaa^2\ccc,\aaa^2\bbb^3\}$, $(x={\mathrm e}^{2\mathrm{i} \ddd},y={\mathrm e}^{\frac{4\mathrm{i} \pi}{3f}})$}

$x^{2} y^{15}+2 x y^{14}-x^{2} y^{12}-x y^{13}-x y^{12}+2 x^{2} y^{10}+x y^{11}-y^{12}-2 x^{2} y^{9}-x y^{10}-x^{2} y^{8}-2 x y^{9}-2 x y^{8}+2 x y^{7}+2 x y^{6}+y^{7}+x y^{5}+2 y^{6}+x^{2} y^{3}-x y^{4}-2 y^{5}+x y^{3}+x y^{2}+y^{3}-2 x y-1=0$.

None.

\subsubsection*{Case $\{\bbb\ddd\eee,\aaa^2\ccc,\aaa\bbb^4\}$, $(x={\mathrm e}^{2\mathrm{i} \ddd},y={\mathrm e}^{\frac{\mathrm{i} \pi}{f}})$}

$2 \zeta_{8}^{3} x y^{19}-\zeta_{8} x^{2} y^{17}-\zeta_{8}^{2} x y^{18}-\zeta_{8}^{2} y^{18}+x y^{16}+\zeta_{8}^{3} x y^{15}+\zeta_{8} x^{2} y^{13}-\zeta_{8}^{2} x y^{14}-2 x^{2} y^{12}-\zeta_{8}^{3} x^{2} y^{11}+x y^{12}-2 \zeta_{8}^{3} x y^{11}+2 \zeta_{8} x^{2} y^{9}+\zeta_{8}^{2} x y^{10}+\zeta_{8} x y^{9}+2 \zeta_{8}^{2} y^{10}-2 x y^{8}+\zeta_{8}^{3} x y^{7}-y^{8}-2 \zeta_{8}^{3} y^{7}-\zeta_{8} x y^{5}+\zeta_{8}^{2} y^{6}+x y^{4}+\zeta_{8}^{3} x y^{3}-\zeta_{8} x^{2} y-\zeta_{8} x y-\zeta_{8}^{2} y^{2}+2 x=0$.

None.

\subsubsection*{Case $\{\bbb\ddd\eee,\aaa^2\ccc,\aaa\bbb^2\ccc^2\}$, $(x={\mathrm e}^{2\mathrm{i} \ddd},y={\mathrm e}^{\frac{2\mathrm{i} \pi}{f}})$}

$x y^{19}+y^{19}+x y^{17}-2 \zeta_{4} x y^{16}-x y^{15}-\zeta_{4} x y^{14}-2 y^{15}+\zeta_{4} x^{2} y^{12}+2 \zeta_{4} x y^{12}-y^{13}-x y^{11}+2 x^{2} y^{9}-\zeta_{4} x y^{10}+y^{11}+\zeta_{4} x^{2} y^{8}-x y^{9}+2 \zeta_{4} y^{10}-\zeta_{4} x y^{8}-\zeta_{4} x^{2} y^{6}+2 x y^{7}+y^{7}-2 \zeta_{4} x^{2} y^{4}-x y^{5}-\zeta_{4} x y^{4}-2 x y^{3}+\zeta_{4} x y^{2}+\zeta_{4} x^{2}+\zeta_{4} x=0$.

None.

\subsubsection*{Case $\{\bbb\ddd\eee,\aaa^2\ccc,\aaa\bbb^2\ddd^2\}$, $(x={\mathrm e}^{\mathrm{i} \ddd+\frac{2\mathrm{i} \pi}{f}},y={\mathrm e}^{\mathrm{i} \ddd-\frac{2\mathrm{i} \pi}{f}})$}

$\zeta_{4} x^{4} y^{3}+x^{5} y+3 x^{4} y^{2}-x^{3} y^{3}-\zeta_{4} x^{5}-2 \zeta_{4} x^{4} y-\zeta_{4} x^{2} y^{3}-x^{3} y-2 x y^{3}-y^{4}-\zeta_{4} x^{2} y+3 \zeta_{4} x y^{2}+\zeta_{4} y^{3}+x y=0$.

None.

\subsubsection*{Case $\{\bbb\ddd\eee,\aaa^2\ccc,\aaa\bbb^2\eee^2\}$, $(x={\mathrm e}^{\mathrm{i} \eee+\frac{2\mathrm{i} \pi}{f}},y={\mathrm e}^{\mathrm{i} \eee-\frac{2\mathrm{i} \pi}{f}})$}

$(x+y) (x^{6}-2 x^{5} y-x^{4} y^{2}+x^{3} y^{3}+2 \zeta_{4} x^{5}-\zeta_{4} x^{4} y+\zeta_{4} x^{3} y^{2}+\zeta_{4} x^{2} y^{3}+x^{4}+x^{3} y-x^{2} y^{2}+2 x y^{3}+\zeta_{4} x^{3}-\zeta_{4} x^{2} y-2 \zeta_{4} x y^{2}+\zeta_{4} y^{3})=0$.

None.

\subsubsection*{Case $\{\bbb\ddd\eee,\aaa^2\ccc,\aaa\bbb\ccc^3\}$, $(x={\mathrm e}^{2\mathrm{i} \ddd},y={\mathrm e}^{\frac{4\mathrm{i} \pi}{f}})$}

$(y+1) (x y^{18}-x y^{17}+y^{18}-2 y^{17}+x y^{15}+y^{16}-2 x y^{14}-y^{15}+5 x y^{13}+2 y^{14}-x^{2} y^{11}-7 x y^{12}-3 y^{13}+2 x^{2} y^{10}+7 x y^{11}+4 y^{12}-3 x^{2} y^{9}-6 x y^{10}-5 y^{11}+4 x^{2} y^{8}+6 x y^{9}+4 y^{10}-5 x^{2} y^{7}-6 x y^{8}-3 y^{9}+4 x^{2} y^{6}+7 x y^{7}+2 y^{8}-3 x^{2} y^{5}-7 x y^{6}-y^{7}+2 x^{2} y^{4}+5 x y^{5}-x^{2} y^{3}-2 x y^{4}+x^{2} y^{2}+x y^{3}-2 x^{2} y+x^{2}-x y+x)=0$.

None.

\subsubsection*{Case $\{\bbb\ddd\eee,\aaa^2\ccc,\aaa\bbb\ccc\ddd^2\}$, $(x={\mathrm e}^{2\mathrm{i} \ddd},y={\mathrm e}^{\frac{4\mathrm{i} \pi}{f}})$}

$(y-1) (y+1) (x^{3} y^{5}+x^{2} y^{5}+x^{2} y^{4}-x^{2} y^{3}+x y^{2}-x y-x-1)=0$.

None.

\subsubsection*{Case $\{\bbb\ddd\eee,\aaa^2\ccc,\aaa\bbb\ccc\eee^2\}$, $(x={\mathrm e}^{2\mathrm{i} \ddd},y={\mathrm e}^{\frac{4\mathrm{i} \pi}{f}})$}

$(y-1) (y+1) (x+y) (x y^{2}-1) =0$.

None.

\subsubsection*{Case $\{\bbb\ddd\eee,\aaa^2\ccc,\aaa\ccc^2\ddd^2\}$, $(x={\mathrm e}^{\mathrm{i} \eee+\frac{2\mathrm{i} \pi}{f}},y={\mathrm e}^{\mathrm{i} \eee-\frac{2\mathrm{i} \pi}{f}})$}

$(x+y) (\zeta_{4} x^{7} y-\zeta_{4} x^{6} y^{2}-\zeta_{4} x^{5} y^{3}-\zeta_{4} x^{4} y^{4}+\zeta_{4} x^{3} y^{5}+4 x^{5} y^{2}-3 x^{4} y^{3}+x^{3} y^{4}+x y^{6}+\zeta_{4} x^{6}+\zeta_{4} x^{4} y^{2}-3 \zeta_{4} x^{3} y^{3}+4 \zeta_{4} x^{2} y^{4}+x^{4} y-x^{3} y^{2}-x^{2} y^{3}-x y^{4}+y^{5})=0$.

None.

\subsubsection*{Case $\{\bbb\ddd\eee,\aaa^2\ccc,\aaa\ccc^2\ddd\eee\}$, $(x={\mathrm e}^{2\mathrm{i} \ddd},y={\mathrm e}^{\frac{4\mathrm{i} \pi}{f}})$}

$(y+1) (x y^{2}-1) (x y^{8}-3 x y^{7}+3 x y^{6}-3 x y^{5}+2 x y^{4}+y^{5}-x y^{3}-2 y^{4}+3 y^{3}-3 y^{2}+3 y-1) =0$.

None.

\subsubsection*{Case $\{\bbb\ddd\eee,\aaa^2\ccc,\aaa\ccc^2\eee^2\}$, $(x={\mathrm e}^{\mathrm{i} \ddd+\frac{2\mathrm{i} \pi}{f}},y={\mathrm e}^{\mathrm{i} \ddd-\frac{2\mathrm{i} \pi}{f}})$}

$x^{5} y^{5}+2 \zeta_{4} x^{6} y^{3}+\zeta_{4} x^{5} y^{4}+\zeta_{4} x^{4} y^{5}-\zeta_{4} x^{3} y^{6}+x^{8}-x^{7} y-4 x^{6} y^{2}-x^{5} y^{3}+3 x^{4} y^{4}-x^{3} y^{5}-x^{2} y^{6}-\zeta_{4} x^{6} y-\zeta_{4} x^{5} y^{2}+3 \zeta_{4} x^{4} y^{3}-\zeta_{4} x^{3} y^{4}-4 \zeta_{4} x^{2} y^{5}-\zeta_{4} x y^{6}+\zeta_{4} y^{7}-x^{5} y+x^{4} y^{2}+x^{3} y^{3}+2 x^{2} y^{4}+\zeta_{4} x^{3} y^{2}=0$.

None.

\subsubsection*{Case $\{\bbb\ddd\eee,\aaa^2\ccc,\bbb^2\ccc\ddd^2\}$, $(x={\mathrm e}^{\mathrm{i} \ddd},y={\mathrm e}^{\frac{4\mathrm{i} \pi}{f}})$}

$x^{4} y^{8}+x^{4} y^{7}-x^{5} y^{5}-x^{4} y^{6}+2 x^{3} y^{7}+x^{2} y^{8}+x^{3} y^{6}-2 x^{4} y^{4}-2 x^{3} y^{5}-2 x^{2} y^{6}+2 x^{3} y^{4}+2 x^{2} y^{4}-2 x^{3} y^{2}-2 x^{2} y^{3}-2 x y^{4}+x^{2} y^{2}+x^{3}+2 x^{2} y-x y^{2}-y^{3}+x y+x=0$.

None.

\subsubsection*{Case $\{\bbb\ddd\eee,\aaa^2\ccc,\ccc^3\ddd^2\}$, $(x={\mathrm e}^{\mathrm{i} \eee},y={\mathrm e}^{\frac{4\mathrm{i} \pi}{f}})$}

$(y-1) (y+1) (x y^{12}-x^{3} y^{9}+2 x^{2} y^{9}+x y^{10}+x^{3} y^{7}+x^{3} y^{6}+2 x^{2} y^{7}+y^{9}-2 x^{2} y^{6}-x y^{7}-y^{8}+x^{3} y^{4}+x^{2} y^{5}+2 x y^{6}-x^{3} y^{3}-2 x y^{5}-y^{6}-y^{5}-x^{2} y^{2}-2 x y^{3}+y^{3}-x^{2})=0$.

None.

\subsubsection*{Case $\{\bbb\ddd\eee,\aaa^2\ccc,\aaa\ccc^2\,(\text{or}\,\ccc^3)\}$, $(x={\mathrm e}^{\mathrm{i} \ddd},y={\mathrm e}^{\mathrm{i} \eee})$}

$x^{5} y^{2}-(\zeta_{3}^{2}-\zeta_{3}+1) x^{4} y^{3}-2 x^{4} y-(\zeta_{3}-2) x^{3} y^{2}+(\zeta_{3}^{2}+2 \zeta_{3}-2) x^{2} y^{3}+(\zeta_{3}^{2}-2 \zeta_{3}+2) x^{3}+(2 \zeta_{3}-1) x^{2} y-2 \zeta_{3} x y^{2}-(\zeta_{3}^{2}+\zeta_{3}-1) x+\zeta_{3} y=0$.

$f=12,(\aaa,\bbb,\ccc,\ddd,\eee)=(4,3,4,6,3)/6$.

\subsubsection*{Case $\{\bbb\ddd\eee,\aaa^2\ccc,\bbb\ccc^2\}$, $(x={\mathrm e}^{2\mathrm{i} \ddd},y={\mathrm e}^{\frac{4\mathrm{i} \pi}{f}})$}

$x y^{17}+x y^{16}+y^{17}-x y^{15}-2 y^{15}-2 x y^{13}-y^{14}+x^{2} y^{11}-x y^{12}+y^{13}+2 x y^{11}-x y^{10}-x y^{9}+2 y^{10}-2 x^{2} y^{7}+x y^{8}+x y^{7}-2 x y^{6}-x^{2} y^{4}+x y^{5}-y^{6}+x^{2} y^{3}+2 x y^{4}+2 x^{2} y^{2}+x y^{2}-x^{2}-x y-x=0$.

$f=20,(\aaa,\bbb,\ccc,\ddd,\eee)=(8,12,4,7,1)/10$.

$f=20,(\aaa,\bbb,\ccc,\ddd,\eee)=(8,12,4,1,7)/10$.

\subsubsection*{Case $\{\bbb\ddd\eee,\aaa^2\ccc,\aaa\ccc^3\}$, $(x={\mathrm e}^{\mathrm{i} \ddd},y={\mathrm e}^{\mathrm{i} \eee})$}

$x^{5} y^{2}-(\zeta_{5}^{3}-\zeta_{5}+1) x^{4} y^{3}-2 x^{4} y-(\zeta_{5}^{4}-\zeta_{5}^{3}-\zeta_{5}^{2}+2 \zeta_{5}-2) x^{3} y^{2}+(\zeta_{5}^{4}+\zeta_{5}^{2}+\zeta_{5}-2) x^{2} y^{3}+(\zeta_{5}^{3}-2 \zeta_{5}^{2}+\zeta_{5}+1) x^{3}+(\zeta_{5}^{4}-\zeta_{5}^{3}+2 \zeta_{5}^{2}-2 \zeta_{5}+1) x^{2} y-2 \zeta_{5}^{2} x y^{2}-(\zeta_{5}^{4}+\zeta_{5}^{2}-\zeta_{5}) x+\zeta_{5}^{2} y=0$.

$f=20,(\aaa,\bbb,\ccc,\ddd,\eee)=(8,12,4,7,1)/10$.

$f=20,(\aaa,\bbb,\ccc,\ddd,\eee)=(8,12,4,1,7)/10$.

\subsubsection*{Case $\{\bbb\ddd\eee,\aaa^2\ccc,\ccc^3\ddd\eee\}$, $(x={\mathrm e}^{2\mathrm{i} \ddd},y={\mathrm e}^{\frac{4\mathrm{i} \pi}{f}})$}

$x^{2} y^{21}-x^{2} y^{20}-2 x^{2} y^{19}-x y^{19}+x^{2} y^{17}+x y^{18}+x y^{17}+2 x^{2} y^{14}-x y^{15}-x y^{14}+2 x y^{13}-x y^{12}+y^{13}-2 x y^{11}+2 x y^{10}-x^{2} y^{8}+x y^{9}-2 x y^{8}+x y^{7}+x y^{6}-2 y^{7}-x y^{4}-x y^{3}-y^{4}+x y^{2}+2 y^{2}+y-1=0$.

$f=20,(\aaa,\bbb,\ccc,\ddd,\eee)=(8,12,4,7,1)/10$.

$f=20,(\aaa,\bbb,\ccc,\ddd,\eee)=(8,12,4,1,7)/10$.

\subsubsection*{Case $\{\bbb\ddd\eee,\aaa^2\ccc,\aaa\bbb^2\}$, $(x={\mathrm e}^{2\mathrm{i} \ddd},y={\mathrm e}^{\frac{2\mathrm{i} \pi}{f}})$}

$\zeta_{4} x^{2} y^{9}-2 \zeta_{4} x y^{9}-\zeta_{4} x^{2} y^{7}+2 x^{2} y^{6}-\zeta_{4} x y^{7}+y^{8}-\zeta_{4} x^{2} y^{5}+\zeta_{4} x y^{5}+x y^{4}-y^{4}+\zeta_{4} x^{2} y-x y^{2}+2 \zeta_{4} y^{3}-y^{2}-2 x+1=0$.

$f=28,(\aaa,\bbb,\ccc,\ddd,\eee)=(6,4,2,3,7)/7$.

\subsubsection*{Case $\{\bbb\ddd\eee,\aaa^2\ccc,\aaa\bbb\ccc^2\}$, $(x={\mathrm e}^{2\mathrm{i} \ddd},y={\mathrm e}^{\frac{4\mathrm{i} \pi}{f}})$}

$(y+1) (x y^{12}-x y^{11}+y^{12}-2 y^{11}+3 x y^{9}+y^{10}-x^{2} y^{7}-5 x y^{8}-y^{9}+2 x^{2} y^{6}+5 x y^{7}+2 y^{8}-3 x^{2} y^{5}-4 x y^{6}-3 y^{7}+2 x^{2} y^{4}+5 x y^{5}+2 y^{6}-x^{2} y^{3}-5 x y^{4}-y^{5}+x^{2} y^{2}+3 x y^{3}-2 x^{2} y+x^{2}-x y+x)=0$.

$f=28,(\aaa,\bbb,\ccc,\ddd,\eee)=(6,4,2,3,7)/7$.

\subsubsection*{Case $\{\bbb\ddd\eee,\aaa^2\ccc,\bbb^2\ddd^2\}$, $(x={\mathrm e}^{\mathrm{i} \ddd},y={\mathrm e}^{\frac{4\mathrm{i} \pi}{f}})$}

$(x+1) (x^{2} y^{5}-x^{3} y^{3}+x^{2} y^{4}-x^{2} y^{3}-x y^{2}+x y-y^{2}+x)=0$.

$f=28,(\aaa,\bbb,\ccc,\ddd,\eee)=(6,4,2,3,7)/7$.

\subsubsection*{Case $\{\bbb\ddd\eee,\aaa^2\ccc,\ccc^2\ddd\eee\}$, $(x={\mathrm e}^{2\mathrm{i} \ddd},y={\mathrm e}^{\frac{4\mathrm{i} \pi}{f}})$}

$x^{2} y^{15}-x^{2} y^{14}-2 x^{2} y^{13}-x y^{13}+x^{2} y^{11}+x y^{12}+2 x^{2} y^{10}-x y^{10}+2 x y^{9}+x y^{8}+y^{9}-x^{2} y^{6}-x y^{7}-2 x y^{6}+x y^{5}-2 y^{5}-x y^{3}-y^{4}+x y^{2}+2 y^{2}+y-1=0$.

$f=28,(\aaa,\bbb,\ccc,\ddd,\eee)=(6,4,2,3,7)/7$.

\subsubsection*{Case $\{\bbb\ddd\eee,\aaa^2\ccc,\aaa\ccc^4\}$, $(x={\mathrm e}^{\mathrm{i} \ddd},y={\mathrm e}^{\mathrm{i} \eee})$}

$x^{5} y^{2}-(\zeta_{7}^{4}-\zeta_{7}+1) x^{4} y^{3}-2 x^{4} y-(\zeta_{7}^{6}-\zeta_{7}^{4}-\zeta_{7}^{3}+2 \zeta_{7}-2) x^{3} y^{2}+(\zeta_{7}^{6}+\zeta_{7}^{3}+\zeta_{7}-2) x^{2} y^{3}+(\zeta_{7}^{4}-2 \zeta_{7}^{3}+\zeta_{7}^{2}+1) x^{3}+(\zeta_{7}^{6}-\zeta_{7}^{4}+2 \zeta_{7}^{3}-2 \zeta_{7}^{2}+1) x^{2} y-2 \zeta_{7}^{3} x y^{2}-(\zeta_{7}^{6}+\zeta_{7}^{3}-\zeta_{7}^{2}) x+\zeta_{7}^{3} y=0$.

$f=28,(\aaa,\bbb,\ccc,\ddd,\eee)=(6,4,2,3,7)/7$.

\subsubsection*{Case $\{\bbb\ddd\eee,\aaa^2\ccc,\bbb\ccc^2\ddd^2\}$, $(x={\mathrm e}^{2\mathrm{i} \ddd},y={\mathrm e}^{\frac{4\mathrm{i} \pi}{f}})$}

$x^{3} y^{17}+x^{3} y^{16}+x^{2} y^{17}-x^{3} y^{15}-2 x^{2} y^{15}-x^{2} y^{14}-x^{2} y^{13}+x^{3} y^{11}-x^{2} y^{12}+2 x^{2} y^{11}-x^{2} y^{10}-x^{2} y^{9}+2 x y^{10}-2 x^{2} y^{7}+x y^{8}+x y^{7}-2 x y^{6}+x y^{5}-y^{6}+x y^{4}+x y^{3}+2 x y^{2}+y^{2}-x-y-1=0$.

$f=28,(\aaa,\bbb,\ccc,\ddd,\eee)=(6,4,2,3,7)/7$.

\subsubsection*{Case $\{\bbb\ddd\eee,\aaa^2\ccc,\aaa\bbb^3\ccc\}$, $(x={\mathrm e}^{2\mathrm{i} \ddd},y={\mathrm e}^{\frac{4\mathrm{i} \pi}{3f}})$}

$\zeta_{12}^{2} x y^{18}-2 \zeta_{12}^{4} x y^{17}+\zeta_{12}^{2} y^{18}+\zeta_{12}^{2} x y^{15}-\zeta_{12}^{4} x y^{14}+\zeta_{12}^{2} x^{2} y^{12}-x y^{13}+\zeta_{12}^{4} x^{2} y^{11}-\zeta_{12}^{2} x y^{12}+2 x^{2} y^{10}+2 \zeta_{12}^{4} x y^{11}-2 \zeta_{12}^{2} y^{12}-\zeta_{12}^{2} x^{2} y^{9}-x y^{10}-\zeta_{12}^{4} x y^{8}-\zeta_{12}^{2} y^{9}-2 \zeta_{12}^{2} x^{2} y^{6}+2 x y^{7}+2 \zeta_{12}^{4} y^{8}-\zeta_{12}^{2} x y^{6}+y^{7}-\zeta_{12}^{4} x y^{5}+\zeta_{12}^{2} y^{6}-x y^{4}+\zeta_{12}^{2} x y^{3}+\zeta_{12}^{2} x^{2}-2 x y+\zeta_{12}^{2} x=0$.

$f=40,(\aaa,\bbb,\ccc,\ddd,\eee)=(18,6,4,31,3)/20$.

$f=40,(\aaa,\bbb,\ccc,\ddd,\eee)=(18,6,4,11,23)/20$.

\subsubsection*{Case $\{\bbb\ddd\eee,\aaa^2\ccc,\bbb^3\ddd^2\}$, $(x={\mathrm e}^{\frac{2\mathrm{i} \ddd}{3}},y={\mathrm e}^{\frac{4\mathrm{i} \pi}{f}})$}

$\zeta_{3} x^{7} y^{3}+x^{5} y^{5}-2 \zeta_{3}^{2} x^{6} y^{3}-2 \zeta_{3} x^{4} y^{5}-\zeta_{3} x^{4} y^{4}+2 x^{5} y^{2}+2 \zeta_{3} x^{4} y^{3}+\zeta_{3}^{2} x^{3} y^{4}+x^{2} y^{5}-\zeta_{3} x^{4} y^{2}+\zeta_{3}^{2} x^{3} y^{3}-x^{5}-\zeta_{3} x^{4} y-2 \zeta_{3}^{2} x^{3} y^{2}-2 x^{2} y^{3}+\zeta_{3}^{2} x^{3} y+2 \zeta_{3}^{2} x^{3}+2 \zeta_{3} x y^{2}-x^{2}-\zeta_{3}^{2} y^{2}=0$.

$f=40,(\aaa,\bbb,\ccc,\ddd,\eee)=(18,6,4,11,23)/20$.

\subsubsection*{Case $\{\bbb\ddd\eee,\aaa\bbb\ccc,\aaa^3\}$, $(x={\mathrm e}^{2\mathrm{i} \ddd},y={\mathrm e}^{\frac{4\mathrm{i} \pi}{f}})$}

$\zeta_{3} x y^{5}-x y^{5}-\zeta_{3} x^{2} y^{3}-2 x y^{4}-y^{5}-\zeta_{3} x^{2} y^{2}-\zeta_{3}^{2} x y^{3}-\zeta_{3} x y^{3}+\zeta_{3}^{2} y^{4}-2 \zeta_{3} y^{4}+2 \zeta_{3}^{2} x^{2} y-\zeta_{3} x^{2} y+\zeta_{3}^{2} x y^{2}+\zeta_{3} x y^{2}+\zeta_{3}^{2} y^{3}+x^{2}+2 x y+\zeta_{3}^{2} y^{2}-\zeta_{3}^{2} x+x=0$.

None.

\subsubsection*{Case $\{\bbb\ddd\eee,\aaa\bbb\ccc,\aaa^4\}$, $(x={\mathrm e}^{2\mathrm{i} \ddd},y={\mathrm e}^{\frac{4\mathrm{i} \pi}{f}})$}

$\zeta_{4} x y^{5}-x y^{5}-\zeta_{4} x^{2} y^{3}-2 x y^{4}-y^{5}-\zeta_{4} x^{2} y^{2}-\zeta_{4} x y^{3}+x y^{3}-2 \zeta_{4} y^{4}-y^{4}-\zeta_{4} x^{2} y-2 x^{2} y+\zeta_{4} x y^{2}-x y^{2}-y^{3}-\zeta_{4} x^{2}-2 \zeta_{4} x y-y^{2}-\zeta_{4} x+x=0$.

None.

\subsubsection*{Case $\{\bbb\ddd\eee,\aaa\bbb\ccc,\aaa^3\bbb\}$, $(x={\mathrm e}^{2\mathrm{i} \ddd},y={\mathrm e}^{\frac{4\mathrm{i} \pi}{3f}})$}

$\zeta_{3} x y^{14}+\zeta_{3}^{2} x y^{13}+\zeta_{3}^{2} y^{13}-y^{12}+2 \zeta_{3}^{2} x y^{10}-2 \zeta_{3} y^{11}-\zeta_{3} x^{2} y^{8}+x y^{9}-\zeta_{3} x y^{8}-y^{9}-\zeta_{3} x^{2} y^{5}-x y^{6}+\zeta_{3} x y^{5}-y^{6}-2 x^{2} y^{3}+2 \zeta_{3}^{2} x y^{4}-\zeta_{3} x^{2} y^{2}+\zeta_{3}^{2} x^{2} y+\zeta_{3}^{2} x y+x=0$.

None.

\subsubsection*{Case $\{\bbb\ddd\eee,\aaa\bbb\ccc,\aaa^3\ccc\,(\text{or}\,\aaa^2\ddd\eee)\}$, $(x={\mathrm e}^{2\mathrm{i} \ddd},y={\mathrm e}^{\frac{2\mathrm{i} \pi}{f}})$}

$x y^{13}-\zeta_{4} x y^{12}+y^{13}+2 x y^{11}+\zeta_{4} x^{2} y^{8}+2 \zeta_{4} y^{10}+\zeta_{4} x y^{8}+y^{9}+\zeta_{4} x^{2} y^{6}-x y^{7}-\zeta_{4} x y^{6}+y^{7}+\zeta_{4} x^{2} y^{4}+x y^{5}+2 x^{2} y^{3}+y^{5}+2 \zeta_{4} x y^{2}+\zeta_{4} x^{2}-x y+\zeta_{4} x=0$.

None.

\subsubsection*{Case $\{\bbb\ddd\eee,\aaa\bbb\ccc,\aaa^2\ddd^2\}$, $(x={\mathrm e}^{\mathrm{i} \ddd},y={\mathrm e}^{\frac{4\mathrm{i} \pi}{f}})$}

$(x^{2} y^{2}+x^{2} y+x y^{2}+x^{2}+y^{2}+x+y+1) (x^{2} y^{3}-x^{3} y+x^{2} y^{2}-x^{2} y-x y^{2}+x y-y^{2}+x)=0$.

None.

\subsubsection*{Case $\{\bbb\ddd\eee,\aaa\bbb\ccc,\aaa\ccc\ddd^2\,(\text{or}\,\ddd^3\eee)\}$, $(x={\mathrm e}^{\mathrm{i} \ccc},y={\mathrm e}^{\frac{4\mathrm{i} \pi}{f}})$}

$(y+1) (x^{3} y^{4}-2 x^{3} y^{3}-3 x^{2} y^{4}-x y^{5}+3 x^{2} y^{3}-2 x^{2} y^{2}-2 x y^{3}+3 x y^{2}-x^{2}-3 x y-2 y^{2}+y)=0$.

None.

\subsubsection*{Case $\{\bbb\ddd\eee,\aaa\bbb\ccc,\aaa^2\bbb\ddd^2\}$, $(x={\mathrm e}^{\mathrm{i} \ddd+\frac{2\mathrm{i} \pi}{f}},y={\mathrm e}^{\mathrm{i} \ddd-\frac{2\mathrm{i} \pi}{f}})$}

$\zeta_{4} x^{4} y^{4}+\zeta_{4} x^{3} y^{5}+\zeta_{4} x^{2} y^{6}+x^{5} y^{2}+2 x^{4} y^{3}+2 x^{2} y^{5}-\zeta_{4} x^{5} y+\zeta_{4} x^{3} y^{3}-\zeta_{4} x^{2} y^{4}+\zeta_{4} x y^{5}+x^{4} y-x^{3} y^{2}+x^{2} y^{3}-y^{5}+2 \zeta_{4} x^{3} y+2 \zeta_{4} x y^{3}+\zeta_{4} y^{4}+x^{3}+x^{2} y+x y^{2}=0$.

None.
\subsubsection*{Case $\{\bbb\ddd\eee,\aaa\bbb\ccc,\aaa^5\}$, $(x={\mathrm e}^{2\mathrm{i} \ddd},y={\mathrm e}^{\frac{4\mathrm{i} \pi}{f}})$}

$(\zeta_{5}^{4}-\zeta_{5}^{3}) x y^{5}-\zeta_{5}^{4} x^{2} y^{3}-2 \zeta_{5}^{3} x y^{4}-\zeta_{5}^{3} y^{5}-\zeta_{5}^{4} x^{2} y^{2}-(\zeta_{5}^{4}+1) x y^{3}-(2 \zeta_{5}^{4}-1) y^{4}-(\zeta_{5}^{4}-2) x^{2} y+(\zeta_{5}^{4}+1) x y^{2}+y^{3}+\zeta_{5} x^{2}+2 \zeta_{5} x y+y^{2}+(\zeta_{5}-1) x=0$.

None.

\subsubsection*{Case $\{\bbb\ddd\eee,\aaa\bbb\ccc,\aaa^4\bbb\}$, $(x={\mathrm e}^{2\mathrm{i} \ddd},y={\mathrm e}^{\frac{\mathrm{i} \pi}{f}})$}

$\zeta_{8}^{3} x y^{19}-\zeta_{8}^{2} x y^{18}-\zeta_{8}^{2} y^{18}-y^{16}-2 \zeta_{8}^{2} x y^{14}-2 \zeta_{8}^{3} y^{15}-\zeta_{8}^{3} x^{2} y^{11}+x y^{12}-\zeta_{8}^{3} x y^{11}-y^{12}-\zeta_{8}^{3} x^{2} y^{7}-x y^{8}+\zeta_{8}^{3} x y^{7}-y^{8}-2 x^{2} y^{4}-2 \zeta_{8} x y^{5}-\zeta_{8}^{3} x^{2} y^{3}-\zeta_{8} x^{2} y-\zeta_{8} x y+x=0$.

None.

\subsubsection*{Case $\{\bbb\ddd\eee,\aaa\bbb\ccc,\aaa^3\bbb^2\}$, $(x={\mathrm e}^{2\mathrm{i} \ddd},y={\mathrm e}^{\frac{4\mathrm{i} \pi}{3f}})$}

$x y^{13}-x y^{11}+y^{12}-y^{11}-x y^{9}-2 y^{10}-x^{2} y^{7}-2 x y^{8}+y^{9}-x y^{7}+x y^{6}-x^{2} y^{4}+2 x y^{5}+y^{6}+2 x^{2} y^{3}+x y^{4}+x^{2} y^{2}-x^{2} y+x y^{2}-x=0$.

None.

\subsubsection*{Case $\{\bbb\ddd\eee,\aaa\bbb\ccc,\aaa^3\ddd^2\}$, $(x={\mathrm e}^{\frac{2\mathrm{i} \ddd}{3}},y={\mathrm e}^{\frac{4\mathrm{i} \pi}{f}})$}

$\zeta_{3} x^{7} y^{3}+x^{5} y^{5}+\zeta_{3} x^{7} y^{2}+2 x^{5} y^{4}-\zeta_{3} x^{4} y^{5}+\zeta_{3} x^{7} y-2 \zeta_{3}^{2} x^{6} y+\zeta_{3} x^{4} y^{3}+x^{2} y^{5}-\zeta_{3} x^{4} y^{2}+\zeta_{3}^{2} x^{3} y^{3}-x^{5}-\zeta_{3}^{2} x^{3} y^{2}+2 \zeta_{3} x y^{4}-\zeta_{3}^{2} y^{4}+\zeta_{3}^{2} x^{3}-2 x^{2} y-\zeta_{3}^{2} y^{3}-x^{2}-\zeta_{3}^{2} y^{2}=0$.

None.

\subsubsection*{Case $\{\bbb\ddd\eee,\aaa\bbb\ccc,\aaa^2\ccc\eee^2\}$, $(x={\mathrm e}^{2\mathrm{i} \ddd},y={\mathrm e}^{\frac{4\mathrm{i} \pi}{f}})$}

$x y^{14}+2 x y^{13}+y^{14}+x^{2} y^{11}-x^{3} y^{9}-x^{3} y^{8}-x^{2} y^{9}-2 x y^{10}-x^{3} y^{7}+x^{2} y^{8}+x^{3} y^{6}-x^{2} y^{7}-2 x^{4} y^{4}-x^{3} y^{5}-x^{2} y^{6}-x^{2} y^{5}+x^{3} y^{3}+x^{5}+2 x^{4} y+x^{4}=0$.

None.

\subsubsection*{Case $\{\bbb\ddd\eee,\aaa\bbb\ccc,\aaa\ccc^2\ddd^2\,(\text{or}\,\ccc\ddd^3\eee)\}$, $(x={\mathrm e}^{2\mathrm{i} \ddd},y={\mathrm e}^{\frac{4\mathrm{i} \pi}{f}})$}

$x^{5} y^{3}+2 x^{4} y^{4}+x^{4} y^{3}-x^{3} y^{4}+x^{2} y^{5}+2 x^{4} y^{2}+x^{3} y^{3}+2 x^{2} y^{4}-x^{3} y^{2}+x^{2} y^{3}-2 x^{3} y-x^{2} y^{2}-2 x y^{3}-x^{3}+x^{2} y-x y^{2}-2 x y-y^{2}=0$.

None.

\subsubsection*{Case $\{\bbb\ddd\eee,\aaa\bbb\ccc,\aaa^2\bbb\eee^2\}$, $(x={\mathrm e}^{\mathrm{i} \ddd+\frac{2\mathrm{i} \pi}{f}},y={\mathrm e}^{\mathrm{i} \ddd-\frac{2\mathrm{i} \pi}{f}})$}

$x y^{10}+2 \zeta_{4} x^{2} y^{8}+x^{4} y^{5}+x^{3} y^{6}+x^{2} y^{7}+2 x y^{8}+y^{9}-\zeta_{4} x^{3} y^{5}+\zeta_{4} x^{2} y^{6}-\zeta_{4} y^{8}-x^{5} y^{2}+x^{3} y^{4}-x^{2} y^{5}+\zeta_{4} x^{5} y+2 \zeta_{4} x^{4} y^{2}+\zeta_{4} x^{3} y^{3}+\zeta_{4} x^{2} y^{4}+\zeta_{4} x y^{5}+2 x^{3} y^{2}+\zeta_{4} x^{4}=0$.

None.

\subsubsection*{Case $\{\bbb\ddd\eee,\aaa\bbb\ccc,\aaa^4\ccc\,(\text{or}\,\aaa^3\ddd\eee)\}$, $(x={\mathrm e}^{2\mathrm{i} \ddd},y={\mathrm e}^{\frac{4\mathrm{i} \pi}{3f}})$}

$\zeta_{3}^{2} x y^{18}+\zeta_{3} x y^{17}+\zeta_{3}^{2} y^{18}+2 \zeta_{3}^{2} x y^{15}-2 \zeta_{3} y^{14}-\zeta_{3} x^{2} y^{11}-y^{13}-\zeta_{3} x y^{11}+x y^{10}-\zeta_{3} x^{2} y^{8}-y^{10}+\zeta_{3} x y^{8}-x y^{7}-\zeta_{3} x^{2} y^{5}-y^{7}-2 x^{2} y^{4}+2 \zeta_{3}^{2} x y^{3}+\zeta_{3}^{2} x^{2}+x y+\zeta_{3}^{2} x=0$.

$f=16,(\aaa,\bbb,\ccc,\ddd,\eee)=(1,3,4,1,4)/4$.

$f=16,(\aaa,\bbb,\ccc,\ddd,\eee)=(1,3,4,2,3)/4$.

\subsubsection*{Case $\{\bbb\ddd\eee,\aaa\bbb\ccc,\aaa^2\ccc\ddd^2\,(\text{or}\,\aaa\ddd^3\eee)\}$, $(x={\mathrm e}^{\mathrm{i} \ddd},y={\mathrm e}^{\frac{4\mathrm{i} \pi}{f}})$}

$(y-1) (y+1) (x^{6} y^{6}+2 x^{6} y^{5}+x^{6} y^{4}+x^{4} y^{6}+x^{6} y^{3}+x^{4} y^{5}+x^{4} y^{4}+2 x^{4} y^{2}-2 x^{2} y^{4}-x^{2} y^{2}-x^{2} y-y^{3}-x^{2}-y^{2}-2 y-1)=0$.

$f=16,(\aaa,\bbb,\ccc,\ddd,\eee)=(1,3,4,1,4)/4$.

\subsubsection*{Case $\{\bbb\ddd\eee,\aaa\bbb\ccc,\ccc\ddd^4\}$, $(x={\mathrm e}^{2\mathrm{i} \ddd},y={\mathrm e}^{\frac{4\mathrm{i} \pi}{f}})$}

$(x-1) (x^{7} y+2 x^{6} y^{2}+2 x^{6} y+2 x^{5} y^{2}+x^{4} y^{3}-x^{3} y^{4}+x^{4} y^{2}-x^{3} y^{2}+x^{4}-x^{3} y-2 x^{2} y^{2}-2 x y^{3}-2 x y^{2}-y^{3})=0$.

$f=16,(\aaa,\bbb,\ccc,\ddd,\eee)=(1,3,4,1,4)/4$.

\subsubsection*{Case $\{\bbb\ddd\eee,\aaa\bbb\ccc,\aaa^2\eee^2\}$, $(x={\mathrm e}^{\mathrm{i} \ddd},y={\mathrm e}^{\frac{4\mathrm{i} \pi}{f}})$}

$x^{5} y^{5}-x^{3} y^{7}+x^{2} y^{8}+x^{5} y^{4}+2 x^{2} y^{7}-2 x^{6} y^{2}+x^{5} y^{3}+x^{4} y^{4}+x^{3} y^{5}+y^{8}-x^{7}-x^{4} y^{3}-x^{3} y^{4}-x^{2} y^{5}+2 x y^{6}-2 x^{5} y-x^{2} y^{4}-x^{5}+x^{4} y-x^{2} y^{3}=0$.

$f=16,(\aaa,\bbb,\ccc,\ddd,\eee)=(1,3,4,2,3)/4$.

\subsubsection*{Case $\{\bbb\ddd\eee,\aaa\bbb\ccc,\aaa^3\eee^2\}$, $(x={\mathrm e}^{\frac{2\mathrm{i} \eee}{3}},y={\mathrm e}^{\frac{4\mathrm{i} \pi}{f}})$}

$x^{9} y^{3}+x^{6} y^{5}+2 \zeta_{3} x^{8} y^{2}+2 x^{6} y^{4}-\zeta_{3} x^{5} y^{5}-\zeta_{3}^{2} x^{7} y^{2}-\zeta_{3}^{2} x^{7} y+\zeta_{3} x^{5} y^{3}-\zeta_{3}^{2} x^{7}-\zeta_{3} x^{5} y^{2}+\zeta_{3}^{2} x^{4} y^{3}+\zeta_{3} x^{2} y^{5}-\zeta_{3}^{2} x^{4} y^{2}+\zeta_{3} x^{2} y^{4}+\zeta_{3} x^{2} y^{3}+\zeta_{3}^{2} x^{4}-2 x^{3} y-2 \zeta_{3}^{2} x y^{3}-x^{3}-y^{2}=0$.

$f=20,(\aaa,\bbb,\ccc,\ddd,\eee)=(2,8,10,5,7)/10$.

$f=28,(\aaa,\bbb,\ccc,\ddd,\eee)=(2,12,14,5,11)/14$.

\subsubsection*{Case $\{\bbb\ddd\eee,\aaa\ccc^2,\aaa\ccc\ddd\eee\,(\text{or}\,\aaa\ddd^2\eee^2)\}$, $(x={\mathrm e}^{2\mathrm{i} \ddd},y={\mathrm e}^{\frac{4\mathrm{i} \pi}{f}})$}

$x^{2} y^{7}+2 x^{2} y^{6}+x y^{7}+x^{2} y^{5}+3 x y^{6}-x^{2} y^{4}+x y^{5}-2 x^{2} y^{3}-x y^{4}+y^{5}-x^{2} y^{2}+x y^{3}+2 y^{4}-x y^{2}+y^{3}-3 x y -y^{2}-x -2 y -1=0$.

None.

\subsubsection*{Case $\{\bbb\ddd\eee,\aaa\ccc^2,\aaa\ccc\eee^2\}$, $(x={\mathrm e}^{\mathrm{i} \ddd+\frac{2\mathrm{i} \pi}{f}},y={\mathrm e}^{\mathrm{i} \ddd-\frac{2\mathrm{i} \pi}{f}})$}

$\zeta_{4} x^{3} y^{6}-3 x^{3} y^{5}+x^{2} y^{6}+x y^{7}+\zeta_{4} x^{6} y+2 \zeta_{4} x^{5} y^{2}+2 \zeta_{4} x^{4} y^{3}-\zeta_{4} x^{3} y^{4}+2 \zeta_{4} x^{2} y^{5}+2 \zeta_{4} x y^{6}+2 x^{5} y+2 x^{4} y^{2}-x^{3} y^{3}+2 x^{2} y^{4}+2 x y^{5}+y^{6}+\zeta_{4} x^{5}+\zeta_{4} x^{4} y-3 \zeta_{4} x^{3} y^{2}+x^{3} y=0$.

None.

\subsubsection*{Case $\{\bbb\ddd\eee,\aaa\ccc^2,\ddd^2\eee^2\}$, $(x={\mathrm e}^{2\mathrm{i} \ddd},y={\mathrm e}^{\frac{4\mathrm{i} \pi}{f}})$}

$x^{2} y^{4}+2 x^{2} y^{3}+3 x y^{4}+x^{2} y^{2}+2 x y^{3}-4 x^{2} y+4 y^{3}-2 x y-y^{2}-3 x-2 y-1=0$.

None.
\subsubsection*{Case $\{\bbb\ddd\eee,\aaa\ccc^2,\aaa^3\bbb\ccc\}$, $(x={\mathrm e}^{2\mathrm{i} \ddd},y={\mathrm e}^{\frac{4\mathrm{i} \pi}{f}})$}

$x y^{16}+y^{17}-x y^{15}-y^{16}-x y^{14}-2 y^{15}-y^{14}-x y^{12}-x^{2} y^{10}+2 y^{12}-2 x y^{10}-2 x y^{9}+2 x y^{8}+2 x y^{7}-2 x^{2} y^{5}+y^{7}+x y^{5}+x^{2} y^{3}+2 x^{2} y^{2}+x y^{3}+x^{2} y+x y^{2}-x^{2}-x y=0$.

None.

\subsubsection*{Case $\{\bbb\ddd\eee,\aaa\ccc^2,\aaa^2\bbb^2\ccc\}$, $(x={\mathrm e}^{2\mathrm{i} \ddd},y={\mathrm e}^{\frac{2\mathrm{i} \pi}{f}})$}

$\zeta_{4} y^{13}+\zeta_{4} x y^{11}-x y^{10}-\zeta_{4} y^{11}-3 \zeta_{4} x y^{9}+2 y^{10}-2 \zeta_{4} y^{9}-2 x^{2} y^{6}-3 \zeta_{4} x y^{7}-3 x y^{6}-2 \zeta_{4} y^{7}-2 x^{2} y^{4}+2 \zeta_{4} x^{2} y^{3}-3 x y^{4}-x^{2} y^{2}-\zeta_{4} x y^{3}+x y^{2}+x^{2}=0$.

None.

\subsubsection*{Case $\{\bbb\ddd\eee,\aaa\ccc^2,\aaa\bbb\ccc\eee^2\}$, $(x={\mathrm e}^{2\mathrm{i} \ddd},y={\mathrm e}^{\frac{4\mathrm{i} \pi}{f}})$}

$x y^{7}+4 x y^{6}+y^{7}-x^{2} y^{4}+3 x y^{5}+y^{6}-x^{2} y^{3}+x y^{4}-x^{3} y -3 x^{2} y^{2}+x y^{3}-x^{3}-4 x^{2} y -x^{2}=0$

None.

\subsubsection*{Case $\{\bbb\ddd\eee,\aaa\ccc^2,\aaa^2\ccc\ddd^2\}$, $(x={\mathrm e}^{\mathrm{i} \eee+\frac{2\mathrm{i} \pi}{f}},y={\mathrm e}^{\mathrm{i} \eee-\frac{2\mathrm{i} \pi}{f}})$}

$(x+y) (x^{6} y^{4}-2 x^{5} y^{5}-x^{3} y^{7}+x y^{9}+2 \zeta_{4} x^{7} y^{2}-2 \zeta_{4} x^{6} y^{3}+2 \zeta_{4} x^{5} y^{4}-\zeta_{4} x^{4} y^{5}+\zeta_{4} x^{3} y^{6}+\zeta_{4} x^{2} y^{7}+\zeta_{4} x y^{8}-\zeta_{4} y^{9}-x^{8}+x^{7} y+x^{6} y^{2}+x^{5} y^{3}-x^{4} y^{4}+2 x^{3} y^{5}-2 x^{2} y^{6}+2 x y^{7}+\zeta_{4} x^{7}-\zeta_{4} x^{5} y^{2}-2 \zeta_{4} x^{3} y^{4}+\zeta_{4} x^{2} y^{5})=0$.

None.

\subsubsection*{Case $\{\bbb\ddd\eee,\aaa\ccc^2,\aaa\ddd\eee^3\}$, $(x={\mathrm e}^{2\mathrm{i} \eee},y={\mathrm e}^{\frac{4\mathrm{i} \pi}{f}})$}

$x^{7} y^{12}-2 x^{6} y^{10}-x^{5} y^{11}-x^{5} y^{10}+x^{5} y^{9}+x^{5} y^{8}+2 x^{4} y^{9}-x^{5} y^{7}+2 x^{4} y^{8}-2 x^{5} y^{6}-x^{5} y^{5}+x^{4} y^{6}+x^{3} y^{6}-x^{2} y^{7}-2 x^{2} y^{6}+2 x^{3} y^{4}-x^{2} y^{5}+2 x^{3} y^{3}+x^{2} y^{4}+x^{2} y^{3}-x^{2} y^{2}-x^{2} y-2 x y^{2}+1=0$.

None.

\subsubsection*{Case $\{\bbb\ddd\eee,\aaa\ccc^2,\ccc\ddd^3\eee\}$, $(x={\mathrm e}^{2\mathrm{i} \ddd},y={\mathrm e}^{\frac{4\mathrm{i} \pi}{f}})$}

$x^{5} y^{3}+2 x^{5} y^{2}+x^{4} y^{3}+2 x^{3} y^{4}+x^{5} y+x^{4} y^{2}+2 x^{3} y^{3}-x^{2} y^{4}-x^{5}-3 x^{4} y-x^{3} y^{2}+x^{2} y^{3}+3 x y^{4}+y^{5}+x^{3} y-2 x^{2} y^{2}-x y^{3}-y^{4}-2 x^{2} y-x y^{2}-2 y^{3}-y^{2}=0$.

None.

\subsubsection*{Case $\{\bbb\ddd\eee,\aaa\ccc^2,\ccc\ddd\eee^3\}$, $(x={\mathrm e}^{2\mathrm{i} \eee},y={\mathrm e}^{\frac{4\mathrm{i} \pi}{f}})$}

$x^{7} y+2 x^{6} y^{2}+2 x^{5} y^{3}+2 x^{4} y^{4}+2 x^{4} y^{3}-x^{3} y^{4}+x^{2} y^{5}-3 x^{5} y+3 x^{2} y^{4}-x^{5}+x^{4} y-2 x^{3} y^{2}-2 x^{3} y-2 x^{2} y^{2}-2 x y^{3}-y^{4}=0$.

None.

\subsubsection*{Case $\{\bbb\ddd\eee,\aaa\ccc^2,\aaa\ddd^3\eee\}$, $(x={\mathrm e}^{2\mathrm{i} \ddd},y={\mathrm e}^{\frac{4\mathrm{i} \pi}{f}})$}

$(x y^{2}-1) (x^{4} y^{8}+2 x^{4} y^{7}+x^{3} y^{8}+x^{4} y^{6}+x^{3} y^{7}-x^{4} y^{5}+2 x^{3} y^{5}-x^{2} y^{6}+x^{3} y^{4}-x^{2} y^{5}+x^{3} y^{3}-x y^{5}+x^{2} y^{3}-x y^{4}+x^{2} y^{2}-2 x y^{3}+y^{3}-x y-y^{2}-x-2 y-1)=0$

None.

\subsubsection*{Case $\{\bbb\ddd\eee,\aaa\ccc^2,\aaa^3\}$, $(x={\mathrm e}^{\mathrm{i} \ddd},y={\mathrm e}^{\mathrm{i} \eee})$}

$x^{5} y^{2}-(\zeta_{3}^{2}-\zeta_{3}+1) x^{4} y^{3}-2 x^{4} y-(\zeta_{3}-2) x^{3} y^{2}+(\zeta_{3}^{2}+2 \zeta_{3}-2) x^{2} y^{3}+(\zeta_{3}^{2}-2 \zeta_{3}+2) x^{3}+(2 \zeta_{3}-1) x^{2} y-2 \zeta_{3} x y^{2}-(\zeta_{3}^{2}+\zeta_{3}-1) x+\zeta_{3} y=0$.

$f=12,(\aaa,\bbb,\ccc,\ddd,\eee)=(4,3,4,6,3)/6$.

\subsubsection*{Case $\{\bbb\ddd\eee,\aaa\ccc^2,\aaa\ccc\ddd^2\}$, $(x={\mathrm e}^{\mathrm{i} \eee+\frac{2\mathrm{i} \pi}{f}},y={\mathrm e}^{\mathrm{i} \eee-\frac{2\mathrm{i} \pi}{f}})$}

$(x+y) (x^{3} y^{3}-2 x^{2} y^{4}-x y^{5}+y^{6}+2 \zeta_{4} x^{3} y^{2}-\zeta_{4} x^{2} y^{3}+\zeta_{4} x y^{4}+\zeta_{4} y^{5}+x^{3} y+x^{2} y^{2}-x y^{3}+2 y^{4}+\zeta_{4} x^{3}-\zeta_{4} x^{2} y-2 \zeta_{4} x y^{2}+\zeta_{4} y^{3})=0$.

$f=14,(\aaa,\bbb,\ccc,\ddd,\eee)=(8,6,10,5,17)/14$.

$f=28,(\aaa,\bbb,\ccc,\ddd,\eee)=(2,10,6,3,1)/7$.

$f=42,(\aaa,\bbb,\ccc,\ddd,\eee)=(8,32,38,19,33)/42$.

\subsubsection*{Case $\{\bbb\ddd\eee,\aaa\ccc^2,\aaa^4\ccc\}$, $(x={\mathrm e}^{\mathrm{i} \ddd},y={\mathrm e}^{\mathrm{i} \eee})$}

$\zeta_{7}^{4} x^{5} y^{2}-(\zeta_{7}^{4}+\zeta_{7}-1) x^{4} y^{3}-2 \zeta_{7}^{4} x^{4} y+(\zeta_{7}^{5}+2 \zeta_{7}^{4}-2) x^{3} y^{2}+(\zeta_{7}^{5}-2 \zeta_{7}^{4}+\zeta_{7}+1) x^{2} y^{3}-(2 \zeta_{7}^{5}-\zeta_{7}^{4}-\zeta_{7}^{2}-\zeta_{7}) x^{3}+(2 \zeta_{7}^{5}+\zeta_{7}^{4}-2 \zeta_{7}^{2}) x^{2} y-2 \zeta_{7}^{5} x y^{2}-(\zeta_{7}^{5}-\zeta_{7}^{2}+\zeta_{7}) x+\zeta_{7}^{5} y=0$.

$f=28,(\aaa,\bbb,\ccc,\ddd,\eee)=(2,10,6,3,1)/7$.

$f=28,(\aaa,\bbb,\ccc,\ddd,\eee)=(2,4,6,4,6)/7$.

\subsubsection*{Case $\{\bbb\ddd\eee,\aaa\ccc^2,\aaa^2\ccc\ddd\eee\}$, $(x={\mathrm e}^{2\mathrm{i} \ddd},y={\mathrm e}^{\frac{4\mathrm{i} \pi}{f}})$}

$x^{2} y^{13}+2 x^{2} y^{12}+x y^{13}+x^{2} y^{11}+x y^{12}-x^{2} y^{10}-x y^{11}+2 x y^{10}+2 x y^{9}-2 x^{2} y^{7}+y^{9}+x y^{7}-x y^{6}-x^{2} y^{4}+2 y^{6}-2 x y^{4}-2 x y^{3}+x y^{2}+y^{3}-x y-y^{2}-x-2 y-1=0$.

$f=20,(\aaa,\bbb,\ccc,\ddd,\eee)=(12,48,24,1,11)/30$.

$f=20,(\aaa,\bbb,\ccc,\ddd,\eee)=(12,48,24,11,1)/30$.

$f=28,(\aaa,\bbb,\ccc,\ddd,\eee)=(2,10,6,3,1)/7$.

\subsubsection*{Case $\{\bbb\ddd\eee,\aaa\ccc^2,\ccc\ddd^2\eee^2\}$, $(x={\mathrm e}^{2\mathrm{i} \ddd},y={\mathrm e}^{\frac{2\mathrm{i} \pi}{f}})$}

$\zeta_{4} x^{2} y^{6}+2 x y^{7}+2 \zeta_{4} x y^{6}-y^{7}+2 \zeta_{4} x^{2} y^{4}+x y^{5}-2 \zeta_{4} y^{6}+\zeta_{4} x y^{4}+2 y^{5}+2 \zeta_{4} x^{2} y^{2}+x y^{3}-2 x^{2} y+\zeta_{4} x y^{2}+2 y^{3}-\zeta_{4} x^{2}+2 x y+2 \zeta_{4} x+y=0$.

$f=28,(\aaa,\bbb,\ccc,\ddd,\eee)=(2,10,6,3,1)/7$.

\subsubsection*{Case $\{\bbb\ddd\eee,\aaa\ccc^2,\ddd\eee^3\}$, $(x={\mathrm e}^{\frac{2\mathrm{i} \ddd}{3}},y={\mathrm e}^{\frac{4\mathrm{i} \pi}{f}})$}

$\zeta_{3}^{2} x^{5} y^{4}-\zeta_{3} x^{7} y+2 \zeta_{3}^{2} x^{5} y^{3}+2 x^{6} y-\zeta_{3} x^{4} y^{3}-2 x^{3} y^{4}-2 x^{3} y^{3}+\zeta_{3}^{2} x^{2} y^{4}+\zeta_{3}^{2} x^{5}-2 \zeta_{3} x^{4} y-2 \zeta_{3} x^{4}-x^{3} y+2 \zeta_{3} x y^{3}+2 \zeta_{3}^{2} x^{2} y-y^{3}+\zeta_{3}^{2} x^{2}=0$.

$f=14,(\aaa,\bbb,\ccc,\ddd,\eee)=(8,6,10,19,3)/14$.

$f=16,(\aaa,\bbb,\ccc,\ddd,\eee)=(4,2,6,13,1)/8$.

\subsubsection*{Case $\{\bbb\ddd\eee,\aaa\ccc^2,\aaa^2\bbb\ccc\}$, $(x={\mathrm e}^{2\mathrm{i} \ddd},y={\mathrm e}^{\frac{4\mathrm{i} \pi}{f}})$}

$(x+y) (y^{10}-y^{9}-2 y^{8}-x y^{6}+y^{7}-x y^{3}+y^{4}+2 x y^{2}+x y-x)=0$.

$f=14,(\aaa,\bbb,\ccc,\ddd,\eee)=(8,2,10,23,3)/14$.

$f=16,(\aaa,\bbb,\ccc,\ddd,\eee)=(4,2,6,13,1)/8$.

$f=18,(\aaa,\bbb,\ccc,\ddd,\eee)=(8,6,14,29,1)/18$.

$f=20,(\aaa,\bbb,\ccc,\ddd,\eee)=(4,4,8,7,9)/10$.

$f=28,(\aaa,\bbb,\ccc,\ddd,\eee)=(2,4,6,4,6)/7$.

$f=36,(\aaa,\bbb,\ccc,\ddd,\eee)=(2,6,8,5,7)/9$.

$(\aaa,\bbb,\ccc,\ddd,\eee)=(16,2f-24,2f-8,f+4,f+20)/2f$.

\subsubsection*{Case $\{\bbb\ddd\eee,\aaa\ccc^2,\aaa^2\ccc\eee^2\}$, $(x={\mathrm e}^{\mathrm{i} \ddd+\frac{2\mathrm{i} \pi}{f}},y={\mathrm e}^{\mathrm{i} \ddd-\frac{2\mathrm{i} \pi}{f}})$}

$\zeta_{4} x^{4} y^{8}-2 x^{5} y^{6}-x^{3} y^{8}+x^{2} y^{9}+x y^{10}+\zeta_{4} x^{9} y+2 \zeta_{4} x^{8} y^{2}+\zeta_{4} x^{7} y^{3}-\zeta_{4} x^{6} y^{4}+\zeta_{4} x^{5} y^{5}+2 \zeta_{4} x^{3} y^{7}+2 \zeta_{4} x^{2} y^{8}+2 x^{7} y^{2}+2 x^{6} y^{3}+x^{4} y^{5}-x^{3} y^{6}+x^{2} y^{7}+2 x y^{8}+y^{9}+\zeta_{4} x^{8}+\zeta_{4} x^{7} y-\zeta_{4} x^{6} y^{2}-2 \zeta_{4} x^{4} y^{4}+x^{5} y^{2}=0$.

$f=16,(\aaa,\bbb,\ccc,\ddd,\eee)=(4,2,6,13,1)/8$.

\subsubsection*{Case $\{\bbb\ddd\eee,\aaa\ccc^2,\aaa^3\ccc\}$, $(x={\mathrm e}^{\mathrm{i} \ddd},y={\mathrm e}^{\mathrm{i} \eee})$}

$\zeta_{5}^{3} x^{5} y^{2}-(\zeta_{5}^{3}+\zeta_{5}-1) x^{4} y^{3}-2 \zeta_{5}^{3} x^{4} y+(\zeta_{5}^{4}+2 \zeta_{5}^{3}-2) x^{3} y^{2}+(\zeta_{5}^{4}-2 \zeta_{5}^{3}+\zeta_{5}+1) x^{2} y^{3}-(2 \zeta_{5}^{4}-\zeta_{5}^{3}-\zeta_{5}^{2}-\zeta_{5}) x^{3}+(2 \zeta_{5}^{4}+\zeta_{5}^{3}-2 \zeta_{5}^{2}) x^{2} y-2 \zeta_{5}^{4} x y^{2}-(\zeta_{5}^{4}-\zeta_{5}^{2}+\zeta_{5}) x+\zeta_{5}^{4} y=0$.

$f=20,(\aaa,\bbb,\ccc,\ddd,\eee)=(12,48,24,1,11)/30$.

$f=20,(\aaa,\bbb,\ccc,\ddd,\eee)=(12,48,24,11,1)/30$.

$f=20,(\aaa,\bbb,\ccc,\ddd,\eee)=(2,2,4,3,5)/5$.

$f=20,(\aaa,\bbb,\ccc,\ddd,\eee)=(4,4,8,7,9)/10$.

\subsubsection*{Case $\{\bbb\ddd\eee,\aaa\ccc^2,\aaa^2\ddd\eee\}$, $(x={\mathrm e}^{2\mathrm{i} \ddd},y={\mathrm e}^{\frac{4\mathrm{i} \pi}{f}})$}

$(y+1) (x^{2} y^{14}+x y^{14}-x y^{13}-x^{2} y^{11}+2 x^{2} y^{10}+x y^{11}-3 x^{2} y^{9}-4 x y^{10}+4 x^{2} y^{8}+5 x y^{9}-3 x^{2} y^{7}-6 x y^{8}-y^{9}+2 x^{2} y^{6}+6 x y^{7}+2 y^{8}-x^{2} y^{5}-6 x y^{6}-3 y^{7}+5 x y^{5}+4 y^{6}-4 x y^{4}-3 y^{5}+x y^{3}+2 y^{4}-y^{3}-x y+x+1)=0$.

$f=28,(\aaa,\bbb,\ccc,\ddd,\eee)=(2,4,6,4,6)/7$.

\subsubsection*{Case $\{\bbb\ddd\eee,\aaa\ccc^2,\aaa^3\ddd\eee\}$, $(x={\mathrm e}^{2\mathrm{i} \ddd},y={\mathrm e}^{\frac{4\mathrm{i} \pi}{f}})$}

$(y+1) (x^{2} y^{20}+x y^{20}-x y^{19}-x^{2} y^{17}+2 x^{2} y^{16}+x y^{17}-3 x^{2} y^{15}-2 y^{16} x+4 x^{2} y^{14}+3 x y^{15}-5 x^{2} y^{13}-6 x y^{14}+6 x^{2} y^{12}+7 x y^{13}-5 x^{2} y^{11}-8 x y^{12}-y^{13}+4 x^{2} y^{10}+8 x y^{11}+2 y^{12}-3 x^{2} y^{9}-8 x y^{10}-3 y^{11}+2 x^{2} y^{8}+8 x y^{9}+4 y^{10}-x^{2} y^{7}-8 x y^{8}-5 y^{9}+7 x y^{7}+6 y^{8}-6 x y^{6}-5 y^{7}+3 x y^{5}+4 y^{6}-2 x y^{4}-3 y^{5}+x y^{3}+2 y^{4}-y^{3}-x y+x+1)=0$.

$f=36,(\aaa,\bbb,\ccc,\ddd,\eee)=(2,6,8,5,7)/9$.

\subsection*{Rational angles solution for all cases with a unique $a^2b$-vertex type  $\aaa\ddd\eee$: $28$ distinct cases}

\subsubsection*{Case $\{\aaa\ddd\eee,\aaa^3,\aaa^2\bbb\,(\text{or}\,\aaa\bbb^2,\bbb^3)\}$, $(x={\mathrm e}^{2\mathrm{i} \ddd},y={\mathrm e}^{\frac{4\mathrm{i} \pi}{f}})$}

$x^{2} y^{2}+(2 \zeta_{3}^{2}+\zeta_{3}-1) x^{2} y+(\zeta_{3}-3) x y^{2}+\zeta_{3}^{2} x^{2}-\zeta_{3}^{2} x y+\zeta_{3}^{2} y^{2}-(3 \zeta_{3}-1) x+(2 \zeta_{3}^{2}-\zeta_{3}+1) y+\zeta_{3}=0$.

None.

\subsubsection*{Case $\{\aaa\ddd\eee,\aaa^3,\aaa\bbb^3\}$, $(x={\mathrm e}^{2\mathrm{i} \ddd},y={\mathrm e}^{\frac{4\mathrm{i} \pi}{f}})$}

$\zeta_{9}^{6} x^{2} y^{2}+(2 \zeta_{9}^{8}-\zeta_{9}^{4}+2 \zeta_{9}^{2}-1) x^{2} y-(2 \zeta_{9}^{7}+\zeta_{9}^{6}-\zeta_{9}^{4}+\zeta_{9}^{3}-\zeta_{9}) x y^{2}+\zeta_{9} x^{2}+(\zeta_{9}^{7}+\zeta_{9}^{6}+\zeta_{9}^{3}-2 \zeta_{9}^{2}+\zeta_{9}) x y+\zeta_{9}^{3} y^{2}-(\zeta_{9}^{7}+2 \zeta_{9}^{6}-\zeta_{9}^{3}+\zeta_{9}-1) x+(2 \zeta_{9}^{5}-\zeta_{9}^{4}+2 \zeta_{9}^{2}-1) y+\zeta_{9}^{7}=0$.

None.

\subsubsection*{Case $\{\aaa\ddd\eee,\aaa^3,\bbb^4\}$, $(x={\mathrm e}^{2\mathrm{i} \ddd},y={\mathrm e}^{\frac{4\mathrm{i} \pi}{f}})$}

$\zeta_{12} x^{2} y^{2}+(2 \zeta_{12}^{4}-2 \zeta_{12}^{2}) x^{2} y-(\zeta_{12}^{4}-\zeta_{12}^{3}+2 \zeta_{12}^{2}+\zeta_{12}+1) x y^{2}-\zeta_{12} x^{2}+2 \zeta_{12}^{2} x y-\zeta_{12}^{3} y^{2}-(\zeta_{12}^{4}+\zeta_{12}^{3}+2 \zeta_{12}^{2}-\zeta_{12}+1) x-(2 \zeta_{12}^{2}-2) y+\zeta_{12}^{3}=0$.

None.

%%\subsubsection*{Case $\{\aaa\ddd\eee,\aaa^3,\bbb\ccc^2\}$, $(x={\mathrm e}^{2\mathrm{i} \ddd},y={\mathrm e}^{\frac{4\mathrm{i} \pi}{f}})$}

%%$\zeta_{3}^{2} x^{2} y^{5}+(2 \zeta_{3}^{2}-\zeta_{3}-1) x y^{6}-\zeta_{3} x^{2} y^{4}-(\zeta_{3}+1) x y^{5}-(2 \zeta_{3}^{2}+2 \zeta_{3}) x^{2} y^{3}+(\zeta_{3}+1) x y^{4}+\zeta_{3}^{2} y^{5}-\zeta_{3} x^{2} y^{2}+2 \zeta_{3}^{2} x y^{3}-y^{4}+\zeta_{3}^{2} x^{2} y+(\zeta_{3}+1) x y^{2}-(2 \zeta_{3}^{2}+2) y^{3}-(\zeta_{3}+1) x y-y^{2}+(2 \zeta_{3}^{2}-\zeta_{3}-1) x+\zeta_{3}^{2} y=0$.

%%$f=20,(\aaa,\bbb,\ccc,\ddd,\eee)=(10,6,12,15,5)/15$.

\subsubsection*{Case $\{\aaa\ddd\eee,\aaa^3,\bbb^3\ccc\}$, $(x={\mathrm e}^{2\mathrm{i} \ddd},y={\mathrm e}^{\frac{2\mathrm{i} \pi}{f}})$}

$\zeta_{12}^{2} x^{2} y^{7}-(\zeta_{12}^{5}+2 \zeta_{12}^{3}+\zeta_{12}) x y^{8}+(\zeta_{12}^{4}-\zeta_{12}^{2}) x y^{7}+x^{2} y^{5}-\zeta_{12}^{4} y^{7}+(2 \zeta_{12}^{5}-2 \zeta_{12}^{3}) x^{2} y^{4}-(\zeta_{12}^{4}-\zeta_{12}^{2}) x y^{5}+\zeta_{12}^{6} x^{2} y^{3}+2 \zeta_{12}^{3} x y^{4}+y^{5}+(\zeta_{12}^{4}-\zeta_{12}^{2}) x y^{3}-(2 \zeta_{12}^{3}-2 \zeta_{12}) y^{4}-\zeta_{12}^{2} x^{2} y+\zeta_{12}^{6} y^{3}-(\zeta_{12}^{4}-\zeta_{12}^{2}) x y-(\zeta_{12}^{5}+2 \zeta_{12}^{3}+\zeta_{12}) x+\zeta_{12}^{4} y=0$.

$f=20,(\aaa,\bbb,\ccc,\ddd,\eee)=(10,6,12,15,5)/15$.

\subsubsection*{Case $\{\aaa\ddd\eee,\aaa^3,\bbb^5\}$, $(x={\mathrm e}^{2\mathrm{i} \ddd},y={\mathrm e}^{\frac{4\mathrm{i} \pi}{f}})$}

$\zeta_{15}^{3} x^{2} y^{2}-(\zeta_{15}^{14}-2 \zeta_{15}^{11}+\zeta_{15}^{8}-2 \zeta_{15}^{6}) x^{2} y-(\zeta_{15}^{13}-\zeta_{15}^{10}+2 \zeta_{15}^{5}+\zeta_{15}^{3}-1) x y^{2}+\zeta_{15}^{9} x^{2}+(\zeta_{15}^{13}-2 \zeta_{15}^{11}+\zeta_{15}^{9}+\zeta_{15}^{4}+\zeta_{15}^{3}) x y+\zeta_{15}^{13} y^{2}+(\zeta_{15}^{12}-\zeta_{15}^{9}+\zeta_{15}^{7}-\zeta_{15}^{4}-2 \zeta_{15}^{2}) x-(\zeta_{15}^{14}-2 \zeta_{15}^{11}+\zeta_{15}^{8}-2 \zeta_{15}) y+\zeta_{15}^{4}=0$.

$f=20,(\aaa,\bbb,\ccc,\ddd,\eee)=(10,6,12,15,5)/15$.

\subsubsection*{Case $\{\aaa\ddd\eee,\aaa^2\ccc,\aaa\bbb^2\}$, $(x={\mathrm e}^{2\mathrm{i} \ddd},y={\mathrm e}^{\frac{4\mathrm{i} \pi}{5f}})$}

$\zeta_{5}^{3} y^{16}-x y^{14}-\zeta_{5} x y^{13}+3 \zeta_{5}^{3} x y^{11}-2 \zeta_{5}^{2} y^{12}-x^{2} y^{9}+\zeta_{5}^{4} x y^{10}-\zeta_{5} x^{2} y^{8}-\zeta_{5}^{2} x^{2} y^{7}+4 \zeta_{5} x y^{8}-y^{9}-\zeta_{5} y^{8}+\zeta_{5}^{3} x y^{6}-\zeta_{5}^{2} y^{7}-2 x^{2} y^{4}+3 \zeta_{5}^{4} x y^{5}-\zeta_{5} x y^{3}-\zeta_{5}^{2} x y^{2}+\zeta_{5}^{4} x^{2}=0$.

None.

\subsubsection*{Case $\{\aaa\ddd\eee,\aaa^2\ccc,\aaa\bbb\ccc\}$, $(x={\mathrm e}^{2\mathrm{i} \ddd},y={\mathrm e}^{\frac{4\mathrm{i} \pi}{f}})$}

$(y+1) (x y^{6}+x y^{5}+y^{6}-x^{2} y^{3}-3 x y^{4}-2 y^{5}+x^{2} y^{2}+4 x y^{3}+y^{4}-2 x^{2} y-3 x y^{2}-y^{3}+x^{2}+x y+x)=0$.

None.

\subsubsection*{Case $\{\aaa\ddd\eee,\aaa^2\ccc,\bbb^3\}$, $(x={\mathrm e}^{2\mathrm{i} \ddd},y={\mathrm e}^{\frac{2\mathrm{i} \pi}{f}})$}

$(\zeta_{12}^{2}+1) x y^{7}+\zeta_{12} y^{8}+2 \zeta_{12} x y^{6}+(\zeta_{12}^{5}+\zeta_{12}^{3}+\zeta_{12}) x^{2} y^{4}+(\zeta_{12}^{4}+\zeta_{12}^{2}) x y^{5}+2 \zeta_{12}^{5} y^{6}-4 \zeta_{12}^{3} x y^{4}+2 \zeta_{12} x^{2} y^{2}+(\zeta_{12}^{4}+\zeta_{12}^{2}) x y^{3}+(\zeta_{12}^{5}+\zeta_{12}^{3}+\zeta_{12}) y^{4}+2 \zeta_{12}^{5} x y^{2}+\zeta_{12}^{5} x^{2}+(\zeta_{12}^{6}+\zeta_{12}^{4}) x y=0$.

None.

%%\subsubsection*{Case $\{\aaa\ddd\eee,\aaa^2\ccc,\bbb\ccc^2\}$, $(x={\mathrm e}^{2\mathrm{i} \ddd},y={\mathrm e}^{\frac{2\mathrm{i} \pi}{f}})$}

%%$(y^{2}+1) (x y^{12}+2 \zeta_{4} x y^{11}+2 \zeta_{4} x^{2} y^{9}-x y^{10}-2 \zeta_{4} x y^{9}-x y^{8}+\zeta_{4} y^{9}+2 \zeta_{4} x y^{7}-\zeta_{4} x^{2} y^{5}-\zeta_{4} y^{7}+2 \zeta_{4} x y^{5}+\zeta_{4} x^{2} y^{3}-\zeta_{4}^{2} x y^{4}-2 \zeta_{4} x y^{3}-\zeta_{4}^{2} x y^{2}+2 \zeta_{4} y^{3}+2 \zeta_{4} x y+\zeta_{4}^{2} x)=0$.

%%None.

\subsubsection*{Case $\{\aaa\ddd\eee,\aaa^2\ccc,\aaa\bbb^3\}$, $(x={\mathrm e}^{2\mathrm{i} \ddd},y={\mathrm e}^{\frac{4\mathrm{i} \pi}{7f}})$}

$\zeta_{7}^{2} y^{24}+\zeta_{7} x y^{21}+\zeta_{7}^{3} x y^{20}+2 \zeta_{7}^{2} x y^{17}-2 y^{18}-\zeta_{7}^{4} x y^{16}-\zeta_{7}^{6} x y^{15}-\zeta_{7}^{3} x^{2} y^{13}-\zeta_{7}^{5} x^{2} y^{12}-x^{2} y^{11}+4 \zeta_{7}^{5} x y^{12}-\zeta_{7}^{3} y^{13}-\zeta_{7}^{5} y^{12}-y^{11}-\zeta_{7}^{4} x y^{9}-\zeta_{7}^{6} x y^{8}-2 \zeta_{7}^{3} x^{2} y^{6}+2 \zeta_{7} x y^{7}+x y^{4}+\zeta_{7}^{2} x y^{3}+\zeta_{7} x^{2}=0$.

None.

\subsubsection*{Case $\{\aaa\ddd\eee,\aaa^2\ccc,\bbb^4\}$, $(x={\mathrm e}^{2\mathrm{i} \ddd},y={\mathrm e}^{\frac{2\mathrm{i} \pi}{f}})$}

$(\zeta_{8}^{3}+\zeta_{8}) x y^{7}+\zeta_{8}^{2} y^{8}+2 \zeta_{8}^{2} x y^{6}+\zeta_{8}^{2} x^{2} y^{4}+(\zeta_{8}^{3}+\zeta_{8}) x y^{5}-2 y^{6}-4 \zeta_{8}^{2} x y^{4}-2 \zeta_{8}^{4} x^{2} y^{2}+(\zeta_{8}^{3}+\zeta_{8}) x y^{3}+\zeta_{8}^{2} y^{4}+2 \zeta_{8}^{2} x y^{2}+\zeta_{8}^{2} x^{2}+(\zeta_{8}^{3}+\zeta_{8}) x y=0$.

None.

\subsubsection*{Case $\{\aaa\ddd\eee,\aaa^2\ccc,\bbb^3\ccc\}$, $(x={\mathrm e}^{2\mathrm{i} \ddd},y={\mathrm e}^{\frac{\mathrm{i} \pi}{f}})$}

$\zeta_{8}^{3} x y^{23}+\zeta_{8}^{2} y^{24}+\zeta_{8} x y^{21}+2 \zeta_{8}^{2} x y^{20}-2 y^{18}-x^{2} y^{14}+\zeta_{8}^{3} x y^{15}+\zeta_{8}^{2} x^{2} y^{12}+\zeta_{8} x y^{13}-y^{14}-4 \zeta_{8}^{2} x y^{12}-\zeta_{8}^{4} x^{2} y^{10}+\zeta_{8}^{3} x y^{11}+\zeta_{8}^{2} y^{12}+\zeta_{8} x y^{9}-\zeta_{8}^{4} y^{10}-2 \zeta_{8}^{4} x^{2} y^{6}+2 \zeta_{8}^{2} x y^{4}+\zeta_{8}^{3} x y^{3}+\zeta_{8}^{2} x^{2}+\zeta_{8} x y=0$.

None.

\subsubsection*{Case $\{\aaa\ddd\eee,\aaa^2\ccc,\bbb^5\}$, $(x={\mathrm e}^{2\mathrm{i} \ddd},y={\mathrm e}^{\frac{2\mathrm{i} \pi}{f}})$}

$(\zeta_{20}^{4}-1) x y^{7}+\zeta_{20}^{7} y^{8}+2 \zeta_{20}^{7} x y^{6}-(\zeta_{20}^{9}-\zeta_{20}^{5}+\zeta_{20}) x^{2} y^{4}+(\zeta_{20}^{8}+\zeta_{20}^{2}) x y^{5}-2 \zeta_{20} y^{6}-4 \zeta_{20}^{5} x y^{4}-2 \zeta_{20}^{9} x^{2} y^{2}+(\zeta_{20}^{8}+\zeta_{20}^{2}) x y^{3}-(\zeta_{20}^{9}-\zeta_{20}^{5}+\zeta_{20}) y^{4}+2 \zeta_{20}^{3} x y^{2}+\zeta_{20}^{3} x^{2}-(\zeta_{20}^{10}-\zeta_{20}^{6}) x y=0$.

None.

\subsubsection*{Case $\{\aaa\ddd\eee,\aaa\ccc^2,\aaa^2\bbb\}$, $(x={\mathrm e}^{2\mathrm{i} \ddd},y={\mathrm e}^{\frac{4\mathrm{i} \pi}{5f}})$}

$x y^{12}-\zeta_{5}^{2} x^{2} y^{10}+\zeta_{5} x y^{11}+\zeta_{5} y^{11}-3 \zeta_{5}^{3} x y^{9}+\zeta_{5}^{2} y^{10}-\zeta_{5}^{4} x y^{8}+\zeta_{5}^{3} y^{9}+2 \zeta_{5} x^{2} y^{6}-4 \zeta_{5} x y^{6}+2 \zeta_{5} y^{6}+\zeta_{5}^{4} x^{2} y^{3}-\zeta_{5}^{3} x y^{4}+x^{2} y^{2}-3 \zeta_{5}^{4} x y^{3}+\zeta_{5} x^{2} y+\zeta_{5} x y-y^{2}+\zeta_{5}^{2} x=0$.

None.

\subsubsection*{Case $\{\aaa\ddd\eee,\aaa\ccc^2,\bbb^3\}$, $(x={\mathrm e}^{2\mathrm{i} \ddd},y={\mathrm e}^{\frac{4\mathrm{i} \pi}{f}})$}

$y^{10}-(\zeta_{3}^{2}-\zeta_{3}) x y^{8}+2 \zeta_{3} y^{9}+(\zeta_{3}^{2}+\zeta_{3}-2) x y^{7}+\zeta_{3}^{2} y^{8}+(2 \zeta_{3}^{2}-\zeta_{3}-1) x^{2} y^{5}+(\zeta_{3}^{2}-\zeta_{3}-2) x y^{6}-4 \zeta_{3}^{2} x y^{5}+(\zeta_{3}^{2}-2 \zeta_{3}-1) x y^{4}+(2 \zeta_{3}^{2}-\zeta_{3}-1) y^{5}+\zeta_{3}^{2} x^{2} y^{2}+(\zeta_{3}^{2}-2 \zeta_{3}+1) x y^{3}+2 x^{2} y-(\zeta_{3}^{2}-1) x y^{2}+\zeta_{3} x^{2}=0$.

None.

\subsubsection*{Case $\{\aaa\ddd\eee,\aaa\ccc^2,\aaa\bbb^3\}$, $(x={\mathrm e}^{2\mathrm{i} \ddd},y={\mathrm e}^{\frac{4\mathrm{i} \pi}{5f}})$}

$\zeta_{5} y^{30}+2 \zeta_{5}^{2} y^{27}-\zeta_{5}^{3} x y^{24}+x y^{23}+\zeta_{5}^{3} y^{24}+\zeta_{5}^{2} x y^{22}-2 \zeta_{5}^{4} x y^{21}+\zeta_{5} x y^{20}+\zeta_{5}^{3} x y^{19}-\zeta_{5}^{2} x^{2} y^{17}-x y^{18}+2 \zeta_{5} x^{2} y^{15}-2 \zeta_{5}^{4} x y^{16}-\zeta_{5}^{2} y^{17}-4 \zeta_{5} x y^{15}-x^{2} y^{13}-2 \zeta_{5}^{3} x y^{14}+2 \zeta_{5} y^{15}-\zeta_{5}^{2} x y^{12}-y^{13}+\zeta_{5}^{4} x y^{11}+\zeta_{5} x y^{10}-2 \zeta_{5}^{3} x y^{9}+x y^{8}+\zeta_{5}^{4} x^{2} y^{6}+\zeta_{5}^{2} x y^{7}-\zeta_{5}^{4} x y^{6}+2 x^{2} y^{3}+\zeta_{5} x^{2}=0$.

None.

\subsubsection*{Case $\{\aaa\ddd\eee,\aaa\ccc^2,\bbb^4\}$, $(x={\mathrm e}^{2\mathrm{i} \ddd},y={\mathrm e}^{\frac{4\mathrm{i} \pi}{f}})$}

$y^{10}+(\zeta_{4}+1) x y^{8}+2 \zeta_{4} y^{9}+2 \zeta_{4} x y^{7}-y^{8}+2 \zeta_{4} x^{2} y^{5}+(3 \zeta_{4}-1) x y^{6}-4 \zeta_{4} x y^{5}-(\zeta_{4}^{2}-3 \zeta_{4}) x y^{4}+2 \zeta_{4} y^{5}-\zeta_{4}^{2} x^{2} y^{2}+2 \zeta_{4} x y^{3}+2 \zeta_{4} x^{2} y+(\zeta_{4}^{2}+\zeta_{4}) x y^{2}+\zeta_{4}^{2} x^{2}=0$.

None.

%%\subsubsection*{Case $\{\aaa\ddd\eee,\aaa\ccc^2,\bbb\ccc^2\}$, $(x={\mathrm e}^{2\mathrm{i} \ddd},y={\mathrm e}^{\frac{4\mathrm{i} \pi}{f}})$}

%%$(y+1) (x^{2} y^{8}+x y^{8}-x y^{7}-x^{2} y^{5}-2 x y^{6}+2 x^{2} y^{4}+3 x y^{5}-x^{2} y^{3}-4 x y^{4}-y^{5}+3 x y^{3}+2 y^{4}-2 x y^{2}-y^{3}-x y+x+1)=0$.

%%$f=20,(\aaa,\bbb,\ccc,\ddd,\eee)=(2,2,4,3,5)/5$.

%%$f=20,(\aaa,\bbb,\ccc,\ddd,\eee)=(4,4,8,7,9)/10$.

\subsubsection*{Case $\{\aaa\ddd\eee,\aaa\ccc^2,\bbb^3\ccc\}$, $(x={\mathrm e}^{2\mathrm{i} \ddd},y={\mathrm e}^{\frac{2\mathrm{i} \pi}{f}})$}

$y^{30}+2 \zeta_{4} y^{27}+\zeta_{4} x y^{25}+x y^{24}-y^{24}-x y^{22}+2 \zeta_{4} x y^{21}+x y^{20}+2 \zeta_{4} x y^{19}-x y^{18}-x^{2} y^{16}+\zeta_{4} x y^{17}+2 \zeta_{4} x^{2} y^{15}-\zeta_{4}^{2} x^{2} y^{14}-4 \zeta_{4} x y^{15}-y^{16}+2 \zeta_{4} y^{15}+\zeta_{4} x y^{13}-\zeta_{4}^{2} y^{14}-\zeta_{4}^{2} x y^{12}+2 \zeta_{4} x y^{11}+\zeta_{4}^{2} x y^{10}+2 \zeta_{4} x y^{9}-\zeta_{4}^{2} x y^{8}-\zeta_{4}^{2} x^{2} y^{6}+\zeta_{4}^{2} x y^{6}+\zeta_{4} x y^{5}+2 \zeta_{4} x^{2} y^{3}+\zeta_{4}^{2} x^{2}=0$.

$f=20,(\aaa,\bbb,\ccc,\ddd,\eee)=(2,2,4,3,5)/5$.

$f=20,(\aaa,\bbb,\ccc,\ddd,\eee)=(4,4,8,7,9)/10$.

\subsubsection*{Case $\{\aaa\ddd\eee,\aaa\ccc^2,\bbb^5\}$, $(x={\mathrm e}^{2\mathrm{i} \ddd},y={\mathrm e}^{\frac{4\mathrm{i} \pi}{f}})$}

$\zeta_{5} y^{10}-(\zeta_{5}^{3}-\zeta_{5}^{2}) x y^{8}+2 \zeta_{5}^{2} y^{9}-(2 \zeta_{5}^{4}-\zeta_{5}^{3}-1) x y^{7}+\zeta_{5}^{3} y^{8}-(\zeta_{5}^{2}-2 \zeta_{5}+1) x^{2} y^{5}-(2 \zeta_{5}^{4}-\zeta_{5}+1) x y^{6}-4 \zeta_{5} x y^{5}-(2 \zeta_{5}^{3}+\zeta_{5}^{2}-\zeta_{5}) x y^{4}-(\zeta_{5}^{2}-2 \zeta_{5}+1) y^{5}+\zeta_{5}^{4} x^{2} y^{2}+(\zeta_{5}^{4}-2 \zeta_{5}^{3}+\zeta_{5}^{2}) x y^{3}+2 x^{2} y-(\zeta_{5}^{4}-1) x y^{2}+\zeta_{5} x^{2}=0$.

$f=20,(\aaa,\bbb,\ccc,\ddd,\eee)=(2,2,4,3,5)/5$.

$f=20,(\aaa,\bbb,\ccc,\ddd,\eee)=(4,4,8,7,9)/10$.

\subsubsection*{Case $\{\aaa\ddd\eee,\ccc^3,\aaa^2\bbb\}$, $(x={\mathrm e}^{2\mathrm{i} \ddd},y={\mathrm e}^{\frac{2\mathrm{i} \pi}{f}})$}

$\zeta_{12}^{5} x^{2} y^{5}-(\zeta_{12}^{2}+1) x y^{6}-2 \zeta_{12} x y^{5}-2 \zeta_{12}^{3} x^{2} y^{3}-(\zeta_{12}^{4}+\zeta_{12}^{2}) x y^{4}-(\zeta_{12}^{5}+\zeta_{12}^{3}-\zeta_{12}) y^{5}+4 \zeta_{12}^{3} x y^{3}+(\zeta_{12}^{5}-\zeta_{12}^{3}-\zeta_{12}) x^{2} y-(\zeta_{12}^{4}+\zeta_{12}^{2}) x y^{2}-2 \zeta_{12}^{3} y^{3}-2 \zeta_{12}^{5} x y-(\zeta_{12}^{6}+\zeta_{12}^{4}) x+\zeta_{12} y=0$.

None.

\subsubsection*{Case $\{\aaa\ddd\eee,\ccc^3,\aaa\bbb\ccc\}$, $(x={\mathrm e}^{2\mathrm{i} \ddd},y={\mathrm e}^{\frac{4\mathrm{i} \pi}{f}})$}

$\zeta_{3} x^{2} y^{3}+(\zeta_{3}-1) x y^{4}+2 \zeta_{3}^{2} x^{2} y^{2}+(2 \zeta_{3}^{2}-\zeta_{3}-3) x y^{3}+(2 \zeta_{3}^{2}-\zeta_{3}-1) y^{4}+x^{2} y-4 \zeta_{3}^{2} x y^{2}+\zeta_{3} y^{3}+(2 \zeta_{3}^{2}-\zeta_{3}-1) x^{2}+(2 \zeta_{3}^{2}-3 \zeta_{3}-1) x y+2 \zeta_{3}^{2} y^{2}-(\zeta_{3}-1) x+y=0$.

None.

\subsubsection*{Case $\{\aaa\ddd\eee,\ccc^3,\aaa\bbb^2\}$, $(x={\mathrm e}^{2\mathrm{i} \ddd},y={\mathrm e}^{\frac{4\mathrm{i} \pi}{f}})$}

$(\zeta_{3}^{2}-\zeta_{3}) x y^{7}+(\zeta_{3}^{2}-2 \zeta_{3}+1) y^{8}-\zeta_{3} x^{2} y^{5}-(\zeta_{3}^{2}+\zeta_{3}-2) x y^{6}-2 \zeta_{3}^{2} x^{2} y^{4}-(\zeta_{3}^{2}-\zeta_{3}-2) x y^{5}-x^{2} y^{3}+4 \zeta_{3}^{2} x y^{4}-\zeta_{3} y^{5}-(\zeta_{3}^{2}-2 \zeta_{3}-1) x y^{3}-2 \zeta_{3}^{2} y^{4}-(\zeta_{3}^{2}-2 \zeta_{3}+1) x y^{2}-y^{3}+(\zeta_{3}^{2}+\zeta_{3}-2) x^{2}+(\zeta_{3}^{2}-1) x y=0$.

None.

\subsubsection*{Case $\{\aaa\ddd\eee,\ccc^3,\aaa\bbb^3\}$, $(x={\mathrm e}^{2\mathrm{i} \ddd},y={\mathrm e}^{\frac{4\mathrm{i} \pi}{f}})$}

$(2 \zeta_{3}^{2}-\zeta_{3}-1) y^{12}+(\zeta_{3}-1) x y^{10}+(2 \zeta_{3}^{2}-\zeta_{3}-1) x y^{9}+\zeta_{3} x^{2} y^{7}-(\zeta_{3}-1) x y^{8}+2 \zeta_{3}^{2} x^{2} y^{6}-2 x y^{7}+x^{2} y^{5}-4 \zeta_{3}^{2} x y^{6}+\zeta_{3} y^{7}-2 \zeta_{3} x y^{5}+2 \zeta_{3}^{2} y^{6}+(\zeta_{3}-1) x y^{4}+y^{5}+(2 \zeta_{3}^{2}-\zeta_{3}-1) x y^{3}-(\zeta_{3}-1) x y^{2}+(2 \zeta_{3}^{2}-\zeta_{3}-1) x^{2}=0$.

None.

\subsubsection*{Case $\{\aaa\ddd\eee,\ccc^3,\bbb^4\}$, $(x={\mathrm e}^{\mathrm{i} \ddd+\mathrm{i} \eee},y={\mathrm e}^{\mathrm{i} \ddd-\mathrm{i} \eee})$}

$2 \zeta_{12} x^{2} y^{2}-(\zeta_{12}^{2}+3 \zeta_{12}+1) x^{2} y+3 \zeta_{12} x^{2}+(2 \zeta_{12}^{2}-4 \zeta_{12}+2) x y+3 \zeta_{12} y^{2}-(\zeta_{12}^{2}+3 \zeta_{12}+1) y+2 \zeta_{12}=0$.

None.

\subsubsection*{Case $\{\aaa\ddd\eee,\ccc^3,\bbb^3\ccc\}$, $(x={\mathrm e}^{\mathrm{i} \ddd},y={\mathrm e}^{\mathrm{i} \eee})$}

$\zeta_{7}^{4} x^{5} y^{2}-(\zeta_{7}^{4}+\zeta_{7}-1) x^{4} y^{3}-2 \zeta_{7}^{4} x^{4} y+(\zeta_{7}^{5}+2 \zeta_{7}^{4}-2) x^{3} y^{2}+(\zeta_{7}^{5}-2 \zeta_{7}^{4}+\zeta_{7}+1) x^{2} y^{3}-(2 \zeta_{7}^{5}-\zeta_{7}^{4}-\zeta_{7}^{2}-\zeta_{7}) x^{3}+(2 \zeta_{7}^{5}+\zeta_{7}^{4}-2 \zeta_{7}^{2}) x^{2} y-2 \zeta_{7}^{5} x y^{2}-(\zeta_{7}^{5}-\zeta_{7}^{2}+\zeta_{7}) x+\zeta_{7}^{5} y=0$.

None.

\subsubsection*{Case $\{\aaa\ddd\eee,\ccc^3,\bbb^5\}$, $(x={\mathrm e}^{\mathrm{i} \ddd+\mathrm{i} \eee},y={\mathrm e}^{\mathrm{i} \ddd-\mathrm{i} \eee})$}

$(\zeta_{15}^{6}-2 \zeta_{15}^{3}+1) x^{2} y^{2}-(\zeta_{15}^{6}-\zeta_{15}^{5}-3 \zeta_{15}^{3}-\zeta_{15}+1) x^{2} y-3 \zeta_{15}^{3} x^{2}-(2 \zeta_{15}^{5}-4 \zeta_{15}^{3}+2 \zeta_{15}) x y-3 \zeta_{15}^{3} y^{2}-(\zeta_{15}^{6}-\zeta_{15}^{5}-3 \zeta_{15}^{3}-\zeta_{15}+1) y+\zeta_{15}^{6}-2 \zeta_{15}^{3}+1=0$.

None.

\subsection*{Case $\{\aaa\ddd\eee,\bbb^2\ccc\}$, $(x={\mathrm e}^{\frac{4\mathrm{i} \pi}{f}}, y={\mathrm e}^{\mathrm{i} \ddd+\mathrm{i} \eee},z={\mathrm e}^{\mathrm{i} \ddd-\mathrm{i} \eee})$}

This is  the only case inducing a three variable polynomial equation.

$f(x,y,z)=x^{4} y^{2} z-2 x^{4} y z-x^{3} y^{2} z+x^{2} y^{2} z^{2}+x^{4} z-x^{3} y^{2}+4 x^{3} y z-x^{3} z^{2}-x^{3} z+2 x^{2} y^{2}-6 x^{2} y z+2 x^{2} z^{2}-x y^{2} z-x y^{2}+4 x y z-x z^{2}+y^{2} z+x^{2}-x z-2 y z+z=0$.

\subsubsection*{Subcase 1. resultant$(f(x,y,z),f(-x,y,z),z)$}

$(y-1) (x^{2}+1) (x^{4} y+x^{3} y^{2}+x^{2} y^{2}+x^{3}+x^{2} y-x y^{2}+x^{2}-x+y) (x^{4} y-x^{3} y^{2}+x^{2} y^{2}-x^{3}+x^{2} y+x y^{2}+x^{2}+x+y)=0$.

None.

\subsubsection*{Subcase 2. resultant$(f(x,y,z),f(x,-y,z),z)$}

$(x^{2}-x+1) (x y^{2}-x^{2}+2 x-1) (x^{2} y^{2}-2 x y^{2}+y^{2}-x)=0$.

None.

\subsubsection*{Subcase 3. resultant$(f(x,y,z),f(x,y,-z),z)$}

$(x y^{2}-x^{2}+2 x-1) (x^{2} y^{2}-2 x y^{2}+y^{2}-x) (x^{4} y^{2}-2 x^{4} y-x^{3} y^{2}+x^{4}+4 x^{3} y-x^{3}-6 x^{2} y-x y^{2}+4 x y+y^{2}-x-2 y+1) =0$.

None.

\subsubsection*{Subcase 4.  resultant$(f(\tilde{x},\tilde{y},z),f(-\tilde{x},-\tilde{y},z),z), (x=\tilde{x}\tilde{y}, y=\tfrac{\tilde{x}}{\tilde{y}})$}

$x^{8} y^{6}+2 x^{7} y^{7}+x^{6} y^{8}-2 x^{8} y^{5}-2 x^{7} y^{6}-2 x^{6} y^{7}-2 x^{5} y^{8}-x^{8} y^{4}-4 x^{7} y^{5}-2 x^{6} y^{6}-4 x^{5} y^{7}-x^{4} y^{8}-4 x^{7} y^{4}-8 x^{6} y^{5}-8 x^{5} y^{6}-4 x^{4} y^{7}+x^{7} y^{3}+8 x^{6} y^{4}+22 x^{5} y^{5}+8 x^{4} y^{6}+x^{3} y^{7}+10 x^{6} y^{3}+34 x^{5} y^{4}+34 x^{4} y^{5}+10 x^{3} y^{6}+4 x^{6} y^{2}+34 x^{5} y^{3}+72 x^{4} y^{4}+34 x^{3} y^{5}+4 x^{2} y^{6}+10 x^{5} y^{2}+34 x^{4} y^{3}+34 x^{3} y^{4}+10 x^{2} y^{5}+x^{5} y+8 x^{4} y^{2}+22 x^{3} y^{3}+8 x^{2} y^{4}+x y^{5}-4 x^{4} y-8 x^{3} y^{2}-8 x^{2} y^{3}-4 x y^{4}-x^{4}-4 x^{3} y-2 x^{2} y^{2}-4 x y^{3}-y^{4}-2 x^{3}-2 x^{2} y-2 x y^{2}-2 y^{3}+x^{2}+2 x y+y^{2}=0$.

None.

\subsubsection*{Subcase 5.  resultant$(f(\tilde{x},y,z),f(-\tilde{x},y,-z),z), (x=\tilde{x}^2)$}

$x^{4} y^{8}-2 x^{4} y^{7}-x^{3} y^{8}+4 x^{4} y^{6}-6 x^{3} y^{7}+4 x^{2} y^{8}-6 x^{4} y^{5}+8 x^{2} y^{7}-x y^{8}+10 x^{4} y^{4}-26 x^{3} y^{5}+36 x^{2} y^{6}-6 x y^{7}+y^{8}-6 x^{4} y^{3}+10 x^{3} y^{4}+16 x^{2} y^{5}-2 y^{7}+4 x^{4} y^{2}-26 x^{3} y^{3}+76 x^{2} y^{4}-26 x y^{5}+4 y^{6}-2 x^{4} y+16 x^{2} y^{3}+10 x y^{4}-6 y^{5}+x^{4}-6 x^{3} y+36 x^{2} y^{2}-26 x y^{3}+10 y^{4}-x^{3}+8 x^{2} y-6 y^{3}+4 x^{2}-6 x y+4 y^{2}-x-2 y+1=0$.

None.

\subsubsection*{Subcase 6.  resultant$(f(x,y,z),f(x,-y,-z),z)$}

$(x-1)  (y^{2}+1) (x^{2}+x+1) (x y^{2}-x^{2}+2 x-1) (x^{2} y^{2}-2 x y^{2}+y^{2}-x)=0$.

$\ddd+\eee=\frac{\pi}2$ or $\frac{3\pi}2$.
None.

\subsubsection*{Subcase 7.  resultant$(f(\tilde{x},\tilde{y},z),f(-\tilde{x},-\tilde{y},-z),z), (x=\tilde{x}\tilde{y}, y=\tfrac{\tilde{x}}{\tilde{y}})$}

$4 x^{7} y^{7}-2 x^{8} y^{5}-6 x^{7} y^{6}-6 x^{6} y^{7}-2 x^{5} y^{8}+x^{8} y^{4}-4 x^{7} y^{5}+6 x^{6} y^{6}-4 x^{5} y^{7}+x^{4} y^{8}+4 x^{7} y^{4}+20 x^{6} y^{5}+20 x^{5} y^{6}+4 x^{4} y^{7}-x^{7} y^{3}-16 x^{6} y^{4}+2 x^{5} y^{5}-16 x^{4} y^{6}-x^{3} y^{7}-6 x^{6} y^{3}-2 x^{5} y^{4}-2 x^{4} y^{5}-6 x^{3} y^{6}+4 x^{6} y^{2}+12 x^{5} y^{3}+56 x^{4} y^{4}+12 x^{3} y^{5}+4 x^{2} y^{6}-6 x^{5} y^{2}-2 x^{4} y^{3}-2 x^{3} y^{4}-6 x^{2} y^{5}-x^{5} y-16 x^{4} y^{2}+2 x^{3} y^{3}-16 x^{2} y^{4}-x y^{5}+4 x^{4} y+20 x^{3} y^{2}+20 x^{2} y^{3}+4 x y^{4}+x^{4}-4 x^{3} y+6 x^{2} y^{2}-4 x y^{3}+y^{4}-2 x^{3}-6 x^{2} y-6 x y^{2}-2 y^{3}+4 x y=0$.

None.

\subsubsection*{Subcase 8.  resultant$(f(x,y,z),f(x^2,y^2,z^2),z)$}

$(x-1) (y-1) (2 x^{18} y^{11}-x^{17} y^{12}-2 x^{19} y^{9}+x^{18} y^{10}-2 x^{17} y^{11}-4 x^{19} y^{8}+8 x^{18} y^{9}+4 x^{17} y^{10}-4 x^{16} y^{11}+x^{15} y^{12}+4 x^{20} y^{6}-6 x^{19} y^{7}+12 x^{18} y^{8}-6 x^{17} y^{9}-23 x^{16} y^{10}+10 x^{15} y^{11}-16 x^{19} y^{6}+22 x^{18} y^{7}-7 x^{17} y^{8}-18 x^{16} y^{9}+28 x^{15} y^{10}-2 x^{14} y^{11}-6 x^{19} y^{5}+34 x^{18} y^{6}-16 x^{17} y^{7}-24 x^{16} y^{8}+10 x^{15} y^{9}-3 x^{14} y^{10}-4 x^{19} y^{4}+22 x^{18} y^{5}-16 x^{17} y^{6}-26 x^{16} y^{7}-x^{15} y^{8}+36 x^{14} y^{9}-8 x^{13} y^{10}-2 x^{19} y^{3}+12 x^{18} y^{4}-16 x^{17} y^{5}-38 x^{16} y^{6}+12 x^{15} y^{7}+76 x^{14} y^{8}-34 x^{13} y^{9}+2 x^{12} y^{10}-8 x^{11} y^{11}+4 x^{10} y^{12}+8 x^{18} y^{3}-7 x^{17} y^{4}-26 x^{16} y^{5}+24 x^{15} y^{6}+78 x^{14} y^{7}-52 x^{13} y^{8}-28 x^{12} y^{9}-32 x^{11} y^{10}+24 x^{10} y^{11}+x^{18} y^{2}-6 x^{17} y^{3}-24 x^{16} y^{4}+12 x^{15} y^{5}+86 x^{14} y^{6}-46 x^{13} y^{7}-56 x^{12} y^{8}-12 x^{11} y^{9}+70 x^{10} y^{10}-8 x^{9} y^{11}+2 x^{18} y+4 x^{17} y^{2}-18 x^{16} y^{3}-x^{15} y^{4}+78 x^{14} y^{5}-56 x^{13} y^{6}-68 x^{12} y^{7}-8 x^{11} y^{8}+68 x^{10} y^{9}-32 x^{9} y^{10}-2 x^{17} y-23 x^{16} y^{2}+10 x^{15} y^{3}+76 x^{14} y^{4}-46 x^{13} y^{5}-32 x^{12} y^{6}-36 x^{11} y^{7}+108 x^{10} y^{8}-12 x^{9} y^{9}+2 x^{8} y^{10}-x^{17}-4 x^{16} y+28 x^{15} y^{2}+36 x^{14} y^{3}-52 x^{13} y^{4}-68 x^{12} y^{5}-120 x^{11} y^{6}+180 x^{10} y^{7}-8 x^{9} y^{8}-28 x^{8} y^{9}-8 x^{7} y^{10}-2 x^{6} y^{11}+x^{5} y^{12}+10 x^{15} y-3 x^{14} y^{2}-34 x^{13} y^{3}-56 x^{12} y^{4}-36 x^{11} y^{5}+312 x^{10} y^{6}-36 x^{9} y^{7}-56 x^{8} y^{8}-34 x^{7} y^{9}-3 x^{6} y^{10}+10 x^{5} y^{11}+x^{15}-2 x^{14} y-8 x^{13} y^{2}-28 x^{12} y^{3}-8 x^{11} y^{4}+180 x^{10} y^{5}-120 x^{9} y^{6}-68 x^{8} y^{7}-52 x^{7} y^{8}+36 x^{6} y^{9}+28 x^{5} y^{10}-4 x^{4} y^{11}-x^{3} y^{12}+2 x^{12} y^{2}-12 x^{11} y^{3}+108 x^{10} y^{4}-36 x^{9} y^{5}-32 x^{8} y^{6}-46 x^{7} y^{7}+76 x^{6} y^{8}+10 x^{5} y^{9}-23 x^{4} y^{10}-2 x^{3} y^{11}-32 x^{11} y^{2}+68 x^{10} y^{3}-8 x^{9} y^{4}-68 x^{8} y^{5}-56 x^{7} y^{6}+78 x^{6} y^{7}-x^{5} y^{8}-18 x^{4} y^{9}+4 x^{3} y^{10}+2 x^{2} y^{11}-8 x^{11} y+70 x^{10} y^{2}-12 x^{9} y^{3}-56 x^{8} y^{4}-46 x^{7} y^{5}+86 x^{6} y^{6}+12 x^{5} y^{7}-24 x^{4} y^{8}-6 x^{3} y^{9}+x^{2} y^{10}+24 x^{10} y-32 x^{9} y^{2}-28 x^{8} y^{3}-52 x^{7} y^{4}+78 x^{6} y^{5}+24 x^{5} y^{6}-26 x^{4} y^{7}-7 x^{3} y^{8}+8 x^{2} y^{9}+4 x^{10}-8 x^{9} y+2 x^{8} y^{2}-34 x^{7} y^{3}+76 x^{6} y^{4}+12 x^{5} y^{5}-38 x^{4} y^{6}-16 x^{3} y^{7}+12 x^{2} y^{8}-2 x y^{9}-8 x^{7} y^{2}+36 x^{6} y^{3}-x^{5} y^{4}-26 x^{4} y^{5}-16 x^{3} y^{6}+22 x^{2} y^{7}-4 x y^{8}-3 x^{6} y^{2}+10 x^{5} y^{3}-24 x^{4} y^{4}-16 x^{3} y^{5}+34 x^{2} y^{6}-6 x y^{7}-2 x^{6} y+28 x^{5} y^{2}-18 x^{4} y^{3}-7 x^{3} y^{4}+22 x^{2} y^{5}-16 x y^{6}+10 x^{5} y-23 x^{4} y^{2}-6 x^{3} y^{3}+12 x^{2} y^{4}-6 x y^{5}+4 y^{6}+x^{5}-4 x^{4} y+4 x^{3} y^{2}+8 x^{2} y^{3}-4 x y^{4}-2 x^{3} y+x^{2} y^{2}-2 x y^{3}-x^{3}+2 x^{2} y)=0$.

None.

\subsubsection*{Subcase 9.  resultant$(f(x,y,z),f(-x^2,y^2,z^2),z)$}

$(y-1) (x^{2}-x+1) (x^{14} y^{12}-x^{15} y^{10}+x^{16} y^{8}-2 x^{15} y^{9}+3 x^{14} y^{10}+4 x^{13} y^{11}-x^{12} y^{12}-4 x^{15} y^{8}+2 x^{14} y^{9}-8 x^{12} y^{11}-2 x^{11} y^{12}+2 x^{16} y^{6}-2 x^{15} y^{7}+5 x^{14} y^{8}+8 x^{13} y^{9}-x^{12} y^{10}-x^{10} y^{12}-6 x^{15} y^{6}-2 x^{14} y^{7}+8 x^{13} y^{8}-14 x^{12} y^{9}-28 x^{11} y^{10}-8 x^{10} y^{11}+2 x^{9} y^{12}+x^{16} y^{4}-2 x^{15} y^{5}+10 x^{14} y^{6}+24 x^{13} y^{7}-14 x^{12} y^{8}-16 x^{11} y^{9}+13 x^{10} y^{10}+12 x^{9} y^{11}+6 x^{8} y^{12}-4 x^{15} y^{4}-2 x^{14} y^{5}-26 x^{12} y^{7}-30 x^{11} y^{8}-26 x^{10} y^{9}-3 x^{9} y^{10}+16 x^{8} y^{11}+2 x^{7} y^{12}-2 x^{15} y^{3}+5 x^{14} y^{4}+24 x^{13} y^{5}+20 x^{12} y^{6}-24 x^{11} y^{7}-37 x^{10} y^{8}+66 x^{9} y^{9}+82 x^{8} y^{10}+12 x^{7} y^{11}-x^{6} y^{12}-x^{15} y^{2}+2 x^{14} y^{3}+8 x^{13} y^{4}-26 x^{12} y^{5}-104 x^{11} y^{6}-46 x^{10} y^{7}+106 x^{9} y^{8}+76 x^{8} y^{9}-3 x^{7} y^{10}-8 x^{6} y^{11}-2 x^{5} y^{12}+3 x^{14} y^{2}+8 x^{13} y^{3}-14 x^{12} y^{4}-24 x^{11} y^{5}+58 x^{10} y^{6}+130 x^{9} y^{7}+122 x^{8} y^{8}+66 x^{7} y^{9}+13 x^{6} y^{10}-x^{4} y^{12}-14 x^{12} y^{3}-30 x^{11} y^{4}-46 x^{10} y^{5}+14 x^{9} y^{6}+164 x^{8} y^{7}+106 x^{7} y^{8}-26 x^{6} y^{9}-28 x^{5} y^{10}-8 x^{4} y^{11}+x^{14}+4 x^{13} y-x^{12} y^{2}-16 x^{11} y^{3}-37 x^{10} y^{4}+130 x^{9} y^{5}+324 x^{8} y^{6}+130 x^{7} y^{7}-37 x^{6} y^{8}-16 x^{5} y^{9}-x^{4} y^{10}+4 x^{3} y^{11}+x^{2} y^{12}-8 x^{12} y-28 x^{11} y^{2}-26 x^{10} y^{3}+106 x^{9} y^{4}+164 x^{8} y^{5}+14 x^{7} y^{6}-46 x^{6} y^{7}-30 x^{5} y^{8}-14 x^{4} y^{9}-x^{12}+13 x^{10} y^{2}+66 x^{9} y^{3}+122 x^{8} y^{4}+130 x^{7} y^{5}+58 x^{6} y^{6}-24 x^{5} y^{7}-14 x^{4} y^{8}+8 x^{3} y^{9}+3 x^{2} y^{10}-2 x^{11}-8 x^{10} y-3 x^{9} y^{2}+76 x^{8} y^{3}+106 x^{7} y^{4}-46 x^{6} y^{5}-104 x^{5} y^{6}-26 x^{4} y^{7}+8 x^{3} y^{8}+2 x^{2} y^{9}-x y^{10}-x^{10}+12 x^{9} y+82 x^{8} y^{2}+66 x^{7} y^{3}-37 x^{6} y^{4}-24 x^{5} y^{5}+20 x^{4} y^{6}+24 x^{3} y^{7}+5 x^{2} y^{8}-2 x y^{9}+2 x^{9}+16 x^{8} y-3 x^{7} y^{2}-26 x^{6} y^{3}-30 x^{5} y^{4}-26 x^{4} y^{5}-2 x^{2} y^{7}-4 x y^{8}+6 x^{8}+12 x^{7} y+13 x^{6} y^{2}-16 x^{5} y^{3}-14 x^{4} y^{4}+24 x^{3} y^{5}+10 x^{2} y^{6}-2 x y^{7}+y^{8}+2 x^{7}-8 x^{6} y-28 x^{5} y^{2}-14 x^{4} y^{3}+8 x^{3} y^{4}-2 x^{2} y^{5}-6 x y^{6}-x^{6}-x^{4} y^{2}+8 x^{3} y^{3}+5 x^{2} y^{4}-2 x y^{5}+2 y^{6}-2 x^{5}-8 x^{4} y+2 x^{2} y^{3}-4 x y^{4}-x^{4}+4 x^{3} y+3 x^{2} y^{2}-2 x y^{3}+y^{4}-x y^{2}+x^{2})=0$.

$f=20,(\aaa,\bbb,\ccc,\ddd,\eee)=(2,4,2,5,3)/5$.

$f=20,(\aaa,\bbb,\ccc,\ddd,\eee)=(4,8,4,9,7)/10$.

$f=20,(\aaa,\bbb,\ccc,\ddd,\eee)=(16,8,4,1,3)/10$.

$f=20,(\aaa,\bbb,\ccc,\ddd,\eee)=(10,12,6,5,15)/15$.

\subsubsection*{Subcase 10.  resultant$(f(x,y,z),f(x^2,-y^2,z^2),z)$}

$2 x^{22} y^{15}-x^{21} y^{16}-2 x^{23} y^{13}-5 x^{22} y^{14}-6 x^{21} y^{15}+4 x^{20} y^{16}+2 x^{23} y^{12}+24 x^{22} y^{13}+22 x^{21} y^{14}-5 x^{19} y^{16}+4 x^{24} y^{10}-10 x^{23} y^{11}-22 x^{22} y^{12}-90 x^{21} y^{13}-51 x^{20} y^{14}+26 x^{19} y^{15}-8 x^{24} y^{9}-16 x^{23} y^{10}+68 x^{22} y^{11}+86 x^{21} y^{12}+194 x^{20} y^{13}+86 x^{19} y^{14}-56 x^{18} y^{15}+5 x^{17} y^{16}+12 x^{24} y^{8}+36 x^{23} y^{9}+17 x^{22} y^{10}-174 x^{21} y^{11}-238 x^{20} y^{12}-332 x^{19} y^{13}-84 x^{18} y^{14}+62 x^{17} y^{15}-4 x^{16} y^{16}-8 x^{24} y^{7}-68 x^{23} y^{8}-70 x^{22} y^{9}+18 x^{21} y^{10}+218 x^{20} y^{11}+538 x^{19} y^{12}+512 x^{18} y^{13}-14 x^{17} y^{14}-40 x^{16} y^{15}+x^{15} y^{16}+4 x^{24} y^{6}+36 x^{23} y^{7}+188 x^{22} y^{8}+126 x^{21} y^{9}-101 x^{20} y^{10}-184 x^{19} y^{11}-864 x^{18} y^{12}-620 x^{17} y^{13}+139 x^{16} y^{14}+6 x^{15} y^{15}+4 x^{14} y^{16}-16 x^{23} y^{6}-70 x^{22} y^{7}-346 x^{21} y^{8}-284 x^{20} y^{9}+194 x^{19} y^{10}+280 x^{18} y^{11}+894 x^{17} y^{12}+406 x^{16} y^{13}-134 x^{15} y^{14}+38 x^{14} y^{15}-16 x^{13} y^{16}-10 x^{23} y^{5}+17 x^{22} y^{6}+126 x^{21} y^{7}+408 x^{20} y^{8}+394 x^{19} y^{9}+24 x^{18} y^{10}-496 x^{17} y^{11}-586 x^{16} y^{12}+186 x^{15} y^{13}+29 x^{14} y^{14}-88 x^{13} y^{15}+24 x^{12} y^{16}+2 x^{23} y^{4}+68 x^{22} y^{5}+18 x^{21} y^{6}-284 x^{20} y^{7}-186 x^{19} y^{8}-88 x^{18} y^{9}-786 x^{17} y^{10}+350 x^{16} y^{11}+98 x^{15} y^{12}-816 x^{14} y^{13}+8 x^{13} y^{14}+112 x^{12} y^{15}-16 x^{11} y^{16}-2 x^{23} y^{3}-22 x^{22} y^{4}-174 x^{21} y^{5}-101 x^{20} y^{6}+394 x^{19} y^{7}-80 x^{18} y^{8}-362 x^{17} y^{9}+1269 x^{16} y^{10}+366 x^{15} y^{11}+846 x^{14} y^{12}+1146 x^{13} y^{13}+8 x^{12} y^{14}-88 x^{11} y^{15}+4 x^{10} y^{16}+24 x^{22} y^{3}+86 x^{21} y^{4}+218 x^{20} y^{5}+194 x^{19} y^{6}-88 x^{18} y^{7}-102 x^{17} y^{8}-140 x^{16} y^{9}-618 x^{15} y^{10}-1020 x^{14} y^{11}-2282 x^{13} y^{12}-1216 x^{12} y^{13}+8 x^{11} y^{14}+38 x^{10} y^{15}+x^{9} y^{16}-5 x^{22} y^{2}-90 x^{21} y^{3}-238 x^{20} y^{4}-184 x^{19} y^{5}+24 x^{18} y^{6}-362 x^{17} y^{7}+452 x^{16} y^{8}+1914 x^{15} y^{9}-x^{14} y^{10}+1074 x^{13} y^{11}+3056 x^{12} y^{12}+1146 x^{11} y^{13}+29 x^{10} y^{14}+6 x^{9} y^{15}-4 x^{8} y^{16}+2 x^{22} y+22 x^{21} y^{2}+194 x^{20} y^{3}+538 x^{19} y^{4}+280 x^{18} y^{5}-786 x^{17} y^{6}-140 x^{16} y^{7}-902 x^{15} y^{8}-3586 x^{14} y^{9}-768 x^{13} y^{10}-944 x^{12} y^{11}-2282 x^{11} y^{12}-816 x^{10} y^{13}-134 x^{9} y^{14}-40 x^{8} y^{15}+5 x^{7} y^{16}-6 x^{21} y-51 x^{20} y^{2}-332 x^{19} y^{3}-864 x^{18} y^{4}-496 x^{17} y^{5}+1269 x^{16} y^{6}+1914 x^{15} y^{7}+3180 x^{14} y^{8}+4012 x^{13} y^{9}+1528 x^{12} y^{10}+1074 x^{11} y^{11}+846 x^{10} y^{12}+186 x^{9} y^{13}+139 x^{8} y^{14}+62 x^{7} y^{15}-x^{21}+86 x^{19} y^{2}+512 x^{18} y^{3}+894 x^{17} y^{4}+350 x^{16} y^{5}-618 x^{15} y^{6}-3586 x^{14} y^{7}-7612 x^{13} y^{8}-3888 x^{12} y^{9}-768 x^{11} y^{10}-1020 x^{10} y^{11}+98 x^{9} y^{12}+406 x^{8} y^{13}-14 x^{7} y^{14}-56 x^{6} y^{15}-5 x^{5} y^{16}+4 x^{20}+26 x^{19} y-84 x^{18} y^{2}-620 x^{17} y^{3}-586 x^{16} y^{4}+366 x^{15} y^{5}-x^{14} y^{6}+4012 x^{13} y^{7}+10116 x^{12} y^{8}+4012 x^{11} y^{9}-x^{10} y^{10}+366 x^{9} y^{11}-586 x^{8} y^{12}-620 x^{7} y^{13}-84 x^{6} y^{14}+26 x^{5} y^{15}+4 x^{4} y^{16}-5 x^{19}-56 x^{18} y-14 x^{17} y^{2}+406 x^{16} y^{3}+98 x^{15} y^{4}-1020 x^{14} y^{5}-768 x^{13} y^{6}-3888 x^{12} y^{7}-7612 x^{11} y^{8}-3586 x^{10} y^{9}-618 x^{9} y^{10}+350 x^{8} y^{11}+894 x^{7} y^{12}+512 x^{6} y^{13}+86 x^{5} y^{14}-x^{3} y^{16}+62 x^{17} y+139 x^{16} y^{2}+186 x^{15} y^{3}+846 x^{14} y^{4}+1074 x^{13} y^{5}+1528 x^{12} y^{6}+4012 x^{11} y^{7}+3180 x^{10} y^{8}+1914 x^{9} y^{9}+1269 x^{8} y^{10}-496 x^{7} y^{11}-864 x^{6} y^{12}-332 x^{5} y^{13}-51 x^{4} y^{14}-6 x^{3} y^{15}+5 x^{17}-40 x^{16} y-134 x^{15} y^{2}-816 x^{14} y^{3}-2282 x^{13} y^{4}-944 x^{12} y^{5}-768 x^{11} y^{6}-3586 x^{10} y^{7}-902 x^{9} y^{8}-140 x^{8} y^{9}-786 x^{7} y^{10}+280 x^{6} y^{11}+538 x^{5} y^{12}+194 x^{4} y^{13}+22 x^{3} y^{14}+2 x^{2} y^{15}-4 x^{16}+6 x^{15} y+29 x^{14} y^{2}+1146 x^{13} y^{3}+3056 x^{12} y^{4}+1074 x^{11} y^{5}-x^{10} y^{6}+1914 x^{9} y^{7}+452 x^{8} y^{8}-362 x^{7} y^{9}+24 x^{6} y^{10}-184 x^{5} y^{11}-238 x^{4} y^{12}-90 x^{3} y^{13}-5 x^{2} y^{14}+x^{15}+38 x^{14} y+8 x^{13} y^{2}-1216 x^{12} y^{3}-2282 x^{11} y^{4}-1020 x^{10} y^{5}-618 x^{9} y^{6}-140 x^{8} y^{7}-102 x^{7} y^{8}-88 x^{6} y^{9}+194 x^{5} y^{10}+218 x^{4} y^{11}+86 x^{3} y^{12}+24 x^{2} y^{13}+4 x^{14}-88 x^{13} y+8 x^{12} y^{2}+1146 x^{11} y^{3}+846 x^{10} y^{4}+366 x^{9} y^{5}+1269 x^{8} y^{6}-362 x^{7} y^{7}-80 x^{6} y^{8}+394 x^{5} y^{9}-101 x^{4} y^{10}-174 x^{3} y^{11}-22 x^{2} y^{12}-2 x y^{13}-16 x^{13}+112 x^{12} y+8 x^{11} y^{2}-816 x^{10} y^{3}+98 x^{9} y^{4}+350 x^{8} y^{5}-786 x^{7} y^{6}-88 x^{6} y^{7}-186 x^{5} y^{8}-284 x^{4} y^{9}+18 x^{3} y^{10}+68 x^{2} y^{11}+2 x y^{12}+24 x^{12}-88 x^{11} y+29 x^{10} y^{2}+186 x^{9} y^{3}-586 x^{8} y^{4}-496 x^{7} y^{5}+24 x^{6} y^{6}+394 x^{5} y^{7}+408 x^{4} y^{8}+126 x^{3} y^{9}+17 x^{2} y^{10}-10 x y^{11}-16 x^{11}+38 x^{10} y-134 x^{9} y^{2}+406 x^{8} y^{3}+894 x^{7} y^{4}+280 x^{6} y^{5}+194 x^{5} y^{6}-284 x^{4} y^{7}-346 x^{3} y^{8}-70 x^{2} y^{9}-16 x y^{10}+4 x^{10}+6 x^{9} y+139 x^{8} y^{2}-620 x^{7} y^{3}-864 x^{6} y^{4}-184 x^{5} y^{5}-101 x^{4} y^{6}+126 x^{3} y^{7}+188 x^{2} y^{8}+36 x y^{9}+4 y^{10}+x^{9}-40 x^{8} y-14 x^{7} y^{2}+512 x^{6} y^{3}+538 x^{5} y^{4}+218 x^{4} y^{5}+18 x^{3} y^{6}-70 x^{2} y^{7}-68 x y^{8}-8 y^{9}-4 x^{8}+62 x^{7} y-84 x^{6} y^{2}-332 x^{5} y^{3}-238 x^{4} y^{4}-174 x^{3} y^{5}+17 x^{2} y^{6}+36 x y^{7}+12 y^{8}+5 x^{7}-56 x^{6} y+86 x^{5} y^{2}+194 x^{4} y^{3}+86 x^{3} y^{4}+68 x^{2} y^{5}-16 x y^{6}-8 y^{7}+26 x^{5} y-51 x^{4} y^{2}-90 x^{3} y^{3}-22 x^{2} y^{4}-10 x y^{5}+4 y^{6}-5 x^{5}+22 x^{3} y^{2}+24 x^{2} y^{3}+2 x y^{4}+4 x^{4}-6 x^{3} y-5 x^{2} y^{2}-2 x y^{3}-x^{3}+2 x^{2} y=0$.

None.

\subsubsection*{Subcase 11.  resultant$(f(x,y,z),f(x^2,y^2,-z^2),z)$}

$x^{23} y^{14}-2 x^{22} y^{15}-x^{21} y^{16}-x^{24} y^{12}-2 x^{23} y^{13}+3 x^{22} y^{14}+14 x^{21} y^{15}+4 x^{20} y^{16}+4 x^{24} y^{11}+6 x^{23} y^{12}-4 x^{22} y^{13}-40 x^{21} y^{14}-40 x^{20} y^{15}-5 x^{19} y^{16}-8 x^{24} y^{10}-22 x^{23} y^{11}+2 x^{22} y^{12}+66 x^{21} y^{13}+133 x^{20} y^{14}+62 x^{19} y^{15}+12 x^{24} y^{9}+47 x^{23} y^{10}+32 x^{22} y^{11}-78 x^{21} y^{12}-238 x^{20} y^{13}-249 x^{19} y^{14}-56 x^{18} y^{15}+5 x^{17} y^{16}-14 x^{24} y^{8}-72 x^{23} y^{9}-107 x^{22} y^{10}+42 x^{21} y^{11}+266 x^{20} y^{12}+504 x^{19} y^{13}+316 x^{18} y^{14}+26 x^{17} y^{15}-4 x^{16} y^{16}+12 x^{24} y^{7}+84 x^{23} y^{8}+190 x^{22} y^{9}+104 x^{21} y^{10}-238 x^{20} y^{11}-578 x^{19} y^{12}-752 x^{18} y^{13}-271 x^{17} y^{14}+x^{15} y^{16}-8 x^{24} y^{6}-72 x^{23} y^{7}-228 x^{22} y^{8}-290 x^{21} y^{9}+35 x^{20} y^{10}+504 x^{19} y^{11}+996 x^{18} y^{12}+800 x^{17} y^{13}+91 x^{16} y^{14}-14 x^{15} y^{15}-4 x^{14} y^{16}+4 x^{24} y^{5}+47 x^{23} y^{6}+190 x^{22} y^{7}+366 x^{21} y^{8}+260 x^{20} y^{9}-287 x^{19} y^{10}-768 x^{18} y^{11}-1302 x^{17} y^{12}-490 x^{16} y^{13}+152 x^{15} y^{14}+42 x^{14} y^{15}+16 x^{13} y^{16}-x^{24} y^{4}-22 x^{23} y^{5}-107 x^{22} y^{6}-290 x^{21} y^{7}-368 x^{20} y^{8}+10 x^{19} y^{9}+516 x^{18} y^{10}+880 x^{17} y^{11}+1070 x^{16} y^{12}-162 x^{15} y^{13}-339 x^{14} y^{14}-88 x^{13} y^{15}-24 x^{12} y^{16}+6 x^{23} y^{4}+32 x^{22} y^{5}+104 x^{21} y^{6}+260 x^{20} y^{7}+126 x^{19} y^{8}-528 x^{18} y^{9}-465 x^{17} y^{10}-570 x^{16} y^{11}-106 x^{15} y^{12}+924 x^{14} y^{13}+439 x^{13} y^{14}+112 x^{12} y^{15}+16 x^{11} y^{16}-2 x^{23} y^{3}+2 x^{22} y^{4}+42 x^{21} y^{5}+35 x^{20} y^{6}+10 x^{19} y^{7}+344 x^{18} y^{8}+838 x^{17} y^{9}+45 x^{16} y^{10}-234 x^{15} y^{11}-1338 x^{14} y^{12}-1494 x^{13} y^{13}-472 x^{12} y^{14}-88 x^{11} y^{15}-4 x^{10} y^{16}+x^{23} y^{2}-4 x^{22} y^{3}-78 x^{21} y^{4}-238 x^{20} y^{5}-287 x^{19} y^{6}-528 x^{18} y^{7}-734 x^{17} y^{8}-124 x^{16} y^{9}+496 x^{15} y^{10}+1160 x^{14} y^{11}+2722 x^{13} y^{12}+1696 x^{12} y^{13}+439 x^{11} y^{14}+42 x^{10} y^{15}+x^{9} y^{16}+3 x^{22} y^{2}+66 x^{21} y^{3}+266 x^{20} y^{4}+504 x^{19} y^{5}+516 x^{18} y^{6}+838 x^{17} y^{7}+532 x^{16} y^{8}-1750 x^{15} y^{9}-1157 x^{14} y^{10}-1746 x^{13} y^{11}-3318 x^{12} y^{12}-1494 x^{11} y^{13}-339 x^{10} y^{14}-14 x^{9} y^{15}-4 x^{8} y^{16}-2 x^{22} y-40 x^{21} y^{2}-238 x^{20} y^{3}-578 x^{19} y^{4}-768 x^{18} y^{5}-465 x^{17} y^{6}-124 x^{16} y^{7}+882 x^{15} y^{8}+3754 x^{14} y^{9}+2081 x^{13} y^{10}+1912 x^{12} y^{11}+2722 x^{11} y^{12}+924 x^{10} y^{13}+152 x^{9} y^{14}+5 x^{7} y^{16}+14 x^{21} y+133 x^{20} y^{2}+504 x^{19} y^{3}+996 x^{18} y^{4}+880 x^{17} y^{5}+45 x^{16} y^{6}-1750 x^{15} y^{7}-4164 x^{14} y^{8}-4856 x^{13} y^{9}-2600 x^{12} y^{10}-1746 x^{11} y^{11}-1338 x^{10} y^{12}-162 x^{9} y^{13}+91 x^{8} y^{14}+26 x^{7} y^{15}-x^{21}-40 x^{20} y-249 x^{19} y^{2}-752 x^{18} y^{3}-1302 x^{17} y^{4}-570 x^{16} y^{5}+496 x^{15} y^{6}+3754 x^{14} y^{7}+8492 x^{13} y^{8}+5112 x^{12} y^{9}+2081 x^{11} y^{10}+1160 x^{10} y^{11}-106 x^{9} y^{12}-490 x^{8} y^{13}-271 x^{7} y^{14}-56 x^{6} y^{15}-5 x^{5} y^{16}+4 x^{20}+62 x^{19} y+316 x^{18} y^{2}+800 x^{17} y^{3}+1070 x^{16} y^{4}-234 x^{15} y^{5}-1157 x^{14} y^{6}-4856 x^{13} y^{7}-10640 x^{12} y^{8}-4856 x^{11} y^{9}-1157 x^{10} y^{10}-234 x^{9} y^{11}+1070 x^{8} y^{12}+800 x^{7} y^{13}+316 x^{6} y^{14}+62 x^{5} y^{15}+4 x^{4} y^{16}-5 x^{19}-56 x^{18} y-271 x^{17} y^{2}-490 x^{16} y^{3}-106 x^{15} y^{4}+1160 x^{14} y^{5}+2081 x^{13} y^{6}+5112 x^{12} y^{7}+8492 x^{11} y^{8}+3754 x^{10} y^{9}+496 x^{9} y^{10}-570 x^{8} y^{11}-1302 x^{7} y^{12}-752 x^{6} y^{13}-249 x^{5} y^{14}-40 x^{4} y^{15}-x^{3} y^{16}+26 x^{17} y+91 x^{16} y^{2}-162 x^{15} y^{3}-1338 x^{14} y^{4}-1746 x^{13} y^{5}-2600 x^{12} y^{6}-4856 x^{11} y^{7}-4164 x^{10} y^{8}-1750 x^{9} y^{9}+45 x^{8} y^{10}+880 x^{7} y^{11}+996 x^{6} y^{12}+504 x^{5} y^{13}+133 x^{4} y^{14}+14 x^{3} y^{15}+5 x^{17}+152 x^{15} y^{2}+924 x^{14} y^{3}+2722 x^{13} y^{4}+1912 x^{12} y^{5}+2081 x^{11} y^{6}+3754 x^{10} y^{7}+882 x^{9} y^{8}-124 x^{8} y^{9}-465 x^{7} y^{10}-768 x^{6} y^{11}-578 x^{5} y^{12}-238 x^{4} y^{13}-40 x^{3} y^{14}-2 x^{2} y^{15}-4 x^{16}-14 x^{15} y-339 x^{14} y^{2}-1494 x^{13} y^{3}-3318 x^{12} y^{4}-1746 x^{11} y^{5}-1157 x^{10} y^{6}-1750 x^{9} y^{7}+532 x^{8} y^{8}+838 x^{7} y^{9}+516 x^{6} y^{10}+504 x^{5} y^{11}+266 x^{4} y^{12}+66 x^{3} y^{13}+3 x^{2} y^{14}+x^{15}+42 x^{14} y+439 x^{13} y^{2}+1696 x^{12} y^{3}+2722 x^{11} y^{4}+1160 x^{10} y^{5}+496 x^{9} y^{6}-124 x^{8} y^{7}-734 x^{7} y^{8}-528 x^{6} y^{9}-287 x^{5} y^{10}-238 x^{4} y^{11}-78 x^{3} y^{12}-4 x^{2} y^{13}+x y^{14}-4 x^{14}-88 x^{13} y-472 x^{12} y^{2}-1494 x^{11} y^{3}-1338 x^{10} y^{4}-234 x^{9} y^{5}+45 x^{8} y^{6}+838 x^{7} y^{7}+344 x^{6} y^{8}+10 x^{5} y^{9}+35 x^{4} y^{10}+42 x^{3} y^{11}+2 x^{2} y^{12}-2 x y^{13}+16 x^{13}+112 x^{12} y+439 x^{11} y^{2}+924 x^{10} y^{3}-106 x^{9} y^{4}-570 x^{8} y^{5}-465 x^{7} y^{6}-528 x^{6} y^{7}+126 x^{5} y^{8}+260 x^{4} y^{9}+104 x^{3} y^{10}+32 x^{2} y^{11}+6 x y^{12}-24 x^{12}-88 x^{11} y-339 x^{10} y^{2}-162 x^{9} y^{3}+1070 x^{8} y^{4}+880 x^{7} y^{5}+516 x^{6} y^{6}+10 x^{5} y^{7}-368 x^{4} y^{8}-290 x^{3} y^{9}-107 x^{2} y^{10}-22 x y^{11}-y^{12}+16 x^{11}+42 x^{10} y+152 x^{9} y^{2}-490 x^{8} y^{3}-1302 x^{7} y^{4}-768 x^{6} y^{5}-287 x^{5} y^{6}+260 x^{4} y^{7}+366 x^{3} y^{8}+190 x^{2} y^{9}+47 x y^{10}+4 y^{11}-4 x^{10}-14 x^{9} y+91 x^{8} y^{2}+800 x^{7} y^{3}+996 x^{6} y^{4}+504 x^{5} y^{5}+35 x^{4} y^{6}-290 x^{3} y^{7}-228 x^{2} y^{8}-72 x y^{9}-8 y^{10}+x^{9}-271 x^{7} y^{2}-752 x^{6} y^{3}-578 x^{5} y^{4}-238 x^{4} y^{5}+104 x^{3} y^{6}+190 x^{2} y^{7}+84 x y^{8}+12 y^{9}-4 x^{8}+26 x^{7} y+316 x^{6} y^{2}+504 x^{5} y^{3}+266 x^{4} y^{4}+42 x^{3} y^{5}-107 x^{2} y^{6}-72 x y^{7}-14 y^{8}+5 x^{7}-56 x^{6} y-249 x^{5} y^{2}-238 x^{4} y^{3}-78 x^{3} y^{4}+32 x^{2} y^{5}+47 x y^{6}+12 y^{7}+62 x^{5} y+133 x^{4} y^{2}+66 x^{3} y^{3}+2 x^{2} y^{4}-22 x y^{5}-8 y^{6}-5 x^{5}-40 x^{4} y-40 x^{3} y^{2}-4 x^{2} y^{3}+6 x y^{4}+4 y^{5}+4 x^{4}+14 x^{3} y+3 x^{2} y^{2}-2 x y^{3}-y^{4}-x^{3}-2 x^{2} y+x y^{2}=0$.

None.

\subsubsection*{Subcase 12.   resultant$(f(x,y,z),f(-x^2,-y^2,z^2),z)$}

$(x^{2}-x+1)(x^{14} y^{16}-x^{15} y^{14}-4 x^{14} y^{15}+x^{16} y^{12}+2 x^{15} y^{13}+13 x^{14} y^{14}+4 x^{13} y^{15}-x^{12} y^{16}-4 x^{16} y^{11}-8 x^{15} y^{12}-22 x^{14} y^{13}-16 x^{13} y^{14}-4 x^{12} y^{15}-2 x^{11} y^{16}+12 x^{16} y^{10}+14 x^{15} y^{11}+36 x^{14} y^{12}+40 x^{13} y^{13}+25 x^{12} y^{14}+8 x^{11} y^{15}-x^{10} y^{16}-20 x^{16} y^{9}-39 x^{15} y^{10}-38 x^{14} y^{11}-62 x^{13} y^{12}-78 x^{12} y^{13}-48 x^{11} y^{14}-4 x^{10} y^{15}+2 x^{9} y^{16}+26 x^{16} y^{8}+56 x^{15} y^{9}+83 x^{14} y^{10}+76 x^{13} y^{11}+123 x^{12} y^{12}+120 x^{11} y^{13}+39 x^{10} y^{14}+4 x^{9} y^{15}+6 x^{8} y^{16}-20 x^{16} y^{7}-80 x^{15} y^{8}-120 x^{14} y^{9}-144 x^{13} y^{10}-130 x^{12} y^{11}-208 x^{11} y^{12}-162 x^{10} y^{13}-39 x^{9} y^{14}-8 x^{8} y^{15}+2 x^{7} y^{16}+12 x^{16} y^{6}+56 x^{15} y^{7}+182 x^{14} y^{8}+248 x^{13} y^{9}+219 x^{12} y^{10}+160 x^{11} y^{11}+272 x^{10} y^{12}+198 x^{9} y^{13}+62 x^{8} y^{14}+4 x^{7} y^{15}-x^{6} y^{16}-4 x^{16} y^{5}-39 x^{15} y^{6}-120 x^{14} y^{7}-356 x^{13} y^{8}-404 x^{12} y^{9}-296 x^{11} y^{10}-234 x^{10} y^{11}-322 x^{9} y^{12}-196 x^{8} y^{13}-39 x^{7} y^{14}-4 x^{6} y^{15}-2 x^{5} y^{16}+x^{16} y^{4}+14 x^{15} y^{5}+83 x^{14} y^{6}+248 x^{13} y^{7}+604 x^{12} y^{8}+536 x^{11} y^{9}+337 x^{10} y^{10}+278 x^{9} y^{11}+380 x^{8} y^{12}+198 x^{7} y^{13}+39 x^{6} y^{14}+8 x^{5} y^{15}-x^{4} y^{16}-8 x^{15} y^{4}-38 x^{14} y^{5}-144 x^{13} y^{6}-404 x^{12} y^{7}-892 x^{11} y^{8}-736 x^{10} y^{9}-353 x^{9} y^{10}-244 x^{8} y^{11}-322 x^{7} y^{12}-162 x^{6} y^{13}-48 x^{5} y^{14}-4 x^{4} y^{15}+2 x^{15} y^{3}+36 x^{14} y^{4}+76 x^{13} y^{5}+219 x^{12} y^{6}+536 x^{11} y^{7}+1170 x^{10} y^{8}+856 x^{9} y^{9}+410 x^{8} y^{10}+278 x^{7} y^{11}+272 x^{6} y^{12}+120 x^{5} y^{13}+25 x^{4} y^{14}+4 x^{3} y^{15}+x^{2} y^{16}-x^{15} y^{2}-22 x^{14} y^{3}-62 x^{13} y^{4}-130 x^{12} y^{5}-296 x^{11} y^{6}-736 x^{10} y^{7}-1344 x^{9} y^{8}-832 x^{8} y^{9}-353 x^{7} y^{10}-234 x^{6} y^{11}-208 x^{5} y^{12}-78 x^{4} y^{13}-16 x^{3} y^{14}-4 x^{2} y^{15}+13 x^{14} y^{2}+40 x^{13} y^{3}+123 x^{12} y^{4}+160 x^{11} y^{5}+337 x^{10} y^{6}+856 x^{9} y^{7}+1496 x^{8} y^{8}+856 x^{7} y^{9}+337 x^{6} y^{10}+160 x^{5} y^{11}+123 x^{4} y^{12}+40 x^{3} y^{13}+13 x^{2} y^{14}-4 x^{14} y-16 x^{13} y^{2}-78 x^{12} y^{3}-208 x^{11} y^{4}-234 x^{10} y^{5}-353 x^{9} y^{6}-832 x^{8} y^{7}-1344 x^{7} y^{8}-736 x^{6} y^{9}-296 x^{5} y^{10}-130 x^{4} y^{11}-62 x^{3} y^{12}-22 x^{2} y^{13}-x y^{14}+x^{14}+4 x^{13} y+25 x^{12} y^{2}+120 x^{11} y^{3}+272 x^{10} y^{4}+278 x^{9} y^{5}+410 x^{8} y^{6}+856 x^{7} y^{7}+1170 x^{6} y^{8}+536 x^{5} y^{9}+219 x^{4} y^{10}+76 x^{3} y^{11}+36 x^{2} y^{12}+2 x y^{13}-4 x^{12} y-48 x^{11} y^{2}-162 x^{10} y^{3}-322 x^{9} y^{4}-244 x^{8} y^{5}-353 x^{7} y^{6}-736 x^{6} y^{7}-892 x^{5} y^{8}-404 x^{4} y^{9}-144 x^{3} y^{10}-38 x^{2} y^{11}-8 x y^{12}-x^{12}+8 x^{11} y+39 x^{10} y^{2}+198 x^{9} y^{3}+380 x^{8} y^{4}+278 x^{7} y^{5}+337 x^{6} y^{6}+536 x^{5} y^{7}+604 x^{4} y^{8}+248 x^{3} y^{9}+83 x^{2} y^{10}+14 x y^{11}+y^{12}-2 x^{11}-4 x^{10} y-39 x^{9} y^{2}-196 x^{8} y^{3}-322 x^{7} y^{4}-234 x^{6} y^{5}-296 x^{5} y^{6}-404 x^{4} y^{7}-356 x^{3} y^{8}-120 x^{2} y^{9}-39 x y^{10}-4 y^{11}-x^{10}+4 x^{9} y+62 x^{8} y^{2}+198 x^{7} y^{3}+272 x^{6} y^{4}+160 x^{5} y^{5}+219 x^{4} y^{6}+248 x^{3} y^{7}+182 x^{2} y^{8}+56 x y^{9}+12 y^{10}+2 x^{9}-8 x^{8} y-39 x^{7} y^{2}-162 x^{6} y^{3}-208 x^{5} y^{4}-130 x^{4} y^{5}-144 x^{3} y^{6}-120 x^{2} y^{7}-80 x y^{8}-20 y^{9}+6 x^{8}+4 x^{7} y+39 x^{6} y^{2}+120 x^{5} y^{3}+123 x^{4} y^{4}+76 x^{3} y^{5}+83 x^{2} y^{6}+56 x y^{7}+26 y^{8}+2 x^{7}-4 x^{6} y-48 x^{5} y^{2}-78 x^{4} y^{3}-62 x^{3} y^{4}-38 x^{2} y^{5}-39 x y^{6}-20 y^{7}-x^{6}+8 x^{5} y+25 x^{4} y^{2}+40 x^{3} y^{3}+36 x^{2} y^{4}+14 x y^{5}+12 y^{6}-2 x^{5}-4 x^{4} y-16 x^{3} y^{2}-22 x^{2} y^{3}-8 x y^{4}-4 y^{5}-x^{4}+4 x^{3} y+13 x^{2} y^{2}+2 x y^{3}+y^{4}-4 x^{2} y-x y^{2}+x^{2})=0$.

None.

\subsubsection*{Subcase 13.  resultant$(f(x,y,z),f(-x^2,y^2,-z^2),z)$}

$2 x^{23} y^{13}-5 x^{22} y^{14}+x^{20} y^{16}-4 x^{23} y^{12}+6 x^{22} y^{13}+28 x^{21} y^{14}-4 x^{20} y^{15}-2 x^{19} y^{16}-4 x^{24} y^{10}+2 x^{23} y^{11}-2 x^{22} y^{12}-64 x^{21} y^{13}-87 x^{20} y^{14}+12 x^{19} y^{15}+2 x^{18} y^{16}+16 x^{24} y^{9}+24 x^{23} y^{10}+6 x^{22} y^{11}+70 x^{21} y^{12}+270 x^{20} y^{13}+210 x^{19} y^{14}-24 x^{18} y^{15}-2 x^{17} y^{16}-24 x^{24} y^{8}-100 x^{23} y^{9}-95 x^{22} y^{10}-72 x^{21} y^{11}-346 x^{20} y^{12}-782 x^{19} y^{13}-420 x^{18} y^{14}+36 x^{17} y^{15}+x^{16} y^{16}+16 x^{24} y^{7}+152 x^{23} y^{8}+396 x^{22} y^{9}+300 x^{21} y^{10}+310 x^{20} y^{11}+1072 x^{19} y^{12}+1792 x^{18} y^{13}+730 x^{17} y^{14}-44 x^{16} y^{15}-4 x^{24} y^{6}-100 x^{23} y^{7}-612 x^{22} y^{8}-1184 x^{21} y^{9}-769 x^{20} y^{10}-910 x^{19} y^{11}-2580 x^{18} y^{12}-3426 x^{17} y^{13}-1133 x^{16} y^{14}+48 x^{15} y^{15}-2 x^{14} y^{16}+24 x^{23} y^{6}+396 x^{22} y^{7}+1844 x^{21} y^{8}+2912 x^{20} y^{9}+1710 x^{19} y^{10}+2096 x^{18} y^{11}+5120 x^{17} y^{12}+5658 x^{16} y^{13}+1596 x^{15} y^{14}-48 x^{14} y^{15}+4 x^{13} y^{16}+2 x^{23} y^{5}-95 x^{22} y^{6}-1184 x^{21} y^{7}-4578 x^{20} y^{8}-6168 x^{19} y^{9}-3312 x^{18} y^{10}-4010 x^{17} y^{11}-8754 x^{16} y^{12}-8216 x^{15} y^{13}-2031 x^{14} y^{14}+48 x^{13} y^{15}-8 x^{12} y^{16}-4 x^{23} y^{4}+6 x^{22} y^{5}+300 x^{21} y^{6}+2912 x^{20} y^{7}+9764 x^{19} y^{8}+11504 x^{18} y^{9}+5694 x^{17} y^{10}+6626 x^{16} y^{11}+13106 x^{15} y^{12}+10626 x^{14} y^{13}+2356 x^{13} y^{14}-48 x^{12} y^{15}+4 x^{11} y^{16}+2 x^{23} y^{3}-2 x^{22} y^{4}-72 x^{21} y^{5}-769 x^{20} y^{6}-6168 x^{19} y^{7}-18372 x^{18} y^{8}-19144 x^{17} y^{9}-8763 x^{16} y^{10}-9624 x^{15} y^{11}-17382 x^{14} y^{12}-12354 x^{13} y^{13}-2480 x^{12} y^{14}+48 x^{11} y^{15}-2 x^{10} y^{16}+6 x^{22} y^{3}+70 x^{21} y^{4}+310 x^{20} y^{5}+1710 x^{19} y^{6}+11504 x^{18} y^{7}+30868 x^{17} y^{8}+28720 x^{16} y^{9}+12212 x^{15} y^{10}+12450 x^{14} y^{11}+20548 x^{13} y^{12}+12976 x^{12} y^{13}+2356 x^{11} y^{14}-48 x^{10} y^{15}-5 x^{22} y^{2}-64 x^{21} y^{3}-346 x^{20} y^{4}-910 x^{19} y^{5}-3312 x^{18} y^{6}-19144 x^{17} y^{7}-46838 x^{16} y^{8}-39088 x^{15} y^{9}-15413 x^{14} y^{10}-14474 x^{13} y^{11}-21740 x^{12} y^{12}-12354 x^{11} y^{13}-2031 x^{10} y^{14}+48 x^{9} y^{15}+x^{8} y^{16}+28 x^{21} y^{2}+270 x^{20} y^{3}+1072 x^{19} y^{4}+2096 x^{18} y^{5}+5694 x^{17} y^{6}+28720 x^{16} y^{7}+64492 x^{15} y^{8}+48516 x^{14} y^{9}+17740 x^{13} y^{10}+15200 x^{12} y^{11}+20548 x^{11} y^{12}+10626 x^{10} y^{13}+1596 x^{9} y^{14}-44 x^{8} y^{15}-2 x^{7} y^{16}-87 x^{20} y^{2}-782 x^{19} y^{3}-2580 x^{18} y^{4}-4010 x^{17} y^{5}-8763 x^{16} y^{6}-39088 x^{15} y^{7}-80920 x^{14} y^{8}-55132 x^{13} y^{9}-18600 x^{12} y^{10}-14474 x^{11} y^{11}-17382 x^{10} y^{12}-8216 x^{9} y^{13}-1133 x^{8} y^{14}+36 x^{7} y^{15}+2 x^{6} y^{16}-4 x^{20} y+210 x^{19} y^{2}+1792 x^{18} y^{3}+5120 x^{17} y^{4}+6626 x^{16} y^{5}+12212 x^{15} y^{6}+48516 x^{14} y^{7}+92656 x^{13} y^{8}+57504 x^{12} y^{9}+17740 x^{11} y^{10}+12450 x^{10} y^{11}+13106 x^{9} y^{12}+5658 x^{8} y^{13}+730 x^{7} y^{14}-24 x^{6} y^{15}-2 x^{5} y^{16}+x^{20}+12 x^{19} y-420 x^{18} y^{2}-3426 x^{17} y^{3}-8754 x^{16} y^{4}-9624 x^{15} y^{5}-15413 x^{14} y^{6}-55132 x^{13} y^{7}-96980 x^{12} y^{8}-55132 x^{11} y^{9}-15413 x^{10} y^{10}-9624 x^{9} y^{11}-8754 x^{8} y^{12}-3426 x^{7} y^{13}-420 x^{6} y^{14}+12 x^{5} y^{15}+x^{4} y^{16}-2 x^{19}-24 x^{18} y+730 x^{17} y^{2}+5658 x^{16} y^{3}+13106 x^{15} y^{4}+12450 x^{14} y^{5}+17740 x^{13} y^{6}+57504 x^{12} y^{7}+92656 x^{11} y^{8}+48516 x^{10} y^{9}+12212 x^{9} y^{10}+6626 x^{8} y^{11}+5120 x^{7} y^{12}+1792 x^{6} y^{13}+210 x^{5} y^{14}-4 x^{4} y^{15}+2 x^{18}+36 x^{17} y-1133 x^{16} y^{2}-8216 x^{15} y^{3}-17382 x^{14} y^{4}-14474 x^{13} y^{5}-18600 x^{12} y^{6}-55132 x^{11} y^{7}-80920 x^{10} y^{8}-39088 x^{9} y^{9}-8763 x^{8} y^{10}-4010 x^{7} y^{11}-2580 x^{6} y^{12}-782 x^{5} y^{13}-87 x^{4} y^{14}-2 x^{17}-44 x^{16} y+1596 x^{15} y^{2}+10626 x^{14} y^{3}+20548 x^{13} y^{4}+15200 x^{12} y^{5}+17740 x^{11} y^{6}+48516 x^{10} y^{7}+64492 x^{9} y^{8}+28720 x^{8} y^{9}+5694 x^{7} y^{10}+2096 x^{6} y^{11}+1072 x^{5} y^{12}+270 x^{4} y^{13}+28 x^{3} y^{14}+x^{16}+48 x^{15} y-2031 x^{14} y^{2}-12354 x^{13} y^{3}-21740 x^{12} y^{4}-14474 x^{11} y^{5}-15413 x^{10} y^{6}-39088 x^{9} y^{7}-46838 x^{8} y^{8}-19144 x^{7} y^{9}-3312 x^{6} y^{10}-910 x^{5} y^{11}-346 x^{4} y^{12}-64 x^{3} y^{13}-5 x^{2} y^{14}-48 x^{14} y+2356 x^{13} y^{2}+12976 x^{12} y^{3}+20548 x^{11} y^{4}+12450 x^{10} y^{5}+12212 x^{9} y^{6}+28720 x^{8} y^{7}+30868 x^{7} y^{8}+11504 x^{6} y^{9}+1710 x^{5} y^{10}+310 x^{4} y^{11}+70 x^{3} y^{12}+6 x^{2} y^{13}-2 x^{14}+48 x^{13} y-2480 x^{12} y^{2}-12354 x^{11} y^{3}-17382 x^{10} y^{4}-9624 x^{9} y^{5}-8763 x^{8} y^{6}-19144 x^{7} y^{7}-18372 x^{6} y^{8}-6168 x^{5} y^{9}-769 x^{4} y^{10}-72 x^{3} y^{11}-2 x^{2} y^{12}+2 x y^{13}+4 x^{13}-48 x^{12} y+2356 x^{11} y^{2}+10626 x^{10} y^{3}+13106 x^{9} y^{4}+6626 x^{8} y^{5}+5694 x^{7} y^{6}+11504 x^{6} y^{7}+9764 x^{5} y^{8}+2912 x^{4} y^{9}+300 x^{3} y^{10}+6 x^{2} y^{11}-4 x y^{12}-8 x^{12}+48 x^{11} y-2031 x^{10} y^{2}-8216 x^{9} y^{3}-8754 x^{8} y^{4}-4010 x^{7} y^{5}-3312 x^{6} y^{6}-6168 x^{5} y^{7}-4578 x^{4} y^{8}-1184 x^{3} y^{9}-95 x^{2} y^{10}+2 x y^{11}+4 x^{11}-48 x^{10} y+1596 x^{9} y^{2}+5658 x^{8} y^{3}+5120 x^{7} y^{4}+2096 x^{6} y^{5}+1710 x^{5} y^{6}+2912 x^{4} y^{7}+1844 x^{3} y^{8}+396 x^{2} y^{9}+24 x y^{10}-2 x^{10}+48 x^{9} y-1133 x^{8} y^{2}-3426 x^{7} y^{3}-2580 x^{6} y^{4}-910 x^{5} y^{5}-769 x^{4} y^{6}-1184 x^{3} y^{7}-612 x^{2} y^{8}-100 x y^{9}-4 y^{10}-44 x^{8} y+730 x^{7} y^{2}+1792 x^{6} y^{3}+1072 x^{5} y^{4}+310 x^{4} y^{5}+300 x^{3} y^{6}+396 x^{2} y^{7}+152 x y^{8}+16 y^{9}+x^{8}+36 x^{7} y-420 x^{6} y^{2}-782 x^{5} y^{3}-346 x^{4} y^{4}-72 x^{3} y^{5}-95 x^{2} y^{6}-100 x y^{7}-24 y^{8}-2 x^{7}-24 x^{6} y+210 x^{5} y^{2}+270 x^{4} y^{3}+70 x^{3} y^{4}+6 x^{2} y^{5}+24 x y^{6}+16 y^{7}+2 x^{6}+12 x^{5} y-87 x^{4} y^{2}-64 x^{3} y^{3}-2 x^{2} y^{4}+2 x y^{5}-4 y^{6}-2 x^{5}-4 x^{4} y+28 x^{3} y^{2}+6 x^{2} y^{3}-4 x y^{4}+x^{4}-5 x^{2} y^{2}+2 x y^{3}=0$.

$f=20,(\aaa,\bbb,\ccc,\ddd,\eee)=(2,4,2,5,3)/5$.

$f=20,(\aaa,\bbb,\ccc,\ddd,\eee)=(4,8,4,9,7)/10$.

$f=20,(\aaa,\bbb,\ccc,\ddd,\eee)=(16,8,4,1,3)/10$.

$f=20,(\aaa,\bbb,\ccc,\ddd,\eee)=(10,12,6,5,15)/15$.

\subsubsection*{Subcase 14.  resultant$(f(x,y,z),f(x^2,-y^2,-z^2),z)$}

$(x-1) (x^{19} y^{14}-2 x^{18} y^{15}-x^{17} y^{16}-x^{20} y^{12}-2 x^{19} y^{13}+5 x^{18} y^{14}+6 x^{17} y^{15}+4 x^{20} y^{11}+8 x^{19} y^{12}-12 x^{18} y^{13}-22 x^{17} y^{14}-4 x^{16} y^{15}+x^{15} y^{16}-12 x^{20} y^{10}-14 x^{19} y^{11}+8 x^{18} y^{12}+46 x^{17} y^{13}+29 x^{16} y^{14}+2 x^{15} y^{15}+20 x^{20} y^{9}+39 x^{19} y^{10}-8 x^{18} y^{11}-68 x^{17} y^{12}-70 x^{16} y^{13}-14 x^{15} y^{14}+2 x^{14} y^{15}-26 x^{20} y^{8}-56 x^{19} y^{9}-25 x^{18} y^{10}+70 x^{17} y^{11}+123 x^{16} y^{12}+62 x^{15} y^{13}+x^{14} y^{14}+20 x^{20} y^{7}+80 x^{19} y^{8}+38 x^{18} y^{9}-74 x^{17} y^{10}-98 x^{16} y^{11}-96 x^{15} y^{12}-32 x^{14} y^{13}+3 x^{13} y^{14}-12 x^{20} y^{6}-56 x^{19} y^{7}-80 x^{18} y^{8}+62 x^{17} y^{9}+135 x^{16} y^{10}+50 x^{15} y^{11}+58 x^{14} y^{12}+6 x^{13} y^{13}-14 x^{12} y^{14}-8 x^{11} y^{15}-4 x^{10} y^{16}+4 x^{20} y^{5}+39 x^{19} y^{6}+38 x^{18} y^{7}-54 x^{17} y^{8}-108 x^{16} y^{9}-18 x^{15} y^{10}-4 x^{14} y^{11}-44 x^{13} y^{12}+28 x^{12} y^{13}+28 x^{11} y^{14}+8 x^{10} y^{15}-x^{20} y^{4}-14 x^{19} y^{5}-25 x^{18} y^{6}+62 x^{17} y^{7}+166 x^{16} y^{8}+30 x^{15} y^{9}-97 x^{14} y^{10}+14 x^{13} y^{11}-25 x^{12} y^{12}-100 x^{11} y^{13}-50 x^{10} y^{14}-8 x^{9} y^{15}+8 x^{19} y^{4}-8 x^{18} y^{5}-74 x^{17} y^{6}-108 x^{16} y^{7}-66 x^{15} y^{8}+66 x^{14} y^{9}+21 x^{13} y^{10}-16 x^{12} y^{11}+156 x^{11} y^{12}+124 x^{10} y^{13}+28 x^{9} y^{14}-2 x^{19} y^{3}+8 x^{18} y^{4}+70 x^{17} y^{5}+135 x^{16} y^{6}+30 x^{15} y^{7}-84 x^{14} y^{8}-68 x^{13} y^{9}+34 x^{12} y^{10}-120 x^{11} y^{11}-254 x^{10} y^{12}-100 x^{9} y^{13}-14 x^{8} y^{14}+x^{19} y^{2}-12 x^{18} y^{3}-68 x^{17} y^{4}-98 x^{16} y^{5}-18 x^{15} y^{6}+66 x^{14} y^{7}+56 x^{13} y^{8}+108 x^{12} y^{9}+100 x^{11} y^{10}+172 x^{10} y^{11}+156 x^{9} y^{12}+28 x^{8} y^{13}+3 x^{7} y^{14}+2 x^{6} y^{15}+x^{5} y^{16}+5 x^{18} y^{2}+46 x^{17} y^{3}+123 x^{16} y^{4}+50 x^{15} y^{5}-97 x^{14} y^{6}-68 x^{13} y^{7}-154 x^{12} y^{8}-340 x^{11} y^{9}-258 x^{10} y^{10}-120 x^{9} y^{11}-25 x^{8} y^{12}+6 x^{7} y^{13}+x^{6} y^{14}+2 x^{5} y^{15}-2 x^{18} y-22 x^{17} y^{2}-70 x^{16} y^{3}-96 x^{15} y^{4}-4 x^{14} y^{5}+21 x^{13} y^{6}+108 x^{12} y^{7}+520 x^{11} y^{8}+448 x^{10} y^{9}+100 x^{9} y^{10}-16 x^{8} y^{11}-44 x^{7} y^{12}-32 x^{6} y^{13}-14 x^{5} y^{14}-4 x^{4} y^{15}-x^{3} y^{16}+6 x^{17} y+29 x^{16} y^{2}+62 x^{15} y^{3}+58 x^{14} y^{4}+14 x^{13} y^{5}+34 x^{12} y^{6}-340 x^{11} y^{7}-828 x^{10} y^{8}-340 x^{9} y^{9}+34 x^{8} y^{10}+14 x^{7} y^{11}+58 x^{6} y^{12}+62 x^{5} y^{13}+29 x^{4} y^{14}+6 x^{3} y^{15}-x^{17}-4 x^{16} y-14 x^{15} y^{2}-32 x^{14} y^{3}-44 x^{13} y^{4}-16 x^{12} y^{5}+100 x^{11} y^{6}+448 x^{10} y^{7}+520 x^{9} y^{8}+108 x^{8} y^{9}+21 x^{7} y^{10}-4 x^{6} y^{11}-96 x^{5} y^{12}-70 x^{4} y^{13}-22 x^{3} y^{14}-2 x^{2} y^{15}+2 x^{15} y+x^{14} y^{2}+6 x^{13} y^{3}-25 x^{12} y^{4}-120 x^{11} y^{5}-258 x^{10} y^{6}-340 x^{9} y^{7}-154 x^{8} y^{8}-68 x^{7} y^{9}-97 x^{6} y^{10}+50 x^{5} y^{11}+123 x^{4} y^{12}+46 x^{3} y^{13}+5 x^{2} y^{14}+x^{15}+2 x^{14} y+3 x^{13} y^{2}+28 x^{12} y^{3}+156 x^{11} y^{4}+172 x^{10} y^{5}+100 x^{9} y^{6}+108 x^{8} y^{7}+56 x^{7} y^{8}+66 x^{6} y^{9}-18 x^{5} y^{10}-98 x^{4} y^{11}-68 x^{3} y^{12}-12 x^{2} y^{13}+x y^{14}-14 x^{12} y^{2}-100 x^{11} y^{3}-254 x^{10} y^{4}-120 x^{9} y^{5}+34 x^{8} y^{6}-68 x^{7} y^{7}-84 x^{6} y^{8}+30 x^{5} y^{9}+135 x^{4} y^{10}+70 x^{3} y^{11}+8 x^{2} y^{12}-2 x y^{13}+28 x^{11} y^{2}+124 x^{10} y^{3}+156 x^{9} y^{4}-16 x^{8} y^{5}+21 x^{7} y^{6}+66 x^{6} y^{7}-66 x^{5} y^{8}-108 x^{4} y^{9}-74 x^{3} y^{10}-8 x^{2} y^{11}+8 x y^{12}-8 x^{11} y-50 x^{10} y^{2}-100 x^{9} y^{3}-25 x^{8} y^{4}+14 x^{7} y^{5}-97 x^{6} y^{6}+30 x^{5} y^{7}+166 x^{4} y^{8}+62 x^{3} y^{9}-25 x^{2} y^{10}-14 x y^{11}-y^{12}+8 x^{10} y+28 x^{9} y^{2}+28 x^{8} y^{3}-44 x^{7} y^{4}-4 x^{6} y^{5}-18 x^{5} y^{6}-108 x^{4} y^{7}-54 x^{3} y^{8}+38 x^{2} y^{9}+39 x y^{10}+4 y^{11}-4 x^{10}-8 x^{9} y-14 x^{8} y^{2}+6 x^{7} y^{3}+58 x^{6} y^{4}+50 x^{5} y^{5}+135 x^{4} y^{6}+62 x^{3} y^{7}-80 x^{2} y^{8}-56 x y^{9}-12 y^{10}+3 x^{7} y^{2}-32 x^{6} y^{3}-96 x^{5} y^{4}-98 x^{4} y^{5}-74 x^{3} y^{6}+38 x^{2} y^{7}+80 x y^{8}+20 y^{9}+x^{6} y^{2}+62 x^{5} y^{3}+123 x^{4} y^{4}+70 x^{3} y^{5}-25 x^{2} y^{6}-56 x y^{7}-26 y^{8}+2 x^{6} y-14 x^{5} y^{2}-70 x^{4} y^{3}-68 x^{3} y^{4}-8 x^{2} y^{5}+39 x y^{6}+20 y^{7}+2 x^{5} y+29 x^{4} y^{2}+46 x^{3} y^{3}+8 x^{2} y^{4}-14 x y^{5}-12 y^{6}+x^{5}-4 x^{4} y-22 x^{3} y^{2}-12 x^{2} y^{3}+8 x y^{4}+4 y^{5}+6 x^{3} y+5 x^{2} y^{2}-2 x y^{3}-y^{4}-x^{3}-2 x^{2} y+x y^{2})=0$.

None.

\subsubsection*{Subcase 15.   resultant$(f(x,y,z),f(-x^2,-y^2,-z^2),z)$}

$2 x^{23} y^{13}-5 x^{22} y^{14}+x^{20} y^{16}-2 x^{23} y^{12}-2 x^{22} y^{13}+28 x^{21} y^{14}-4 x^{20} y^{15}-2 x^{19} y^{16}-4 x^{24} y^{10}+10 x^{23} y^{11}-6 x^{22} y^{12}-24 x^{21} y^{13}-83 x^{20} y^{14}+12 x^{19} y^{15}+2 x^{18} y^{16}+8 x^{24} y^{9}+16 x^{23} y^{10}-42 x^{22} y^{11}+60 x^{21} y^{12}+134 x^{20} y^{13}+202 x^{19} y^{14}-24 x^{18} y^{15}-2 x^{17} y^{16}-12 x^{24} y^{8}-36 x^{23} y^{9}-63 x^{22} y^{10}+104 x^{21} y^{11}-226 x^{20} y^{12}-462 x^{19} y^{13}-408 x^{18} y^{14}+36 x^{17} y^{15}+x^{16} y^{16}+8 x^{24} y^{7}+68 x^{23} y^{8}+140 x^{22} y^{9}+172 x^{21} y^{10}-178 x^{20} y^{11}+674 x^{19} y^{12}+1168 x^{18} y^{13}+714 x^{17} y^{14}-44 x^{16} y^{15}-4 x^{24} y^{6}-36 x^{23} y^{7}-284 x^{22} y^{8}-440 x^{21} y^{9}-349 x^{20} y^{10}+154 x^{19} y^{11}-1624 x^{18} y^{12}-2386 x^{17} y^{13}-1113 x^{16} y^{14}+48 x^{15} y^{15}-2 x^{14} y^{16}+16 x^{23} y^{6}+140 x^{22} y^{7}+872 x^{21} y^{8}+1168 x^{20} y^{9}+646 x^{19} y^{10}+120 x^{18} y^{11}+3302 x^{17} y^{12}+4098 x^{16} y^{13}+1572 x^{15} y^{14}-48 x^{14} y^{15}+4 x^{13} y^{16}+10 x^{23} y^{5}-63 x^{22} y^{6}-440 x^{21} y^{7}-2174 x^{20} y^{8}-2704 x^{19} y^{9}-1116 x^{18} y^{10}-802 x^{17} y^{11}-5786 x^{16} y^{12}-6096 x^{15} y^{13}-2007 x^{14} y^{14}+48 x^{13} y^{15}-8 x^{12} y^{16}-2 x^{23} y^{4}-42 x^{22} y^{5}+172 x^{21} y^{6}+1168 x^{20} y^{7}+4736 x^{19} y^{8}+5456 x^{18} y^{9}+1798 x^{17} y^{10}+1914 x^{16} y^{11}+8868 x^{15} y^{12}+8010 x^{14} y^{13}+2332 x^{13} y^{14}-48 x^{12} y^{15}+4 x^{11} y^{16}+2 x^{23} y^{3}-6 x^{22} y^{4}+104 x^{21} y^{5}-349 x^{20} y^{6}-2704 x^{19} y^{7}-9180 x^{18} y^{8}-9688 x^{17} y^{9}-2719 x^{16} y^{10}-3320 x^{15} y^{11}-11946 x^{14} y^{12}-9394 x^{13} y^{13}-2456 x^{12} y^{14}+48 x^{11} y^{15}-2 x^{10} y^{16}-2 x^{22} y^{3}+60 x^{21} y^{4}-178 x^{20} y^{5}+646 x^{19} y^{6}+5456 x^{18} y^{7}+15928 x^{17} y^{8}+15248 x^{16} y^{9}+3796 x^{15} y^{10}+4738 x^{14} y^{11}+14266 x^{13} y^{12}+9904 x^{12} y^{13}+2332 x^{11} y^{14}-48 x^{10} y^{15}-5 x^{22} y^{2}-24 x^{21} y^{3}-226 x^{20} y^{4}+154 x^{19} y^{5}-1116 x^{18} y^{6}-9688 x^{17} y^{7}-24926 x^{16} y^{8}-21472 x^{15} y^{9}-4845 x^{14} y^{10}-5794 x^{13} y^{11}-15132 x^{12} y^{12}-9394 x^{11} y^{13}-2007 x^{10} y^{14}+48 x^{9} y^{15}+x^{8} y^{16}+28 x^{21} y^{2}+134 x^{20} y^{3}+674 x^{19} y^{4}+120 x^{18} y^{5}+1798 x^{17} y^{6}+15248 x^{16} y^{7}+35224 x^{15} y^{8}+27268 x^{14} y^{9}+5620 x^{13} y^{10}+6192 x^{12} y^{11}+14266 x^{11} y^{12}+8010 x^{10} y^{13}+1572 x^{9} y^{14}-44 x^{8} y^{15}-2 x^{7} y^{16}-83 x^{20} y^{2}-462 x^{19} y^{3}-1624 x^{18} y^{4}-802 x^{17} y^{5}-2719 x^{16} y^{6}-21472 x^{15} y^{7}-45072 x^{14} y^{8}-31396 x^{13} y^{9}-5928 x^{12} y^{10}-5794 x^{11} y^{11}-11946 x^{10} y^{12}-6096 x^{9} y^{13}-1113 x^{8} y^{14}+36 x^{7} y^{15}+2 x^{6} y^{16}-4 x^{20} y+202 x^{19} y^{2}+1168 x^{18} y^{3}+3302 x^{17} y^{4}+1914 x^{16} y^{5}+3796 x^{15} y^{6}+27268 x^{14} y^{7}+52228 x^{13} y^{8}+32896 x^{12} y^{9}+5620 x^{11} y^{10}+4738 x^{10} y^{11}+8868 x^{9} y^{12}+4098 x^{8} y^{13}+714 x^{7} y^{14}-24 x^{6} y^{15}-2 x^{5} y^{16}+x^{20}+12 x^{19} y-408 x^{18} y^{2}-2386 x^{17} y^{3}-5786 x^{16} y^{4}-3320 x^{15} y^{5}-4845 x^{14} y^{6}-31396 x^{13} y^{7}-54896 x^{12} y^{8}-31396 x^{11} y^{9}-4845 x^{10} y^{10}-3320 x^{9} y^{11}-5786 x^{8} y^{12}-2386 x^{7} y^{13}-408 x^{6} y^{14}+12 x^{5} y^{15}+x^{4} y^{16}-2 x^{19}-24 x^{18} y+714 x^{17} y^{2}+4098 x^{16} y^{3}+8868 x^{15} y^{4}+4738 x^{14} y^{5}+5620 x^{13} y^{6}+32896 x^{12} y^{7}+52228 x^{11} y^{8}+27268 x^{10} y^{9}+3796 x^{9} y^{10}+1914 x^{8} y^{11}+3302 x^{7} y^{12}+1168 x^{6} y^{13}+202 x^{5} y^{14}-4 x^{4} y^{15}+2 x^{18}+36 x^{17} y-1113 x^{16} y^{2}-6096 x^{15} y^{3}-11946 x^{14} y^{4}-5794 x^{13} y^{5}-5928 x^{12} y^{6}-31396 x^{11} y^{7}-45072 x^{10} y^{8}-21472 x^{9} y^{9}-2719 x^{8} y^{10}-802 x^{7} y^{11}-1624 x^{6} y^{12}-462 x^{5} y^{13}-83 x^{4} y^{14}-2 x^{17}-44 x^{16} y+1572 x^{15} y^{2}+8010 x^{14} y^{3}+14266 x^{13} y^{4}+6192 x^{12} y^{5}+5620 x^{11} y^{6}+27268 x^{10} y^{7}+35224 x^{9} y^{8}+15248 x^{8} y^{9}+1798 x^{7} y^{10}+120 x^{6} y^{11}+674 x^{5} y^{12}+134 x^{4} y^{13}+28 x^{3} y^{14}+x^{16}+48 x^{15} y-2007 x^{14} y^{2}-9394 x^{13} y^{3}-15132 x^{12} y^{4}-5794 x^{11} y^{5}-4845 x^{10} y^{6}-21472 x^{9} y^{7}-24926 x^{8} y^{8}-9688 x^{7} y^{9}-1116 x^{6} y^{10}+154 x^{5} y^{11}-226 x^{4} y^{12}-24 x^{3} y^{13}-5 x^{2} y^{14}-48 x^{14} y+2332 x^{13} y^{2}+9904 x^{12} y^{3}+14266 x^{11} y^{4}+4738 x^{10} y^{5}+3796 x^{9} y^{6}+15248 x^{8} y^{7}+15928 x^{7} y^{8}+5456 x^{6} y^{9}+646 x^{5} y^{10}-178 x^{4} y^{11}+60 x^{3} y^{12}-2 x^{2} y^{13}-2 x^{14}+48 x^{13} y-2456 x^{12} y^{2}-9394 x^{11} y^{3}-11946 x^{10} y^{4}-3320 x^{9} y^{5}-2719 x^{8} y^{6}-9688 x^{7} y^{7}-9180 x^{6} y^{8}-2704 x^{5} y^{9}-349 x^{4} y^{10}+104 x^{3} y^{11}-6 x^{2} y^{12}+2 x y^{13}+4 x^{13}-48 x^{12} y+2332 x^{11} y^{2}+8010 x^{10} y^{3}+8868 x^{9} y^{4}+1914 x^{8} y^{5}+1798 x^{7} y^{6}+5456 x^{6} y^{7}+4736 x^{5} y^{8}+1168 x^{4} y^{9}+172 x^{3} y^{10}-42 x^{2} y^{11}-2 x y^{12}-8 x^{12}+48 x^{11} y-2007 x^{10} y^{2}-6096 x^{9} y^{3}-5786 x^{8} y^{4}-802 x^{7} y^{5}-1116 x^{6} y^{6}-2704 x^{5} y^{7}-2174 x^{4} y^{8}-440 x^{3} y^{9}-63 x^{2} y^{10}+10 x y^{11}+4 x^{11}-48 x^{10} y+1572 x^{9} y^{2}+4098 x^{8} y^{3}+3302 x^{7} y^{4}+120 x^{6} y^{5}+646 x^{5} y^{6}+1168 x^{4} y^{7}+872 x^{3} y^{8}+140 x^{2} y^{9}+16 x y^{10}-2 x^{10}+48 x^{9} y-1113 x^{8} y^{2}-2386 x^{7} y^{3}-1624 x^{6} y^{4}+154 x^{5} y^{5}-349 x^{4} y^{6}-440 x^{3} y^{7}-284 x^{2} y^{8}-36 x y^{9}-4 y^{10}-44 x^{8} y+714 x^{7} y^{2}+1168 x^{6} y^{3}+674 x^{5} y^{4}-178 x^{4} y^{5}+172 x^{3} y^{6}+140 x^{2} y^{7}+68 x y^{8}+8 y^{9}+x^{8}+36 x^{7} y-408 x^{6} y^{2}-462 x^{5} y^{3}-226 x^{4} y^{4}+104 x^{3} y^{5}-63 x^{2} y^{6}-36 x y^{7}-12 y^{8}-2 x^{7}-24 x^{6} y+202 x^{5} y^{2}+134 x^{4} y^{3}+60 x^{3} y^{4}-42 x^{2} y^{5}+16 x y^{6}+8 y^{7}+2 x^{6}+12 x^{5} y-83 x^{4} y^{2}-24 x^{3} y^{3}-6 x^{2} y^{4}+10 x y^{5}-4 y^{6}-2 x^{5}-4 x^{4} y+28 x^{3} y^{2}-2 x^{2} y^{3}-2 x y^{4}+x^{4}-5 x^{2} y^{2}+2 x y^{3}=0$.

None.
\end{document}